%% file: root.tex
\documentclass[12pt]{article}
\usepackage{graphicx}
\usepackage{setspace}
\usepackage{bigints}
\pdfoutput=1
\usepackage{graphics} 
\usepackage{epstopdf}
\usepackage{epsfig} 
\usepackage{mathptmx} 
\usepackage{times} 
\usepackage{amsmath} 
\usepackage{amsthm}
\usepackage{amssymb}  
\usepackage{algorithmicx}
\usepackage{algorithm}
\usepackage{bm}
\usepackage[noend]{algpseudocode}
\usepackage{cite}
\usepackage{subfig}
\usepackage{gensymb}
\usepackage{booktabs}
\usepackage{url}
\usepackage{tikz,pgfplots}
\pgfplotsset{compat=newest} 
\pgfplotsset{plot coordinates/math parser=false} 

\DeclareMathOperator*{\minimize}{minimize}

\DeclareMathOperator*{\subjectto}{subject\:to}

\newtheorem{remark}{Remark}

\DeclareMathAlphabet{\pazocal}{OMS}{zplm}{m}{n}

\newlength\figureheight 
\newlength\figurewidth  

\makeatletter
\def\my@tag@font{\normalsize}
\def\maketag@@@#1{\hbox{\m@th\normalfont\my@tag@font#1}}
\let\amsmath@eqref\eqref
\renewcommand\eqref[1]{{\let\my@tag@font\relax\amsmath@eqref{#1}}}
\makeatother

\title{\LARGE \bf Estimation-aware model predictive path-following control for a general 2-trailer with a car-like tractor} 

\author
{\small Oskar~Ljungqvist$^1$, Daniel~Axehill$^1$, Henrik~Pettersson$^2$, Johan L\"ofberg$^1$\\
	\small{$^1$Department of Automatic Control, Link\"oping University, Link\"oping, Sweden.} \\ 
	\small{E-mail: \texttt{\{oskar.ljungqvist, daniel.axehill, johan.lofberg\}@liu.se.}} \\
	\small{$^2$ Scania CV, S\"odert\"alje, Sweden.} \\ 
	\small{E-mail: \texttt{henrik\_x.pettersson@scania.com.}}
}

\date{}

\begin{document} 
		
\baselineskip 16pt
	
\maketitle 

\begin{abstract} 
The design of the path-following controller is crucial for reliable autonomous vehicle operation. This design problem is especially challenging for a general 2-trailer with a car-like tractor due to the vehicle's unstable joint-angle kinematics in backward motion. Additionally, advanced sensors placed in the rear of the tractor have been proposed to solve the joint-angle estimation problem. Since these sensors typically have a limited field of view, the estimation solution introduces restrictions on the joint-angle configurations that can be estimated with high accuracy. To explicitly consider these constraints in the controller, a model predictive path-following control approach is proposed. Two approaches with different computation complexity and performance are presented. In the first approach, the joint-angle constraints are modeled as a union of convex polytopes, making it necessary to incorporate binary decision variables. The second approach avoids binary variables at the expense of a more conservative controller. In simulations and field experiments, the performance of the proposed path-following control approach is compared with a previously proposed control strategy.  
\end{abstract}


\section{Introduction}
Autonomous transport solutions and advanced driver assistance systems are experiencing massive interest in order to increase efficiency and safety, and to reduce the environmental impact related to freight transport.  
Today, autonomous driving in urban areas still faces many unsolved problems and legislation changes are needed to drive autonomously on public roads. 
In contrast, autonomous driving in closed areas such as mines and harbors are predicted to be more suitable for initial deployment of such systems.
Within these sites, different tractor-trailer combinations are frequently used for transportation of goods. 
These vehicles are composed of a car-like tractor and several passive trailers that are interconnected though hitches that are of off-axle or on-axle type. 
When the connections are of mixed hitching types, the tractor-trailer vehicle is called a general N-trailer (GNT) and when only pure on-axle hitching is present, it is referred to as a standard N-trailer (SNT).
Due to the specific kinematic properties of tractor-trailer vehicles~\cite{sordalen1993conversion,altafini1998general,CascadeNtrailernonmin}, the feedback-control problem is in general very difficult. The feedback-control problems that have been investigated in the literature for various tractor-trailer combinations are mainly path following (see e.g.,~\cite{hybridcontrol2001,rimmer2017implementation,Cascade-nSNT,bolzern1998path,virtualMorales2013,sampei1995arbitrary,minimumSweep,SamsonChainedform1995,altafini2003path,astolfi2004path,michalek2014highly,werling2014reversing}), and trajectory tracking and set-point stabilization (see e.g.,~\cite{CascadeNtrailer,michalek2018forward,kayacan2014robust,pradalier2008robust,hillary,divelbiss1997trajectory}). Additionally, to aid human drivers when reversing several advanced drivers-assistance systems concepts have been proposed (see e.g.,~\cite{michalek2014concept,evestedtLjungqvist2016,balcerak2004maneuver,michalek2013helping,hafner2017control}).

\begin{figure}[t]
	\centering
	\includegraphics[width=0.9\linewidth]{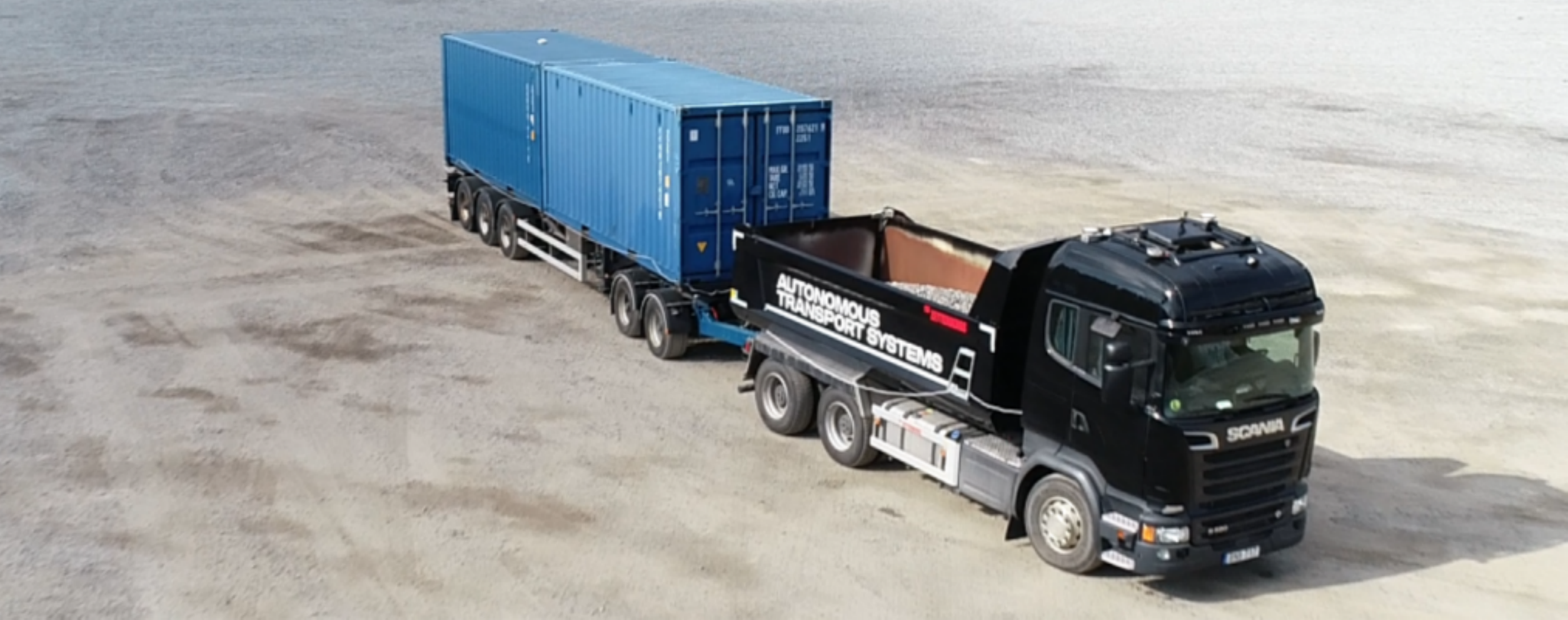}   
	\caption{The test vehicle that is used as research platform. The tractor is a modified version of a Scania R580 6x4, whereas the semitrailer nor the dolly is equipped with any sensor.}	\label{j2:fig:truck_scania} 
\end{figure}

Most of the control approaches presented in the literature address the problem of tracking a geometric trajectory or path defined in the position and orientation of the last trailer's axle.
In this work, the path-following control problem during low-speed maneuvers for a full-scale G2T with a car-like tractor (see Figure~\ref{j2:fig:truck_scania}) is considered for the case when the nominal path contains full state and control information, i.e., it is designed to operate in series with a motion planner as in~\cite{li2019tractor,evestedtLjungqvist2016planning,LjungqvistJFR2019}. 
In such an architecture, it is crucial that all nominal vehicle states are followed to avoid collision with surrounding obstacles.

The feedback-control problem considered in this work is challenging due to the tractor-trailer vehicle's structurally unstable joint-angle kinematics in backward motion and the tractor's limited curvature and curvature rate. 
Thus, if the vehicle is not accurately steered, these system properties can cause the vehicle segments to fold and enter a jackknife state. 
Additionally, nonlinear observers together with advanced sensors such as cameras, LIDARs or RADARs mounted in the rear of the car-like tractor have been proposed to solve the problem of estimating the trailer pose and the joint angles~\cite{CameraSolSaxe,caup2013video,LjungqvistJFR2019,Daniel2018,Patrik2016}. 
These solutions are promising because the system becomes independent of any trailer sensor. 
However, because such sensors typically have a limited field of view (FOV), it is important that the vehicle is controlled such that its joint angles remain in the region where high-accuracy state estimates can be computed by the used estimation solution. 

The contribution of this work is a path-following control approach for a G2T with a car-like tractor where the vehicle's physical constraints and the rear-view sensor's sensing limitations are modeled and incorporated as constraints in the controller. 
It is done by proposing a path-following control approach that is based on the framework of model predictive control (MPC)~\cite{faulwasser2015nonlinear,garcia1989model,mayne2000constrained,pedroSpatial,falcone2007predictive}. 

A preliminary version of the framework has been presented in~\cite{Ljungqvist2020ICRA}. 
Although the controller in~\cite{Ljungqvist2020ICRA} is shown to yield a significant performance enhancement compared to previous works, it is restricted to only use a single convex polytope to model the joint-angle constraints. 
On the expense of a conservative controller, the resulting MPC formulation could in our preliminary version be represented as a quadratic programming (QP) problem. 
However, due to the nonlinear mapping from the rear-view sensor's FOV to the joint-angle space, the allowed joint-angle region may in many applications become non-convex~\cite{Ljungqvist2020ICRA}. 
To alleviate the controller's conservativeness and thus enhance its performance, our previously presented preliminary results are in this work extended by a more complex modeling of the constraints on the joint angles as a union of convex polytopes, where a single polytope as in our previous work is a special case. 
By incorporating binary decision variables together with big-M modeling strategies~\cite{williams2013model}, the proposed extension results in an MPC formulation that can be cast as a mixed-integer quadratic programming (MIQP) problem. The extension significantly extends the usefulness of the method to even more advanced sensors.

The performance and computation complexity of the proposed MPC approaches are evaluated in a simulation campaign.
In the simulations and in field experiments on the full-scale test vehicle shown in Fig.~\ref{j2:fig:truck_scania}, the performance of the proposed predictive path-following control approach in terms of suppressing disturbances and recovering from nontrivial initial states, is compared to the proposed path-following controller in~\cite{LjungqvistJFR2019} where the vehicle's physical constraints and the estimation solution's sensing-limitations are neglected. 

The remainder of the paper is organized as follows. In Section~\ref{j2:sec:Modeling}, the path-following error model is derived and in Section~\ref{j2:sec:MPC}, the proposed MPC approach is presented. 
The control design as well as the modeling of the vehicle's physical and sensing limiting constraints are presented in Section~\ref{j2:sec:controller_synthesis}. 
Results from simulations and field experiments on a full-scale test vehicle are presented in Section~\ref{j2:sec:results}. 
Finally, the paper is concluded in Section~\ref{j2:sec:conclusions} by summarizing the contributions and a discussion of directions for future work.

\section{Vehicle model}
\label{j2:sec:Modeling}
\begin{figure}[t]
	\begin{center}
		\includegraphics[width=0.8\linewidth]{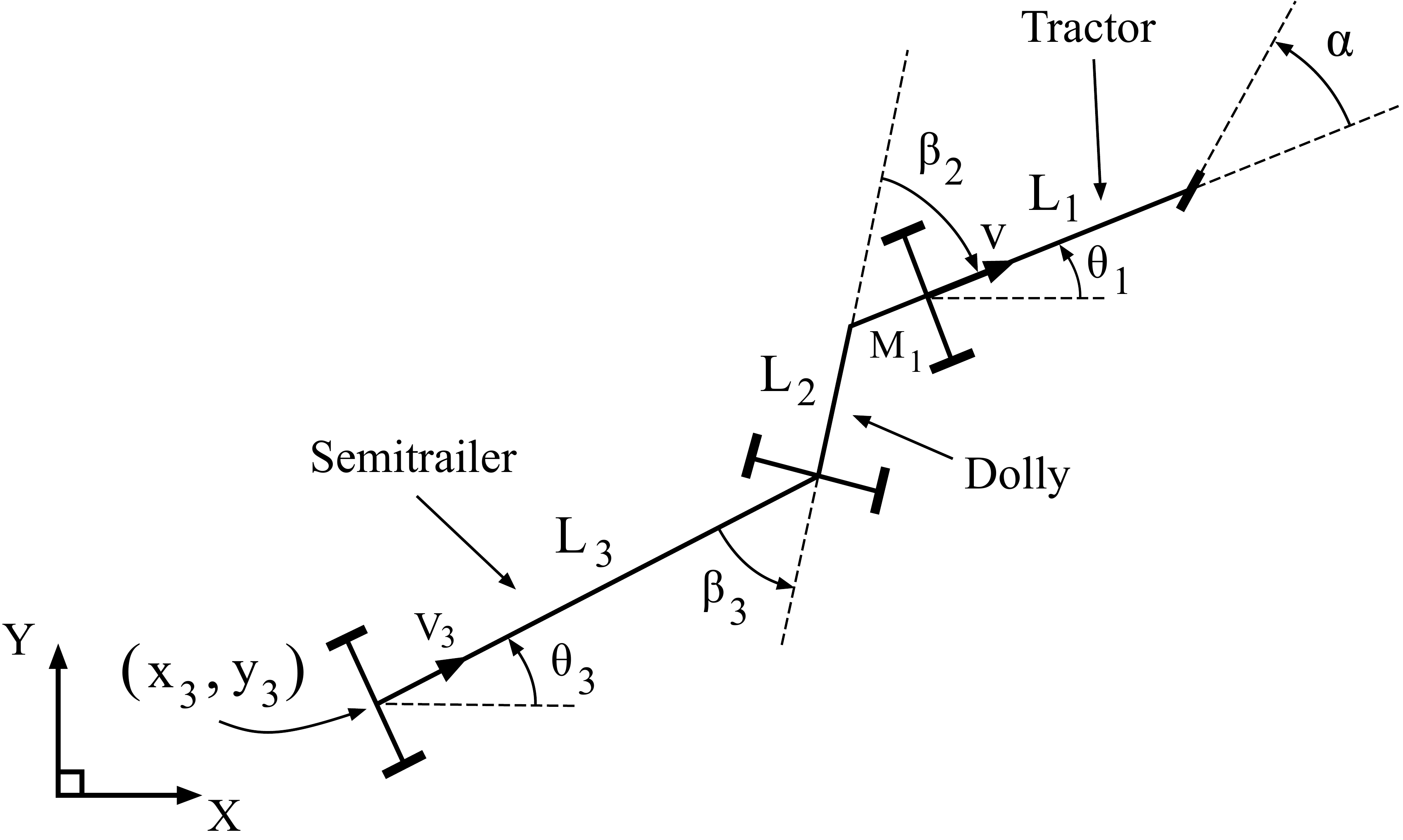}   
		\caption{A schematic description of the geometric lengths, states and control inputs that are of relevance for modeling the G2T with a car-like tractor.} 
		\label{j2:fig:schematic_model_description}
	\end{center} 
\end{figure}

The G2T with a car-like tractor is schematically illustrated in Figure~\ref{j2:fig:schematic_model_description}. 
The vehicle is composed of three interconnected vehicle segments including a car-like tractor, a dolly and a semitrailer. 
It is assumed that the vehicle has an off-axle hitch connection between the car-like tractor and the dolly, and an on-axle connection between the dolly and the semitrailer. 
The state vector $x=[x_3\hspace{5pt}y_3 \hspace{5pt} \theta_3 \hspace{5pt} \beta_3 \hspace{5pt} \beta_2]^T$ is used to represent a configuration of the vehicle, where $(x_3,y_3)$ is the position of the center of the semitrailer's axle, $\theta_3$ is the orientation of the semitrailer, $\beta_3$ is the joint angle between the semitrailer and the dolly, and $\beta_2$ is the joint angle between the dolly and the car-like tractor. 
The length $L_3$ represents the distance between the semitrailer's axle and the dolly's axle, $L_2$ is the distance between the dolly's axle and the off-axle hitching connection at the car-like tractor, $M_1$ is the signed hitching offset at the tractor (positive in Figure~\ref{j2:fig:schematic_model_description}), and $L_1$ denotes the wheelbase of the car-like tractor. 
The car-like tractor is assumed to be front-wheeled steered with perfect Ackermann steering geometry, where $\alpha$ denotes its steering angle. 
The control inputs are the car-like tractor's curvature $u = \frac{\tan \alpha}{L_1}$ and the longitudinal velocity $v$ of its rear axle. 
Since low-speed maneuvers are considered, a kinematic model
is used to describe the vehicle. The model has been presented in, e.g.,~\cite{hybridcontrol2001,LjungqvistJFR2019} and is derived based on various assumptions including rolling without slipping of the wheels and that the vehicle is operating on a flat surface. The kinematic vehicle model is given by
\begin{subequations}
	\label{j2:eq:model}
	\begin{align} 
	\dot{x}_3 &= v_3 \cos \theta_3,  \label{j2:eq:model1}\\
	\dot{y}_3 &= v_3 \sin \theta_3,  \label{j2:eq:model2} \\
	\dot{\theta}_3 &= v_3 \frac{\tan \beta_3 }{L_3}, \label{j2:eq:model3}\\
	\dot{\beta}_3 &= v_3\left(\frac{\sin\beta_2 - M_1\cos\beta_2u}{L_2  C_1(\beta_2,\beta_3,u)} - \frac{\tan\beta_3}{L_3}\right), \label{j2:eq:model4}\\
	\dot{\beta}_2 &= v_3 \left( \frac{u - \frac{\sin \beta_2}{L_2} + \frac{M_1}{L_2}\cos \beta_2 u}{C_1(\beta_2,\beta_3,u)}\right), \label{j2:eq:model5}
	\end{align}
\end{subequations}
where $C_1(\beta_2,\beta_3,u)$ is defined as 
\begin{align}
C_1(\beta_2,\beta_3,u) = \cos\beta_3\left(\cos{\beta_2} + M_1\sin\beta_2u\right),
\label{j2:eq:C1}
\end{align}
which describes the relationship, $v_3 = vC_1(\beta_2,\beta_3,u)$, between the longitudinal velocity of the semitrailer's axle, $v_3$, and the velocity of the car-like tractor's rear axle, $v$. 
When $C_1(\beta_2,\beta_3,u)=0$, the system in~\eqref{j2:eq:model} is singular~\cite{altafini1998general} and therefore fundamentally difficult to control. 
It is therefore further assumed that \mbox{$C_1(\beta_2,\beta_3,u)>0$}. 

The model in~\eqref{j2:eq:model} is compactly represented as $\dot x = v_3f(x,u)$. 
Since $v_3$ enters bilinear in~\eqref{j2:eq:model}, time-scaling~\cite{sampei1986time,LjungqvistJFR2019} can be applied to eliminate the longitudinal speed dependence and it is therefore without loss of generality further assumed that the velocity of the tractor is restricted to $v\in\{-1,1\}$, where $v=1$ denotes forward motion and $v=-1$ backward motion. 
The direction of motion is essential for the stability of the system~\eqref{j2:eq:model}, where the joint-angle kinematics are structurally unstable in backward motion ($v < 0$), where it risks to fold and enter a jack-knife state~\cite{hybridcontrol2001}. 
In forward motion ($v > 0$), these modes are stable but in case of positive off-axle hitching ($M_1>0$), some of the system's output channels poses non-minimum phase properties (see, e.g.,~\cite{CascadeNtrailernonmin} for an analysis).

\subsection{Constraints}
The car-like tractor is assumed to have physical bounds on its curvature $u$ and curvature rate $\dot u$, which are modeled as box constraints
\begin{align}
\label{j2:control_constraints}
|u|\leq u_{\text{max}},\quad |\dot u|\leq \dot u_{\text{max}},
\end{align} 
where the positive constants $u_{\text{max}}$ and $\dot u_{\text{max}}$ denote maximum curvature and curvature rate, respectively. 
In practice, advanced sensors mounted in the rear of the car-like tractor have been proposed as solutions to the joint-angle estimation problem~\cite{CameraSolSaxe,caup2013video,LjungqvistJFR2019}. 
Such sensors typically have a limited FOV which enforce non-convex restrictions on the set of joint angles that can be estimated with high accuracy. 
Additionally, constraints on the joint angles that prevent the vehicle from entering a jack-knife state should also be considered. 
The above mentioned restrictions on $\beta_2$ and $\beta_3$ are assumed to be described by a union of $n\in\mathbb Z_{+}$ convex polytopes $\mathbb P_i$ in the form
\begin{align}
\label{j2:joint_angle_constaints}
\mathbb X = \bigcup\limits_{i=1}^{n}\mathbb P_i =  \bigcup\limits_{i=1}^{n}\left\{(\beta_3,\beta_2)\in\mathbb R^2 \middle| H_i\left(\begin{matrix}
\beta_3 & \beta_2
\end{matrix}\right)^T \leq h_i \right\} ,
\end{align}
where $H_i\in\mathbb R^{m_i \times 2}$, $h_i\in\mathbb R^{m_i}$ and $m_i\in\mathbb Z_{+}$ for $i=1,\hdots,n$. 
The set $\mathbb X$ is assumed to be closed, compact and contain the origin $(\beta_2,\beta_3)=(0,0)$ in its interior. 
For compactness, the constraint in \eqref{j2:joint_angle_constaints} is represented as $(\beta_{3},\beta_{2})\in\mathbb X$. 
Note that if $n\geq 2$, the set in~\eqref{j2:joint_angle_constaints} is in general a non-convex set and if $n=1$, it is a set of linear inequality constraints. 

Even though the constraints in~\eqref{j2:control_constraints} and~\eqref{j2:joint_angle_constaints} have been considered by the motion planner while computing the motion plan, disturbances are always present during plan execution making it important to also consider them in the controller.

\subsection{Path-following error model}
Given a nominal trajectory $(x_r(\cdotp),u_r(\cdotp),v_{3r}(\cdotp))$ satisfying the model of the G2T with a car-like tractor~\eqref{j2:eq:model}:
\begin{align}
\dot{x}_r = v_{3r}f( x_r,u_r), \label{j2:eq:tray}
\end{align}
which is assumed to satisfy the constraints in~\eqref{j2:control_constraints} and~\eqref{j2:joint_angle_constaints}. 
Except from satisfying the constraints, the objective of the path-following controller is to control the tractor's curvature $u$ such that the motion plan obtained from the motion planner is followed with a small and bounded path-following error. 
Similar to~\cite{LjungqvistJFR2019}, this is performed by first deriving a path-following error model that describes the vehicle in terms of deviation from a nominal path. 
Given the vehicle's current state $x(t)$, define $s(t)$ as the distance traveled by the position of the semitrailer's axle onto its projection to its nominal path $(x_{3r}(\cdotp), y_{3r}(\cdotp))$ up to time $t$.  
Since $\frac{\text dx_r}{\text dt}=\frac{\text dx_r}{\text ds}|v_{3r}|$, the nominal trajectory in~\eqref{j2:eq:tray} can instead be interpreted as a nominal path
\begin{align}
\frac{\text d x_r}{\text ds} = \bar v_{3r}f(x_r,u_r), \label{j2:eq:tray_s}
\end{align}
where $\bar v_{3r}=\text{sign}(v_{3r})\in\{-1,1\}$ represents the nominal motion direction, i.e., $\bar v_{3r}=1$ represents forward motion and $\bar v_{3r}=-1$ backward motion.

\begin{figure}[t!]
	\centering
	\includegraphics[width=0.8\linewidth]{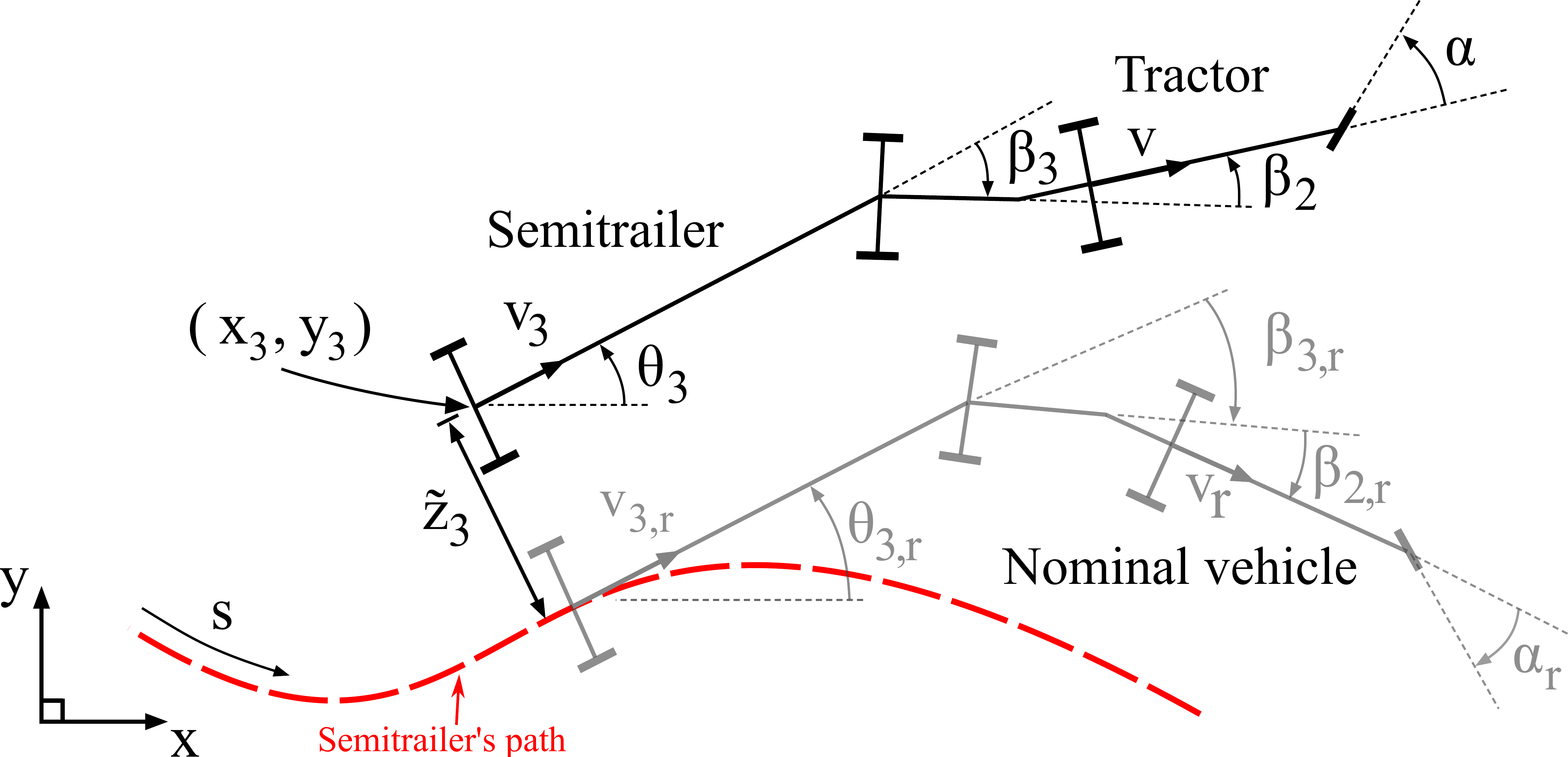}
	\caption{The G2T with a car-like tractor in the Frenet frame path coordinate system.}
	\label{j2:fig:frenet}%
\end{figure}

To exploit that a nominal path~\eqref{j2:eq:tray_s} is provided, the G2T with car-like tractor in~\eqref{j2:eq:model} is first modeled in terms of deviation from this nominal path, as illustrated in Figure~\ref{j2:fig:frenet}.
Denote $\tilde z_3(t)$ as the signed lateral distance between the center of the semitrailer's axle $(x_3(t),y_3(t))$ and its projection to its nominal path in $(x_{3r}(\cdotp), y_{3r}(\cdotp))$. 
Define the orientation error of the semitrailer as $\tilde{\theta}_3(t)=\theta_3(t)-\theta_{3r}(s(t))$, and define joint-angle errors as $\tilde{\beta}_3(t)=\beta_3(t)-\beta_{3r}(s(t))$ and $\tilde{\beta}_2(t)=\beta_2(t)-\beta_{2r}(s(t))$. 
Finally, define the controlled curvature deviation as $\tilde{u}(t) =u(t) - u_r(s(t))$.
Using the Frenet-frame transformation together with the chain rule (see~\cite{LjungqvistJFR2019} for details), the G2T with a car-like tractor~\eqref{j2:eq:model} can be modeled in terms of deviation from the nominal path~\eqref{j2:eq:tray_s} as
\begingroup\makeatletter\def\f@size{11}\check@mathfonts
\begin{subequations}
	\label{j2:eq:model_frenet_frame}
	\begin{align}
	\dot s &= v_3 \frac{\bar v_{3r}\cos \tilde \theta_3}{1-\kappa_{3r} \tilde z_3}, 
	\label{j2:eq:model_s1}
	\\ 
	\dot{\tilde z}_3 &= v_3 \sin \tilde \theta_3, 
	\label{j2:eq:model_s2}
	\\
	\dot{\tilde \theta}_3 &= v_3 \left( \frac{\tan(\tilde{\beta}_3+\beta_{3r})}{L_3} - \frac{\kappa_{3r}\cos \tilde \theta_3}{1-\kappa_{3r} \tilde z_3} \right), \label{j2:eq:model_s3}
	\\
	\dot{\tilde \beta}_3 &= v_3\left(\frac{\sin(\tilde \beta_2+\beta_{2r})-M_1\cos(\tilde \beta_2+\beta_{2r})(\tilde u+ u_r)}{L_2 C_1(\tilde \beta_2+\beta_{2r},\tilde \beta_3+\beta_{3r}, \tilde u + u_r)} - \frac{\tan(\tilde \beta_3+\beta_{3r})}{L_3} \nonumber  \right.  \\ &\left. -\frac{\cos{\tilde{\theta}_3}}{1-\kappa_{3r}\tilde z_3}\left[\frac{\sin\beta_{2r} - M_1\cos\beta_{2r} u_r}{L_2C_1(\beta_{2r},\beta_{3r},u_r)}-\kappa_{3r}\right]\right), 
	\label{j2:eq:model_s4} \\
	\dot{\tilde \beta}_2 &=v_3\left( \left[ \frac{\tilde u+ u_r - \frac{\sin(\tilde \beta_2+\beta_{2r})}{L_2} + \frac{M_1}{L_2}\cos(\tilde \beta_2+\beta_{2r})(\tilde u+ u_r)}{C_1(\tilde \beta_2+\beta_{2r},\beta_3+\beta_{3r}, \tilde u+ u_r)}\right]  -\frac{\cos{\tilde{\theta}_3}}{1-\kappa_{3r}\tilde z_3}\left[ \frac{u_r - \frac{\sin \beta_{2r}}{L_2} + \frac{M_1}{ L_2}\cos \beta_{2r}u_r}{C_1(\beta_{2r},\beta_{3r}, u_r)}\right]\right), \label{j2:eq:model_s5}
	\end{align}
\end{subequations}
\endgroup
where
\begin{equation}
\kappa_{3r}(s)=\frac{\text d\theta_{3r}}{\text ds}=\frac{\tan \beta_{3r}(s)}{L_3},
\end{equation}
is the nominal curvature for the center of the semitrailer's axle. The transformation to the Frenet frame path-coordinate system is valid as long as $\tilde z_3$ and $\tilde\theta_3$ satisfy 
\begin{align}\label{j2:frenet_frame_transformation_constraints}
\kappa_{3r}(s)\tilde z_3<1, \quad |\tilde\theta_3|<\pi/2.
\end{align}
The former constraint can be approximated as $|\tilde z_3|\leq\tilde z_3^{\text{max}}$, where the constant $\tilde z_3^{\text{max}} < 1/\kappa_{3r}^{\text{max}}$, where $\kappa_{3r}^{\text{max}}$ is the maximum curvature of the nominal path in $(x_{3r}(\cdotp),y_{3r}(\cdotp))$. 
Note that $\bar v_{3r}$ is included in~\eqref{j2:eq:model_s1} to make $\dot s> 0$ as long as the constraints in~\eqref{j2:frenet_frame_transformation_constraints} are satisfied, and the semitrailer's velocity $v_3$ and its nominal motion direction $\bar v_{3r}$ have the same sign. 
Since it is required that $C_1(\beta_2,\beta_3,u)>0$ and the relationship $v = v_3C_1(\beta_2,\beta_3,u)$ holds, this is equivalent to that tractor's velocity $v$ is selected such that $\text{sign}(v)=\bar v_{3r}$. 

Define the path-following error state $\tilde x = [\tilde z_3 \hspace{5pt} \tilde\theta_3\hspace{5pt} \tilde\beta_3 \hspace{5pt} \tilde\beta_2]^T$, where its model is given by~\eqref{j2:eq:model_s2}--\eqref{j2:eq:model_s5}. 
It is easily verified that the origin $(\tilde x,\tilde u)=(0,0)$ to this system is an equilibrium point for all $s(t)$. 
Moreover, since the velocity of the car-like tractor $v$ is selected such that $\dot s(t)>0$, it is possible to perform time-scaling~\cite{sampei1986time} and eliminate the time-dependency presented in~\eqref{j2:eq:model_s2}--\eqref{j2:eq:model_s5}. 
Using the chain rule, it holds that $\frac{\text d\tilde x}{\text  ds}=\frac{\text d\tilde x}{\text  dt}\frac{1}{\dot s}$, and the spatial version of the path-following error model~\eqref{j2:eq:model_s2}--\eqref{j2:eq:model_s5} becomes 
\begin{subequations}
	\label{j2:eq:MPC_spatial_path_following_error_model}
	\begin{align}
	&\frac{\text d \tilde z_3 }{\text d s} = \bar v_{3r}   (1-\kappa_{3r} \tilde z_3)  \tan \tilde \theta_3, \\
	&\frac{\text d \tilde \theta_3 }{\text d s} = \bar v_{3r}\left( \left[1-\kappa_{3r} \tilde z_3\right]\frac{\tan(\tilde{\beta}_3+\beta_{3r})}{L_3\cos \tilde \theta_3} - \kappa_{3r} \right), \\
	&\frac{\text d \tilde \beta_3 }{\text d s} = \bar v_{3r}  \left( \frac{1-\kappa_{3r} \tilde z_3}{\cos \tilde \theta_3}\left[\frac{\sin(\tilde \beta_2+\beta_{2r})-M_1\cos(\tilde \beta_2+\beta_{2r}) (\tilde u+ u_r)}{L_2 C_1(\beta_{2r}+\tilde\beta_2,\beta_{3r}+\tilde\beta_3,u_r+\tilde u)} \nonumber \right. \right. \\ &\left. \left.- \frac{\tan(\tilde \beta_3+\beta_{3r})}{L_3 }\right] 
	-\left[\frac{\sin\beta_{2r} -M_1 \cos\beta_{2r}u_r}{L_2 C_1(\beta_{2r},\beta_{3r},u_r)}-\kappa_{3r}\right]\right), \\
	&\frac{\text d \tilde \beta_2 }{\text d s} =  \bar v_{3r} \left( \frac{1-\kappa_{3r}\tilde z_3}{\cos{\tilde{\theta}_3}} \left[\frac{\tilde u+ u_r  + \frac{M_1}{L_2}\cos(\tilde \beta_2+\beta_{2r})(\tilde u+ u_r)}{C_1(\beta_{2r}+\tilde\beta_2,\beta_{3r}+\tilde\beta_3,u_r+\tilde u)}  \nonumber \right. \right. \\
	&\left. \left. -\frac{\sin(\tilde \beta_2+\beta_{2r})}{L_2{C_1(\beta_{2r}+\tilde\beta_2,\beta_{3r}+\tilde\beta_3,u_r+\tilde u)}}\right]  -\left[ \frac{u_r - \frac{\sin \beta_{2r}}{L_2} + \frac{M_1}{L_2}\cos \beta_{2r}u_r}{C_1(\beta_{2r},\beta_{3r},u_r)}\right]\right),
	\end{align}
\end{subequations}
which can compactly be represented as $\frac{\text d\tilde x}{\text ds}=\tilde f(s,\tilde x,\tilde u)$, where the origin $(\tilde x,\tilde u)=(0,0)$ is an equilibrium point for all $s$. 

\section{Model predictive path-following controller}
\label{j2:sec:MPC}
The task of the MPC controller is to control the tractor's curvature $u$ such that the path-following error is minimized while the vehicle's constraints~\eqref{j2:control_constraints}--\eqref{j2:joint_angle_constaints} are satisfied for all time instances. 
Since the vehicle's joint angles are restricted to a union of convex polytopes~\eqref{j2:joint_angle_constaints}, binary decision variables are in general needed to incorporate these constraints within the MPC controller. 
Therefore, due to the combination of binary and continuous variables, the resulting MPC problem will be in the form of a mixed-integer programming (MIP) problem. 
To obtain an MIP problem that has the potential of being solved at a sufficiently high rate, each subproblem when the binary variables are fixed should be simple to solve. 
As a consequence, the aim is to derived an MPC formulation that can be converted into the form of an MIQP problem, where efficient state-of-the-art commercial MIQP solvers exist, such as Gurobi~\cite{gurobi} and CPLEX~\cite{cplex}.  

First, the nonlinear path-following error model~\eqref{j2:eq:MPC_spatial_path_following_error_model} is linearized around the origin $(\tilde x,\tilde u) = (0,0)$:
\begin{align}
\label{j2:linear_cont_model}
\frac{\text d \tilde x}{\text ds} = A(s)\tilde x +  B(s)\tilde u,
\end{align}
where the distance-varying matrices $A(s)$ and $B(s)$ are presented in Appendix A. 
With the sampling distance $\Delta_s$, Euler-forward discretization yields a discrete approximation of~\eqref{j2:linear_cont_model} in the form
\begin{align}
\label{j2:d_lin_sys}
\tilde x_{k+1} = F_k\tilde x_k + G_k\tilde u_k,
\end{align}  
where
\begin{align}
\label{j2:d_system_matrix}
F_k = I +\Delta_s A_k, \quad
G_k = \Delta_sB_k.
\end{align}    
Since the tractor's curvature is \mbox{$u_k=\tilde u_k+u_{r,k}$}, the deviation in the curvature is bounded as
\begin{align}
\label{j2:curvature_constaint_MPC}
-u_{\text max} \leq \tilde u_k + u_{r,k}  \leq u_{\text max}.
\end{align}  
Moreover, since $\dot s>0$ the constraint on the tractor's curvature rate \mbox{$|\dot u|\leq \dot u_{\text{max}}$} can be described in $s$ using the chain rule as
\begin{align}
\label{j2:rate-limit-constraints-ds}
\left|\frac{\text du}{\text ds}\right|\leq\frac{\dot u_{\text{max}}}{\dot s} = \frac{1-\kappa_{3r} \tilde z_3}{|v| C_1(\beta_2,\beta_3,u)\cos \tilde \theta_3}\dot u_{\text{max}},
\end{align}
since $v_3 = v C_1(\beta_2,\beta_3,u)$. Locally around the nominal path $(\tilde x,\tilde u)=(0,0)$, it holds that $\cos\tilde\theta_3\approx 1$ and $\kappa_{3r} \tilde z_3\approx 0$. 
Thus, to avoid coupling between $\tilde x$ and $\tilde u$, the constraint in~\eqref{j2:rate-limit-constraints-ds} is approximated as
\begin{align}
\label{j2:rate-limit-constraints-ds2}
\left|\frac{\text du}{\text ds}\right| \leq \frac{\dot u_{\text{max}}}{|v|C_1(\beta_{2r}(s),\beta_{3r}(s),u_r(s))} \triangleq c_{\text{max}}(s).
\end{align} 
By discretizing~\eqref{j2:rate-limit-constraints-ds2} using Euler forward, the rate limit on the controlled curvature deviation can be described by the following slew-rate constraint
\begin{align}
\label{j2:curvature_rate_constaint_MPC}
-c_{\text{max},k}\Delta_s \leq \tilde u_{k} - \tilde u_{k-1} - \bar u_{r,k} \leq c_{\text{max},k}\Delta_s,   
\end{align}
where $\bar u_{r,k} = u_{r,k} - u_{r,k-1}$. Denote the linear inequality constraints in~\eqref{j2:curvature_constaint_MPC} and~\eqref{j2:curvature_rate_constaint_MPC} as \mbox{$\tilde u_k\in\tilde {\mathbb U}_k$}. 
What remains is to describe the non-convex joint-angle constraint~\eqref{j2:joint_angle_constaints} as a function of $\tilde x$ and to model it in an MIQP representable form. 
Since the joint angles are \mbox{$\beta_{j,k}=\beta_{jr,k}+\tilde\beta_{j,k}$}, $j=2,3$, each convex polytope $\mathbb P_i$, $i=1,\hdots,n$ that is used to model the constraint on the joint angles~\eqref{j2:joint_angle_constaints} can be written as
\begin{align}
\label{j2:joint_angle_constaints_error_states}
H_i\left[\begin{matrix}
\tilde\beta_{3,k} \\ \tilde\beta_{2,k}
\end{matrix}\right] \leq h_i - H_i\left[\begin{matrix}
\beta_{3r,k} \\ \beta_{2r,k}
\end{matrix}\right] = \bar h_{i,k},
\end{align}
which is denoted as $(\tilde\beta_{3,k},\tilde\beta_{2,k})\in\tilde{\mathbb P}_{i,k}$. 
Now, at each sample $k$ along the future prediction horizon, a binary decision variable is introduced $\delta_{i,k}=\{0,1\}$ to indicate if $(\tilde\beta_{3,k},\tilde\beta_{2,k}) \in \tilde{\mathbb P}_{i,k}$. 
Using logical implications, this can be enforced using hybrid logic on the constraint on the joint angles with
\begin{align}
\delta_{i,k} \rightarrow (\tilde\beta_{3,k},\tilde\beta_{2,k}) \in \tilde{\mathbb P}_{i,k}, \quad i=1,\hdots,n.
\end{align}
Define $\bar\delta_k=[\delta_{1,k}\hspace{5pt}\delta_{2,k}\hspace{5pt}\hdots\hspace{5pt}\delta_{n,k}]^T\in\{0,1\}^n$ as a binary vector and by adding the constraint $\sum_{i=1}^{n}\delta_{i,k}=1$, the logical model ensures that the joint angles are in at least one of the convex polytopes at sample $k$. 
The logical implications are easily converted to linear inequality constraints using \mbox{big-M} modeling strategies, but the details are omitted, as the employed modeling tool YALMIP~\cite{lofberg2004yalmip} will do this step automatically given the logical model.

Given the current path-following error state $\tilde x(s(t))$, the MPC problem with prediction horizon $N$ is defined as 
\begin{subequations}
	\label{j2:MPC_problem}
	\begin{align}
	\minimize_{\mathbf{\tilde x},\hspace{0.1ex} \bm{\tilde u},\hspace{0.1ex} \bm{\delta}} \hspace{2ex} &V_N(\mathbf{\tilde x},\bm{\tilde u}) = V_f(\tilde x_N) + \sum_{k=0}^{N-1} l(\tilde x_k,\tilde u_k) \\
	\subjectto \hspace{1ex} 
	&\tilde x_{k+1} = F_k\tilde x_k + G_k\tilde u_k, \quad k = 0,1,\hdots,N-1, \\
	&\delta_{i,k} \rightarrow (\tilde\beta_{3,k},\tilde\beta_{2,k})\in \tilde{\mathbb P}_{i,k}, \hspace{2pt} i=1,\hdots,n,\quad k = 0,1,\hdots,N-1, \\
	&\sum_{i=1}^{n}\delta_{i,k}=1, \quad \tilde u_k\in\tilde{\mathbb U}_k, \quad k = 0,1,\hdots,N-1, \\
	&|\tilde z_{3,k}|\leq \tilde z_{3}^\text{max}, \quad |\tilde\theta_{3,k}|\leq\tilde\theta_{3}^\text{max}, \quad k = 0,1,\hdots,N-1,	\label{j2:MPC_problem_constraint_lat_head}\\
	&\tilde x_0  = \tilde x(s(t)) \text{ given,}
	\end{align}
\end{subequations}
where $\mathbf{\tilde x}^T = [\tilde x_0^T \hspace{5pt} \tilde x_1^T \hspace{5pt} \hdots \hspace{5pt} \tilde x_N^T]$ is the predicted path-following error state-vector sequence, $\mathbf{\tilde u} = [\tilde u_0 \hspace{5pt} \tilde u_1 \hspace{5pt} \hdots \hspace{5pt} \tilde u_{N-1}]^T$ is the curvature-deviation sequence and $\bm{\delta}^T = [\bar \delta_0^T\hspace{5pt}\bar\delta_1^T\hspace{5pt}\hdots\hspace{5pt}\bar\delta_{N-1}^T]$ is the binary-vector sequence. 
In~\eqref{j2:MPC_problem_constraint_lat_head}, the constants $\tilde z_{3}^\text{max}$ and $\tilde\theta_{3}^\text{max}<\pi/2$ are the semitrailer's maximum lateral and orientation error, respectively. 
The stage-cost is chosen as $l(\tilde x_k,\tilde u_k) = \tilde x_k^T Q \tilde x_k + \tilde u_k^2$ and the terminal-cost $V_f(\tilde x_N)=\tilde x_N^T P_N \tilde x_N$, where $Q\succeq 0$ and $P_N\succ 0$ are design matrices of appropriate dimensions. 
Since the cost function $V_N$ is quadratic and there are only linear equality and inequality constraints for a fixed binary-vector sequence, the MPC problem~\eqref{j2:MPC_problem} can be converted into an MIQP problem. 

\begin{remark}
	If the constraint on the joint angles~\eqref{j2:joint_angle_constaints} is modeled as a single polytope ($n=1$), binary variables are not needed and the MPC problem~\eqref{j2:MPC_problem} simplifies to a QP problem. 
\end{remark}

At each sampling instance, the MPC problem~\eqref{j2:MPC_problem} is solved to obtain the optimal open-loop controlled curvature deviation sequence $\bm{\tilde u}^*$. 
As in standard receding horizon control, only the first control input $\tilde u_0^*$ is used
\begin{align}
u(t) = u_r(s(t)) + \tilde u_0^*, 
\end{align}
and the optimizing problem~\eqref{j2:MPC_problem} is repeatedly solved at each sampling instance using new state information. 
Note that the MPC controller only computes the feedback part of the control signal, as the optimal feed forward $u_r(s(t))$ already is provided by the motion planner.

\section{Estimation-aware controller design}
\label{j2:sec:controller_synthesis}
In this section, the control design of the proposed model predictive path-following controller~\eqref{j2:MPC_problem} is presented. It is tailored for the case when a rear-view LIDAR with a limited FOV is used as part of the estimation solution~\cite{LjungqvistJFR2019}. However, note that similar techniques can also be used if a rear-view camera or RADAR is used as sensor.

\begin{figure}[t]
	\centering
	\includegraphics[width=0.9\linewidth]{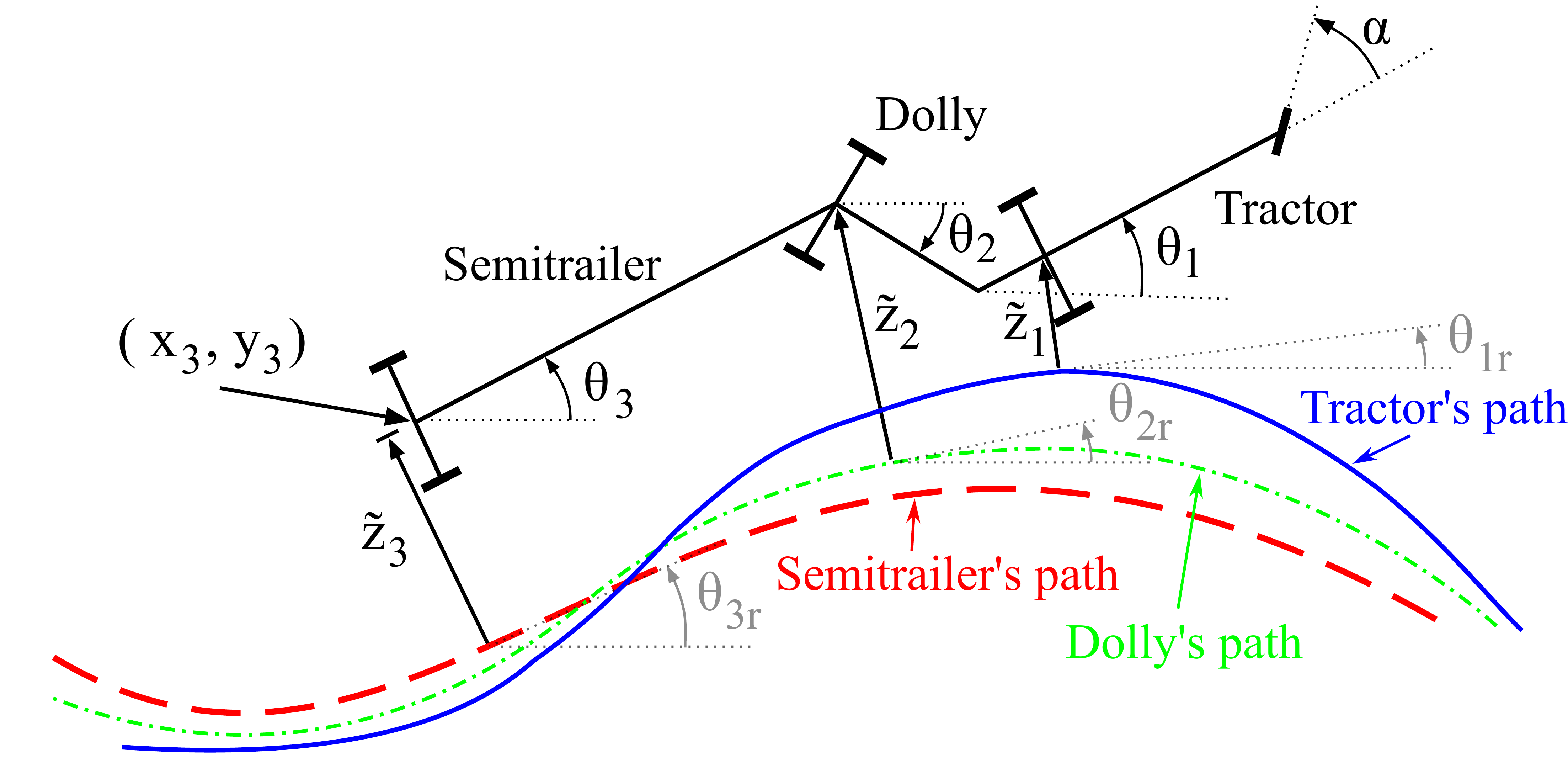}
	\caption{Illustration of the additional path-following error states that are used in the design of the MPC controller.}
	\label{j2:fig:truck_frenet_advanced}
\end{figure}

\subsection{Design of cost function}
Since a nominal path that satisfies the vehicle model~\eqref{j2:eq:tray_s} is provided, it is possible to calculate the nominal paths for the position and orientation of the dolly as well as the car-like tractor using the nominal state path $x_r(\cdotp)$ together with holonomic relationships~\cite{altafini1998general}. 
To reduce the risk of colliding with surrounding obstacles, it is desired that the MPC controller is tuned such that the transient response of all path-following errors are prioritized. 
Denote $\tilde z_1$ and $\tilde z_2$ as the signed lateral distances of the axle of the dolly and the car-like tractor onto their nominal paths, see Figure~\ref{j2:fig:truck_frenet_advanced}. 
Moreover, define their corresponding orientation errors as $\tilde \theta_1=\theta_1 -\theta_{1r}$ and $\tilde \theta_2 = \theta_2 -\theta_{2r}$, respectively. 
In general, it is not possible to derive a closed-form expression to describe these additional path-following error states as a function of the modeled path-following error states $\tilde x$ (see~\cite{altafini2003path} for details). 
However, for the special case of a straight nominal path, closed-form expressions exist and the signed lateral errors $\tilde z_2$ and $\tilde z_1$ can be described as 
\begin{subequations}
	\label{j2:eq:add_lateral_errors}
	\begin{align}
	\tilde z_2 &= \tilde z_3 + L_3\sin\tilde\theta_3, \\
	\tilde z_1 &= \tilde z_2 + L_2\sin(\tilde\theta_3+\tilde\beta_3) + M_1\sin(\tilde\theta_3+\tilde\beta_3+\tilde\beta_2),
	\end{align}
\end{subequations}
and their corresponding orientation errors $\tilde \theta_2$ and $\tilde \theta_1$ as
\begin{subequations}
	\begin{align}
	\tilde \theta_2 &= \tilde\theta_3+\tilde\beta_3, \\
	\tilde \theta_1 &= \tilde\theta_3+\tilde\beta_3+\tilde\beta_2.
	\end{align}
\end{subequations}
Using these approximate relationships also for arbitrary nominal paths the control-measure vector is defined as 
\begin{align}
z=\begin{bmatrix}
\tilde z_1 & \tilde\theta_1 & \tilde z_2 & \tilde\theta_2 & \tilde\beta_2 & \tilde z_3 & \tilde\theta_3 & \tilde\beta_3\end{bmatrix}^T\triangleq h(\tilde x).
\end{align}
The function $h(\tilde x)$  is nonlinear and its Jacobian linearization around the origin yields $z= \frac{\partial h(0)}{\partial \tilde x}\tilde x \triangleq M\tilde x$, where 
\begin{align}
M = \begin{bmatrix}
1 & L_3+L_2+M_1 & L_2+M_1 & M_1 \\
0 & 1 & 1 & 1 \\
1 & L_3 & 0 & 0 \\
0 & 1 & 1 & 0 \\
0 & 0 & 0 & 1 \\
1 & 0 & 0 & 0 \\
0 & 1 & 0 & 0 \\
0 & 0 & 1 & 0
\end{bmatrix}.
\end{align}
By selecting the weight matrix for the quadratic stage-cost as $Q=M^T\bar QM$, where $\bar Q\succeq 0$ is a diagonal matrix, each diagonal element in $\bar Q$ corresponds to penalizing a specific control objective in $z$. 

After the weight matrix $Q$ has been selected, the weight matrix $P_N\succ 0$ for the terminal cost is selected as the solution to the discrete-time algebraic Riccati equation (DARE):
\begin{align}
\label{j2:eq:DARE}
F^T P_N F - P_N - F^T P_N GK + Q = 0,
\end{align}  
where $K = (1+G^TP_NG)^{-1}G^TP_N F$ is the LQ feedback gain, and the matrices $F=I+\Delta_s A$ and $G=\Delta_s B$ are the discrete-time system matrices~\eqref{j2:d_system_matrix} for the linearized path-following error model~\eqref{j2:d_lin_sys} around a straight nominal path. In this case, the matrices $A$ and $B$ are given by
\begin{align}
A = \bar v_{3r}\begin{bmatrix}
0 & 1 & 0              & 0 \\ 
0 & 0 & \frac{1}{L_3}  & 0 \\ 
0 & 0 & -\frac{1}{L_3} & \frac{1}{L_2} \\ 
0 & 0 & 0 & -\frac{1}{L_2}
\end{bmatrix},\hspace{5pt}  B&= \bar v_{3r}\begin{bmatrix}
0 \\ 0 \\ \frac{-M_1}{L2} \\ \frac{L_2+M_1}{L_2}
\end{bmatrix},
\end{align} 
where $v_{3r}\in\{-1,1\}$ specifies the nominal motion direction. Thus, since the system's stability properties vary depending on the direction of motion, different terminal costs are used during backward and forward motion tasks~\cite{LjungqvistJFR2019}. 

\subsection{Modeling of the constraint on the joint angles}  
The constraint on the joint angles~\eqref{j2:joint_angle_constaints} is intended to be selected such that the system avoids jackknifing, but also to restrict the joint angles to remain in the region where the used state-estimation solution is able to compute reliable and accurate state estimates of the trailer pose and the joint angles. 

\begin{figure}[t!]
	\centering
	\includegraphics[width=0.82\linewidth]{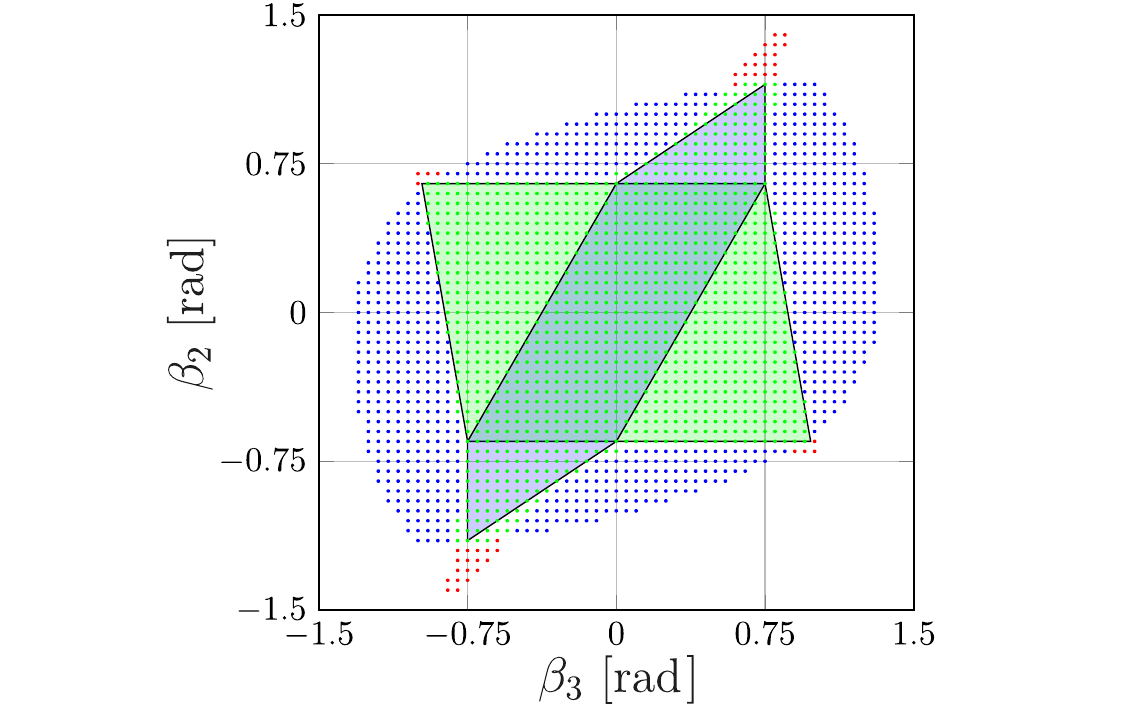}
	\caption{The simulated stability region for different joint angles (blue and green dots) around a straight nominal path and the region where the LIDAR-based estimation solution can compute estimates of the joint angles (red and green dots). The union of the blue and green polytopes is used to model the constraint for the joint angles~\eqref{j2:joint_angle_constaints}.}
	\label{j2:fig:LIDAR_FOV_stability_region}
\end{figure}

Since the tractor's curvature is limited by~\eqref{j2:curvature_constaint_MPC} and \eqref{j2:curvature_rate_constaint_MPC}, it is not possible to globally stabilize the path-following error system~\eqref{j2:eq:MPC_spatial_path_following_error_model} since for sufficiently large joint angles, jackknifing is impossible to prevent by only driving backwards~\cite{hybridcontrol2001}. 
This limit is possible to calculate analytically for the single-trailer case, but approximate methods are most often utilized when more than one trailer is present~\cite{hybridcontrol2001}. 
Given a straight nominal path, the vehicle parameters presented in Table~\ref{j2:tab:vehicle_parameters}, the MPC controller in~\eqref{j2:MPC_problem} with no joint-angle constraints and the design parameters in Table~\ref{j2:tab:design_parameters}. 
This system is simulated from different initial joint angles in backward motion with $v=-1$ m/s and is checked for convergence onto the straight nominal path. 
Using this technique, it is possible to numerically evaluate from which joint angles the closed-loop system is able to recover from and which will lead to jackknifing. The simulated stability region for the system is illustrated by the blue and green dots in Figure~\ref{j2:fig:LIDAR_FOV_stability_region}.

\begin{figure}[t!]
	\centering
	\includegraphics[width=0.8\linewidth]{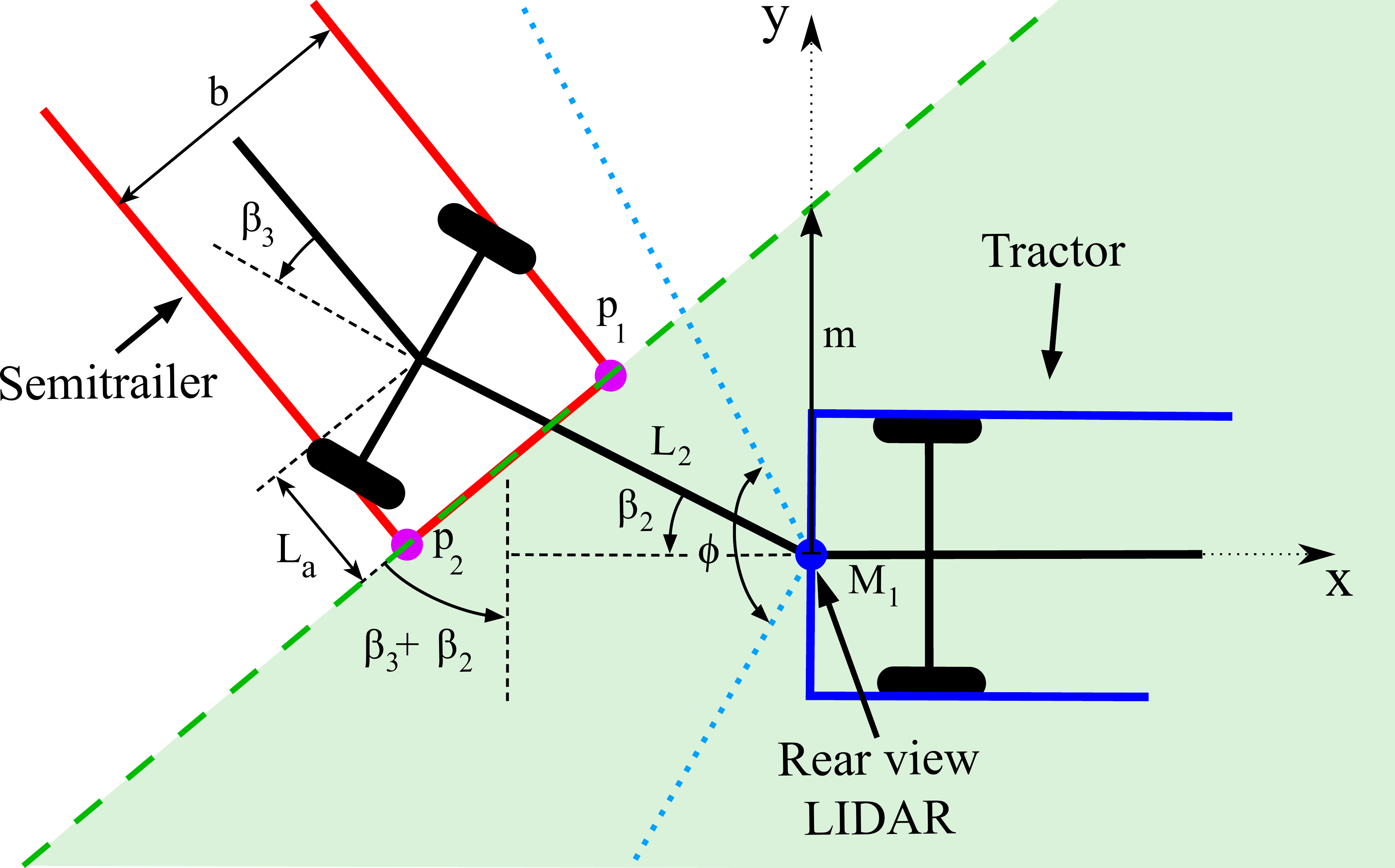}
	\caption{Illustration of the placement of the rear-view LIDAR sensor, its FOV (blue dotted lines) and relevant physical quantities. }
	\label{j2:fig:LIDAR_FOV_illustration}
\end{figure}

The sensor placement used in~\cite{LjungqvistJFR2019} is illustrated in Figure~\ref{j2:fig:LIDAR_FOV_illustration} together with a definition of relevant physical quantities that are also explained in Table~\ref{j2:tab:vehicle_parameters}. 
As long as the entire front of the semitrailer's body is visible from the LIDAR's point cloud, the LIDAR-based estimation technique presented in~\cite{LjungqvistJFR2019} computes accurate estimates of the joint angles as well as the remaining states of the semitrailer. 
Geometrically, this is equivalent to the following two conditions. First, the two corners, $p_1$ and $p_2$, of the semitrailer's front are inside the LIDAR's FOV. Second, the placement of the LIDAR is inside the half-space ahead of the semitrailer's front (green area in Figure~\ref{j2:fig:LIDAR_FOV_illustration}).

To model these conditions, a local coordinate system is introduced that is aligned with the tractor's orientation and has its origin at the tractor's off-axle hitch connection. Using basic trigonometry, the positions for the two corners of the semitrailer's front $p_1=(x_{p_1},y_{p_1})$ and $p_2=(x_{p_2},y_{p_2})$, can be expressed as a function of $\beta_2$ and $\beta_3$, and the vehicle parameters $L_2$, $L_a$ and $b$. 
Define angle of the sensor's horizontal scan field as $\phi$, then the following inequality constraints have to be satisfied for the two features\footnote{Note that is it straightforward to use the proposed design strategy for other types of advanced sensors (e.g., camera or RADAR) and if other features on the semitrailer's body are important to keep in the sensor's FOV.}, $p_1$ and $p_2$, to lie within the LIDAR's FOV:
\begin{subequations}
	\label{j2:eq:fov_cont_p1p2}
	\begin{align}
	&\cos{\left(\frac{\phi}{2}\right)} y_{p_i} + \sin{\left(\frac{\phi}{2}\right)} x_{p_i} \leq 0, \\
	&\cos{\left(\frac{\phi}{2}\right)} y_{p_i} - \sin{\left(\frac{\phi}{2}\right)} x_{p_i} \geq 0,\quad i=1,2.
	\end{align}
\end{subequations}
Furthermore, the border of the half-space that is ahead of the semitrailer's front (green dashed line in Figure~\ref{j2:fig:LIDAR_FOV_illustration}) can be expressed as
\begin{align} 
\label{j2:eq:line}
\sin{(\beta_2+\beta_3)}y-\cos{(\beta_2+\beta_3)}x = \bar m,   
\end{align}
where  $\bar m = L_2\cos{\beta_3} - L_a$. When $\beta_2+\beta_3\neq 0$, the line in~\eqref{j2:eq:line} intersects the y-axis in the tractor's local coordinate system at $m=\bar m/\sin{(\beta_2+\beta_3)}$. 
Hence, for the placement of the rear-view LIDAR sensor to strictly lie inside the half-space ahead of the semitrailer's front, one of the following three conditions has to be satisfied
\begin{subequations}
	\label{j2:eq:fov_line_intersection}
	\begin{align}
	\beta_2+\beta_3>0&\text{  and  } m\leq  \epsilon, \\
	\beta_2+\beta_3<0&\text{  and  } m \geq -\epsilon,\\
	\beta_2+\beta_3=0&\text{  and  } L_a -L_2\cos{\beta_2} \leq -\epsilon,
	\end{align}   
\end{subequations}
where $\epsilon>0$ is a constant that is used to enforce a certain robustness margin. 
Given the vehicle parameters in Table~\ref{j2:tab:vehicle_parameters} and the robustness margin $\epsilon=1$ m, the joint angles that satisfy the above mentioned constraints are illustrated by the green and red dots in Figure~\ref{j2:fig:LIDAR_FOV_stability_region}. 
As can be seen, the numerically computed stability region for the closed-loop system in backward motion (blue and green dots in Figure~\ref{j2:fig:LIDAR_FOV_stability_region}) is almost completely covering the region where accurate state estimates can be computed. 

The constraint on the joint angles~\eqref{j2:joint_angle_constaints} is now modeled as an inner approximation of the intersection of these two sets which is represented by the green dots in Figure~\ref{j2:fig:LIDAR_FOV_stability_region}. 
Based on the method used to represent this inner approximation, two alternative MPC controllers are proposed:
\begin{itemize}
	\item The union of two convex polytopes (The blue and green polytopes in Figure~\ref{j2:fig:LIDAR_FOV_stability_region}) is used as inner approximation. Two binary decision variables ($n=2$) are therefore introduced in the MIQP-MPC controller~\eqref{j2:MPC_problem}. 
	\item A single convex polytope (e.g. green polytope in Figure~\ref{j2:fig:LIDAR_FOV_stability_region}) is used to obtain an inner approximation. In this case, no binary decision variables are needed in the QP-MPC controller~\eqref{j2:MPC_problem}. 
\end{itemize}
The computational complexity of these two MPC controllers differs substantially because the MIQP-MPC controller requires that an MIQP problem is solved at each sampling instance, whereas the QP-MPC controller only needs to solve a QP problem. 
However, since the QP-MPC controller is restricted to use a single convex polytope to model the constraint on the joint angles, its performance may become noticeably suboptimal unless the polytope is carefully selected. 
To analyze the impact of this choice, different convex polytopes will be evaluated in the result section. 
Finally, since the vehicle model used in the MPC controllers only is an approximation, the hard constraints on the joint angles as well as the constraints on the semitrailer's lateral and orientation errors are replaced with soft constraints using standard techniques~\cite{mayne2000constrained,garcia1989model,falcone2007predictive}. 

\begin{table}[b!]
	\centering
	\caption{The vehicle parameters for the full-scale G2T with a car-like tractor.}
	\begin{tabular}{l l}
		\hline \noalign{\smallskip} Vehicle parameter  & Value   \\  \hline \noalign{\smallskip}	
		The tractor's wheelbase $L_1$            &   4.62 m  \\ 
		Maximum curvature $u_{\text{max}}$ & $0.18$ m$^{-1}$ \\
		Maximum curvature rate $\dot u_{\text{max}}$ & $0.13$ m$^{-1}$s$^{-1}$ \\
		Length of the off-hitch $M_1$      &   1.66 m  \\  
		Length of the dolly $L_2$          &   3.87 m  \\    
		Length of the semitrailer $L_3$    &   8.00 m  \\
		Length of the semitrailer's overhang $L_a$ &   1.73 m  \\
		Width of the semitrailer's front $b$   &   2.45 m  \\
		Angle of the horizontal scan field for the LIDAR $\phi$ & $140\times\pi/180$ rad \\
		\hline \noalign{\smallskip}
	\end{tabular}
	\label{j2:tab:vehicle_parameters}
\end{table}

\begin{table}[t!]
	\vspace{-15pt}
	\centering
	\caption{Design parameters for the MPC controllers.}
	\begin{tabular}{l l}
		\hline \noalign{\smallskip} Design parameter  & Value   \\  \hline \noalign{\smallskip}	
		Prediction horizon $N_{\text{QP}}$   & 40 \\ 
		Prediction horizon $N_{\text{MIQP}}$ & 30 \\ 
		Weight matrix $\bar Q$       &$1/35\times$diag$([0.5, 1, 0.5, 1, 4, 0.5, 1,4])$ \\  
		Sampling distance $\Delta_s$ &0.2 m\\
		Controller frequency  $f_s$      &10 Hz \\
		\hline \noalign{\smallskip}
	\end{tabular}
	\label{j2:tab:design_parameters}
\end{table}

\section{Results}
\label{j2:sec:results}
The performance of the proposed estimation-aware MPC approach is first evaluated in a simulation study and then in field experiments using the full-scale tractor-trailer vehicle shown in Figure~\ref{j2:fig:truck_scania}. 
Due to the extensive work-load required to interface an MIQP solver on the test vehicle, only the QP-MPC controller is experimentally validated in the field experiments. The implementation details for the simulations study and the field experiments are thoroughly explained in Section~\ref{j2:sec:simlation_setup} and Section~\ref{j2:sec:experimental_setup}, respectively.         

\subsection{Simulation setup}
\label{j2:sec:simlation_setup}
The performance of the MIQP-MPC controller and the QP-MPC controller are evaluated in a simulation study of backward and forward tracking a straight and a figure-eight nominal path. 
To enable rapid prototyping, the MPC controllers have been implemented in Matlab using YALMIP~\cite{lofberg2004yalmip}, where the state-of-the-art commercial MIQP solver Gurobi 8.1.1~\cite{gurobi} is used to solve~\eqref{j2:MPC_problem} at each sampling instance for both MPC controllers. All simulations have been performed on a standard laptop computer with an Intel Core i7-4600U@2.1GHz CPU. 
The design parameters for the MPC controllers are listed in Table~\ref{j2:tab:design_parameters} and the vehicle parameters in Table~\ref{j2:tab:vehicle_parameters}. 
The vehicle parameters are selected to coincide with the full-scale test vehicle shown in Figure~\ref{j2:fig:truck_scania}. 
As previously mentioned, since binary decision variables are not needed in the QP-MPC formulation it can be represented as a QP problem. 
This enables the use of a longer prediction horizon $N_{\text{QP}}=40$ for the QP-MPC controller compared to a prediction horizon $N_{\text{MIQP}}=30$ for the MIQP-MPC controller. Moreover, default settings are used in Gurobi with the exceptions that it is specified to perform warm starts and for MIQP-MPC to terminate once a solution with a relative suboptimality gap $\delta = 0.02$ is found.   

The performance of the proposed MPC controllers are benchmarked with the LQ controller presented in~\cite{LjungqvistJFR2019}. 
The weight matrix $Q$ is tuned equivalently for all controllers and the LQ feedback gains $K_{\text{fwd}}$ and $K_{\text{rev}}$ are computed offline by solving the DARE in~\eqref{j2:eq:DARE} in Matlab with $v_{3r}=\pm 1$. 
It would be natural to compare the proposed approach with other path-following approaches~\cite{hybridcontrol2001,rimmer2017implementation,Cascade-nSNT,bolzern1998path,virtualMorales2013}. However, since these approaches are not designed to follow nominal paths with full state and control-input information, a comparison is not included as the results would become misleading.    

\subsection{Simulation results}
The first set of simulations involve backward tracking of a straight nominal path aligned with the $x$-axis ($y_{3r}=0$), where the longitudinal velocity of the tractor is selected as \mbox{$v=-1$ m/s}. 
In this case, the nominal $\beta_{3r}=\beta_{2r}=u_{r}=0$ and therefore $\tilde\beta_3 =\beta_3$, $\tilde\beta_2=\beta_2$ and $\tilde u= u$. In the simulations, the initial path-following error state $\tilde x(0)$ is perturbed to compare how the different controllers handle disturbance rejection while satisfying the constraint on the joint angles. 
The initial state is $\tilde x(0)=[0 \hspace{5pt} 0 \hspace{5pt} \beta_3^i \hspace{5pt} \beta_2^i]^T$, where the initial joint-angle errors $\tilde\beta_3^i$ and $\tilde \beta_2^i$ are perturbed to various degrees. 
In this set of simulations, the green polytope shown in Figure~\ref{j2:fig:LIDAR_FOV_stability_region} is used as the joint-angle constraint in the QP-MPC controller.

\begin{figure}
	\centering
	\setlength\figureheight{0.27\linewidth}
	\setlength\figurewidth{0.8\linewidth} 
	\begin{tikzpicture}
	\node[anchor=south west] (myplot) at (0,0) {
		\input{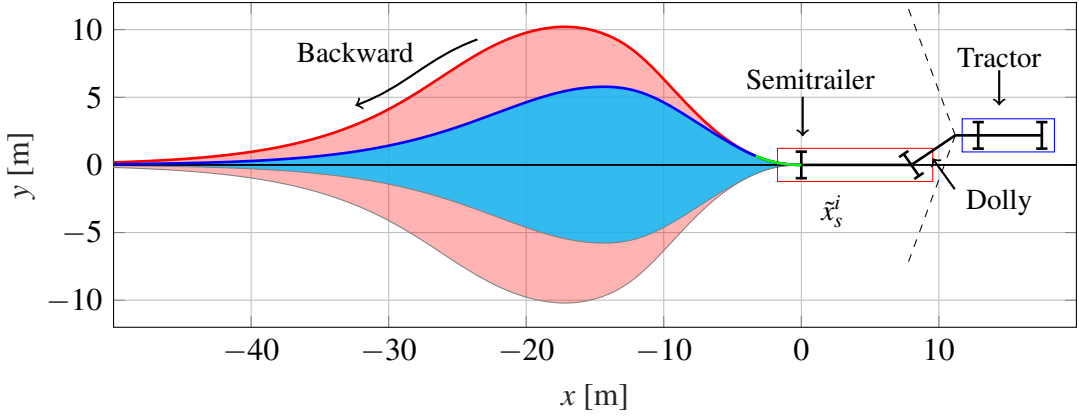}
	};
	\begin{scope}[x={(myplot.south east)}, y={(myplot.north west)}]
	\node[text=black] at (0.77,0.5) {\small $\tilde x^i_s$};
	\draw[<-, thick] (0.92,0.74) -- (0.92,0.82); 
	\node[text=black] at (0.92,0.85) {\small Tractor};
	\draw[<-, thick] (0.743,0.67) -- (0.743,0.76); 
	\node[text=black] at (0.75,0.8) {\small Semitrailer};
	\draw[->, thick] (0.88,0.55) -- (0.86,0.62); 
	\node[text=black] at (0.92,0.52) {\small Dolly};
	\draw[<-, thick] (0.34,0.74) to [out=20,in=200] (0.45,0.89);  
	\node at (0.34,0.86) {\small Backward};
	\end{scope}
	\end{tikzpicture}
	\vspace{-10pt}
	\caption{Envelopes of the trajectories for the axle of the semitrailer $(x_3(\cdotp),y_3(\cdotp))$ during backward tracking of a straight nominal path (black line), using MIQP-MPC (blue set) and QP-MPC (red set). For the high-lighted initial state $\tilde x^i_{s}$, the LQ controller leads to jackknifing (see Figure~\ref{j2:fig:traj_b23_straight}.). }
	\label{j2:fig:xy_path_straight}
\end{figure}

\begin{figure}[t!]
	\centering
	\setlength\figureheight{0.28\linewidth}
	\setlength\figurewidth{0.28\linewidth}
	\captionsetup[subfloat]{captionskip=-2pt}
	\vspace*{-10pt}
	\subfloat[][The trajectories for the joint angles from initial state $\tilde x^i_{s}$ in Figure~\ref{j2:fig:xy_path_straight}. Initial (desired) state denoted by a black (blue) star.]{
		\begin{tikzpicture}
		\node[anchor=south west] (myplot) at (0,0) {
			\input{path_b23_straight_MIQP.tex}
		};
		\end{tikzpicture}
		\label{j2:fig:traj_b23_straight}
	}
	~
	\setlength\figureheight{0.21\linewidth}
	\setlength\figurewidth{0.3\linewidth} 
	\subfloat[][The car-like tractor's curvatures from $\tilde x^i_{s}$.]
	{
		\begin{tikzpicture}
		\node[anchor=south west] (myplot) at (0,0) {
			\input{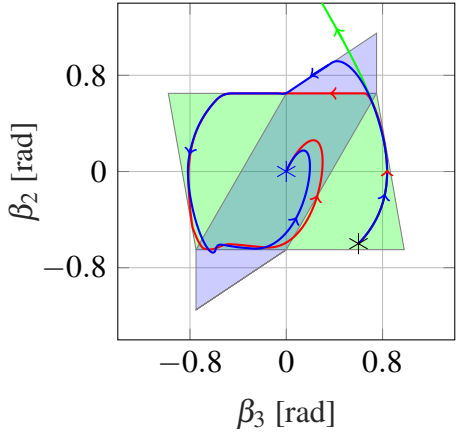}
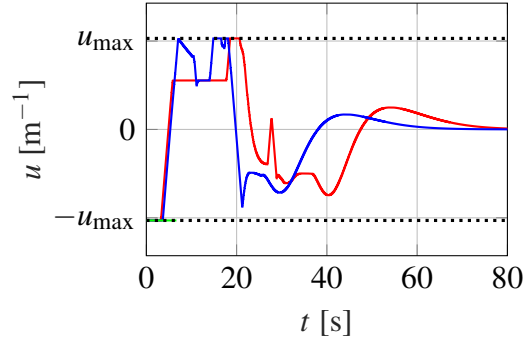
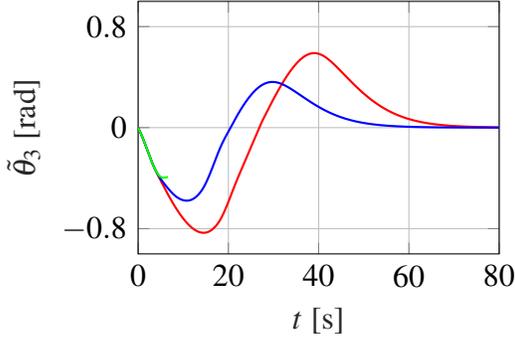
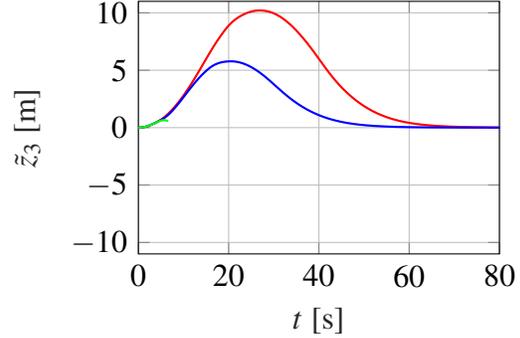
		};
		\end{tikzpicture}
		\label{j2:fig:straight_controls}
	}
	\quad
	\setlength\figureheight{0.21\linewidth}
	\setlength\figurewidth{0.3\linewidth}
	\subfloat[][The orientation errors from initial state $\tilde x^i_{s}$.]{
		\begin{tikzpicture}
		\node[anchor=south west] (myplot) at (0,0) {
			\input{path_th3_straight_MIQP.tex}
		};
		\end{tikzpicture}
		\label{j2:fig:traj_th3_sim_straight}
	}
	~ 
	\subfloat[][The lateral errors from initial state $\tilde x^i_{s}$.]{
		\begin{tikzpicture}
		\node[anchor=south west] (myplot) at (0,0) {
			\input{path_z3_straight_MIQP.tex}
		};
		\end{tikzpicture}
		\label{j2:fig:traj_z3_sim_straight}
	}
	\vspace{-7pt}
	\caption{Simulation results from path following of a straight nominal path ($y_{3r}=0$) in backward motion from perturbed initial joint-angle errors $\tilde\beta^i_3,\tilde\beta^i_2\in[-0.6,\hspace{1pt}0.6]$ rad, using MIQP-MPC (blue lines), QP-MPC (red lines) and LQ (green lines). High-lighted initial state $\tilde x^i_{s}$ in Figure~\ref{j2:fig:xy_path_straight} is from $(\tilde\beta_3^i,\tilde\beta_2^i)=(0.6,-0.6)$ rad.}
	\label{j2:fig:traj_sim_straight}
\end{figure}

The results from the simulations are presented in Figure~\ref{j2:fig:xy_path_straight}--\ref{j2:fig:traj_sim_straight}. 
The envelopes of the trajectories for the position of the semitrailer ($x_3(\cdotp),y_3(\cdotp)$) using MIQP-MPC and QP-MPC from initial joint-angle errors $\tilde\beta_3^i,\tilde\beta_3^i\in[-0.6,0.6]$ rad are presented in Figure~\ref{j2:fig:xy_path_straight}\footnote{The envelope for ($x_3(\cdotp),y_3(\cdotp)$) using the LQ controller is not presented in Figure~\ref{j2:fig:xy_path_straight} because the vehicle enters a jack-knife state for some $\tilde\beta_3^i,\tilde\beta_3^i\in[-0.6,0.6]$ rad, e.g., from initial state $\tilde x_s^i$ with $(\beta_3^i,\beta_2^i)=(0.6,-0.6)$ rad.}. 
The results show that the maximum transient in the lateral error of the semitrailer using the QP-MPC is 10.1 m in comparison to 5.7 m using the MIQP-MPC. 
The maximum convergence time to the straight nominal path is also longer using the QP-MPC compared to the MIQP-MPC. 
The reason why the MIQP-MPC outperforms the QP-MPC is because the MIQP-MPC exploits its larger joint-angle region (blue and green polytope in Figure~\ref{j2:fig:traj_b23_straight}) in comparison to the QP-MPC controller's more restrictive region (green polytope in Figure~\ref{j2:fig:traj_b23_straight}). 
This can be observed in Figure~\ref{j2:fig:traj_b23_straight} where the trajectories for the joint angles are provided from initial state $\tilde x_s^i$ with initial joint-angle errors $(\tilde\beta_3^i,\tilde\beta_2^i)=(0.6,-0.6)$ rad. 
The results show that MIQP-MPC improves the convergence time and transient response by steering the vehicle such that the joint angles purely enters into the blue polytope for a small portion of the maneuver. 
This additional feasible joint-angle region is also used by the MIQP-MPC controller to reduce the maximum transient error in the semitrailer's orientation error from initial state $\tilde x_s^i$ (see Figure~\ref{j2:fig:traj_th3_sim_straight}) which is 0.57 rad in comparison to  0.81 rad for the QP-MPC controller.
Moreover, the LQ controller is not able to stabilize the vehicle from $\tilde x_s^i$. 
This is because the LQ controller saturates the tractor's curvature (green line in Figure~\ref{j2:fig:straight_controls}) and since the system is open-loop unstable jackknifing occurs almost instantly (green line in Figure~\ref{j2:fig:traj_b23_straight}).


The second set of simulations involve path following of a figure-eight nominal path in $(x_{3r}(\cdotp),y_{3r}(\cdotp))$ in both forward \mbox{($v=1$ m/s)} and backward motion $(v=-1$ m/s). 
The figure-eight nominal path has been computed following the steps presented in~\cite{LjungqvistJFR2019}, which can be executed in both forward and backward motion because the system is symmetric~\cite{LjungqvistJFR2019}. 
Based on the results from the tracking of the straight nominal path, the convex polytope representing the joint-angle constraint in the QP-MPC controller is adjusted. 
The new constraint is illustrated by the green polytope in Figure~\ref{j2:fig:traj_b23_eight_rev}--\ref{j2:fig:traj_b23_eight_fwd}. That is, the convex polytope has been rotated to increase the allowed joint-angle region where the joint angles have equal sign. 
Also in this set of simulations, the initial state $\tilde x(0)$ is perturbed to compare the performance of the controllers. 
The initial state is chosen as \mbox{$\tilde x_{\text{r}}^i=[{-4}\text{ m} \hspace{5pt} 0 \hspace{5pt} 0.9\text{ rad} \hspace{5pt} 0.3\text{ rad}]^T$} for the backward simulations and $\tilde x_{\text{f}}^i=[3\text{ m} \hspace{5pt} 0.4\text{ rad} \hspace{5pt} {-1}\text{ rad} \hspace{5pt} {-0.7}\text{ rad}]^T$ for the forward simulations. 

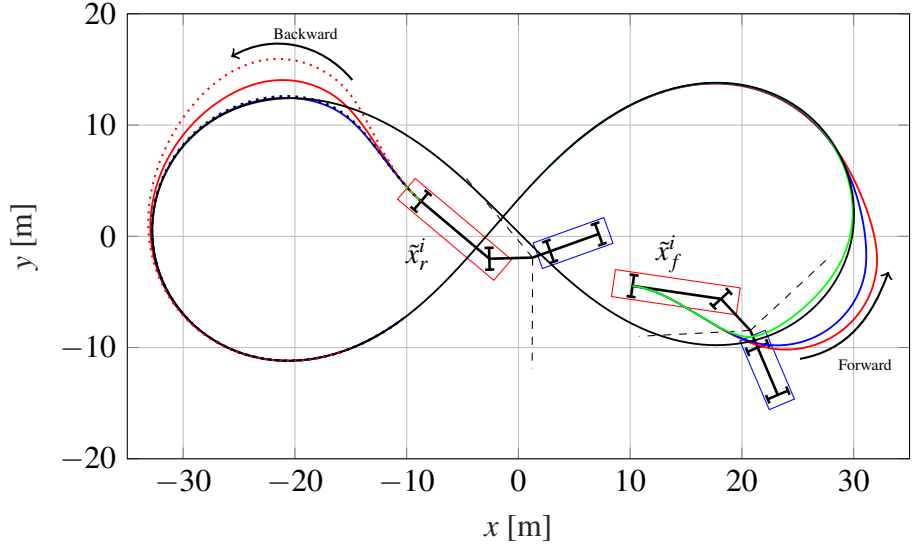
\begin{figure}
	\vspace{-15pt}
	\centering
	\setlength\figureheight{0.37\textwidth}
	\setlength\figurewidth{0.65\textwidth}
	\begin{tikzpicture}
	\node[anchor=south west] (myplot) at (0,0) {
		\input{path_xy_eight_MIQP.tex}
	};
	\begin{scope}[x={(myplot.south east)}, y={(myplot.north west)}]
	\draw[<-, thick] (0.26,0.87) to [out=30,in=130] (0.39,0.83);  
	\node at (0.34,0.91) {\tiny Backward};
	\draw[->, thick] (0.87,0.35) to [out=10,in=250] (0.965, 0.5);
	\node at (0.94,0.34) {\tiny Forward};
	\node[text=black] at (0.73,0.53) {\small $\tilde x^i_f$};
	\node[text=black] at (0.46,0.53) {\small $\tilde x^i_r$};
	\end{scope}
	\end{tikzpicture}
	\vspace{-10pt}
	\caption{Simulation results from path following a figure-eight nominal path in $(x_{3r}(\cdotp),y_{3r}(\cdotp))$ (black line) using MIQP-MPC (blue lines), QP-MPC (red lines) and LQ (green lines) from initial state $\tilde x_r^i$ in backward motion and from $\tilde x_f^i$ in forward motion. From $\tilde x^i_{r}$, the LQ controller leads to jack-knifing (see Figure~\ref{j2:fig:traj_b23_eight_rev}). The red dotted and black dotted lines are the trajectories from $\tilde x^i_{r}$ using QP-MPC tuned as in the first set of simulations and using MIQP-MPC with termination criterion $\delta=0.2$, respectively.}
	\label{j2:fig:traj_sim_eight_xy}
\end{figure}
\begin{figure}[t!]
	\centering
	\captionsetup[subfloat]{captionskip=-1pt}
	\setlength\figureheight{0.25\textwidth}
	\setlength\figurewidth{0.25\textwidth}
	\subfloat[][The joint-angle trajectories from initial state $\tilde x^i_{r}$ in Figure~\ref{j2:fig:traj_sim_eight_xy}. The black line is the nominal path in $(\beta_{3r}(\cdotp),\beta_{2r}(\cdotp))$.]{
		\input{path_b23_eight_rev_MIQP.tex}
		\label{j2:fig:traj_b23_eight_rev}
	}
	~
	\subfloat[][The joint-angle trajectories from initial state $\tilde x^i_{f}$ in Figure~\ref{j2:fig:traj_sim_eight_xy}. The black line is the nominal path in $(\beta_{3r}(\cdotp),\beta_{2r}(\cdotp))$.]{
		\input{path_b23_eight_fwd_MIQP.tex}
		\label{j2:fig:traj_b23_eight_fwd}
	}
	\quad
	\setlength\figureheight{0.18\columnwidth}
	\setlength\figurewidth{0.3\columnwidth}
	\subfloat[][The orientation errors from initial state $\tilde x^i_{r}$.]{
		\begin{tikzpicture}
		\node[anchor=south west] (myplot) at (0,0) {
			\input{path_th3_eight_rev_MIQP.tex}
		};
		\end{tikzpicture}
		\label{j2:fig:traj_th3_sim_eight}
	}
	~
	\subfloat[][The lateral errors from initial state $\tilde x^i_{r}$.]{
		\begin{tikzpicture}
		\node[anchor=south west] (myplot) at (0,0) {
			\input{path_z3_eight_rev_MIQP.tex}
		};
		\end{tikzpicture}
		\label{j2:fig:traj_z3_sim_eight}
	}
	\quad
	\setlength\figureheight{0.18\columnwidth}
	\setlength\figurewidth{0.5\columnwidth} 
	\subfloat[][The car-like tractor's curvatures from initial state $\tilde x^i_{r}$.]
	{
		\begin{tikzpicture}
		\node[anchor=south west] (myplot) at (0,0) {
			\input{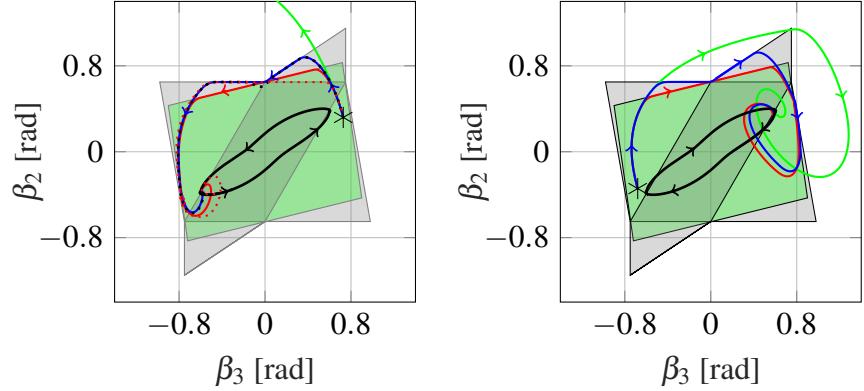}
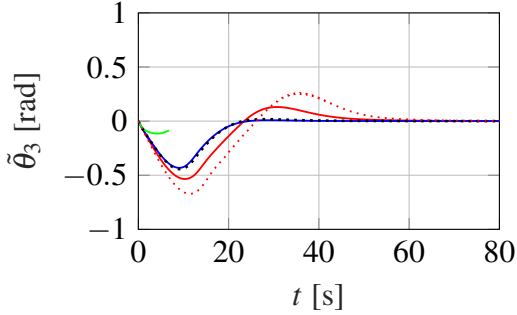
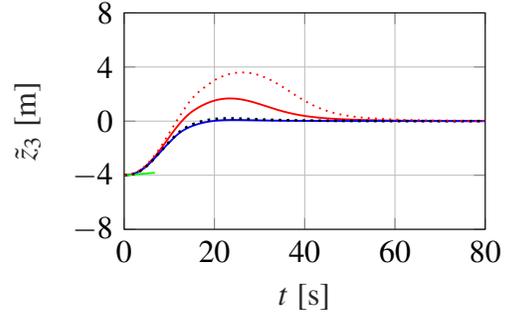
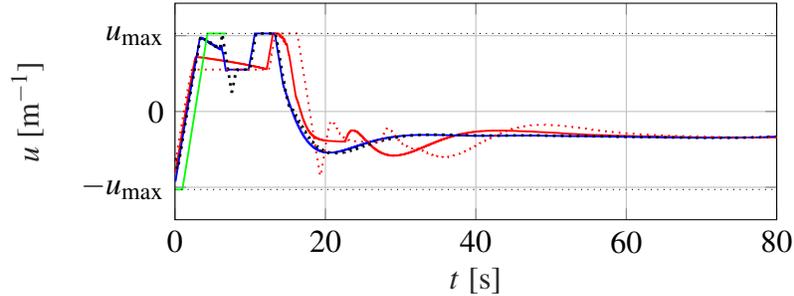
		};
		\end{tikzpicture}
		\label{j2:fig:traj_eight_controls}
	}
	\caption{Simulation results from path following a figure-eight nominal path in $(x_{3r}(\cdotp),y_{3r}(\cdotp))$ using MIQP-MPC (blue lines), QP-MPC (red lines) and LQ (green lines) from initial state $\tilde x_r^i$ in backward motion and from $\tilde x_f^i$ in forward motion. In Figure~\ref{j2:fig:traj_b23_eight_rev}--\ref{j2:fig:traj_b23_eight_fwd}, the initial states for $(\beta_{3}^i,\beta_{2}^i)$ are marked with black stars. The red dotted and black dotted lines are the trajectories from $\tilde x^i_{r}$ using QP-MPC tuned as in the first set of simulations and using MIQP-MPC with termination criterion $\delta=0.2$, respectively.}
	\label{j2:fig:traj_sim_eight}
	\vspace*{-10pt}
\end{figure}

The simulation results are presented in Figure~\ref{j2:fig:traj_sim_eight_xy}--\ref{j2:fig:traj_sim_eight}. 
For the backward simulations from $\tilde x_{\text{r}}^i$, the LQ controller is not able to stabilize the system and jackknifing occurs already after 5 seconds (see Figure~\ref{j2:fig:traj_b23_eight_rev}). 
As previously mentioned, this is because the LQ controller saturates the tractor's curvature as it is not aware of its constraints (see Figure~\ref{j2:fig:traj_eight_controls}). 
The QP-MPC controller converges to the figure-eight nominal path but with overshoots in both lateral error $\tilde z_3$ and orientation error $\tilde\theta_3$ of the semitrailer (red solid lines in Figure~\ref{j2:fig:traj_z3_sim_eight}-\ref{j2:fig:traj_th3_sim_eight}). 
However, by using the rotated polytope as constraint on the joint angles, the overshoots are decreased in comparison to the tuning used for QP-MPC in the first set of simulations (red-dotted lines in Figure~\ref{j2:fig:traj_sim_eight}). 
The MIQP-MPC controller is able to smoothly converge to the nominal path with no overshoot in $\tilde z_3$ nor in $\tilde\theta_3$ (blue lines in Figure~\ref{j2:fig:traj_z3_sim_eight}-\ref{j2:fig:traj_th3_sim_eight}). 
This performance enhancement is because the MIQP-MPC controller exploits its larger feasibility region in the joint angles. 
This can be seen in Figure~\ref{j2:fig:traj_b23_eight_rev} where the joint-angle trajectories $(\beta_2(\cdotp),\beta_3(\cdotp))$ are presented together with its nominal path $(\beta_{3r}(\cdotp),\beta_{2r}(\cdotp))$ (black solid line). 
For the forward cases, all controllers are able to converge to the figure-eight nominal path from initial state $\tilde x_f^i$, where the LQ controller converges fastest with no overshoots in the lateral error nor the heading error of the semitrailer (see Figure~\ref{j2:fig:traj_sim_eight_xy}). 
However, in the LQ case the trajectories for the joint angles (green line in Figure~\ref{j2:fig:traj_b23_eight_fwd}) are exiting the LIDAR-based estimation solution's sensing region which can cause unsatisfactory behaviors in practice. 
This is because the used estimation solution is not guaranteed to compute accurate estimates of the semitrailer's pose nor the joint angles.  

%
%

\begin{figure}[t!]
	\centering
	\setlength\figureheight{0.26\textwidth}
	\setlength\figurewidth{0.33\textwidth}
	\subfloat[][Theoretical (dashed lines) and actual suboptimality gaps (crosses) for each MPC interation of MIQP-MPC.]
	{
		\begin{tikzpicture}
		\node[anchor=south west] (myplot) at (0,0) {
			\input{Eight_suboptimalitygap.tex}
		};
		\end{tikzpicture}
		\label{j2:fig:eight_gap}
	}	
	~
	\hspace{-10pt}
	\subfloat[][Computation times for each MPC interation of MIQP-MPC. ]
	{
		\begin{tikzpicture}
		\node[anchor=south west] (myplot) at (0,0) {
			\input{Eight_computation_times_new_2.tex}
		};
		\end{tikzpicture}
		\label{j2:fig:eight_computation_time}
	}
	\caption{Computation time and relative suboptimality gap for the MIQP-MPC controller's first 200 iterations from $\tilde x_r^i$ in Figure~\ref{j2:fig:traj_sim_eight_xy} using the suboptimality gaps: $\delta=0.02$ (blue crosses), $\delta=0.1$ (green crosses) and $\delta=0.2$ (red crosses). In Figure~\ref{j2:fig:eight_computation_time}, the black crosses denote the computation time of the QP-MPC controller and the black-dashed line is the controller's sampling time.}
	\label{j2:fig:eight_computation_eval}
\end{figure}
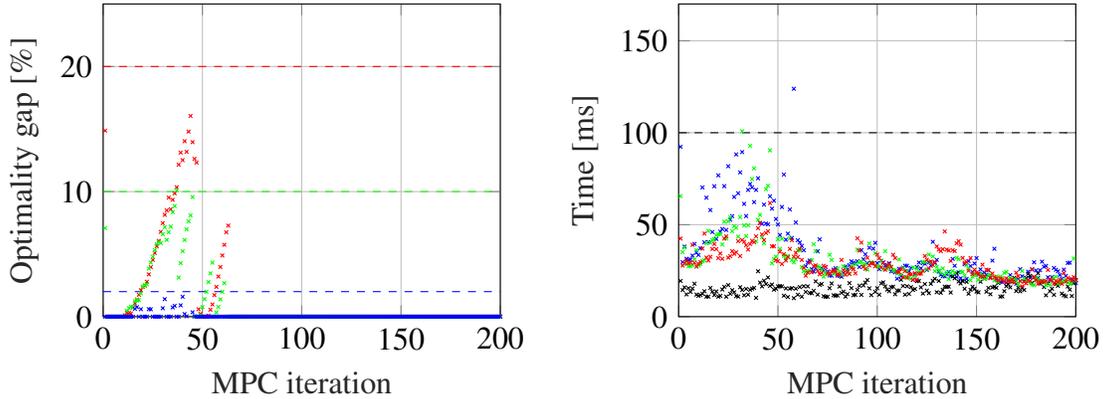

In the simulations, Gurobi's average computation time for the MIQP-MPC controller is \mbox{$40$ ms} (max: 181 ms) in comparison to 17 ms (max: 24 ms) for the QP-MPC controller. 
That is, the found worst-case computation time for MIQP-MPC is larger than the controller's specified sampling time of 100 ms. 
However, global solvers use a fair amount of its computation time only to prove that a found solution satisfies the specified relative suboptimality gap $\delta$. 
Inspired by the findings in~\cite{axehill2014parametric}, the backward tracking simulation from $\tilde x_{\text{r}}^i$ in Figure~\ref{j2:fig:traj_sim_eight_xy} using MIQP-MPC is repeated by relaxing the MIQP solver's termination criterion.
The computation times and the actual suboptimality gaps for the first 200 MPC iterations (20 seconds) using termination criteria for relative suboptimality gap as $\delta=0.02$ (blue), $\delta=0.1$ (green) and $\delta=0.2$ (red) are presented in Figure~\ref{j2:fig:eight_computation_eval}. 
From Figure~\ref{j2:fig:eight_computation_time} it can be concluded that by selecting $\delta =0.2$ the worst-case computation time is reduced to 58 ms (mean: 26 ms), on the expense that the computed control-input sequence is suboptimal at some MPC iterations (see Figure~\ref{j2:fig:eight_gap}). 
It is however guaranteed to never be worse than the specified suboptimality gap. 
Moreover, since the MPC controller operates in a receding horizon fashion, a suboptimal control-input sequence may only have a minor impact on the final trajectory taken by the vehicle. 
These arguments have strong similarities to well established methods for approximate nonlinear MPC, such as real-time iterations~\cite{AcadoRealTime}. 
As an example, the trajectories taken by the vehicle using MIQP-MPC with $\delta =0.2$ from $\tilde x_{\text{r}}^i$ is shown in Figure~\ref{j2:fig:traj_sim_eight} by black-dotted lines. 
The results show that the resulting trajectory is almost identical to the one obtained using MIQP-MPC with termination criterion $\delta =0.02$ (blue lines).

\subsection{Field experiments}
\label{j2:sec:experimental_setup}
The QP-MPC controller has been implemented and experimentally validated on a modified version of a Scania R580 6x4 tractor that is shown in Figure~\ref{j2:fig:truck_scania}. 
The tractor is equipped with additional computation power compared to its commercial version and a servo motor for automated control of its steering column. 
The tractor is also equipped with a localization system (real-time kinematic GPS [RTK-GPS]  and IMUs) and a rear-view LIDAR sensor that is aiming at the semitrailer body (see Figure~\ref{j2:fig:LIDAR_FOV_illustration}). 
Neither the semitrailer nor the dolly is equipped with any sensor. Instead the joint angles and the pose of the semitrailer are estimated using an extended Kalman filter (EKF) with the tractor's pose and virtual measurements of the trailer states as inputs. By running a random sample consensus algorithm~\cite{LjungqvistJFR2019}, the virtual measurements are computed by extracting features of the semitrailer body from the LIDAR's point cloud. 
The vehicle parameters are listed in Table~\ref{j2:tab:vehicle_parameters} and coincide with the ones used in the simulation study. 
For more details of the test platform, including the LIDAR-based estimation solution, the reader is referred to~\cite{LjungqvistJFR2019}.  

The QP-MPC controller is implemented in \texttt{C++} and qpOASES~\cite{qpOASES} is used as QP solver, which is an open-source active-set solver with warm-starts. Based on the results from the tracking of the figure-eight nominal path, the rotated polytope is used to represent the constraint on the joint angles.  
The controllers are operating in serial with the EKF, where the state estimate $\hat x_k$ is used to compute the path-following error $\tilde x_k$ at each sampling instance. 
The performance of the QP-MPC controller is evaluated in a set of real-world experiments of backward tracking a straight and a figure-eight nominal path, where it is also benchmarked with an LQ controller. 
The design parameters for the controllers are presented in Table~\ref{j2:tab:design_parameters}, which are equal to the ones used in the simulations with the exception that the controller frequency is increased to 20 Hz to meet system architectural constraints. 

\begin{figure}[t!]
	\centering
	\captionsetup[subfloat]{captionskip=-2pt}
	\setlength\figureheight{0.2\linewidth}
	\setlength\figurewidth{0.78\linewidth} 
	\vspace*{-7pt}
	\subfloat[][The trajectories taken by the axle of the semitrailer $(x_3(\cdotp),y_3(\cdotp))$ and the nominal path in $(x_{3r}(\cdotp),y_{3r}(\cdotp))$ (black solid line).]{
		\begin{tikzpicture}
		\node[anchor=south west] (myplot) at (0,0) {
			\input{path_xy_straight.tex}
		};
		\begin{scope}[x={(myplot.south east)}, y={(myplot.north west)}]
		\node[text=black] at (0.62,0.89) {\small $\tilde x^i_A$};
		\node[text=black] at (0.95,0.55) {\small $\tilde x^i_B$};
		\node[text=black] at (0.723,0.46) {\small $\tilde x^i_C$};
		\end{scope}
		\end{tikzpicture}
		\label{j2:fig:traj_straight_xy}
	}
	\quad
	\setlength\figureheight{0.25\columnwidth}
	\setlength\figurewidth{0.25\columnwidth}
	\subfloat[][The trajectories for the joint angles where the stars highlight the initial states.]{
		\begin{tikzpicture}
		\node[anchor=south west] (myplot) at (0,0) {
			\input{trajectory_beta23_straight.tex}
		};
		\end{tikzpicture}
		\label{j2:fig:traj_straight_beta23}
	}
	~
	\setlength\figureheight{0.2\columnwidth}
	\setlength\figurewidth{0.4\columnwidth}
	\subfloat[][The car-like tractor's curvatures.]{
		\begin{tikzpicture}
		\node[anchor=south west] (myplot) at (0,0) {
			\input{curvature_straight.tex}
		};
		\end{tikzpicture}
		\label{j2:fig:traj_straight_kappa}
	}
	\vspace*{-5pt}
	\quad
	\setlength\figureheight{0.18\columnwidth}
	\setlength\figurewidth{0.3\columnwidth}
	\subfloat[][The orientation errors of the semitrailer.]{
		\begin{tikzpicture}
		\node[anchor=south west] (myplot) at (0,0) {
			\input{heading_error_straight.tex}
		};
		\end{tikzpicture}
		\label{j2:fig:traj_straight_heading_error}
	}
	~
	\subfloat[][The lateral errors of the semitrailer.]{
		\begin{tikzpicture}
		\node[anchor=south west] (myplot) at (0,0) {
			\input{lateral_error_straight.tex}
		};
		\end{tikzpicture}
		\label{j2:fig:traj_straight_lateral_error}
	}
	\vspace*{-5pt}
	\caption{Experimental results from path following a straight nominal path ($y_{3r}=0$) from the three different initial states $\tilde x^i_A$ (red), $\tilde x^i_B$ (blue) and $\tilde x^i_C$ (green) using the QP-MPC controller (solid lines) and the LQ controller (dashed lines).} 
	\label{j2:fig:traj_straigt}
\end{figure}
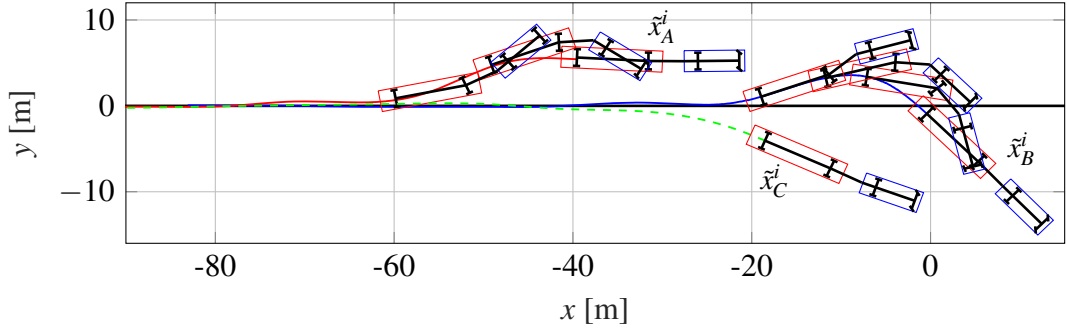
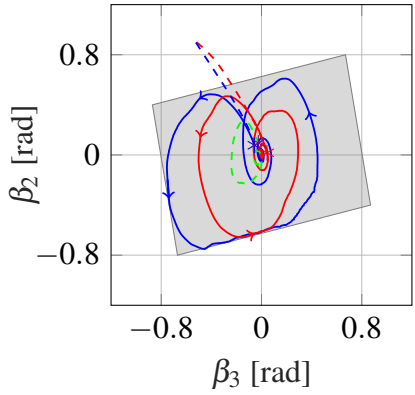
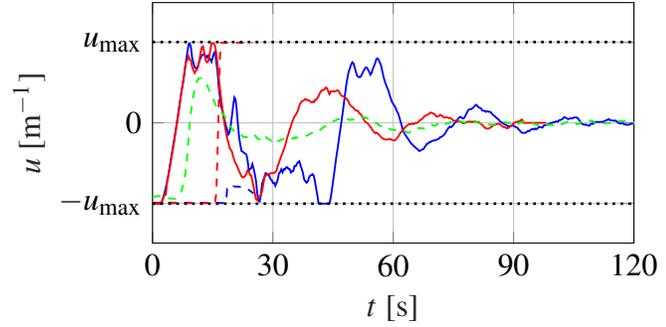
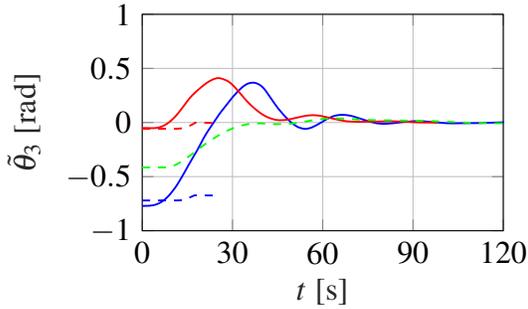
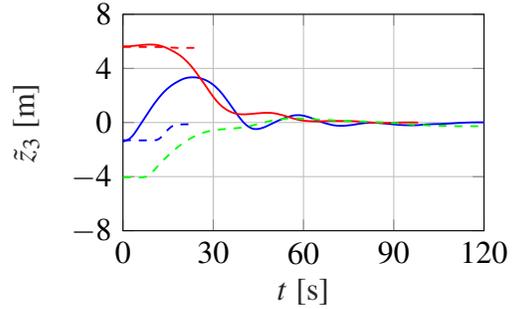

\subsection{Results from field experiments}
The QP-MPC controller is first evaluated by tracking a straight nominal path in backward motion ($v=-1$ m/s). 
As in the simulation study, the initial state $\tilde x(0)$ is perturbed (see Figure~\ref{j2:fig:traj_straight_xy}) to evaluated how the controllers handle disturbance rejection while satisfying the constraint on the joint angles. 
In Experiment A, the initial path-following error is $\tilde x^i_A=[{5.6}\text{ m} \hspace{5pt} 0 \hspace{5pt} 0 \hspace{5pt} 0]^T$, in Experiment B the initial error $\tilde x^i_B=[{-1.2}\text{ m} \hspace{5pt} {-0.77}\text{ rad} \hspace{5pt} 0 \hspace{5pt} 0]^T$, and in Experiment C the initial path-following error is $\tilde x^i_C=[{-4.1}\text{ m} \hspace{5pt} {-0.42}\text{ rad} \hspace{5pt} 0 \hspace{5pt} 0]^T$. 

The results from the experiments are presented in Figure~\ref{j2:fig:traj_straigt}. In Experiment A and B, the LQ controller is not able to stabilize the tractor-trailer vehicle due to the constraints on the tractor's curvature and jackknifing occurs almost instantly. 
This can be seen in Figure~\ref{j2:fig:traj_straight_beta23} where the joint-angle trajectories $(\beta_2(\cdotp),\beta_3(\cdotp))$ are plotted. 
In contrast to this behavior, the QP-MPC controller is able to make the system converge to the straight nominal path while satisfying the constraints on the joint angles. 
In experiment C, both controllers are able to stabilize the vehicle around the nominal path. The reason for this can be seen in Figure~\ref{j2:fig:traj_straight_kappa} where the tractor's curvature is plotted. 
The results show that in Experiment C, the feedback computed by the LQ controller (green dashed line) is not saturating the tractor's curvature, which it does in Experiment A and B. 

\begin{figure}[b!]
	\centering
	\setlength\figureheight{0.39\textwidth}
	\setlength\figurewidth{0.69\textwidth}
	\begin{tikzpicture}
	\node[anchor=south west] (myplot) at (0,0) {
		\input{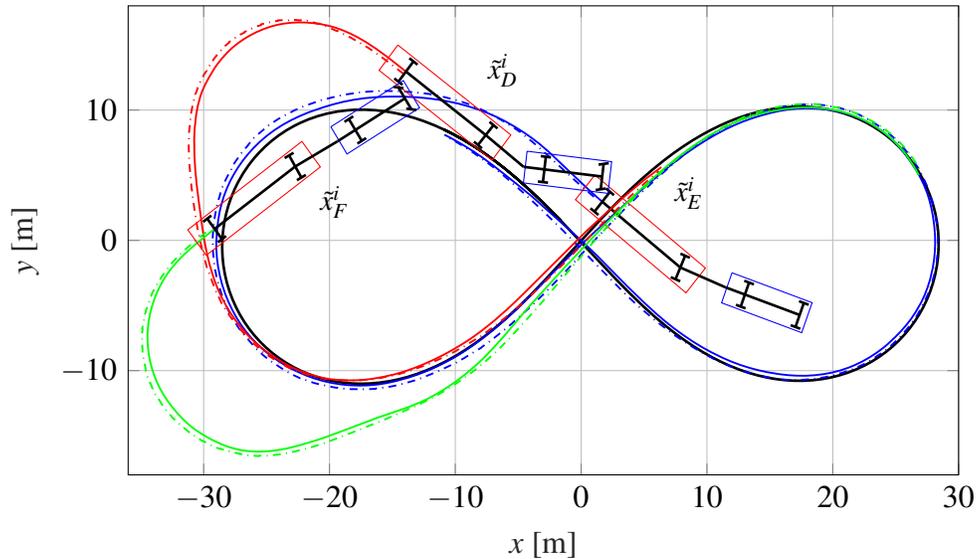}
	};
	\begin{scope}[x={(myplot.south east)}, y={(myplot.north west)}]
	\node[text=black] at (0.51 ,0.87) {\small $\tilde x^i_D$};
	\node[text=black] at (0.695 ,0.66) {\small $\tilde x^i_E$};
	\node[text=black] at (0.343,0.645) {\small $\tilde x^i_F$};
	\end{scope}
	\end{tikzpicture}
	\vspace{-10pt}
	\caption{Experimental results from backward tracking a figure-eight nominal path in $(x_{3r}(\cdotp),y_{3r}(\cdotp))$ (black solid line) from three different initial states $\tilde x^i_D$ (red), $\tilde x^i_E$ (blue) and $\tilde x^i_F$ (green) using the QP-MPC controller (solid lines). The dashed-dotted lines are measured path of the axle of the semitrailer by the external RTK-GPS.}
	\label{j2:fig:traj_eight_xy}
	\vspace{5pt}
\end{figure}

The second set of experiments involve backward tracking of a figure-eight shaped nominal path in $(x_{3r}(\cdotp),y_{3r}(\cdotp))$. Also in this set of experiments, the initial path-following error state $\tilde x(0)$ is perturbed (see Figure~\ref{j2:fig:traj_eight_xy}) to compare the performance of the controllers. In Experiment D, the initial path-following error is \mbox{$\tilde x^i_D=[{3.4}\text{ m} \hspace{5pt} {-0.46}\text{ rad}  \hspace{5pt} 0.46\text{ rad} \hspace{5pt} 0.73\text{ rad}]^T$}, in Experiment E the initial error is $\tilde x^i_E=[{3}\text{ m} \hspace{5pt} {0} \hspace{5pt} 0.26\text{ rad} \hspace{5pt} 0.27\text{ rad}]^T$, and in Experiment F the initial error is  $\tilde x^i_F=[{1.2}\text{ m} \hspace{5pt} {-0.8}\text{ rad} \hspace{5pt} 0.55\text{ rad} \hspace{5pt} 0.44\text{ rad}]^T$.

\begin{figure}[t!]
	\centering
	\captionsetup[subfloat]{captionskip=-2pt}
	\setlength\figureheight{0.27\columnwidth}
	\setlength\figurewidth{0.27\columnwidth}
	\subfloat[][The trajectories for the joint angles where the stars highlight the initial states.]{
		\begin{tikzpicture}
		\node[anchor=south west] (myplot) at (0,0) {
			\input{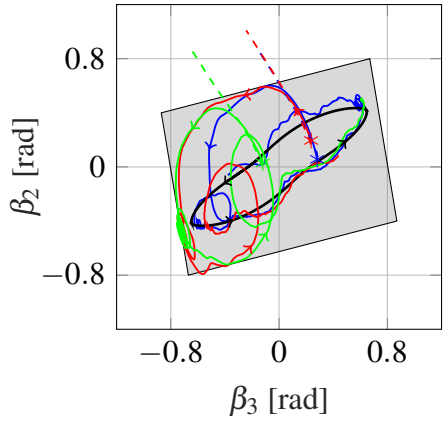}
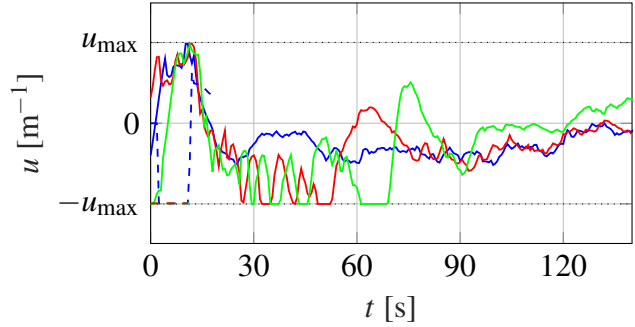
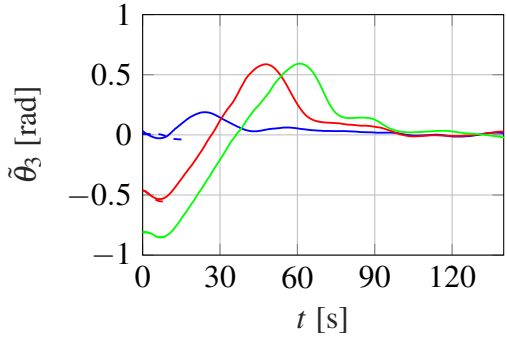
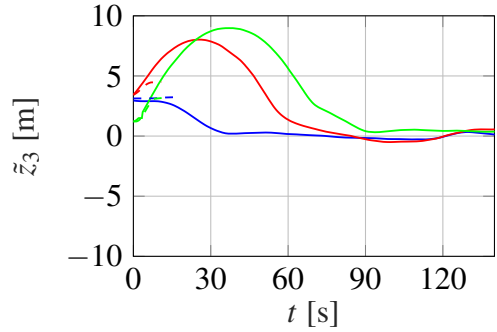
		};
		\end{tikzpicture}
		\label{j2:fig:traj_eight_beta23}
	}
	\hspace{-10pt}
	~
	\setlength\figureheight{0.2\columnwidth}
	\setlength\figurewidth{0.40\columnwidth}
	\subfloat[][The car-like tractor's curvatures.]{
		\begin{tikzpicture}
		\node[anchor=south west] (myplot) at (0,0) {
			\input{curvature_eight.tex}
		};
		\end{tikzpicture}
		\label{j2:fig:traj_eight_kappa}
	}	
	\quad
	\setlength\figureheight{0.2\columnwidth}
	\setlength\figurewidth{0.3\columnwidth}
	\subfloat[][The orientation errors of the semitrailer.]{
		\begin{tikzpicture}
		\node[anchor=south west] (myplot) at (0,0) {
			\input{heading_error_eight.tex}
		};
		\end{tikzpicture}
		\label{j2:fig:traj_eight_heading_error}
	}
	~
	\subfloat[][The lateral errors of the semitrailer.]{
		\begin{tikzpicture}
		\node[anchor=south west] (myplot) at (0,0) {
			\input{lateral_error_eight.tex}
		};
		\end{tikzpicture}
		\label{j2:fig:traj_eight_lateral_error}
	}
	\caption{Experimental results from backward tracking a figure-eight nominal path in $(x_{3r}(\cdotp),y_{3r}(\cdotp))$ from three different initial states $\tilde x^i_D$ (red), $\tilde x^i_E$ (blue) and $\tilde x^i_F$ (green) using the QP-MPC controller (solid lines) and the LQ controller (dashed lines). In Figure~\ref{j2:fig:traj_eight_beta23}, the black line is the nominal path in $(\beta_{3r}(\cdotp),\beta_{2r}(\cdotp))$.} 
	\label{j2:fig:traj_eight}
\end{figure}

The results from the experiments are presented in Figure~\ref{j2:fig:traj_eight_xy}--\ref{j2:fig:traj_eight}. In all three experiments, the LQ controller fails to stabilize the vehicle and jackknifing occurs almost instantly (see Figure~\ref{j2:fig:traj_eight_beta23}). 
In contrast to this behavior, the QP-MPC controller is able to stabilize the vehicle around the figure-eight nominal path in all three experiments while satisfying the constraints on the joint angles during the majority of the maneuvers. 
As can be seen in Figure~\ref{j2:fig:traj_eight_beta23}, the joint-angle trajectories for QP-MPC are violating their soft constraints from $\tilde x^i_D$ and $\tilde x^i_F$. Even though there are only minor violations, the effect on the linear penalty on constraint violation can be observed in Figure~\ref{j2:fig:traj_eight_kappa} where the tractor's curvature is plotted. 
Between 25--50 s, there are minor oscillations in the control input which is a result of model errors, small errors in the joint-angle estimates and the linear cost for violating the soft constraints on the joint angles.

To validate the performance of the used LIDAR-based estimation solution, an external RTK-GPS is mounted on the semitrailer's axle.
As can be seen in Figure~\ref{j2:fig:traj_eight_xy}, even though the maneuvers are advanced the observer is able to track the position of the semitrailer's axle. 
One important reason for this is that the QP-MPC controller is restricted to control the vehicle such that its joint angles remain in the region where high-accuracy state estimates can be computed by the used LIDAR-based estimation solution. 
As a final note, the average computation time for the QP-MPC controller during the experiments was only $5$ ms which is an order of magnitude less than the controller's sampling time of $50$ ms.


\section{Conclusions and future work}
\label{j2:sec:conclusions}
An estimation-aware model predictive path-following control approach for a G2T with a car-like tractor is proposed. 
The approach targets low-speed maneuvers and is designed to follow nominal paths in forward and backward motion that contain full state and control-input information. 
The path-following controller is tailored to operate in series with an estimation solution that uses an advanced sensor with a limited FOV to estimate the joint angles and the semitrailer's pose. 
To ensure that high-accuracy state estimates can be computed, the estimation solution's sensing region is modeled as constraints on the joint angles which are included in the MPC formulation. 
Two modeling approaches are proposed with different computation complexity and performance of the resulting MPC controller. In the first approach that is called MIQP-MPC, the constraints on the joint angles are modeled as a union of convex polytopes, making it necessary to incorporate binary decision variables in the MPC formulation. 
The second approach is called QP-MPC and avoids binary variables using a single convex polytope leading to a more restrictive approximation of the estimation solution's sensing region. 
The proposed MPC approaches are first evaluated in simulations and the QP-MPC is also evaluated in field experiments on a full-scale test vehicle. The computation complexity and performance of the proposed MPC approaches in terms of suppressing disturbances and recovering from non-trivial initial states is benchmarked, and shown to outperform, a previously proposed control strategy where the joint-angle constraints are neglected.

As future work, we would like to investigate if it would be feasible to deploy a suboptimal version of the MIQP-MPC controller on the full-scale test vehicle. 
This could be done either by relaxing the MIQP solver's termination criterion or by employing methods for approximate nonlinear MPC~\cite{AcadoRealTime}.   
We would also like to investigate if there exist other systems with advanced sensors with similar sensing limitations that need to be rigorously respected in the controller.

\section*{Acknowledgment}
The research leading to these results has been founded by Strategic vehicle research and innovation programme (FFI). The authors gratefully acknowledge the Royal Institute of Technology for providing the external RTK-GPS. The authors would also like to express their gratitude to Scania CV for providing hardware, software and technical support.

\section*{Appendix A}
The linearization of the path-following error model in~\eqref{j2:eq:MPC_spatial_path_following_error_model} around the origin $(\tilde x,\tilde u) = (0,0)$ is
\begin{align}
\frac{\text d \tilde x}{\text ds} = A(s)\tilde x +  B(s)\tilde u,
\end{align}
where the matrices $A(s)$ and $B(s)$ are given by
\begin{align}
A(s) = \frac{\partial\tilde f(s,0,0)}{\partial \tilde x} = \bar v_{3r}\begin{bmatrix}
0 & 1 & 0 & 0 \\
a_{21}(s) & 0 & a_{23}(s) & 0 \\
a_{31}(s) & 0 & a_{33}(s) & a_{34}(s) \\
a_{41}(s) & 0 & a_{43}(s) & a_{44}(s) 
\end{bmatrix}, 
\end{align}
and 
\begin{align}
B(s) = \frac{\partial\tilde f(s,0,0)}{\partial \tilde u} = \bar v_{3r}\begin{bmatrix}
0 \\ 0 \\ b_3(s) \\ b_4(s) 
\end{bmatrix}, 
\end{align}
where
\begin{align*}
a_{21}(s) &= -\frac{\tan^2{\beta_{3r}}}{L_3^2}, \\
a_{23}(s) &= \frac{1}{L_3\cos^2{\beta_{3r}}}, \\ 
a_{31}(s) &= \frac{\tan{\beta_{3r}}(u_rM_1\cos\beta_{2r}- \sin\beta_{2r})}{L_3L_2\cos\beta_{3r}(\cos\beta_{2r} + u_rM_1\sin{\beta_{2r}})} + \frac{\tan^2\beta_{3r}}{L^2_3}, \\
a_{33}(s) &= \frac{\sin\beta_{3r}(\sin\beta_{2r}-u_r M_1 \cos\beta_{2r})}{L_2 \cos\beta^2_{3r} (\cos\beta_{2r} + u_r M_1 \sin\beta_{2r})} -\frac{1}{\cos^2\beta_{3r}L_3}, 
\end{align*}
\begin{align*}
a_{34}(s) &= \frac{1+u_r^2M_1^2}{L_2\cos\beta_{3r}(\cos\beta_{2r}+u_rM_1\sin{\beta_{2r}})^2}, \\
a_{41}(s) &= 
-\frac{\tan\beta_{3r}}{L_2L_3}\left(
\frac{u_rL_2 - \sin\beta_{2r} + M_1\cos\beta_{2r}u_r}{\cos\beta_{3r}(\cos\beta_{2r} + M_1u_r \sin\beta_{2r})}\right), \\
a_{43}(s) &= \frac{\tan\beta_{3r}}{L_2}\left(\frac{u_rL_2 + u_r M_1 \cos\beta_{2r} - \sin\beta_{2r}}{\cos\beta_{3r}(\cos\beta_{2r} +u_rM_1\sin\beta_{2r})}\right), \\
a_{44}(s) &= \frac{1 + u_r^2 M_1^2 + u_r^2L_2M_1\cos\beta_{2r} -u_r L_2 \sin\beta_{2r}}{L_2\cos\beta_{3r}(\cos\beta_{2r}+u_r M_1\sin\beta_{2r})^2}, \\
b_3(s) &= - \frac{M_1}{L_2\cos\beta_{3r}(\cos\beta_{2r}+u_r M_1 \sin\beta_{2r})^2},\\
b_4(s) &= \frac{M_1 + L_2\cos\beta_{2r}}{\cos\beta_{3r}L_2(\cos\beta_{2r}+u_r M_1 \sin\beta_{2r})^2}.
\end{align*}

\bibliographystyle{abbrv}
\bibliography{root}

\end{document}

%% file: path_b23_straight_MIQP.tex
%
%
%
\begin{tikzpicture}

\begin{axis}[%
width=\figurewidth,
height=\figureheight,
at={(0\figurewidth,0\figureheight)},
scale only axis,
xmin=-1.4,
xmax=1.4,
xtick={-0.8,0,0.8},
xlabel={$\beta_3$ [rad]},
ymin=-1.4,
ymax=1.4,
ytick={-0.8,0,0.8},
xlabel style={font=\color{white!15!black},at={(axis description cs:0.5,-0.14)},anchor=north},
ylabel style={font=\color{white!15!black},at={(axis description cs:-0.2,.5)},anchor=south},
ylabel={$\beta_2$ [rad]},
axis background/.style={fill=white},
xmajorgrids,
ymajorgrids
]

\addplot[area legend, draw=gray, fill=green, fill opacity=0.3, forget plot]
table[row sep=crcr] {%
x	y\\
-0.98	0.65\\
-0.75	-0.65\\
0.98	-0.65\\
0.75	0.65\\
-0.98	0.65\\
}--cycle;

\addplot[area legend, draw=gray, fill=blue, fill opacity=0.2, forget plot]
table[row sep=crcr] {%
x	y\\
-0.75	-1.15\\
0	-0.65\\
0.75	0.65\\
0.75	1.15\\
0	0.65\\
-0.75	-0.65\\
-0.75	-1.15\\
0	-0.65\\
}--cycle;
\addplot [color=green, line width=0.8pt, forget plot]
  table[row sep=crcr]{%
0.6	-0.6\\
0.627875153427256	-0.569755481734023\\
0.654767770145175	-0.537846617384246\\
0.680538251160882	-0.504169676627281\\
0.705034971322429	-0.468615281592838\\
0.728093829969311	-0.431068557436909\\
0.749537927399505	-0.391409446290993\\
0.769177417808319	-0.349513230155611\\
0.786809602861975	-0.305251316533252\\
0.802219346327474	-0.258492348977963\\
0.815179909285144	-0.209103712547005\\
0.820668045059065	-0.183381957672095\\
0.825454327159076	-0.156953510246163\\
0.829507875290067	-0.129802409577256\\
0.832797473288153	-0.10191308932812\\
0.83529168321192	-0.0732704918030138\\
0.836958977273648	-0.0438601934988114\\
0.837767889217301	-0.0136685422114793\\
0.837687186746891	0.0173171941584951\\
0.83668606656211	0.049108667141653\\
0.834734373456387	0.0817162758386905\\
0.831802844762684	0.115148984190108\\
0.827863381180369	0.149414128210846\\
0.822889344667355	0.18451721437508\\
0.816855883620366	0.220461710728442\\
0.809740284978269	0.257248832734305\\
0.801522352156512	0.294877326320447\\
0.792184806845874	0.333343251072386\\
0.781713711682435	0.372639767006264\\
0.77009890962187	0.412756928830062\\
0.757334474543865	0.453681492046311\\
0.743419166197954	0.495396735637836\\
0.728356881121464	0.537882306383657\\
0.694835247884465	0.625064114725105\\
0.656919359470601	0.714987397645707\\
0.614865400071874	0.807350063110566\\
0.569042691253421	0.901789244633707\\
0.519934497044334	0.997886431178233\\
0.468131471474043	1.09517602794426\\
0.414317005246837	1.19315708994105\\
0.331484604845746	1.34028137725545\\
0.276020108428582	1.43770026579747\\
};
\addplot [color=red, line width=0.8pt, forget plot]
  table[row sep=crcr]{%
0.6	-0.6\\
0.627875153427296	-0.569755481734186\\
0.654767770145312	-0.537846617384785\\
0.680538251161092	-0.504169676628026\\
0.7050349713227	-0.468615281593686\\
0.728093829969672	-0.431068557437952\\
0.74953792739996	-0.391409446292201\\
0.769177417809076	-0.34951323015769\\
0.78680960286291	-0.305251316535614\\
0.802219346328579	-0.258492348980507\\
0.81517990928658	-0.209103712550253\\
0.820668045060697	-0.183381957675738\\
0.825454327160834	-0.156953510249921\\
0.829507875292122	-0.129802409581697\\
0.832797473290402	-0.101913089332829\\
0.835291683214351	-0.0732704918079269\\
0.836958977276278	-0.0438601935039666\\
0.837767889220108	-0.0136685422167753\\
0.837687186749929	0.0173171941529009\\
0.836686066565626	0.0491086671350111\\
0.83473437346021	0.0817162758316241\\
0.831802844766781	0.115148984182756\\
0.823147170924153	0.183725952929981\\
0.812078112810855	0.250417538144296\\
0.799577541225447	0.312719148186548\\
0.785947885144829	0.370392238755343\\
0.77148124860912	0.423215475808756\\
0.756460245497407	0.47098135841383\\
0.74116012457791	0.513492028925957\\
0.725852051808105	0.550554416417023\\
0.710807415129219	0.581974836204645\\
0.696303026527562	0.607553140492141\\
0.68262710933026	0.627076488194199\\
0.670028933413376	0.640487362390194\\
0.664211678100976	0.644886920028457\\
0.658725633103363	0.647844293419815\\
0.653622919276115	0.649290434695161\\
0.639225850519508	0.649698605193875\\
0.382008552094527	0.649731011087024\\
-0.465149822647672	0.649264046097342\\
-0.481662524929206	0.647780718326212\\
-0.497954375264681	0.644640490513647\\
-0.513983487164657	0.639784529827692\\
-0.529678681452172	0.633037832926375\\
-0.544995948060614	0.624375157440367\\
-0.559886405410731	0.613764868078032\\
-0.574278877359907	0.601093929072067\\
-0.601727903404736	0.571188361103068\\
-0.628139855669366	0.539358614082116\\
-0.653432235185451	0.505765706960984\\
-0.677454619345618	0.47030052051002\\
-0.700044111783412	0.432848428288206\\
-0.721025002537605	0.393289602105893\\
-0.740208604157038	0.351499524146174\\
-0.757393328777816	0.307349759159213\\
-0.772365087743854	0.260709048514811\\
-0.784898114771322	0.211444795776427\\
-0.790176682328147	0.185787508682821\\
-0.798884860413258	0.13329225400407\\
-0.805804378981555	0.0816846579917391\\
-0.810517168824188	0.0295815959393944\\
-0.812789820991679	-0.0235218692549074\\
-0.812918426082493	-0.0764679748419776\\
-0.811566352902552	-0.127222253365743\\
-0.807951550955122	-0.196371088424985\\
-0.80240489887762	-0.271987462632347\\
-0.794186915780133	-0.354379476556883\\
-0.786774662759985	-0.412765042843375\\
-0.77965673478227	-0.4586117537403\\
-0.777586903846824	-0.467840500539602\\
-0.768373443244993	-0.492992335043094\\
-0.752316253232099	-0.533009930356742\\
-0.737913222107194	-0.56578151430584\\
-0.723709292153026	-0.593580722872681\\
-0.710242943408559	-0.615403013330488\\
-0.697816346707244	-0.631023187437128\\
-0.69209334722024	-0.636430877631221\\
-0.686754215701178	-0.640192573519111\\
-0.681844282007276	-0.642272374888109\\
-0.668247999825817	-0.643851569338389\\
-0.648938263434609	-0.647163383597908\\
-0.634352987684922	-0.648107185326908\\
-0.619817402238739	-0.647561398373864\\
-0.605580523468202	-0.645031234614119\\
-0.591706063753451	-0.640699785571783\\
-0.564402563703762	-0.629645773702272\\
-0.523074063369474	-0.612622259175804\\
-0.508551088750613	-0.608594268386828\\
-0.493318610451633	-0.60603517406534\\
-0.477293948510006	-0.604860632905883\\
-0.45463867798196	-0.605061834718007\\
-0.418026745218426	-0.607501101015513\\
-0.342637243783601	-0.61479839504315\\
-0.173916095185194	-0.63183354524871\\
-0.145301692705584	-0.632202696746176\\
-0.126048041345683	-0.631151670588101\\
-0.106735709626752	-0.628898501210404\\
-0.0874303517351469	-0.625346541519072\\
-0.0681934958056822	-0.620431067876999\\
-0.0490833292243011	-0.614113388861252\\
-0.0301551775768194	-0.606376529198138\\
-0.0114618358371014	-0.5972217528744\\
0.00694621612048707	-0.586665700947348\\
0.0250206885289029	-0.574738026149054\\
0.0427164006221775	-0.561482777759687\\
0.0599908072164147	-0.546954292680227\\
0.0768029444399582	-0.531212049795158\\
0.0931146311663166	-0.514323511320696\\
0.10889076970374	-0.496363411798465\\
0.131482333137043	-0.467593156966066\\
0.152706520930659	-0.436892377351485\\
0.172491307958974	-0.404579396545472\\
0.190785435889695	-0.370987175478063\\
0.207556191564793	-0.336447803284943\\
0.222793026596385	-0.301300131813505\\
0.236506480370578	-0.265880187184713\\
0.248724290442313	-0.230506992295043\\
0.259489184319512	-0.195477489403325\\
0.271676073052317	-0.149760073414693\\
0.281528851528568	-0.105670951873523\\
0.289219672408688	-0.0636584847536968\\
0.294928347457622	-0.0240631405320817\\
0.298834801948952	0.012875764810272\\
0.301113887627253	0.0470086172162564\\
0.301927003303656	0.0856231916362949\\
0.300757063906706	0.119741358484032\\
0.297879338334033	0.149468094449412\\
0.293544333952389	0.174985188544387\\
0.28797759287688	0.196529547588403\\
0.282773976298936	0.211086226067321\\
0.277005487301691	0.223421445639316\\
0.27075753475132	0.23368657374681\\
0.264106295262699	0.24203796776272\\
0.257122266658027	0.248624596270121\\
0.249870431453889	0.253589125828416\\
0.242409615301015	0.257070193786816\\
0.232870146922186	0.259538091286858\\
0.223179080047839	0.260143218156834\\
0.213417081863506	0.25911639527655\\
0.201707173108857	0.256029749687347\\
0.190099112667557	0.251236553389645\\
0.176799202858646	0.243888187088114\\
0.162053199675965	0.233685660577953\\
0.146192357030896	0.220545272886891\\
0.12962301529235	0.204620311354107\\
0.112805585246788	0.186295343504063\\
0.0949268595345152	0.164492575410268\\
0.075897975419131	0.138628944873034\\
0.0573534302279722	0.110618501239422\\
0.0394971769590377	0.0806843197992487\\
0.0230276123428637	0.0499235400877708\\
0.0084328699257229	0.0192543867146567\\
0.00141429596788945	0.00310221723794202\\
};
\addplot [color=blue, line width=0.8pt, forget plot]
  table[row sep=crcr]{%
0.6	-0.6\\
0.627875153427256	-0.569755481734023\\
0.654767770145175	-0.537846617384246\\
0.680538251160881	-0.504169676627281\\
0.705034971322429	-0.468615281592838\\
0.728093829969311	-0.431068557436909\\
0.749537927399505	-0.391409446290993\\
0.769177417808319	-0.349513230155611\\
0.786809602861975	-0.305251316533252\\
0.802219346327474	-0.258492348977963\\
0.815179909285144	-0.209103712547005\\
0.820668045059065	-0.183381957672095\\
0.825454327159076	-0.156953510246163\\
0.829507875290067	-0.129802409577256\\
0.832797473288153	-0.10191308932812\\
0.83529168321192	-0.0732704918030139\\
0.836958977273648	-0.0438601934988115\\
0.837767889217301	-0.0136685422114792\\
0.837687186746891	0.0173171941584951\\
0.83668606656211	0.0491086671416531\\
0.834734373456387	0.0817162758386906\\
0.831802844762684	0.115148984190108\\
0.827863381180369	0.149414128210846\\
0.822889344667355	0.18451721437508\\
0.816855883620366	0.220461710728442\\
0.809740284978269	0.257248832734305\\
0.801522352156512	0.294877326320447\\
0.782636114336423	0.369925633765866\\
0.761664098302972	0.441413031064266\\
0.739105018511507	0.508481188932556\\
0.715239010277479	0.570926889849422\\
0.690325004725212	0.628570554292035\\
0.664601672777233	0.681252391879329\\
0.638289640673077	0.72882785828207\\
0.611594783783781	0.771162597498725\\
0.5847110625622	0.808131630160924\\
0.571242150762451	0.824614933299632\\
0.557769054125233	0.839797930668109\\
0.544312686857161	0.853668422966931\\
0.530893929507125	0.866214484955584\\
0.517531148869841	0.877433173487036\\
0.504251487480955	0.88729188598393\\
0.491068891614378	0.895806085651076\\
0.478020252019422	0.902911303507963\\
0.465112880527283	0.908652467081045\\
0.452397844971658	0.912918558810373\\
0.439872901638355	0.915798200975148\\
0.427612616175462	0.917100696363715\\
0.415611682897271	0.916941213881956\\
0.403947393401018	0.915129862211628\\
0.392637673923982	0.91171292585839\\
0.381750093298477	0.906556334805196\\
0.350549838144418	0.885850408481532\\
0.0174729547714774	0.661668766107721\\
0.0072993649877291	0.655867054183818\\
-0.00314159134418257	0.651023332645778\\
-0.014055434776443	0.64782708499777\\
-0.0254871546591439	0.646320792482569\\
-0.0374811420523682	0.646546358928078\\
-0.0624379358570214	0.648950532825439\\
-0.100286505501902	0.64963116319765\\
-0.339715631662871	0.649731009390846\\
-0.468500806510388	0.649041406516666\\
-0.484975816538587	0.647250020452347\\
-0.50125622440478	0.643942068678842\\
-0.51724511274664	0.638802789044776\\
-0.532930459585124	0.631933822327525\\
-0.548229481253589	0.623122618170043\\
-0.563100962885088	0.612375291257195\\
-0.577470155753403	0.599562895856135\\
-0.604895436847642	0.569519079038141\\
-0.63128765455561	0.537597159138999\\
-0.656552834849544	0.503906351727712\\
-0.680539852146457	0.468337236354057\\
-0.703085086160699	0.43077489731347\\
-0.724012090708681	0.391099239004064\\
-0.743131441425108	0.34918551060569\\
-0.760240828129349	0.304905093940129\\
-0.775125474250045	0.258126616752312\\
-0.794531073658637	0.188775041522279\\
-0.80476349193992	0.140577973987306\\
-0.8125869933908	0.0905435475649325\\
-0.815355779043017	0.0641799911953741\\
-0.817004090191659	0.0358927353006498\\
-0.817420590109719	0.00553738057590158\\
-0.816926254281059	-0.0256607315256783\\
-0.815503957824068	-0.0576674645144627\\
-0.813123859615429	-0.0904928554994765\\
-0.809757098876918	-0.124145456295704\\
-0.800270430800829	-0.19298733714017\\
-0.788585431671198	-0.259297244561588\\
-0.775313767062897	-0.321726046830142\\
-0.76068876091543	-0.380156629322927\\
-0.744948007862385	-0.434424317169092\\
-0.728564620222472	-0.48368513019744\\
-0.711807648214432	-0.527745062397321\\
-0.694940173452539	-0.566416426796144\\
-0.678223445756343	-0.599512539533201\\
-0.661921944551191	-0.626841950703371\\
-0.646309258436525	-0.648202291512083\\
-0.638850820824978	-0.656576243517573\\
-0.631674693062144	-0.663373782541563\\
-0.624820684141261	-0.668563600054935\\
-0.618330522422629	-0.672112434011602\\
-0.612248109702815	-0.673984842110383\\
-0.606619779713528	-0.674142972359752\\
-0.601494558645973	-0.672546332057218\\
-0.596924424485669	-0.669151556747404\\
-0.592964560890158	-0.663912181311203\\
-0.589673599998479	-0.656778416089735\\
-0.584362241183082	-0.639710000098819\\
-0.578392093491495	-0.626412883238673\\
-0.571450174992136	-0.61727749694417\\
-0.563547046838104	-0.611705048420884\\
-0.554769495248088	-0.608953200211795\\
-0.545264310494595	-0.608188493075871\\
-0.535157775051737	-0.608760913593942\\
-0.507893596709872	-0.612981076834213\\
-0.420314797262789	-0.62728388537442\\
-0.370528440193305	-0.633086687084109\\
-0.26645353796954	-0.641537099117319\\
-0.240702980216235	-0.641136828579566\\
-0.214760437914928	-0.638603973741366\\
-0.197448979645402	-0.635593269827365\\
-0.18018198827058	-0.631467143809832\\
-0.163005343482792	-0.626199453986139\\
-0.145962279718904	-0.619780781103684\\
-0.129093688289695	-0.612215549095111\\
-0.112438189184641	-0.603520174417705\\
-0.0960315678417043	-0.593723375730217\\
-0.0719635537487974	-0.577038658635291\\
-0.0486422924840736	-0.558093603082618\\
-0.0261669823178103	-0.537056857918151\\
-0.00462606303811497	-0.514119964880158\\
0.0159035346946259	-0.489492711367362\\
0.0353571149491855	-0.463398946486956\\
0.0536827946600645	-0.436072565198267\\
0.0708418379162994	-0.40775350086406\\
0.0868084548681972	-0.378682929824803\\
0.10156790851544	-0.34909469193771\\
0.119370354638099	-0.30924492647519\\
0.135058179074103	-0.269426305298876\\
0.148692536026827	-0.230126017642269\\
0.160358816459495	-0.19177503036763\\
0.170160506021318	-0.154741744272257\\
0.179967689813763	-0.110758823631429\\
0.1872837606702	-0.069754882045398\\
0.192355024468003	-0.0320348221175147\\
0.195424292590576	0.00222526870632567\\
0.196724983060157	0.032957589177713\\
0.196477587636182	0.0601816871406076\\
0.194887819720469	0.0839849124826431\\
0.192145869687811	0.104505222754335\\
0.188426344272599	0.121916681168784\\
0.183888343359674	0.136418388509856\\
0.178674169710596	0.148229513030972\\
0.1729147937374	0.157567059142522\\
0.16672855024518	0.164649001928962\\
0.160221105317827	0.169692044511339\\
0.153486253764215	0.172907811022157\\
0.146606720605894	0.174500048678485\\
0.139654957914573	0.17466263332727\\
0.131305531261426	0.173226918609724\\
0.121669400951649	0.169670441905406\\
0.112258542058668	0.164485913411999\\
0.101896945651585	0.15704154482241\\
0.0896805494193189	0.146149324335627\\
0.0773254236492081	0.132980239026569\\
0.0633632150324718	0.115619153620237\\
0.0490732716966528	0.0951519849040132\\
0.0346034727949408	0.0715133776674443\\
0.0205927616666487	0.0454577737794625\\
0.0077489936726205	0.0181011670198291\\
0.000159238796980543	0.000329722991185277\\
};
\addplot [color=black, draw=none, mark size=4.0pt, mark=asterisk, mark options={solid, black}, forget plot]
  table[row sep=crcr]{%
0.6	-0.6\\
};
\addplot [color=blue, draw=none, mark size=4.0pt, mark=asterisk, mark options={solid, blue}, forget plot]
  table[row sep=crcr]{%
0.000159238796980512	0.000329722991185284\\
};

\addplot[->, thick, color=blue] coordinates
{ (0.8151,-0.2091) (0.8206,-0.1833)};

\addplot[->, thick, color=blue] coordinates
{ (0.3505,0.8858) (0.1974, 0.7916)};

\addplot[->, thick, color=blue] coordinates
{ (0.0708,-0.4077) (0.0868, -0.3786)};

\addplot[->, thick, color=blue] coordinates
{ (-0.7945,0.1887) (-0.8047, 0.14057)};

\addplot[->, thick, color=red] coordinates
{ (0.8377,-0.01366) (0.8376,0.0173)};

\addplot[->, thick, color=red] coordinates
{ (0.2487,-0.2305) (0.2594, -0.1954)};

\addplot[->, thick, color=red] coordinates
{ (0.4,0.6496) (0.36, 0.6497)};

\addplot[->, thick, color=green] coordinates
{ (0.4681,1.0951) (0.4143,1.1931)};

\end{axis}
\end{tikzpicture}%

%% file: path_th3_straight_MIQP.tex
%
%
\begin{tikzpicture}

\begin{axis}[%
width=\figurewidth,
height=\figureheight,
at={(0\figurewidth,0\figureheight)},
scale only axis,
xmin=0,
xmax=80,
xtick={0,20,40,60,80},
xlabel={$t$ [s]},
ymin=-1,
ymax=1,
ytick={-0.8,0,0.8},
xlabel style={font=\color{white!15!black},at={(axis description cs:0.5,-0.17)},anchor=north},
ylabel style={font=\color{white!15!black},at={(axis description cs:-0.24,.5)},anchor=south},
ylabel={$\tilde \theta_3$ [rad]},
axis background/.style={fill=white},
xmajorgrids,
ymajorgrids
]
\addplot [color=red, line width=0.8pt, forget plot]
  table[row sep=crcr]{%
0	0\\
0.299999999999997	-0.0225724202462345\\
0.599999999999994	-0.0466388448038373\\
0.900000000000006	-0.072061944488766\\
1.3	-0.10776560405121\\
1.7	-0.145046380529081\\
3	-0.267802000646768\\
3.3	-0.29390963686437\\
3.59999999999999	-0.318466865004538\\
3.90000000000001	-0.341697225381665\\
4.3	-0.371012431038338\\
4.8	-0.405849839851712\\
6.59999999999999	-0.528083807474943\\
7.2	-0.566723147459214\\
7.7	-0.597526281205091\\
8.2	-0.626968095700576\\
8.7	-0.654955383739534\\
9.2	-0.681387737887775\\
9.7	-0.70615708619971\\
10.2	-0.729147248384493\\
10.6	-0.746174866217302\\
11	-0.761916265898307\\
11.4	-0.776299835367794\\
11.8	-0.789250214469718\\
12.2	-0.800688302325142\\
12.6	-0.810531329133298\\
13	-0.818693010088822\\
13.4	-0.825083802333168\\
13.8	-0.829611289373048\\
14.2	-0.83218072110823\\
14.5	-0.832764891838892\\
14.8	-0.832152524981268\\
15.1	-0.830302739715322\\
15.4	-0.827175048914185\\
15.7	-0.822729745009255\\
16	-0.816928347191848\\
16.3	-0.809734114225847\\
16.6	-0.801112626291655\\
16.9	-0.791032437856117\\
17.2	-0.779465801917993\\
17.5	-0.766389463308485\\
17.8	-0.751766303764356\\
18.1	-0.735210612533677\\
18.3	-0.722900168215659\\
18.5	-0.709545946075934\\
18.8	-0.687822149605665\\
19.1	-0.664598103535852\\
19.4	-0.640020550800642\\
19.8	-0.605443634922665\\
20.2	-0.569281839259304\\
21.6	-0.440666463116713\\
22.1	-0.397272928241534\\
22.8	-0.338453925879875\\
26.3	-0.0468768365432624\\
27.2	0.0280363299737871\\
27.5	0.0515060336979332\\
27.8	0.0733929302020329\\
29.7	0.206405562659313\\
30.4	0.253195261478012\\
31.1	0.298363004134785\\
31.7	0.335557805636654\\
32.3	0.371156008829487\\
32.8	0.399380828747667\\
33.3	0.426136419684909\\
33.8	0.451326457809643\\
34.3	0.474862298323103\\
34.8	0.496627311190693\\
35.2	0.512674590952869\\
35.6	0.527426889175104\\
36	0.540804096355728\\
36.4	0.552720785593252\\
36.8	0.56308595057088\\
37.2	0.571828285121939\\
37.6	0.578909984464431\\
38	0.584292274163715\\
38.4	0.587938547438654\\
38.8	0.58982127756903\\
39.2	0.589929006622796\\
39.6	0.588272351515059\\
40	0.584888248067969\\
40.4	0.579842120935623\\
40.8	0.573227611697448\\
41.2	0.565163637068878\\
41.6	0.55578918139453\\
42.1	0.542461061416574\\
42.6	0.527633298314541\\
43.2	0.50828620744457\\
43.9	0.484173568704861\\
44.9	0.44804402570017\\
47.6	0.349583163423446\\
48.6	0.314996338310451\\
49.5	0.285347540594699\\
50.4	0.257310046452517\\
51.2	0.233851385357539\\
52	0.211832417261647\\
52.8	0.191282633726544\\
53.6	0.172206546670267\\
54.4	0.154587596267973\\
55.2	0.138391789421476\\
56.1	0.121812320967777\\
57	0.106888101891528\\
57.9	0.0935216938651706\\
58.8	0.0816073421887467\\
59.8	0.069937912724285\\
60.8	0.0597715476952203\\
61.9	0.0501440902449701\\
63.1	0.0412793318435973\\
64.4	0.0333391307394635\\
65.8	0.0264182349738746\\
67.3	0.0205455156511221\\
69	0.0154288009396311\\
71	0.0110114629792406\\
73.3	0.00748392348090476\\
76.2	0.00462465291569458\\
80	0.00249225046424328\\
80.1	0.00245252030290999\\
};
\addplot [color=blue, line width=0.8pt, forget plot]
  table[row sep=crcr]{%
0	0\\
0.299999999999997	-0.0225724202462487\\
0.599999999999994	-0.0466388448038657\\
0.900000000000006	-0.0720619444887944\\
1.3	-0.107765604051252\\
1.7	-0.145046380529109\\
3	-0.267802000646427\\
3.3	-0.293893045528918\\
3.59999999999999	-0.318021300098579\\
3.90000000000001	-0.33991093382599\\
4.2	-0.35988762359392\\
4.5	-0.378223758964509\\
4.8	-0.395186141176282\\
5.2	-0.416109898048802\\
5.59999999999999	-0.435581590591966\\
6.09999999999999	-0.458407659068158\\
6.7	-0.484261302537448\\
7.3	-0.508895423088759\\
7.7	-0.523639090528903\\
8.09999999999999	-0.536712446295653\\
8.5	-0.548072637752938\\
8.90000000000001	-0.557679348807653\\
9.3	-0.565495439354137\\
9.7	-0.571487069336953\\
10.1	-0.575623804284305\\
10.5	-0.577883966243547\\
10.9	-0.578226825425361\\
11.2	-0.577147008569\\
11.5	-0.57480331143131\\
11.8	-0.571154235019222\\
12.1	-0.56617022946142\\
12.4	-0.559814847008312\\
12.7	-0.552052462421528\\
13	-0.542849811622077\\
13.3	-0.532176848658253\\
13.6	-0.520007529485454\\
13.9	-0.506320637191195\\
14.2	-0.490995875553025\\
14.5	-0.473556108364377\\
14.7	-0.460640777332969\\
15	-0.439406779804813\\
15.3	-0.416557294615444\\
15.6	-0.392310251551649\\
16	-0.358086502125914\\
16.4	-0.322150106127808\\
17.1	-0.257231452044365\\
17.7	-0.201772813946349\\
18	-0.175322109445318\\
18.3	-0.15040337072098\\
18.6	-0.126976343613066\\
18.9	-0.104884865147596\\
19.3	-0.0771128045730052\\
19.8	-0.0442221095773192\\
21.2	0.0462884002469934\\
21.7	0.0796330402346968\\
22.2	0.111139330994803\\
22.7	0.140957114790737\\
23.2	0.169256248340147\\
23.7	0.19604332795646\\
24.2	0.221215505816303\\
24.7	0.244671980977444\\
25.2	0.266291621840537\\
25.6	0.282194702630505\\
26	0.29671757536164\\
26.4	0.309827144281385\\
26.8	0.321508522978263\\
27.2	0.33173933310465\\
27.6	0.340495588489432\\
28	0.347756531689029\\
28.4	0.353509632042943\\
28.8	0.357754691698602\\
29.2	0.360507356956006\\
29.6	0.361801081669199\\
30	0.361687787087106\\
30.4	0.360237173459211\\
30.9	0.356674861653474\\
31.4	0.351351817681817\\
31.9	0.344482557614924\\
32.5	0.334512190675937\\
33.1	0.323014604720484\\
33.8	0.308140037241515\\
34.7	0.287457174782091\\
36.3	0.248768838489582\\
38	0.208115266853\\
39.1	0.183253774099569\\
40.1	0.162071966102076\\
41.1	0.142442047081445\\
42.1	0.124468564693899\\
43	0.109739269181716\\
44	0.0949661953630283\\
45	0.0818151372234297\\
46	0.0702002761363332\\
47.1	0.0590711630972152\\
48.2	0.049514546327714\\
49.4	0.0406889778373909\\
50.7	0.0327724575077752\\
52.1	0.0258733919273055\\
53.6	0.0200299728019218\\
55.3	0.0149572708838406\\
57.3	0.010603210724824\\
59.6	0.00715444371381579\\
62.5	0.00438861046762895\\
66.4	0.00231437961477354\\
72.2	0.000924756981234509\\
80.1	0.000273078343766997\\
};
\addplot [color=green, line width=0.8pt, forget plot]
  table[row sep=crcr]{%
0	0\\
0.3	-0.0225724202462523\\
0.600000000000001	-0.0466388448038702\\
0.9	-0.0720619444888007\\
1.3	-0.107765604051251\\
1.7	-0.145046380529113\\
3	-0.267802000646429\\
3.3	-0.293893045528921\\
3.6	-0.31802130009858\\
3.8	-0.332731857067077\\
4	-0.346151101725312\\
4.2	-0.358131079397014\\
4.4	-0.368548011258219\\
4.6	-0.377312162521712\\
4.8	-0.384377625344602\\
5	-0.389750854806048\\
5.2	-0.393496603898241\\
5.4	-0.395739931850733\\
5.6	-0.396663284377071\\
5.9	-0.396090306425521\\
6.3	-0.39312427097513\\
6.6	-0.390491629490904\\
};
\end{axis}
\end{tikzpicture}%

%% file: path_z3_straight_MIQP.tex
%
%
\begin{tikzpicture}

\begin{axis}[%
width=\figurewidth,
height=\figureheight,
at={(0\figurewidth,0\figureheight)},
scale only axis,
xmin=0,
xmax=80,
xtick={0,20,40,60,80},
xlabel={$t$ [s]},
ymin=-11,
ymax=11,
xlabel style={font=\color{white!15!black},at={(axis description cs:0.5,-0.17)},anchor=north},
ylabel style={font=\color{white!15!black},at={(axis description cs:-0.24,.5)},anchor=south},
ylabel={$\tilde z_3$ [m]},
axis background/.style={fill=white},
xmajorgrids,
ymajorgrids
]
\addplot [color=red, line width=0.8pt, forget plot]
  table[row sep=crcr]{%
0	0\\
0.5	0.00784624814450297\\
1	0.0317050999199893\\
1.5	0.0715444831367194\\
2	0.126672013043859\\
2.59999999999999	0.211080840166005\\
3.3	0.329694176986806\\
4	0.464891044391408\\
4.59999999999999	0.596725451811579\\
5.09999999999999	0.722131648885522\\
5.59999999999999	0.86646068104362\\
6.09999999999999	1.03216382294815\\
6.7	1.25111806557925\\
7.3	1.49013792391648\\
7.90000000000001	1.74885283981831\\
8.5	2.02678832978808\\
9.09999999999999	2.32334457158768\\
9.7	2.6377799486457\\
10.3	2.96919112617732\\
11	3.3758281644544\\
11.7	3.80197509220864\\
12.5	4.30955551411044\\
13.4	4.90101442562145\\
15	5.97855983352051\\
16.1	6.70946327835918\\
16.8	7.15556408482288\\
17.4	7.51849295861034\\
18	7.85980520064408\\
19	8.41119147424773\\
19.5	8.66151690678183\\
19.9	8.84537822646224\\
20.3	9.01417122887882\\
20.7	9.16783271740582\\
21.1	9.30623798537937\\
21.6	9.4594579674647\\
22.1	9.59365181061895\\
22.7	9.73513215956204\\
23.3	9.85861673646943\\
23.9	9.9647916932754\\
24.5	10.0530068000818\\
25.1	10.1222538753114\\
25.7	10.1717778313703\\
26.2	10.1976245301253\\
26.7	10.2092849281732\\
27.3	10.2047138305405\\
27.9	10.1827617270884\\
28.5	10.1437279439208\\
29	10.095070057713\\
29.6	10.0183064021374\\
30.2	9.92147166174018\\
30.7	9.8248997120153\\
31.3	9.69020069670006\\
31.9	9.53558507023592\\
32.5	9.36209112950657\\
33.1	9.17112125766234\\
33.7	8.96341440231136\\
34.3	8.73934467403721\\
35	8.45823959707641\\
35.7	8.15759793879135\\
36.4	7.83971612097797\\
37.2	7.45866931509495\\
38	7.0587086497485\\
38.8	6.6403900284258\\
39.9	6.04441126428684\\
41.6	5.12251451917771\\
42.4	4.70825143956324\\
43.1	4.36266839549155\\
43.8	4.03488827528932\\
44.5	3.72555173397377\\
45.2	3.4344849216495\\
45.9	3.16108901524657\\
46.6	2.90460982524434\\
47.3	2.66428478610055\\
48.1	2.40850176592681\\
48.9	2.17194183061255\\
49.7	1.9537433552882\\
50.5	1.75309504829541\\
51.3	1.56919895246079\\
52.1	1.40124729351892\\
53	1.23032313603234\\
53.9	1.07729703944463\\
54.8	0.940902477210869\\
55.8	0.807290768510939\\
56.8	0.69084999944333\\
57.9	0.580525168399191\\
59	0.486653971668176\\
60.2	0.40056569895988\\
61.5	0.323720225793508\\
63	0.25269145948171\\
64.6	0.193752013359457\\
66.5	0.141263326972762\\
68.7	0.0980680903156497\\
71.4	0.062877100526066\\
74.8	0.0362028831343792\\
79.5	0.0171257877267124\\
80.1	0.0155809376102667\\
};
\addplot [color=blue, line width=0.8pt, forget plot]
  table[row sep=crcr]{%
0	0\\
0.5	0.00784624814451718\\
1	0.0317050999200177\\
1.5	0.071544483136833\\
2	0.126672013044043\\
2.59999999999999	0.211080840166289\\
3.3	0.32966236351281\\
4.3	0.518934068571554\\
5	0.66576567478667\\
5.59999999999999	0.808452713224881\\
6.09999999999999	0.945130310764242\\
6.5	1.06998125634591\\
6.90000000000001	1.21218125052819\\
7.40000000000001	1.41446586028681\\
8	1.67562838085227\\
8.7	2.00115866672124\\
9.40000000000001	2.34509217813873\\
10.3	2.80672449406863\\
11.2	3.26913280374521\\
12.5	3.90234403665566\\
13.2	4.22279668763525\\
13.8	4.47833567408702\\
14.5	4.75301406591494\\
15.2	5.01192554135012\\
15.7	5.17499039290412\\
16.2	5.31658095101675\\
16.6	5.41406575566516\\
17	5.49728886633577\\
17.5	5.58268284840396\\
18	5.64882123491044\\
18.5	5.69783672610053\\
19.1	5.73916949519959\\
19.8	5.76834777418563\\
20.5	5.77843322920623\\
21.1	5.76932256699227\\
21.7	5.74038560850327\\
22.3	5.69350163039509\\
22.9	5.62948764433419\\
23.5	5.54835032923373\\
24.1	5.45047851823324\\
24.7	5.33644600141602\\
25.3	5.20703398215554\\
26	5.03820652410958\\
26.7	4.8518692954319\\
27.4	4.64835064522927\\
28.2	4.3966576612945\\
29.1	4.09416564873017\\
30.5	3.60197063739018\\
32	3.0791055227194\\
33	2.74908310359176\\
33.9	2.47039887603795\\
34.7	2.23887684065502\\
35.5	2.02305846669019\\
36.3	1.82287157877602\\
37.2	1.61590421284731\\
38.1	1.42755467249985\\
39	1.25697431298883\\
39.9	1.10323393224739\\
40.8	0.965343377831644\\
41.8	0.829457641767874\\
42.8	0.710384949343052\\
43.9	0.597034043863388\\
45	0.500222613177897\\
46.2	0.41120003373392\\
47.5	0.331613198761175\\
48.9	0.262406748599631\\
50.5	0.200425332603388\\
52.3	0.147865055738905\\
54.4	0.103759940757925\\
56.9	0.0682905408087322\\
60.1	0.0403177652734286\\
64.5	0.019879211474489\\
71.3	0.00686102116888776\\
80.1	0.00175046532582712\\
};
\addplot [color=green, line width=0.8pt, forget plot]
  table[row sep=crcr]{%
0	0\\
0.5	0.00784624814451185\\
1	0.031705099920023\\
1.5	0.0715444831368268\\
2	0.126672013044048\\
2.6	0.211080840166296\\
3.3	0.329662363512805\\
4.6	0.560731995591348\\
5	0.613210085429941\\
5.3	0.639205214707317\\
5.6	0.650714135419406\\
5.9	0.645849268141509\\
6.2	0.623784628580728\\
6.5	0.584781702836406\\
6.6	0.568195736773119\\
};
\end{axis}
\end{tikzpicture}%

%% file: path_xy_eight_MIQP.tex
%
%
\begin{tikzpicture}

\begin{axis}[%
width=\figurewidth,
height=\figureheight,
at={(0\figurewidth,0\figureheight)},
scale only axis,
xmin=-35,
xmax=35,
xlabel={$x$ [m]},
ymin=-20,
ymax=20,
xlabel style={font=\color{white!15!black},at={(axis description cs:0.5,-0.1)},anchor=north},
ylabel style={font=\color{white!15!black},at={(axis description cs:-0.1,.5)},anchor=south},
ylabel={$y$ [m]},
axis background/.style={fill=white},
xmajorgrids,
ymajorgrids
]
\addplot [color=black, dashed, forget plot]
  table[row sep=crcr]{%
1.26721717753919	-1.91576302677241\\
-5.13205354653223	5.76859343294771\\
};
\addplot [color=black, dashed, forget plot]
  table[row sep=crcr]{%
1.26721717753919	-1.91576302677241\\
1.22993383500526	-11.9156935241494\\
};
\addplot [color=red, forget plot]
  table[row sep=crcr]{%
-9.23925104947098	5.20358011225474\\
-10.8214774154334	3.33300273431678\\
-2.20384196041885	-3.95621726267025\\
-0.621615594456397	-2.0856398847323\\
-9.23925104947098	5.20358011225474\\
};
\addplot [color=black, line width=1.0pt, forget plot]
  table[row sep=crcr]{%
-8.38837990327485	4.16298121241453\\
-7.76485411062886	3.63557242376038\\
};
\addplot [color=black, line width=1.0pt, forget plot]
  table[row sep=crcr]{%
-9.65416099604481	2.66651931006416\\
-9.03063520339882	2.13911052141001\\
};
\addplot [color=black, line width=1.0pt, forget plot]
  table[row sep=crcr]{%
-8.07661700695185	3.89927681808745\\
-9.34239809972182	2.40281491573709\\
};
\addplot [color=black, line width=1.0pt, forget plot]
  table[row sep=crcr]{%
-8.70950755333683	3.15104586691227\\
-2.60149978864146	-2.01540757296513\\
};
\addplot [color=black, line width=1.0pt, forget plot]
  table[row sep=crcr]{%
-2.60149978864146	-2.01540757296513\\
1.26721717753919	-1.91576302677241\\
};
\addplot [color=black, line width=1.0pt, forget plot]
  table[row sep=crcr]{%
-3.03493073174033	-1.04624621969338\\
-2.21853481725777	-1.02521873147872\\
};
\addplot [color=black, line width=1.0pt, forget plot]
  table[row sep=crcr]{%
-2.98446476002516	-3.00559641445154\\
-2.16806884554259	-2.98456892623689\\
};
\addplot [color=black, line width=1.0pt, forget plot]
  table[row sep=crcr]{%
-2.62673277449905	-1.03573247558605\\
-2.57626680278387	-2.99508267034421\\
};
\addplot [color=blue, forget plot]
  table[row sep=crcr]{%
1.32339419902442	-0.593692944569106\\
2.15275418014946	-2.89904801954308\\
8.46787449950726	-0.627160219819185\\
7.63851451838221	1.67819485515479\\
1.32339419902442	-0.593692944569106\\
};
\addplot [color=black, line width=1.0pt, forget plot]
  table[row sep=crcr]{%
6.47271651755048	0.998422691214664\\
7.24116820920847	1.27487601825634\\
};
\addplot [color=black, line width=1.0pt, forget plot]
  table[row sep=crcr]{%
7.13620450245051	-0.845861368764518\\
7.90465619410851	-0.569408041722838\\
};
\addplot [color=black, line width=1.0pt, forget plot]
  table[row sep=crcr]{%
2.12077071275267	-0.56720584458261\\
2.88922240441066	-0.290752517540929\\
};
\addplot [color=black, line width=1.0pt, forget plot]
  table[row sep=crcr]{%
2.7842586976527	-2.41148990456179\\
3.55271038931069	-2.13503657752011\\
};
\addplot [color=black, line width=1.0pt, forget plot]
  table[row sep=crcr]{%
6.85694236337948	1.1366493547355\\
7.52043034827951	-0.707634705243677\\
};
\addplot [color=black, line width=1.0pt, forget plot]
  table[row sep=crcr]{%
2.50499655858166	-0.428979181061769\\
3.1684845434817	-2.27326324104095\\
};
\addplot [color=black, line width=1.0pt, forget plot]
  table[row sep=crcr]{%
1.26721717753919	-1.91576302677241\\
7.18868635582949	0.214507324745913\\
};
\addplot [color=red, line width=0.7pt, forget plot]
  table[row sep=crcr]{%
-8.70950755333683	3.15104586691227\\
-9.12975303169038	3.5285807210705\\
-9.54374458099068	3.94275649385857\\
-9.99090557852492	4.43473441218477\\
-10.4627869808681	5.00163711442302\\
-10.9824317745657	5.67930228462149\\
-11.5487368379279	6.47324272932265\\
-12.2660476686344	7.53675490056922\\
-13.5654846329285	9.47275895680237\\
-14.0997159075967	10.2092591044938\\
-14.5863798951834	10.8262198849164\\
-15.0192341906362	11.3239012841715\\
-15.4114689226537	11.7294917821787\\
-15.8180432976033	12.1034422523783\\
-16.2218537418709	12.4300371271725\\
-16.6184876316911	12.7109115221974\\
-17.0124112182709	12.955108944556\\
-17.4176116027865	13.1744826788022\\
-17.8378471253929	13.371489236267\\
-18.2749141496027	13.5464126300242\\
-18.7258947597806	13.6972583386161\\
-19.1872301230979	13.8222897171016\\
-19.656211076575	13.9204977250944\\
-20.1306194436506	13.9912898088249\\
-20.6085739828217	14.0343305340962\\
-21.0884350307588	14.0494396073929\\
-21.5680061340655	14.0365596902934\\
-22.0376492732865	13.9967784760838\\
-22.5683095890993	13.9196385398707\\
-23.1015039148977	13.8082486696084\\
-23.6420805163215	13.6612590182567\\
-24.1912515728572	13.4772970327545\\
-24.7478809231114	13.2554082798436\\
-25.3088980789203	12.9954836416859\\
-25.8707342911252	12.6980032930382\\
-26.4291305813129	12.3642726003526\\
-26.9799924049111	11.9960483418803\\
-27.5196019232146	11.5953357700001\\
-28.0445845217319	11.1643172232841\\
-28.5519277307272	10.7052618774287\\
-29.0391333879611	10.2202958005579\\
-29.5039926426579	9.71150835627486\\
-29.944587675121	9.18089820483619\\
-30.359269674765	8.63034634180866\\
-30.7465738190894	8.06168506833237\\
-31.1052313992218	7.4766596680961\\
-31.4341567438936	6.87691391932143\\
-31.7323647442624	6.26411005482262\\
-31.9989892957737	5.63990461564097\\
-32.2332965408068	5.00591236901646\\
-32.4346730389611	4.36370296684797\\
-32.6025907977395	3.71487720326643\\
-32.7366042871333	3.06109556776966\\
-32.8363686133006	2.40403254537463\\
-32.901644967391	1.7453224027052\\
-32.9322842693852	1.08659814649839\\
-32.9282245913603	0.429613642790898\\
-32.8895128193863	-0.223894220976426\\
-32.8162917000317	-0.872302386693299\\
-32.7240855055878	-1.43424324865633\\
-32.605845210925	-1.98990912652025\\
-32.461851029823	-2.53816166187401\\
-32.2924157798218	-3.07797049419617\\
-32.0978680252129	-3.6083909835091\\
-31.8786242999125	-4.12834424350353\\
-31.6351761598051	-4.63671874804108\\
-31.3680594694171	-5.13245530480388\\
-31.0778259948638	-5.61458146595956\\
-30.7650811774239	-6.08212623708057\\
-30.4305268192117	-6.5340721986856\\
-30.0749037881568	-6.96945896350023\\
-29.6989638796378	-7.38740955051829\\
-29.3034461635567	-7.78713557388362\\
-28.8891940635309	-8.16780841692714\\
-28.4571148734085	-8.5286231415389\\
-27.9426462110821	-8.91574348175816\\
-27.4073821978775	-9.27505540239475\\
-26.8527058971448	-9.60566863135116\\
-26.2801336466293	-9.90672074322369\\
-25.6911573619289	-10.1774820315597\\
-25.0871687126116	-10.4173676600959\\
-24.4696352231724	-10.6258426979218\\
-23.8401371298545	-10.8024293430017\\
-23.2001928342102	-10.9467744588732\\
-22.5511816005925	-11.0586455647896\\
-21.8945524759006	-11.1378723605776\\
-21.2318307007599	-11.1843658843319\\
-20.5644468105469	-11.1981453473551\\
-19.8936658462757	-11.1793196281016\\
-19.2208064663123	-11.1280734213893\\
-18.5472325603288	-11.0446866496403\\
-17.8741808315191	-10.9295265433\\
-17.2026992361325	-10.7830175331874\\
-16.5338823044055	-10.6056746241495\\
-15.8688309149414	-10.398114767811\\
-15.2084889677524	-10.1610171841128\\
-14.553596517944	-9.89508679575892\\
-13.9049376360296	-9.60113914071059\\
-13.2632568630764	-9.28008672274794\\
-12.5503948518207	-8.88767871284406\\
-12.2367629401808	-8.70299655008945\\
};
\addplot [color=blue, line width=0.7pt, forget plot]
  table[row sep=crcr]{%
-8.70950755333683	3.15104586691227\\
-9.11721361029933	3.5164828442399\\
-9.55312393671919	3.95036987500224\\
-10.0306596997839	4.4701557886134\\
-10.5829349832133	5.11871270041058\\
-11.2089155455231	5.90157700970538\\
-12.1861415368885	7.1801315359242\\
-13.18043722244	8.46671114564315\\
-13.731763978661	9.13134750525476\\
-14.1761029888075	9.62521313605926\\
-14.5946115831144	10.0490461966029\\
-15.0198880285574	10.4357471964856\\
-15.4441174411061	10.7768767777231\\
-15.8626594127284	11.0711352628801\\
-16.2705978662736	11.3205311255285\\
-16.6767361422682	11.5361315743434\\
-17.1505716396059	11.7511294507688\\
-17.6416188348027	11.9370721373073\\
-18.1559098555095	12.0960795107779\\
-18.6956015439607	12.2276509662177\\
-19.2599898694169	12.3300958065637\\
-19.8464793218342	12.4014378685921\\
-20.4513207988924	12.4399079156204\\
-20.9922327616257	12.4455453648631\\
-21.5409092191901	12.4244113099119\\
-22.0944274566562	12.3762072704843\\
-22.6500395124868	12.3008729801676\\
-23.2052107441441	12.1985393064265\\
-23.7576269300952	12.0694851763766\\
-24.3051955092518	11.9140989129262\\
-24.846003768685	11.732856555407\\
-25.3783043327431	11.5263026046588\\
-25.9005063928067	11.2950324117202\\
-26.4111491609232	11.0396854950942\\
-26.908892915211	10.7609356769878\\
-27.3924990848491	10.4594893457542\\
-27.8608158326193	10.1360847770599\\
-28.3127698312008	9.79149047203089\\
-28.7473684996542	9.42649612553628\\
-29.1636839833084	9.04191755087622\\
-29.5608546983216	8.638588415257\\
-29.9380595692078	8.21737985142469\\
-30.2945261326546	7.77919114587797\\
-30.6295273802589	7.32495044942085\\
-30.942384324438	6.85560722491629\\
-31.2324505622045	6.3721505923377\\
-31.4991028235727	5.87563077344547\\
-31.7417597315673	5.36713685518722\\
-31.9598873563565	4.84778797191529\\
-32.1530024489014	4.31872279707881\\
-32.3206609124265	3.78112711223658\\
-32.4624534928043	3.23627024643033\\
-32.5780305560797	2.68545675014783\\
-32.6671114799606	2.12999373231731\\
-32.7294846490927	1.57116137822253\\
-32.764993903562	1.01028941553938\\
-32.7735485363457	0.448776190085511\\
-32.7551365337801	-0.112008334768177\\
-32.7098158191065	-0.670778501874871\\
-32.6376971175709	-1.22632601068216\\
-32.5389622980025	-1.7773344659165\\
-32.4138654818184	-2.32249682888694\\
-32.2627130967148	-2.86059966050532\\
-32.0858317010367	-3.39056378001332\\
-31.8836092119256	-3.9112596058906\\
-31.6565270581937	-4.42147007394959\\
-31.405121155699	-4.92002635944532\\
-31.1299586335932	-5.40584154838945\\
-30.8316306300223	-5.87789681023332\\
-30.5108259270021	-6.33511959030317\\
-30.168292936875	-6.77646956115147\\
-29.8048010091176	-7.20098997923003\\
-29.4211047076284	-7.60782954633871\\
-29.0180158146776	-7.99614611256726\\
-28.5964443935708	-8.36508189636941\\
-28.1573194064625	-8.71385213453895\\
-27.6351275546621	-9.08687652247126\\
-27.0925214479701	-9.4317513386886\\
-26.5310023497126	-9.74755255973352\\
-25.9521047642782	-10.0334603069737\\
-25.3572895415683	-10.2888045390272\\
-24.747967517175	-10.5130247631738\\
-24.1257150039382	-10.7055854541203\\
-23.4921079743663	-10.8660599230706\\
-22.848608353166	-10.9941516568673\\
-22.196606839582	-11.089658249199\\
-21.5376429344401	-11.1524359409873\\
-20.8732165645444	-11.1824473188382\\
-20.2046824257248	-11.1797576435575\\
-19.5333090405836	-11.1445102021832\\
-18.8604849357861	-11.0769358697665\\
-18.187529164831	-10.9773616908375\\
-17.5155835958603	-10.8461866502845\\
-16.8456961438606	-10.6838741583499\\
-16.1790074760253	-10.4909988348808\\
-15.5165543830029	-10.2682163996946\\
-14.8591689767387	-10.0162202310501\\
-14.2075950523372	-9.73576317901549\\
-13.5626326304228	-9.42772505858783\\
-12.8457697497082	-9.04938326665164\\
-12.1387256919445	-8.63866517109056\\
-11.4420557497017	-8.19705316885985\\
-10.7562523920234	-7.72619080268411\\
-10.1559250277305	-7.28442330099264\\
};
\addplot [color=green, line width=0.7pt, forget plot]
  table[row sep=crcr]{%
-8.70950755333683	3.15104586691227\\
-9.13646948842353	3.53427354414599\\
-9.74073460888635	4.12780076990395\\
-10.2044363214775	4.5707608486211\\
};
\addplot [color=black, dashed, forget plot]
  table[row sep=crcr]{%
20.7507152356411	-8.46894872361294\\
28.0630437798902	-1.64766942563953\\
};
\addplot [color=black, dashed, forget plot]
  table[row sep=crcr]{%
20.7507152356411	-8.46894872361294\\
10.7645127731103	-8.99407763858749\\
};
\addplot [color=red, forget plot]
  table[row sep=crcr]{%
8.68160986792619	-2.99120174471033\\
8.32422443350206	-5.41499539336557\\
19.4904933773272	-7.06144820901669\\
19.8478788117514	-4.63765456036145\\
8.68160986792619	-2.99120174471033\\
};
\addplot [color=black, line width=1.0pt, forget plot]
  table[row sep=crcr]{%
9.95340082338419	-3.42637474869038\\
10.7613320396026	-3.54550322683176\\
};
\addplot [color=black, line width=1.0pt, forget plot]
  table[row sep=crcr]{%
9.66749247584488	-5.36540966761457\\
10.4754236920633	-5.48453814575595\\
};
\addplot [color=black, line width=1.0pt, forget plot]
  table[row sep=crcr]{%
10.3573664314934	-3.48593898776107\\
10.0714580839541	-5.42497390668526\\
};
\addplot [color=black, line width=1.0pt, forget plot]
  table[row sep=crcr]{%
10.2144122577237	-4.45545644722316\\
18.1288404982306	-5.62242929432238\\
};
\addplot [color=black, line width=1.0pt, forget plot]
  table[row sep=crcr]{%
18.1288404982306	-5.62242929432238\\
20.7507152356411	-8.46894872361294\\
};
\addplot [color=black, line width=1.0pt, forget plot]
  table[row sep=crcr]{%
18.573024009754	-4.65814867168088\\
19.1263050697588	-5.25883537295236\\
};
\addplot [color=black, line width=1.0pt, forget plot]
  table[row sep=crcr]{%
17.1313759267025	-5.98602321569241\\
17.6846569867073	-6.58670991696388\\
};
\addplot [color=black, line width=1.0pt, forget plot]
  table[row sep=crcr]{%
18.8496645397564	-4.95849202231662\\
17.4080164567049	-6.28636656632814\\
};
\addplot [color=blue, forget plot]
  table[row sep=crcr]{%
22.073852323398	-8.45068928762941\\
19.8187850796451	-9.40838011777924\\
22.4422129787986	-15.585745646906\\
24.6972802225514	-14.6280548167561\\
22.073852323398	-8.45068928762941\\
};
\addplot [color=black, line width=1.0pt, forget plot]
  table[row sep=crcr]{%
23.9530242985117	-13.5023291427383\\
24.2722545752283	-14.2540182239893\\
};
\addplot [color=black, line width=1.0pt, forget plot]
  table[row sep=crcr]{%
22.1489705035094	-14.2684818068582\\
22.468200780226	-15.0201708881091\\
};
\addplot [color=black, line width=1.0pt, forget plot]
  table[row sep=crcr]{%
22.1451385477186	-9.24531444789884\\
22.4643688244352	-9.99700352914978\\
};
\addplot [color=black, line width=1.0pt, forget plot]
  table[row sep=crcr]{%
20.3410847527164	-10.0114671120187\\
20.660315029433	-10.7631561932696\\
};
\addplot [color=black, line width=1.0pt, forget plot]
  table[row sep=crcr]{%
24.11263943687	-13.8781736833638\\
22.3085856418677	-14.6443263474837\\
};
\addplot [color=black, line width=1.0pt, forget plot]
  table[row sep=crcr]{%
22.3047536860769	-9.62115898852431\\
20.5006998910747	-10.3873116526442\\
};
\addplot [color=black, line width=1.0pt, forget plot]
  table[row sep=crcr]{%
20.7507152356411	-8.46894872361294\\
23.2106125393688	-14.2612500154237\\
};
\addplot [color=red, line width=0.7pt, forget plot]
  table[row sep=crcr]{%
10.2144122577237	-4.45545644722316\\
10.7551054181758	-4.5513704862232\\
11.2674087423184	-4.67435765304165\\
11.8021387697465	-4.83740196134681\\
12.3511803843909	-5.04076387670088\\
12.8886905190193	-5.27254076640563\\
13.5056656712271	-5.5729491027714\\
14.2161362689257	-5.95545083486788\\
15.0944129544818	-6.46717715985542\\
16.6073686765191	-7.39464231945882\\
17.8597120932324	-8.14702720074485\\
18.703815621773	-8.61929699585993\\
19.4563142684961	-9.00442240816135\\
20.1190319472575	-9.30916772590992\\
20.7646711304056	-9.57041271799615\\
21.3415023435263	-9.77060288168085\\
21.864893310212	-9.92136146721764\\
22.3650751620589	-10.0343095619152\\
22.8319208932447	-10.1110885715359\\
23.3073046707997	-10.161185176367\\
23.7852177916311	-10.183443043303\\
24.2636259797582	-10.1777185791438\\
24.7408477134897	-10.1440279035004\\
25.2152479826803	-10.0824876274438\\
25.6851842419147	-9.9933142170058\\
26.1490569026598	-9.87681537110651\\
26.6053130258507	-9.73338160189614\\
27.0524318513223	-9.56348571439937\\
27.4888500763474	-9.36771741627145\\
27.9130762799893	-9.14674467725711\\
28.3236279958228	-8.90134157338996\\
28.7190982484739	-8.63235403641951\\
29.0981293542887	-8.34071092156427\\
29.459504583152	-8.02734377618\\
29.8019846883785	-7.69331359945631\\
30.1244087766152	-7.33973760155472\\
30.4256051226872	-6.96789781586706\\
30.7046248960255	-6.5789574724307\\
30.9604581028488	-6.17431053870431\\
31.192269367185	-5.75526239304979\\
31.3992660902373	-5.3232137912232\\
31.5802213330529	-4.88096345193684\\
31.7344629495982	-4.43082567719454\\
31.8612486868954	-3.97630856870142\\
31.9718207774009	-3.45797674829239\\
32.0462462454157	-2.94610392951711\\
32.0835009091806	-2.4874988774174\\
32.0923732869472	-1.99883543795056\\
32.0691353053665	-1.48942096024369\\
32.0140936861035	-0.963091092067479\\
31.9270419154825	-0.408318026237239\\
31.7874948164714	0.267265293737353\\
31.6084048970827	0.973465078532804\\
31.3951138887554	1.69487214607796\\
31.1210736605729	2.50950990576767\\
30.8169466647775	3.31609892568034\\
30.4869393631639	4.10737428876531\\
30.1335479769338	4.87879630748642\\
29.7580963679949	5.62727851637297\\
29.3610986755476	6.35065130678787\\
28.9427310286684	7.04702559865563\\
28.5029736799927	7.71478212423139\\
28.0939977484083	8.28323690822675\\
27.6679014131616	8.8271721437559\\
27.2246190601094	9.34580972962829\\
26.7641034604367	9.83846461454473\\
26.2863052579182	10.3045490893576\\
25.7912499189708	10.7434813832644\\
25.2790174218788	11.1547032847763\\
24.7497750473257	11.5376481538622\\
24.2037583300813	11.8917553915018\\
23.6413177065057	12.2164414733526\\
23.0629172460288	12.5111127059916\\
22.4691592057958	12.7751654419535\\
21.860753067329	13.0080148455365\\
21.2385532625484	13.2090933321839\\
20.6035414816299	13.3778758599669\\
19.9568161906393	13.5138986381182\\
19.2995861523954	13.616775456655\\
18.6331431578839	13.686215846762\\
17.9588478914578	13.7220370340625\\
17.278116774881	13.7241743448568\\
16.5923923074918	13.6926893082337\\
15.9031260546617	13.6277735291671\\
15.2117496663454	13.5297485888161\\
14.5196614752137	13.3990641816044\\
13.828217206702	13.2362960593192\\
13.1387085945374	13.0421390908956\\
12.4523374337925	12.817394964178\\
11.7702120086782	12.562964209441\\
11.0933556460715	12.2798435180236\\
10.4226797139401	11.9691072619163\\
9.67653124808369	11.5879353133803\\
8.94013805910514	11.1749479388375\\
8.21427504031296	10.7319247853164\\
7.49951686587783	10.2606863481491\\
6.79623900409235	9.7630714247076\\
6.02853384984725	9.1814955177861\\
5.27529528796578	8.57219935098022\\
5.12636776770015	8.44723429661776\\
};
\addplot [color=blue, line width=0.7pt, forget plot]
  table[row sep=crcr]{%
10.2144122577237	-4.45545644722316\\
10.7551054181755	-4.55137048622316\\
11.2674087423184	-4.67435765304173\\
11.8021387697487	-4.83740196134769\\
12.3525741879296	-5.04132535023232\\
12.9130849589315	-5.28394358415639\\
13.5359059767559	-5.59002338800844\\
14.2146976349142	-5.95858949809105\\
15.1414929967287	-6.50055301080237\\
17.9328497970683	-8.16165840895679\\
18.7008260629663	-8.57336936770981\\
19.4077178272084	-8.91786671505598\\
20.037102356832	-9.18888552309035\\
20.5736344633128	-9.38571411463887\\
21.0574835632742	-9.53110153446133\\
21.5121465210896	-9.63883191267719\\
21.9780171727164	-9.72083055019085\\
22.4484024926062	-9.77529117794712\\
22.9211418419755	-9.80190894552159\\
23.3945835730551	-9.80058950928711\\
23.8670940680207	-9.77133697346569\\
24.3372583496945	-9.71421849501665\\
24.8027840880429	-9.62953461166278\\
25.2625921716698	-9.51750426858128\\
25.7149726158062	-9.37853082674765\\
26.1583421261908	-9.21309984827922\\
26.5911799067835	-9.02177684343076\\
27.0119638659473	-8.80523270703286\\
27.4192125820679	-8.56423126790387\\
27.8115084611082	-8.29961413810672\\
28.1874711411274	-8.01231653968331\\
28.5457812338097	-7.70335457418767\\
28.8852770258436	-7.37372828841305\\
29.2048506232204	-7.0244877241136\\
29.5033039478328	-6.65691745203161\\
29.7795510958933	-6.2723596483541\\
30.0327594838088	-5.87195402170362\\
30.2608999286571	-5.4592301596312\\
30.4626668047495	-5.0375594329114\\
30.6583411416709	-4.55326001372663\\
30.8164148908447	-4.0728382338846\\
30.9295438505005	-3.63715395257896\\
31.0200911416019	-3.1674947058195\\
31.0831743696309	-2.67071716489638\\
31.1184044709839	-2.14697389321485\\
31.1266161614886	-1.58463869072461\\
31.1072984555796	-0.983754011996659\\
31.051954354188	-0.26519393565221\\
30.9634634348382	0.471000228306096\\
30.8449439077505	1.21184741508578\\
30.6792151407045	2.04010881024502\\
30.4820964847606	2.85519391328084\\
30.2554335317529	3.65252419506715\\
30.0002060284303	4.42921557116586\\
29.7168653471386	5.18319355198238\\
29.4414831241693	5.83311578825995\\
29.1439162093899	6.46302745649711\\
28.8241269691955	7.07205889344613\\
28.4820897536332	7.65940924349962\\
28.1177779798141	8.22437864475683\\
27.7312361329834	8.76625378584222\\
27.3225922583985	9.28430864430911\\
26.8920468478009	9.77782497344265\\
26.4398400687062	10.2461223133978\\
25.9663068666938	10.6884923798144\\
25.4718793128499	11.1042072871997\\
24.957105613067	11.4925135224134\\
24.4226166152837	11.8526675928856\\
23.8691492460552	12.1839251676879\\
23.2975555346171	12.4855482724275\\
22.7833054927348	12.724591717311\\
22.256627263013	12.9399419998662\\
21.7182837485103	13.1311757441959\\
21.1690877731895	13.2979071976336\\
20.5292516939071	13.458017796701\\
19.8777219146154	13.5852483279963\\
19.2158979081182	13.6792663685032\\
18.5452188005329	13.739844228682\\
17.8671484715116	13.766866230854\\
17.1831629062418	13.7603338411447\\
16.4947364586252	13.7203693651534\\
15.8033208821293	13.6472168876878\\
15.1103306277994	13.5412409040909\\
14.4171308725545	13.4029238801782\\
13.7250388003556	13.232865738862\\
13.0352978873682	13.0317759613092\\
12.3490607309262	12.8004635752219\\
11.6673891612191	12.539831514345\\
10.9912613587068	12.2508744228373\\
10.3215430185698	11.9346595975138\\
9.57669356470693	11.5477400554421\\
8.8417856247806	11.1294321500574\\
8.1175476284346	10.6814997721462\\
7.4045090526273	10.2057439865326\\
6.62578451744093	9.64670989968181\\
5.86153208200778	9.05810664996132\\
5.1117190437437	8.44244802549638\\
4.30325493934069	7.73692575298291\\
4.08565491908321	7.53977728206272\\
};
\addplot [color=green, line width=0.7pt, forget plot]
  table[row sep=crcr]{%
10.2144122577237	-4.45545644722316\\
10.7551054181755	-4.55137048622316\\
11.2674087423184	-4.67435765304173\\
11.8021387697487	-4.83740196134769\\
12.3525741879296	-5.04132535023232\\
12.9130849589315	-5.28394358415639\\
13.4779116627643	-5.56012964489242\\
14.1216163927159	-5.90729538222396\\
14.9133316589554	-6.36923593647118\\
17.6502352569354	-8.00313444795441\\
18.2732849745949	-8.32379192658647\\
18.8453045146555	-8.58316572485881\\
19.3541025668712	-8.77886556642866\\
19.7905150444897	-8.91439964685582\\
20.1582619689368	-8.99704996428278\\
20.4563825272961	-9.03198367477184\\
20.7122328906492	-9.03283609674317\\
20.9873934743112	-9.00582785198178\\
21.3213530486103	-8.94104290791375\\
21.6634599179654	-8.84485871997447\\
22.0493394809825	-8.70597385204391\\
22.4745577781796	-8.52066296948553\\
22.9332506515232	-8.28723362573071\\
23.4179780800709	-8.00641822318816\\
23.9199007061585	-7.68155194070895\\
24.4293977159096	-7.31826632040275\\
25.008676394149	-6.86529580963397\\
25.5722663225056	-6.38322738448798\\
26.1094000060644	-5.88319350528873\\
26.6123878234984	-5.37439257004002\\
27.0771625599131	-4.86289225541757\\
27.5026259445688	-4.35163434642152\\
27.8897311727184	-3.84084919008815\\
28.2404053833463	-3.32885576354211\\
28.5566183689612	-2.81299883217259\\
28.8397317545552	-2.29061018396386\\
29.0904508794136	-1.75932230389298\\
29.3089578721549	-1.21708685749038\\
29.4734911800653	-0.732481837901631\\
29.6126850506468	-0.237814540493371\\
29.7259651141743	0.26704488595151\\
29.8127173130981	0.781927500082002\\
29.8723196524367	1.30629220681112\\
29.9041942097797	1.83914557099733\\
29.9078765626022	2.37902758840564\\
29.8830725593536	2.92419902382361\\
29.8296736538181	3.47274207742023\\
29.7477723695178	4.02253476684258\\
29.6376706751342	4.57131967892703\\
29.4998634881981	5.11681373343331\\
29.3349862656514	5.65685781438571\\
29.1437859439723	6.18941479155773\\
28.9271027220577	6.71256223565042\\
28.6858219491342	7.22456006808366\\
28.4208270245899	7.72388327212483\\
28.1329727671267	8.20920745909333\\
27.8230832453348	8.67936513588474\\
27.442959744389	9.19678559264955\\
27.0362058013243	9.69168276710302\\
26.6592526232359	10.1053219814351\\
26.2633364771845	10.5001147924816\\
25.849095993978	10.8753904563124\\
25.417170435799	11.2305099627966\\
24.9682070080019	11.5648620780815\\
24.5028726546217	11.8778575378724\\
24.0218657297871	12.1689263329147\\
23.5259331721366	12.4375145154876\\
23.0158657713562	12.6830950976275\\
22.4925045895737	12.9051714169133\\
21.9567281848955	13.1032885601042\\
21.4094589902569	13.277035193571\\
20.8516560305188	13.426050898775\\
20.2025442344542	13.5656832770111\\
19.542485537889	13.6722617444168\\
18.8730111292449	13.7455408075743\\
18.1956576602082	13.7853886611154\\
17.5119538532424	13.7917892200083\\
16.8234081995409	13.7648427695226\\
16.1314991612156	13.7047656110548\\
15.4376546940875	13.6118874565598\\
14.743245987486	13.4866484398283\\
14.0495858962877	13.3295974606497\\
13.3579249429358	13.1413901716174\\
12.6694248813832	12.9227790379481\\
11.9851559308193	12.6746073450848\\
11.3061057026113	12.3978084754212\\
10.6331692279073	12.0933974979032\\
9.88438726242202	11.7192744693831\\
9.14525963164414	11.3131894096601\\
8.41656454256065	10.8768500883486\\
7.69890189545815	10.4120238486968\\
6.99267511923846	9.92050437451224\\
6.22167221481422	9.34528044247606\\
5.46516500519736	8.74186502072534\\
4.72297779551809	8.1127368287018\\
3.92262568398306	7.39392142654347\\
3.13828782770041	6.65015581521147\\
2.71683050410904	6.23506197998693\\
};
\addplot [color=black, line width=0.7pt, forget plot]
  table[row sep=crcr]{%
-1.63750553683739	2.00474688458309\\
-3.10351278436603	3.44959561945274\\
-4.23124541908783	4.52002401660389\\
-5.2401886247016	5.43690959359547\\
-6.17332794409008	6.24479038514883\\
-7.07551341214065	6.98538092652238\\
-7.94539091006117	7.65914978202599\\
-8.78108970140114	8.26743235252739\\
-9.58042163656931	8.81232735282567\\
-10.3921436193238	9.32792059791066\\
-11.1641376064929	9.7824597686469\\
-11.9463337188449	10.2070174705936\\
-12.685036902578	10.5744708337617\\
-13.4313958238615	10.9125070708038\\
-14.184556791351	11.2196137752676\\
-14.9435347814354	11.4943711032672\\
-15.6525203177429	11.7193787553849\\
-16.3645179439448	11.9143736116807\\
-17.078400440861	12.0784718174661\\
-17.7929689778837	12.2108972587315\\
-18.5069455550903	12.3109830815476\\
-19.218976604442	12.3781768237536\\
-19.8732938626976	12.410632768683\\
-20.523624593263	12.4144481735735\\
-21.1688138883006	12.3894665718936\\
-21.8076875336576	12.3356201833339\\
-22.4390604839112	12.252929061244\\
-23.0617493960349	12.1414984473261\\
-23.6745770020152	12.0015155337036\\
-24.2763750073337	11.8332462030583\\
-24.8659860350596	11.6370321049879\\
-25.442265804764	11.4132880800737\\
-26.0040917417308	11.16249702451\\
-26.550364107404	10.8852065755927\\
-27.0800065850534	10.5820262471372\\
-27.5919635270986	10.2536270841462\\
-28.0851996603385	9.90074144144936\\
-28.5587053514585	9.52415952250233\\
-29.0114985241485	9.12472748840212\\
-29.442626175225	8.7033455844746\\
-29.8511554598696	8.26097753495429\\
-30.2361777750984	7.79865124467877\\
-30.5968186697037	7.31745190064889\\
-30.9322438035879	6.81851555425988\\
-31.2416481307018	6.30304659643998\\
-31.5242505566967	5.77233993924121\\
-31.779321629376	5.22775009805913\\
-32.0061927257312	4.67067821564743\\
-32.2042610708754	4.10255780107674\\
-32.3729690168787	3.52491822416612\\
-32.511823165485	2.93937111556416\\
-32.620418601405	2.3475539000276\\
-32.6984455215868	1.75108378056967\\
-32.7456770275076	1.15158965622607\\
-32.7619678191294	0.550803718831624\\
-32.747276600789	-0.0495600540937673\\
-32.7016584935694	-0.647886158489889\\
-32.6252430484559	-1.24265817676064\\
-32.5182529294185	-1.8322449807428\\
-32.3810042064553	-2.4150406146125\\
-32.2138771398777	-2.98957063149287\\
-32.0172829368996	-3.55450441598448\\
-31.7917390783951	-4.10838768766451\\
-31.5378622000358	-4.64974098236871\\
-31.2563297905171	-5.17716564918487\\
-30.9478480875029	-5.68937629164785\\
-30.613212873337	-6.18507494582964\\
-30.2533230795482	-6.66295754673347\\
-29.8691189190168	-7.12181517850922\\
-29.4615332994225	-7.56057109727154\\
-29.0315538675793	-7.97818644613579\\
-28.580291606346	-8.3736053722423\\
-28.1088852700009	-8.74586409486327\\
-27.6184424551362	-9.09412479831886\\
-27.1100700132385	-9.41762909318646\\
-26.5850091308924	-9.71561168003453\\
-26.0445228425178	-9.98739064613199\\
-25.4898257536731	-10.2323972767902\\
-24.9220710362435	-10.4501598156507\\
-24.3425376344354	-10.6402225368229\\
-23.7525339868623	-10.8022020102487\\
-23.1533084701532	-10.93581354954\\
-22.5460021944082	-11.0408635884408\\
-21.931838294362	-11.1172024576161\\
-21.3120847958934	-11.1647478224467\\
-20.6879400030664	-11.1835020878381\\
-20.0604738390858	-11.1735402214957\\
-19.4307655808707	-11.1349953571255\\
-18.7473776043514	-11.0612231776237\\
-18.0640535210255	-10.954504039157\\
-17.381952293743	-10.8152645015595\\
-16.7022069550605	-10.6440133375895\\
-16.0259986234736	-10.4413725313732\\
-15.3543679438264	-10.2080405072171\\
-14.688137611407	-9.94474891149299\\
-14.0281244315918	-9.65232651666869\\
-13.3252157452738	-9.3059138687392\\
-12.6312228874017	-8.92804126180142\\
-11.9466545168124	-8.51995970184277\\
-11.2720435459148	-8.08309462667246\\
-10.5607271722491	-7.58480984017807\\
-9.86137365673387	-7.05707154797217\\
-9.12854357010236	-6.46392168502266\\
-8.40910140608314	-5.84211504242578\\
-7.65863984776689	-5.15287971081638\\
-6.87862958203782	-4.39453080523275\\
-6.07001963717768	-3.56650561202633\\
-5.14969421129202	-2.5783376524946\\
-4.12001617966882	-1.42617644216261\\
-2.73590163204273	0.172244476764511\\
-0.125530339565607	3.19463046126664\\
0.999095641321027	4.44368447249124\\
1.9762289385085	5.48319817112548\\
2.88721995981314	6.40672136133559\\
3.72912683571262	7.21688249880809\\
4.54363725998338	7.95851732479574\\
5.32889535392092	8.63271587968442\\
6.08265501100122	9.24130276713737\\
6.85089393021619	9.82220359315374\\
7.58423839807747	10.3394242964541\\
8.33004407027472	10.8280427698674\\
9.0370002313569	11.2564869626897\\
9.75391840030474	11.6566127556842\\
10.4801057338721	12.0268330438978\\
11.2147305314141	12.3656386979386\\
11.9035734572614	12.6508637717853\\
12.59791621152	12.906691610139\\
13.2967324608822	13.132130998346\\
13.99891237024	13.3262896687235\\
14.7032539902018	13.4883754538682\\
15.4084678784769	13.6177029062159\\
16.1132042807807	13.7137034864059\\
16.7620982660803	13.7723491605201\\
17.4082600807596	13.8019515190126\\
18.0505269845824	13.8023351375171\\
18.6877206167171	13.7734118858148\\
19.3186560091241	13.7151790990723\\
19.9421453308126	13.6277169910557\\
20.5570031100231	13.5111856632888\\
21.162048229718	13.3658223773188\\
21.7561068379367	13.1919390373401\\
22.3380208739819	12.9899178963627\\
22.9066514176425	12.7602075779249\\
23.4608783656637	12.5033202359928\\
23.9995931030318	12.2198331721152\\
24.5217041916458	11.9103862560223\\
25.0261445112793	11.5756773801024\\
25.5118714669162	11.2164608207616\\
25.9778637758993	10.8335485799363\\
26.4231096933009	10.4278219796876\\
26.8466232139235	10.0002226789915\\
27.2474490412875	9.55174960147196\\
27.6246666416179	9.0834545317376\\
27.9773693598336	8.59647030221049\\
28.3046707803218	8.09201745623495\\
28.6057330146689	7.57137650943251\\
28.879774854414	7.03587736561143\\
29.1260725540322	6.48689637114556\\
29.3439375897604	5.92591494910798\\
29.5327465788738	5.35448280700922\\
29.6919622109139	4.7741730587331\\
29.8211359601842	4.18655756415458\\
29.9198862002417	3.59330233571551\\
29.9879239550544	2.99611416120589\\
30.0250673079913	2.39664899705239\\
30.0312284509057	1.79647517345533\\
30.0064012360706	1.19724402594232\\
29.9506757548231	0.600604835280016\\
29.8642239245795	0.00807966673824012\\
29.747270955767	-0.578933601776725\\
29.6001354816465	-1.15889250171636\\
29.4232277468282	-1.73024377802052\\
29.2170231840119	-2.29152314738995\\
28.9820410789071	-2.84136526219216\\
28.7189063717957	-3.37833737441795\\
28.4283318117967	-3.90102416395427\\
28.1110749343988	-4.40811143806727\\
27.7679060217747	-4.89840527165587\\
27.399701293366	-5.37067760968935\\
27.0074278361483	-5.82371801045981\\
26.5920780344201	-6.25642080565824\\
26.1546207713811	-6.66781298335667\\
25.6960993364575	-7.05693833454546\\
25.2176647306746	-7.42284719103362\\
24.7204774684284	-7.76469019546986\\
24.2056472557591	-8.0817443043538\\
23.6743013596354	-8.37334724790291\\
23.1276921918276	-8.63884422763364\\
22.5670731359632	-8.877671322449\\
21.9936291145492	-9.08937629634421\\
21.4084954222462	-9.27358992653112\\
20.8129372147492	-9.42996561313615\\
20.2082203656025	-9.5582377746782\\
19.5955325934052	-9.65823611741624\\
18.9759586651213	-9.72987129981433\\
18.3506771710249	-9.77310173518037\\
17.7208813790253	-9.78795946620336\\
17.0876790975507	-9.77455703227371\\
16.3989916034846	-9.7283537036511\\
15.7087340088216	-9.64950115877677\\
15.0182213677677	-9.53841225471226\\
14.3286249861291	-9.3955892477473\\
13.6409169968373	-9.22159282001969\\
12.9561211782407	-9.01708953911952\\
12.2752306058807	-8.78285087882237\\
11.599055621795	-8.51970759862705\\
10.8768492977325	-8.20499259489683\\
10.1616586145456	-7.85900878009388\\
9.454179818734	-7.48311462092931\\
8.75481504687261	-7.07868897871089\\
8.0148283444249	-6.6153950738365\\
7.28480721926112	-6.12306061781699\\
6.51706404405928	-5.56818261803822\\
5.71333507475475	-4.94778127487943\\
4.8752649039203	-4.25987594426088\\
4.00379566160176	-3.50339205738547\\
3.05434791038679	-2.63670923187266\\
1.98352358354477	-1.61516364999156\\
0.659266166754939	-0.306125436225926\\
-1.63750553683739	2.00474688458309\\
};
\addplot [color=black, dotted, line width=1.0pt, forget plot]
table[row sep=crcr]{%
	-8.70950755333683	3.15104586691227\\
	-8.99081443681792	3.39893637173331\\
	-9.29468466646252	3.68863496266096\\
	-9.64461063307267	4.04844216976358\\
	-10.0262133731352	4.46907627986654\\
	-10.4039547384158	4.91033168515069\\
	-10.8540275209527	5.46330801052976\\
	-11.3645858973611	6.11860144262673\\
	-12.1212590441379	7.12087564925235\\
	-13.0504961506458	8.34804645591824\\
	-13.4891552236363	8.90036547417132\\
	-13.8710874619555	9.35576098391661\\
	-14.202212904222	9.72650070929019\\
	-14.4960253038462	10.0341084979534\\
	-14.7970846693872	10.3271902034031\\
	-15.0991522292649	10.5982690945158\\
	-15.400635127777	10.8458837438451\\
	-15.6998932343609	11.0694964521437\\
	-15.9949868566592	11.2693036164172\\
	-16.2843617893491	11.4466132735636\\
	-16.6300995310734	11.6363523188595\\
	-16.9829441456598	11.8073869671756\\
	-17.3435632791833	11.9605713626683\\
	-17.7111387900236	12.0959831994765\\
	-18.0887858721306	12.215027740678\\
	-18.4786150574081	12.318199948763\\
	-18.8816812011749	12.4053354323853\\
	-19.298101115846	12.4758750736542\\
	-19.7272672214322	12.5290873542419\\
	-20.1680369538518	12.5642199676429\\
	-20.6188811319997	12.5805963427292\\
	-21.0780564088074	12.5776733968317\\
	-21.5437036187708	12.5550631343268\\
	-22.0139326204286	12.5125357914694\\
	-22.4869346541857	12.4500024841183\\
	-22.9609517059247	12.367502102184\\
	-23.4343687642001	12.2651711483024\\
	-23.9056806469126	12.1432279538894\\
	-24.3735085667297	12.0019514713329\\
	-24.8365848061466	11.8416683546309\\
	-25.2937338267413	11.6627459471981\\
	-25.7438736766241	11.4655817197494\\
	-26.1860117141553	11.2505945653894\\
	-26.6192329632373	11.0182207234857\\
	-27.0426755415452	10.7689203091397\\
	-27.4555528475178	10.5031594382816\\
	-27.8571133213389	10.2214294391805\\
	-28.2466646384406	9.92422847875064\\
	-28.6235622167413	9.61206395705413\\
	-28.9871959335249	9.2854597313036\\
	-29.3369980815113	8.94494651310298\\
	-29.6724348225996	8.59106559953833\\
	-29.99299557542	8.22437783710954\\
	-30.2982000325078	7.84545515346239\\
	-30.5875917722696	7.45488641452575\\
	-30.8607434768909	7.05326893846707\\
	-31.1172521858519	6.6412114648554\\
	-31.3567247918956	6.21935721295846\\
	-31.5787867851801	5.78837523627481\\
	-31.7830894313496	5.34894948929394\\
	-31.9693078727316	4.90178279524962\\
	-32.1371467434949	4.44758242996219\\
	-32.2863320103629	3.98707716052505\\
	-32.4166042360375	3.52104704452887\\
	-32.5277328171038	3.05029665510956\\
	-32.6195216846512	2.57564117799087\\
	-32.6918117722706	2.0978928047505\\
	-32.7444804854837	1.61784306761687\\
	-32.7774326500967	1.13632080383697\\
	-32.7906070538213	0.654198618665852\\
	-32.7839845711382	0.172341881761881\\
	-32.7575864369861	-0.308418042234145\\
	-32.711465773291	-0.787310736245551\\
	-32.6457032825087	-1.26357510043036\\
	-32.5604215326323	-1.73637698822146\\
	-32.4557815571173	-2.20489646954698\\
	-32.3319737478751	-2.66836224759328\\
	-32.18920215702	-3.12607698232796\\
	-32.027683377042	-3.57737897795664\\
	-31.8476916659263	-4.0215150360013\\
	-31.649541631403	-4.45772641835226\\
	-31.433574942618	-4.88528530290233\\
	-31.2001488129857	-5.30350990111936\\
	-30.9496277708433	-5.71176483732879\\
	-30.6824267739171	-6.10938474149591\\
	-30.3990058940405	-6.49569887193793\\
	-30.0998446968384	-6.87007379419866\\
	-29.7854272519612	-7.23192682465903\\
	-29.4562252692446	-7.58073378021655\\
	-29.1127473374282	-7.91596718823075\\
	-28.7555634704323	-8.23708206699681\\
	-28.3852577496448	-8.54357009471876\\
	-28.0024107988819	-8.83497129872559\\
	-27.6075786086603	-9.11088186673047\\
	-27.2013227754008	-9.3709228732389\\
	-26.7842801789996	-9.61469852469924\\
	-26.3571012056212	-9.84184532575948\\
	-25.9204259432712	-10.0520443523233\\
	-25.4748622311307	-10.2450259470862\\
	-25.0209878972492	-10.4205590674823\\
	-24.5594557441931	-10.5784082997754\\
	-24.0909417264527	-10.7183640079754\\
	-23.6161061987658	-10.8402566799417\\
	-23.1355728143034	-10.9439583563909\\
	-22.6499104292837	-11.029378160626\\
	-22.1597359922448	-11.0964385054344\\
	-21.6657072839125	-11.1450852545676\\
	-21.1684606089111	-11.1752999966404\\
	-20.6685934759849	-11.1870983762091\\
	-20.166640145451	-11.1805241905866\\
	-19.6631368709629	-11.155643096966\\
	-19.1586767429133	-11.1125496970189\\
	-18.6538329407834	-11.0513701126216\\
	-18.1491365441053	-10.9722583765207\\
	-17.6450516310303	-10.8753874820281\\
	-17.1420044555273	-10.7609472596709\\
	-16.6404811788283	-10.6291670404318\\
	-16.14095312401	-10.4803066052334\\
	-15.6438465420067	-10.3146481906418\\
	-15.1495198899721	-10.132485094699\\
	-14.6582673503103	-9.93411543730638\\
	-14.1704293353842	-9.71988421707092\\
	-13.6863400725565	-9.49017035219838\\
	-13.2062861366826	-9.24537003795211\\
	-12.7304864372854	-8.98588354648344\\
	-12.1809563452161	-8.66511453219736\\
	-11.6376845236844	-8.32559941387711\\
	-11.1009245179067	-7.96810793747242\\
	-10.5708309275019	-7.59342523641646\\
	-10.3456968256957	-7.42777423604222\\
};

\addplot [color=red, dotted, line width=0.8pt, forget plot]
table[row sep=crcr]{%
	-8.70950755333683	3.15104586691227\\
	-9.09605742962898	3.4968849206668\\
	-9.48050220921502	3.87875537900857\\
	-9.88312778427456	4.32013946875811\\
	-10.2850430283832	4.80515493163228\\
	-10.6754618356643	5.3213562000382\\
	-11.1010217114612	5.93696821596563\\
	-11.5129014937766	6.5862718395061\\
	-11.9564707051295	7.34205971494443\\
	-12.4738884879474	8.28626469394414\\
	-14.0926228988396	11.302345705088\\
	-14.4847863478391	11.9418768323874\\
	-14.8626306514826	12.4965692154491\\
	-15.2133834911363	12.9550443099493\\
	-15.5709848803084	13.36864209565\\
	-15.931375783168	13.7346750710659\\
	-16.2891438324973	14.0532379109964\\
	-16.6932695840936	14.3677988048192\\
	-17.1151015849485	14.6543579225505\\
	-17.5823051258108	14.9307590428063\\
	-18.0399981223346	15.1652273572209\\
	-18.5023932863745	15.3678861112881\\
	-18.945536864337	15.5307604017565\\
	-19.3881168166815	15.6643205019091\\
	-19.84276357617	15.7730298847794\\
	-20.3071993467091	15.8556870759875\\
	-20.7774733204943	15.9111289615752\\
	-21.2509111862666	15.9388484045662\\
	-21.7254316367599	15.9386470711999\\
	-22.1992115679903	15.9104964680021\\
	-22.6705379541544	15.8544877665783\\
	-23.1377495675686	15.7708149256301\\
	-23.5972449608416	15.6603067332646\\
	-24.101386887694	15.5056677812716\\
	-24.6003377335367	15.3174197803632\\
	-25.1012966399517	15.0925146183271\\
	-25.6020544435866	14.8308265382144\\
	-26.1004466537276	14.5325488630726\\
	-26.5925350988751	14.1995787553531\\
	-27.0806065160976	13.8300533777561\\
	-27.5648858475526	13.4230388437299\\
	-28.0439352442121	12.9788199008819\\
	-28.5151746873837	12.4989399651375\\
	-28.9752670322993	11.9861948209001\\
	-29.4211531767191	11.4436755019661\\
	-29.8499722534626	10.8747503692531\\
	-30.2592252338572	10.2827503335489\\
	-30.6468170955601	9.67074103994985\\
	-31.0107923423377	9.0417976089908\\
	-31.3496661728542	8.39836723727596\\
	-31.6620743380451	7.7427475469637\\
	-31.9468257816442	7.07696092550136\\
	-32.2028228114807	6.40291105097619\\
	-32.4291026664716	5.72229537224768\\
	-32.6247637476883	5.03683738717957\\
	-32.7890406552034	4.34811969916053\\
	-32.9212721433093	3.65770477947326\\
	-33.0209085310332	2.96715498202698\\
	-33.0875267209157	2.27799917688053\\
	-33.1208306072879	1.59172754670807\\
	-33.1206443872633	0.909760232396884\\
	-33.0869055918214	0.233612254324164\\
	-33.0196672329044	-0.435317620718415\\
	-32.9190738346407	-1.09573149981137\\
	-32.7853592149658	-1.74630207959979\\
	-32.6188661916982	-2.38561241147058\\
	-32.4200096237433	-3.01234841866035\\
	-32.1892334449941	-3.6253366440241\\
	-31.9270969598738	-4.2232475271312\\
	-31.6342764174799	-4.80467277761199\\
	-31.3115207380337	-5.3682535399306\\
	-30.9596127554198	-5.91272587674304\\
	-30.6284986382959	-6.37240786579373\\
	-30.2764555639684	-6.81552502884125\\
	-29.8494157219513	-7.3005746904629\\
	-29.3970323538939	-7.76175549989632\\
	-28.9204759527033	-8.19792459401881\\
	-28.4210567084276	-8.60793359968627\\
	-27.9000922742083	-8.99076137792013\\
	-27.3588487669249	-9.34553203412811\\
	-26.7987097800284	-9.67138167751061\\
	-26.2211664208657	-9.96749097846771\\
	-25.6276824217212	-10.2331695724094\\
	-25.0196207043879	-10.4678656386678\\
	-24.398453992495	-10.6710626226792\\
	-23.7657412294505	-10.8423120324517\\
	-23.1229731630658	-10.9812881546348\\
	-22.471501881976	-11.0877805692235\\
	-21.8127855582111	-11.1616344464024\\
	-21.1483277129471	-11.2027824211818\\
	-20.4795262340312	-11.2112628484747\\
	-19.8076221094258	-11.1871989377163\\
	-19.1339400205001	-11.1307907606049\\
	-18.4598211491442	-11.0423336736915\\
	-17.7864675173414	-10.9222060337181\\
	-17.1149030789531	-10.7708409427139\\
	-16.4462285313085	-10.588769290457\\
	-15.7815195678879	-10.3766179463157\\
	-15.4509453514351	-10.2594751699117\\
};

\end{axis}
\end{tikzpicture}%

%% file: path_b23_eight_rev_MIQP.tex
%
%
\definecolor{mycolor1}{rgb}{0.50000,0.00000,0.50000}%
\begin{tikzpicture}

\begin{axis}[%
width=\figurewidth,
height=\figureheight,
at={(0\figurewidth,0\figureheight)},
scale only axis,
xmin=-1.4,
xmax=1.4,
xtick={-0.8,0,0.8},
xlabel={$\beta_3$ [rad]},
ymin=-1.4,
ymax=1.4,
ytick={-0.8,0,0.8},
xlabel style={font=\color{white!15!black},at={(axis description cs:0.5,-0.14)},anchor=north},
ylabel style={font=\color{white!15!black},at={(axis description cs:-0.2,.5)},anchor=south},
ylabel={$\beta_2$ [rad]},
axis background/.style={fill=white},
xmajorgrids,
ymajorgrids
]

\addplot[area legend, draw=gray, fill=gray, fill opacity=0.3, forget plot]
table[row sep=crcr] {%
x	y\\
-0.98	0.65\\
-0.75	-0.65\\
0.98	-0.65\\
0.75	0.65\\
-0.98	0.65\\
}--cycle;

\addplot[area legend, draw=gray, fill=green, fill opacity=0.3, forget plot]
table[row sep=crcr] {%
x	y\\
-0.72	-0.83\\
0.9	-0.43\\
0.72	0.83\\
-0.9	0.43\\
-0.72	-0.83\\
}--cycle;

\addplot[area legend, draw=gray, fill=gray, fill opacity=0.3, forget plot]
table[row sep=crcr] {%
x	y\\
-0.75	-1.15\\
0	-0.65\\
0.75	0.65\\
0.75	1.15\\
0	0.65\\
-0.75	-0.65\\
-0.75	-1.15\\
0	-0.65\\
}--cycle;

\addplot [color=green, line width=0.8pt, forget plot]
table[row sep=crcr]{%
	0.727830156269241	0.319586716944193\\
	0.71725399234584	0.358589225955401\\
	0.70554533957607	0.39841628094559\\
	0.692698031649021	0.439055830382303\\
	0.6787105755753	0.480492406975362\\
	0.66358660823495	0.522706963718091\\
	0.629971970880887	0.609375127556064\\
	0.592001930846463	0.698831869871196\\
	0.550212963458651	0.789846919498193\\
	0.528541929200527	0.834018816899305\\
	0.506490180793211	0.87700747461737\\
	0.484098740008008	0.918790500354168\\
	0.461405096492763	0.959348629290502\\
	0.438443203787778	0.998665559569136\\
	0.415243514859805	1.0367277712241\\
	0.39183305226277	1.0735243327219\\
	0.368235507919143	1.10904669893273\\
	0.344471367554335	1.14328850395284\\
	0.320558054984802	1.17624535177234\\
	0.296510091733638	1.20791460734581\\
	0.272339267800566	1.23829519018922\\
	0.248054819822897	1.26738737220934\\
	0.223663613308788	1.29519258107862\\
	0.199170326085541	1.32171321010659\\
	0.174577630568179	1.34695243523119\\
	0.149886372904842	1.37091403946181\\
	0.125095747486478	1.39360224485104\\
	0.100203465712317	1.41502155185127\\
};

\addplot [color=red, line width=0.8pt, forget plot]
table[row sep=crcr]{%
	0.727830156269241	0.319586716944193\\
	0.710753967020006	0.383229857988313\\
	0.692161342186626	0.443157319205843\\
	0.672593672561276	0.498307766484651\\
	0.652288040796649	0.548555743620674\\
	0.631463505583503	0.593798982765052\\
	0.610325911946884	0.633948380126083\\
	0.589071519173149	0.668923705655511\\
	0.567890727397455	0.698650450955035\\
	0.557378252740279	0.711551197003844\\
	0.546973787348848	0.723051112224241\\
	0.53667385380451	0.733232949478233\\
	0.526527973576707	0.742004498260223\\
	0.516541929630102	0.74942179947513\\
	0.506752655962879	0.755437374936629\\
	0.497179275817615	0.760067536583689\\
	0.487834787555999	0.763348025922842\\
	0.478770421389761	0.765191630018651\\
	0.469964407450219	0.765757057642353\\
	0.46148392159747	0.764904625471898\\
	0.436596576690512	0.759312520598148\\
	-0.609478701096805	0.500784965932137\\
	-0.623906409091647	0.494810298817594\\
	-0.637926288764083	0.486994665888755\\
	-0.65147171224241	0.477234619321531\\
	-0.664520541056981	0.465625268465751\\
	-0.676958363860719	0.451894454750697\\
	-0.688743131434138	0.436098221840126\\
	-0.699766006556452	0.418036406287263\\
	-0.719929027651947	0.377695869971108\\
	-0.738206856984308	0.335025368476128\\
	-0.754405437231901	0.289945616083157\\
	-0.768339610671809	0.242443533573027\\
	-0.780135929443664	0.193577744467542\\
	-0.789568270695005	0.14297000837988\\
	-0.799079339691499	0.0641339224617632\\
	-0.80227387183806	0.00999301510643213\\
	-0.803202535471514	-0.0442827448172681\\
	-0.802258074829277	-0.0973218958380951\\
	-0.800091648360854	-0.147392584497877\\
	-0.794506549323384	-0.233246790891414\\
	-0.785567813799615	-0.345743575708506\\
	-0.771256430350331	-0.462725504316287\\
	-0.769479109492471	-0.468118422816066\\
	-0.761401207647457	-0.488105973116323\\
	-0.751590275356006	-0.508438936084575\\
	-0.73965974472736	-0.529208388295086\\
	-0.729691978997809	-0.543851408511362\\
	-0.719108444535997	-0.557065112757998\\
	-0.708046158053127	-0.568628938725759\\
	-0.696720477489625	-0.578183190914573\\
	-0.685215732187992	-0.585772556957823\\
	-0.673659623354816	-0.591355834461071\\
	-0.662197106700354	-0.594819514045807\\
	-0.650830561491771	-0.596516586912196\\
	-0.639682514002526	-0.596397684115998\\
	-0.628822317735561	-0.594567223639304\\
	-0.618265375509928	-0.591282108497974\\
	-0.608065827139973	-0.586647570995587\\
	-0.598269418843095	-0.580780714660388\\
	-0.585899485842054	-0.571255497988138\\
	-0.574395834033128	-0.560027347886376\\
	-0.563728287742774	-0.547642663616645\\
	-0.551683884770884	-0.530726380991885\\
	-0.543019485834479	-0.516580636305758\\
	-0.531788342069732	-0.494443932331455\\
	-0.522525014863194	-0.47216199535051\\
	-0.515177810304064	-0.450212978183398\\
	-0.50881664565868	-0.425828263417732\\
	-0.504681126075716	-0.403120088411913\\
	-0.502520624789099	-0.382535529392716\\
	-0.502165642584787	-0.361890331522502\\
	-0.503629230858703	-0.344702713942011\\
	-0.50614027614288	-0.332399969423572\\
	-0.509556444985962	-0.322565495176381\\
	-0.513631927533043	-0.315145430534104\\
	-0.518234497005263	-0.309755023793356\\
	-0.523916440190932	-0.305876268893751\\
	-0.52983166589574	-0.304154936205063\\
	-0.535806538826346	-0.304236388885137\\
	-0.543214908672155	-0.306183903136532\\
	-0.551878674418645	-0.310292754449399\\
	-0.56130261742213	-0.316959324923232\\
	-0.572166177606968	-0.326653238977167\\
	-0.579731219425782	-0.335604152618678\\
	-0.587873359290003	-0.347848609531143\\
	-0.592862074212407	-0.359309076856625\\
	-0.595469724481263	-0.368196094054055\\
	-0.595868383726694	-0.375106150020284\\
	-0.594357833865465	-0.380126750004087\\
	-0.591020699103667	-0.384549213999533\\
	-0.586228777737009	-0.388204424349342\\
	-0.577265126823078	-0.39223264358345\\
	-0.568769795619375	-0.394520104542376\\
	-0.551175015388858	-0.397306011953089\\
	-0.530418051475847	-0.398261720159319\\
	-0.512239977103418	-0.397940791656398\\
	-0.478417271442365	-0.395253689882893\\
	-0.446492638140234	-0.390471862178298\\
	-0.41958376718348	-0.385291368137108\\
	-0.362962603995708	-0.370761141027002\\
};

\addplot [color=blue, line width=0.8pt, forget plot]
  table[row sep=crcr]{%
0.727830156269241	0.319586716944193\\
0.708310640072438	0.390333824915359\\
0.686896432730548	0.457560187914684\\
0.664133326900392	0.520214237645116\\
0.640284750121428	0.578103847755773\\
0.615584106523079	0.63109019913128\\
0.590216086887905	0.679141661586721\\
0.564357634482674	0.722207277135143\\
0.538176726684821	0.760231004275064\\
0.525013555524548	0.777341670618099\\
0.51183329577909	0.793157767659514\\
0.498645662263721	0.80770498414055\\
0.485483580699918	0.82092933537105\\
0.472348985369059	0.832884032737066\\
0.459287073387616	0.843471632141077\\
0.44629797112426	0.852755809593616\\
0.433422572371539	0.860655616684822\\
0.420669339582074	0.867209866763772\\
0.408072851453236	0.872364188314179\\
0.395653776091271	0.876119638347874\\
0.383435901006447	0.878467894194826\\
0.371459857657622	0.879344721024164\\
0.359727208189088	0.878833463951755\\
0.348327786140605	0.876695350409636\\
0.337241326945874	0.873115111605691\\
0.326571635293524	0.867814617721359\\
0.295648076080698	0.847954826657343\\
0.0440879617214684	0.679537858256379\\
0.0246751881767093	0.664426107035423\\
0.0146013889626706	0.65838410584574\\
0.00424585757003837	0.653372812609511\\
-0.00657763224857455	0.650003876600954\\
-0.0179141410604376	0.648319803654509\\
-0.0541395987173183	0.648902779648706\\
-0.4570641624876	0.649516701689282\\
-0.473624163486398	0.648604279094091\\
-0.490004456980329	0.646178084860623\\
-0.506145285514831	0.642084848688312\\
-0.521999783542135	0.636250972949279\\
-0.537539876117584	0.62870210477943\\
-0.552679744880279	0.619216284294962\\
-0.567396370312526	0.607884243145511\\
-0.581591573725263	0.594469125407653\\
-0.608796670425856	0.563967153606023\\
-0.634980421044265	0.531738428308122\\
-0.660010999296365	0.497721767967588\\
-0.683735109626314	0.461806759687327\\
-0.705986891743708	0.423877538222866\\
-0.726587616800532	0.383813134996361\\
-0.745345573164164	0.341488045397132\\
-0.762056210336529	0.296773069723657\\
-0.776502626766942	0.249536491463107\\
-0.782805812725526	0.22493113267988\\
-0.788456507362374	0.199645664196037\\
-0.793424636097935	0.173663699689808\\
-0.797679639046926	0.146969084041338\\
-0.801190560953372	0.119545990178424\\
-0.803926156728875	0.091379026478751\\
-0.805855014163479	0.0624533551923963\\
-0.806945695422726	0.0327548222516582\\
-0.806957789599287	-0.0275231278810089\\
-0.804442480356314	-0.11111300453169\\
-0.800391165198442	-0.183845804856715\\
-0.79447519605118	-0.26248138922705\\
-0.789352663667793	-0.31160237656601\\
-0.782818109025195	-0.356103079219426\\
-0.775241010171004	-0.395671616500294\\
-0.766882683418035	-0.43035492286164\\
-0.75798144161685	-0.460250398503313\\
-0.748697235588962	-0.485662727564339\\
-0.739187015177649	-0.506853867917701\\
-0.72956230357242	-0.524177926554486\\
-0.7199269999963	-0.537959651817171\\
-0.710367288512099	-0.548525744911308\\
-0.700953792955068	-0.556202060957375\\
-0.691748373371659	-0.561293517408744\\
-0.682800346425521	-0.564095997106199\\
-0.674157182523923	-0.56486773485107\\
-0.66585546739365	-0.56386048193663\\
-0.657921992745713	-0.561315813233238\\
-0.650380803699901	-0.557446406003986\\
-0.640965349506088	-0.550568295184436\\
-0.632303635797702	-0.542116904973678\\
-0.622551881443726	-0.5299186336836\\
-0.61398427886638	-0.51651365851495\\
-0.605218034908516	-0.499600755673999\\
-0.596949958142068	-0.47969926404483\\
-0.589845413872723	-0.457816923148599\\
-0.584439790503406	-0.435388191089537\\
-0.581075656486411	-0.414021035360277\\
-0.579828938015601	-0.395303632001386\\
-0.580486085957971	-0.380529847852964\\
-0.582395925497899	-0.371016723264969\\
-0.585090789967649	-0.365536411676968\\
-0.587601369158814	-0.363884283402521\\
-0.590127104715684	-0.364691190458506\\
-0.59612542927393	-0.374233022848456\\
-0.601780456994572	-0.373649052636217\\
-0.604486230297608	-0.376231190946825\\
-0.605311441103567	-0.381269112346022\\
-0.604085915712375	-0.384903253026174\\
-0.600420462865743	-0.388848319996771\\
-0.592612912358591	-0.393149719470166\\
-0.58156883692057	-0.396522079150992\\
-0.57144865887798	-0.398406864817504\\
-0.54959639620496	-0.400448761303941\\
-0.520537091473998	-0.400390231842122\\
-0.489337287329745	-0.397937434285968\\
-0.466709311633394	-0.394919727746075\\
-0.427202281025939	-0.38788418287424\\
-0.392360954326601	-0.37951997696987\\
-0.36067980981327	-0.370673615237071\\
-0.30080574050209	-0.350118611780636\\
};

\addplot[->, thick, color=blue] coordinates
{ (0.6402,0.5781) (0.6155,0.6310)};
\addplot[->, thick, color=blue] coordinates
{ (0.2956,0.8479) (0.0440,0.6795)};
\addplot[->, thick, color=blue] coordinates
{ (-0.7265,0.3838) (-0.7453,0.3414)};
\addplot[->, thick, color=green] coordinates
{ (0.3682,1.1090) (0.34447,1.1432)};
\addplot[->, thick, color=red] coordinates
{ (-0.365,0.57) (-0.4094,0.5557)};

\addplot [color=black, dotted, line width=1.0pt, forget plot]
table[row sep=crcr]{%
	0.727830156269241	0.319586716944193\\
	0.721223159891246	0.346660789344477\\
	0.713490192415232	0.375357705096438\\
	0.704308614896605	0.406525418728891\\
	0.69389429776452	0.439346295235609\\
	0.682869444002935	0.471849454230265\\
	0.671517494270049	0.503216781371337\\
	0.659872231882201	0.533420669848444\\
	0.647969363107463	0.562425635912161\\
	0.635832449762613	0.59023138112775\\
	0.623494760726583	0.616804616566067\\
	0.61096976830894	0.642171370305781\\
	0.598293235175096	0.666284282889654\\
	0.585474495836368	0.689178037207199\\
	0.572549787968276	0.710798292606667\\
	0.559526330251228	0.731183239006137\\
	0.54643776703684	0.750282383969468\\
	0.53329210124844	0.768128948615313\\
	0.520119844971405	0.784679422368344\\
	0.506931828603803	0.799956288906277\\
	0.493754343941771	0.813928449191963\\
	0.48060348584229	0.826600275089271\\
	0.467499699358075	0.837960883035201\\
	0.454467489336695	0.847986381197805\\
	0.441519151762608	0.856696820143155\\
	0.428691238022605	0.864027862836546\\
	0.415985245681	0.870042886623925\\
	0.403460337625507	0.874599044245226\\
	0.391114056729188	0.877784572140452\\
	0.379000733117948	0.879480725491341\\
	0.367130634153747	0.879739400393108\\
	0.355547599156475	0.878490417913693\\
	0.344281292876005	0.875726265195084\\
	0.333344109346537	0.871505110091063\\
	0.322818100628672	0.865635446680156\\
	0.291975434955248	0.845530564137837\\
	0.150405505617793	0.750696547836934\\
	0.0605367268321255	0.690128509674386\\
	0.0509692379254716	0.681982884211041\\
	0.0419341767149212	0.671986159696694\\
	0.025281736815858	0.647804863321013\\
	0.0167840309373611	0.637092145408007\\
	0.00798220535234706	0.627874475262268\\
	-0.00116968258376482	0.620195821768735\\
	-0.0107177686944101	0.614100028964259\\
	-0.0207082696767947	0.609630801415799\\
	-0.0311873033020575	0.606831531689911\\
	-0.0422006946693788	0.605745195790792\\
	-0.0537937496000551	0.606414259001245\\
	-0.0660109716102359	0.60888052337389\\
	-0.078895723212403	0.613184943268\\
	-0.0924898245798517	0.619367404480068\\
	-0.106707269413893	0.626971812871356\\
	-0.121036388175882	0.633923974162151\\
	-0.135326148235975	0.639669456584338\\
	-0.149482619325993	0.643904151941655\\
	-0.163535951776159	0.646830566314087\\
	-0.177454877266749	0.648400129272593\\
	-0.205207256851336	0.649083219094282\\
	-0.322042202586814	0.649377479261473\\
	-0.449420000447544	0.649522596503176\\
	-0.465971304335404	0.648932754369924\\
	-0.482452204685425	0.647293726485416\\
	-0.498740288895858	0.644135982585849\\
	-0.514758638079973	0.639234118238479\\
	-0.530484797238388	0.632633927828707\\
	-0.545843036569801	0.624150220214879\\
	-0.56080275823398	0.613829290829977\\
	-0.575313084394155	0.601634122240806\\
	-0.589284767006482	0.587379011084452\\
	-0.616275559599864	0.556441435053036\\
	-0.642221997069286	0.523796245126252\\
	-0.654759041810304	0.506800348813981\\
	-0.666980036135016	0.48933730796588\\
	-0.678865058061389	0.471392981218078\\
	-0.690393393817893	0.452952889151402\\
	-0.701543528041324	0.434002224899027\\
	-0.71229313665429	0.41452586793596\\
	-0.722619082892612	0.394508401448702\\
	-0.732497417011281	0.373934133716559\\
	-0.741903380261613	0.352787123969573\\
	-0.750811413802086	0.331051213221052\\
	-0.759195173280898	0.308710060604566\\
	-0.767027549909498	0.285747185774856\\
	-0.774280698932806	0.262146017958204\\
	-0.780926076492806	0.237889952258869\\
	-0.786934485976608	0.21296241384238\\
	-0.792276135036389	0.187346930621748\\
	-0.796920704564579	0.16102721506645\\
	-0.800837431000472	0.133987255733747\\
	-0.803995203430235	0.106211419084226\\
	-0.80636267701636	0.0776845620852654\\
	-0.807908404349038	0.0483921560235657\\
	-0.808751466689293	0.0188081075325347\\
	-0.808955021194553	-0.0391253566227041\\
	-0.807776007045768	-0.0938212451104837\\
	-0.805788797319617	-0.144044328211714\\
	-0.801997436035749	-0.210338581026223\\
	-0.796779786654198	-0.281383972472154\\
	-0.791336063614775	-0.345423153701312\\
	-0.788622505514108	-0.367530342924158\\
	-0.785170539340221	-0.389128002051651\\
	-0.781083597090557	-0.409942298309899\\
	-0.77646460810693	-0.429737262835808\\
	-0.771403331716778	-0.448346527341341\\
	-0.765973406839527	-0.465670946260252\\
	-0.760245126847997	-0.481634582626914\\
	-0.75427992091796	-0.496198303725568\\
	-0.748128787317572	-0.509358860858987\\
	-0.741839598571621	-0.521124225002644\\
	-0.735454844154064	-0.531519079134193\\
	-0.729010514878069	-0.540585753514286\\
	-0.722540002290832	-0.548371006662481\\
	-0.716073161262069	-0.554928453441264\\
	-0.709634373501886	-0.560323209384364\\
	-0.703247672884428	-0.564615357370211\\
	-0.696931721896413	-0.567875546909615\\
	-0.690706916728413	-0.570162840537608\\
	-0.684590414128927	-0.57153980931766\\
	-0.678597354733001	-0.57206994148935\\
	-0.672741320278283	-0.571814156563755\\
	-0.667030445659104	-0.570844185022413\\
	-0.658761649614542	-0.568166749088571\\
	-0.650873155906498	-0.564194883610343\\
	-0.643383379509692	-0.559117525028426\\
	-0.636312227888569	-0.553091669041853\\
	-0.629665402019464	-0.546290790763921\\
	-0.621475107488135	-0.536255064829775\\
	-0.6140501989442	-0.525411424111751\\
	-0.607374677934048	-0.514053034862766\\
	-0.600068562691574	-0.499453540667026\\
	-0.593845236191557	-0.484837896040256\\
	-0.58771902207648	-0.467699005314721\\
	-0.582250311252079	-0.448793827488833\\
	-0.578372235979865	-0.431501557253376\\
	-0.575600216040158	-0.414065906297535\\
	-0.574269715151217	-0.399315860346915\\
	-0.574103040291066	-0.385940872485192\\
	-0.5749512333055	-0.375810446005963\\
	-0.576479470524324	-0.368560647151108\\
	-0.578379693031484	-0.363847011947168\\
	-0.580643694972712	-0.361053391596666\\
	-0.582630358761172	-0.36025017894981\\
	-0.584836009878249	-0.360801644973986\\
	-0.587002195322327	-0.362763965049181\\
	-0.593218727146735	-0.370197609141736\\
	-0.600648140789853	-0.371744512751822\\
	-0.602889187834339	-0.374009291645988\\
	-0.603906250833447	-0.377214650601458\\
	-0.60421181551862	-0.380956404091247\\
	-0.603225448208874	-0.38397954150486\\
	-0.600787403651663	-0.387182339089964\\
	-0.597213572359225	-0.389778895688952\\
	-0.591216190797464	-0.392905440880539\\
	-0.581984531974505	-0.395808202711045\\
	-0.573999199456596	-0.397570793285246\\
	-0.5567017365027	-0.399735883133277\\
	-0.544310889006168	-0.400237283065348\\
	-0.519809606348472	-0.400113453389291\\
	-0.495582121420918	-0.398313246947115\\
	-0.476723696186498	-0.396258182193017\\
	-0.449384655835219	-0.39192923753185\\
	-0.426263353905177	-0.387564357425369\\
	-0.399418200852109	-0.381205491618003\\
	-0.368282076742563	-0.372898265823261\\
	-0.341847390036734	-0.364510492714946\\
	-0.30687750505117	-0.352344662077864\\
};

\addplot [color=black, line width=1.0pt, forget plot]
  table[row sep=crcr]{%
-0.0408598394484972	0.13466520638993\\
-0.142880355317615	0.0425116099110959\\
-0.165790007183586	0.0243363978601857\\
-0.207792446652068	-0.00619805834531528\\
-0.379727274189688	-0.122585643275196\\
-0.421790650072576	-0.153997633105531\\
-0.454559715096812	-0.180436049038848\\
-0.483834279693127	-0.205985574362196\\
-0.509368282939933	-0.230408118394342\\
-0.531115334474125	-0.253330947980684\\
-0.549797922321808	-0.275364585524183\\
-0.565273302455181	-0.296096826582156\\
-0.577507091598239	-0.315180439378931\\
-0.586604540196093	-0.332242623443335\\
-0.59319650337228	-0.347968746384261\\
-0.597203963419761	-0.361807404254658\\
-0.598738317277314	-0.372796482947882\\
-0.598487369328993	-0.379421503086881\\
-0.601013674108382	-0.376391595304695\\
-0.603298229778452	-0.375956963721876\\
-0.605382958825144	-0.377548357263932\\
-0.606062065673226	-0.381878150143506\\
-0.604756500932058	-0.385467104119656\\
-0.600950967176329	-0.389347319813641\\
-0.593036813203388	-0.393550839721977\\
-0.581939992226893	-0.396853746022613\\
-0.57146846982738	-0.398734496993717\\
-0.549923990652537	-0.40067549765775\\
-0.520712345713422	-0.400568426157832\\
-0.48894538650299	-0.39803760523401\\
-0.465253652041426	-0.394800329267211\\
-0.426403064351405	-0.387791457136648\\
-0.391196145447763	-0.379302274861166\\
-0.359422232257052	-0.370335598480062\\
-0.319700204240819	-0.357065001344296\\
-0.288124569844594	-0.345153565668635\\
-0.247324071292396	-0.327760832059443\\
-0.214828939919369	-0.312422909768168\\
-0.176746911581933	-0.2923809485882\\
-0.14320955034361	-0.273162188194853\\
-0.108049113236798	-0.250945858736107\\
-0.0731543065995718	-0.227199976915445\\
-0.0390435046810411	-0.201710425738353\\
-0.00278326072706225	-0.172752919736408\\
0.034772879575632	-0.14018062297753\\
0.12222101782624	-0.0602793486005441\\
0.147837350724528	-0.0385532467579188\\
0.178487104403983	-0.0151200011325642\\
0.212810258784891	0.00964261761365393\\
0.344679669943556	0.0978743471052701\\
0.38580933749421	0.126997562165325\\
0.426399826572588	0.157575129460162\\
0.459944653744835	0.184945781978777\\
0.487766468086906	0.209596864355854\\
0.512879924711121	0.233950222808982\\
0.534165750290406	0.256779768610305\\
0.552351719063648	0.278652838292397\\
0.567314517920431	0.299164072917108\\
0.578637708782478	0.317295628238769\\
0.587640683137094	0.334696596103494\\
0.593750547823124	0.350025366254316\\
0.597413298595257	0.363672869514385\\
0.598389048382896	0.37919286264729\\
0.601347697055104	0.376059034750358\\
0.604012463108397	0.376139481137542\\
0.605509629425276	0.378694277559492\\
0.605912531093133	0.382258224076919\\
0.604390648910618	0.385901578121503\\
0.601074971259114	0.389067324428908\\
0.594758932813508	0.392691231697316\\
0.582576504204521	0.396672003284468\\
0.564403335774116	0.399665331628139\\
0.543587382427411	0.40080958201004\\
0.525906967395448	0.400619845561091\\
0.493525231384274	0.398472747847316\\
0.457289818556047	0.393617338228412\\
0.431254641022287	0.388666524818272\\
0.391411175855567	0.379480134315336\\
0.356982004703935	0.369503052077386\\
0.321403811832181	0.357792774820342\\
0.282259302604925	0.342771841054673\\
0.247536568695296	0.328034523460971\\
0.210148426275164	0.310106423125453\\
0.175221557745807	0.291813367885065\\
0.140041108045668	0.271273132572285\\
0.103106199311675	0.248019373640193\\
0.0716233453926332	0.226124507137795\\
0.0342547183523489	0.198287682157441\\
-0.00172868802770976	0.168826818035807\\
-0.0408598394484972	0.13466520638993\\
};

\addplot [color=red, dotted, line width=0.8pt, forget plot]
table[row sep=crcr]{%
	0.727830156269241	0.319586716944193\\
	0.713739618837556	0.374549077256329\\
	0.69859176762567	0.425581457141168\\
	0.682923852539411	0.471626555430038\\
	0.666975778938442	0.512557393067329\\
	0.650976733376787	0.548265643590978\\
	0.635149181375618	0.578655939971158\\
	0.619720796029353	0.603617020968377\\
	0.604923151400703	0.623037235760165\\
	0.59097183439196	0.636876121920154\\
	0.584366196737218	0.641741152736332\\
	0.578042381069272	0.645238645206063\\
	0.571978057506826	0.647519375680154\\
	0.560642621541119	0.648718103409332\\
	0.502281884062518	0.648917371225839\\
	-0.482655876888779	0.649002831864855\\
	-0.499224100258374	0.646974542968521\\
	-0.515576827988639	0.643345836769781\\
	-0.531634751644887	0.637888691795357\\
	-0.5473930926814	0.630739022229964\\
	-0.562757895834908	0.621631469906699\\
	-0.577673065207956	0.610509553816232\\
	-0.592076542787623	0.597311170686111\\
	-0.619664995793474	0.567004983112814\\
	-0.646239892847336	0.534944177778879\\
	-0.671678150089933	0.501105859222293\\
	-0.695827331249418	0.465380155776681\\
	-0.718522481573592	0.427651719348215\\
	-0.739585813506189	0.387800055872696\\
	-0.758826574870187	0.34570006913132\\
	-0.776041165725568	0.301222872462193\\
	-0.791013586760209	0.254236931278465\\
	-0.797588299655841	0.229761723388978\\
	-0.803516321609709	0.204609606989809\\
	-0.808767649449628	0.178764161233906\\
	-0.813311777611348	0.152209178849434\\
	-0.817117784372452	0.124928760996465\\
	-0.820154433224048	0.0969074225806695\\
	-0.822390290924154	0.0681302085033281\\
	-0.823793863827497	0.0385828212374162\\
	-0.824333754108747	0.00825175999502736\\
	-0.823978837484797	-0.0228755283986104\\
	-0.820641272511179	-0.0870033671385799\\
	-0.814964299510259	-0.149711987696428\\
	-0.80806935804541	-0.207954273744702\\
	-0.800270183579494	-0.261474073788178\\
	-0.791876211417261	-0.310027950365546\\
	-0.783193639363061	-0.353381823985683\\
	-0.770296922935157	-0.408162160418067\\
	-0.758483162723669	-0.449950733749049\\
	-0.748844657862397	-0.477929925667213\\
	-0.742562960706666	-0.491194853923763\\
	-0.741437818537297	-0.49217156350485\\
	-0.74088332439583	-0.491356422931141\\
	-0.741709805976118	-0.484166877606207\\
	-0.745515384269097	-0.469230896060429\\
	-0.752815978773304	-0.446105042196849\\
	-0.76240789721238	-0.419495106805034\\
	-0.771502612209039	-0.397856102556549\\
	-0.779418070449307	-0.382514294968068\\
	-0.785760578756605	-0.373829672207206\\
	-0.788213927948274	-0.372096027557752\\
	-0.790124540569704	-0.372140992775258\\
	-0.791482998340794	-0.373874440235722\\
	-0.792801967330322	-0.38120842381035\\
	-0.791929377297347	-0.398749398976911\\
	-0.78825683318462	-0.42737688275486\\
	-0.781179507078537	-0.470231419643704\\
	-0.772554858143473	-0.49254910879267\\
	-0.761628949324587	-0.515943662560922\\
	-0.75224547337162	-0.532717673218122\\
	-0.742310789604983	-0.547570873405453\\
	-0.732104374757567	-0.559996992224084\\
	-0.721767956245266	-0.569973721408636\\
	-0.711350273551196	-0.57774950805533\\
	-0.697042643413643	-0.5861415771889\\
	-0.68607061129542	-0.59090902149801\\
	-0.671306351368632	-0.595134410787201\\
	-0.614235198614503	-0.606283606809955\\
	-0.601103069917829	-0.607205347916507\\
	-0.583599788347145	-0.606512302560834\\
	-0.570638819216307	-0.604291916124497\\
	-0.557793080588131	-0.600987758007454\\
	-0.541157740979716	-0.594328690277615\\
	-0.525217478621602	-0.585253633553724\\
	-0.510136082409628	-0.57388766209631\\
	-0.499390632607968	-0.564202263903948\\
	-0.48285924365872	-0.54564575779676\\
	-0.467941290820393	-0.525210463162132\\
	-0.457355536373273	-0.507471120574587\\
	-0.447845995516357	-0.489278653386228\\
	-0.437711797659592	-0.465583609442792\\
	-0.430875614882624	-0.446562190417516\\
	-0.423936903086326	-0.422899606038267\\
	-0.418804926749441	-0.399504308297297\\
	-0.415369552251284	-0.376868516580874\\
	-0.413318562361756	-0.351234264610944\\
	-0.413200687296271	-0.331413763028428\\
	-0.414390676194053	-0.313208637303736\\
	-0.416776262469606	-0.296698799083431\\
	-0.420236710670483	-0.281951103076601\\
	-0.424619571019291	-0.26907543994147\\
	-0.429805966102836	-0.258071835651012\\
	-0.435672572956632	-0.248896009936005\\
	-0.443429951027211	-0.240256861156703\\
	-0.450368377897442	-0.234931466428104\\
	-0.459059674820499	-0.230759216701362\\
	-0.468026650366672	-0.228701761025271\\
	-0.477141605162408	-0.22849810085655\\
	-0.487770059804632	-0.230389135129039\\
	-0.496756360970694	-0.233543333330664\\
	-0.506930029231703	-0.238810624749598\\
	-0.519357044298509	-0.247463985931703\\
	-0.53215955125552	-0.258689166101421\\
	-0.542502518205364	-0.270001134160282\\
	-0.554429589242899	-0.285771512201662\\
	-0.564350222129265	-0.301120490625195\\
	-0.57266559614868	-0.317146212420341\\
	-0.582575726919825	-0.340202845857461\\
	-0.587060144211072	-0.355248817259584\\
	-0.588333251057165	-0.36462355244353\\
	-0.587661040404424	-0.372122804829571\\
	-0.585573179926774	-0.377520240184656\\
	-0.581986908906741	-0.382118204141553\\
	-0.576137609191481	-0.386360912160415\\
	-0.568700760545369	-0.389780873714283\\
	-0.55587795821131	-0.393121504396721\\
	-0.543826761773927	-0.394839653748979\\
	-0.521066478093583	-0.396077762352091\\
	-0.500770828439831	-0.395403783530262\\
	-0.472769018243349	-0.393105916478065\\
	-0.440380230680312	-0.388323012972037\\
};

\addplot [color=black, draw=none, mark size=4.0pt, mark=asterisk, mark options={solid, black}, forget plot]
  table[row sep=crcr]{%
0.727830156269241	0.319586716944193\\
};

\addplot[<-, thick, color=black] coordinates
{ (0.4976,0.2181) (0.470,0.193)};
\addplot[<-, thick, color=black] coordinates
{ (-0.3356,-0.3625) (-0.374,-0.3747)};
\addplot[<-, thick, color=black] coordinates
{ (-0.1849,0.010) (-0.1492,0.037)};

\end{axis}
\end{tikzpicture}%

%% file: path_b23_eight_fwd_MIQP.tex
%
%
\begin{tikzpicture}

\begin{axis}[%
width=\figurewidth,
height=\figureheight,
at={(0\figurewidth,0\figureheight)},
scale only axis,
xmin=-1.4,
xmax=1.4,
xtick={-0.8,0,0.8},
xlabel={$\beta_3$ [rad]},
ymin=-1.4,
ymax=1.4,
ytick={-0.8,0,0.8},
xlabel style={font=\color{white!15!black},at={(axis description cs:0.5,-0.14)},anchor=north},
ylabel style={font=\color{white!15!black},at={(axis description cs:-0.2,.5)},anchor=south},
ylabel={$\beta_2$ [rad]},
axis background/.style={fill=white},
xmajorgrids,
ymajorgrids
]

\addplot[area legend, draw=black, fill=gray, fill opacity=0.3, forget plot]
table[row sep=crcr] {%
x	y\\
-0.98	0.65\\
-0.75	-0.65\\
0.98	-0.65\\
0.75	0.65\\
-0.98	0.65\\
}--cycle;

\addplot[area legend, draw=black, fill=green, fill opacity=0.3, forget plot]
table[row sep=crcr] {%
x	y\\
-0.72	-0.83\\
0.9	-0.43\\
0.72	0.83\\
-0.9	0.43\\
-0.72	-0.83\\
}--cycle;

\addplot[area legend, draw=black, fill=gray, fill opacity=0.3, forget plot]
table[row sep=crcr] {%
x	y\\
-0.75	-1.15\\
0	-0.65\\
0.75	0.65\\
0.75	1.15\\
0	0.65\\
-0.75	-0.65\\
-0.75	-1.15\\
0	-0.65\\
}--cycle;

\addplot [color=green, line width=0.8pt, forget plot]
table[row sep=crcr]{%
	-0.680061691709607	-0.342737783036121\\
	-0.690511671959332	-0.304070952187463\\
	-0.699849031974694	-0.266240177245472\\
	-0.708088628920542	-0.229250078113204\\
	-0.715248538420745	-0.193102536497781\\
	-0.721349631126078	-0.157796912186705\\
	-0.726415171776878	-0.123330260253717\\
	-0.730470443502191	-0.0896975461306817\\
	-0.733542399085377	-0.056891855964883\\
	-0.735659340070087	-0.0249046001467899\\
	-0.736850623874994	0.00627429166927285\\
	-0.737146398520057	0.0366561853275367\\
	-0.736577364127986	0.0662535676313467\\
	-0.735174560036773	0.0950798830989461\\
	-0.732969176126988	0.12314938233691\\
	-0.729992386815827	0.150476982121116\\
	-0.72627520608427	0.177078137082487\\
	-0.716742188177655	0.228164929491824\\
	-0.704610682391764	0.276539979185265\\
	-0.690112345348992	0.32233572156966\\
	-0.673468041898154	0.365684499470154\\
	-0.654886627273868	0.406716774659789\\
	-0.634564254270024	0.445559780956649\\
	-0.612684084102249	0.48233654482218\\
	-0.589416301830847	0.517165206309147\\
	-0.564918356636694	0.550158581342644\\
	-0.539335363704232	0.581423914708272\\
	-0.512800618101658	0.611062781163601\\
	-0.485436182151377	0.639171099439261\\
	-0.457353516707321	0.665839230392975\\
	-0.4286541338556	0.691152136180758\\
	-0.399430254159584	0.715189582044775\\
	-0.369765455966924	0.738026366258285\\
	-0.324604802682127	0.770182370108732\\
	-0.278844811329805	0.800011463575605\\
	-0.232676771776062	0.82771093496481\\
	-0.186262829393748	0.853459302248706\\
	-0.139739615647146	0.877418105732492\\
	-0.093221512589492	0.899733577345848\\
	-0.0468035655923218	0.920538175014049\\
	0.014798039377151	0.946133314718182\\
	0.0759376012998771	0.969509085885616\\
	0.136506265202481	0.990889458918033\\
	0.1964259001063	1.01047201600705\\
	0.270333805701473	1.03268624738378\\
	0.343085817730883	1.05266117383343\\
	0.414659813328698	1.07065420996753\\
	0.485062802408126	1.0868883280745\\
	0.568042994463829	1.10431862875456\\
	0.649461518899136	1.119782518108\\
	0.742619686938959	1.13566476938913\\
	0.76890011769597	1.13981997534204\\
	0.782202827624697	1.14100043613467\\
	0.795836744770274	1.14063960513906\\
	0.809767381776034	1.13876531492628\\
	0.82396270789289	1.13540270631668\\
	0.838392872963014	1.13057446572328\\
	0.853029946237398	1.12430105175438\\
	0.867847668730636	1.11660091215371\\
	0.882821217795777	1.10749069219049\\
	0.897926982633137	1.09698543561424\\
	0.913142349517365	1.08509877926758\\
	0.928445495627975	1.07184314240744\\
	0.943815190493276	1.05722991172573\\
	0.959230604201491	1.0412696229849\\
	0.974671121692312	1.02397214009425\\
	0.990116162614647	1.00534683234884\\
	1.00554500641988	0.985402750434826\\
	1.02093662255289	0.964148801672339\\
	1.05152151709047	0.917747264599379\\
	1.08169028484292	0.866217247653284\\
	1.11124748258242	0.809642632517475\\
	1.13997848819767	0.74812164877543\\
	1.16764584757623	0.681772200116318\\
	1.1939864644914	0.610736867748127\\
	1.21870986647155	0.535187376631375\\
	1.23030209154189	0.495972842797711\\
	1.24074697032758	0.457594783787667\\
	1.24994921531893	0.420391412864243\\
	1.2579588153163	0.384287407490342\\
	1.26479347326492	0.349304589225684\\
	1.27050642730576	0.315363422962594\\
	1.27511186471223	0.282501752225778\\
	1.27866221381232	0.250645726945832\\
	1.28117044297883	0.219841147570498\\
	1.28268695466731	0.190022434831498\\
	1.28322438491755	0.161240098630015\\
	1.28283143940601	0.133435268177055\\
	1.28152075465233	0.106661759250308\\
	1.27933976538615	0.0808657027833866\\
	1.27630238268219	0.0560995518427996\\
	1.27243807377365	0.0323700533941618\\
	1.26779268543816	0.00962903074984767\\
	1.26238437337109	-0.0120823494266455\\
	1.25625421247175	-0.0327972148875655\\
	1.24942077055922	-0.0524748345187589\\
	1.24192503059081	-0.0711490983877316\\
	1.23378667627289	-0.0887833827098805\\
	1.2250335367817	-0.105366684961295\\
	1.21570797153649	-0.120941256782327\\
	1.20583278363779	-0.135484490711701\\
	1.19543381073655	-0.148983064455558\\
	1.18454206815016	-0.161443998012833\\
	1.17319806311066	-0.172912445932261\\
	1.16142432425002	-0.183373117235073\\
	1.1492462097217	-0.192820945919362\\
	1.13669304197048	-0.201268261816176\\
	1.12379370657063	-0.20872920363266\\
	1.11058442008591	-0.215249416251815\\
	1.09708858916061	-0.220832493021619\\
	1.08332726300602	-0.225474006059571\\
	1.0693269218264	-0.229192979850563\\
	1.05511354623274	-0.232010218797741\\
	1.04071242496162	-0.233947539887512\\
	1.02614815066981	-0.235027673743677\\
	1.01144462472023	-0.235274181708653\\
	0.981712014696515	-0.233364199092381\\
	0.951692319970273	-0.228419372858065\\
	0.921533151781921	-0.220571732095815\\
	0.891412484508119	-0.210104755880334\\
	0.861474863933147	-0.197224613063115\\
	0.831850510699414	-0.182110200093184\\
	0.802676774854536	-0.16500276249513\\
	0.774055696140164	-0.146036441425341\\
	0.746113311225936	-0.125454563947236\\
	0.718946096040559	-0.103425355148451\\
	0.692646974961041	-0.0801211594088456\\
	0.667306574764954	-0.0557287553987991\\
	0.642991318753345	-0.0303737036115408\\
	0.619789191831032	-0.00427076028631435\\
	0.597768072682028	0.0224043075671172\\
	0.576965006450886	0.0495622891389542\\
	0.557431823116276	0.0770410147006966\\
	0.539206744076248	0.104710816088965\\
	0.522326235580451	0.132415794267181\\
	0.506793733901643	0.160098024765936\\
	0.486110704264846	0.201133115796446\\
	0.474051318031425	0.22806128814774\\
	0.463359038623105	0.25457070584068\\
	0.449866735214186	0.293320455385331\\
	0.439450227541742	0.330308678256391\\
	0.431929368674916	0.365521321338198\\
	0.42720019123601	0.398665309971338\\
	0.425131439583793	0.429481276101725\\
	0.425567775685931	0.457746485186086\\
	0.427191490497414	0.474988011846098\\
	0.429806727984494	0.490906317042792\\
	0.433320304301154	0.50557970477018\\
	0.437669419851791	0.518993562350441\\
	0.442791085397888	0.531137562408397\\
	0.44864601792946	0.5419251220871\\
	0.455162872341966	0.551365095260735\\
	0.462246339589273	0.559571052197408\\
	0.469859986019308	0.566470664691169\\
	0.477936691847678	0.572084258586499\\
	0.486408582888505	0.576458420351082\\
	0.495212234264166	0.579621885751108\\
	0.504283563686161	0.581629320287006\\
	0.518253955675603	0.582588475270863\\
	0.53251329963388	0.581160875292225\\
	0.546863947623954	0.577570718426144\\
	0.561119930431647	0.572066839483531\\
	0.575136040391609	0.564823176971662\\
	0.588740242196197	0.5561444652343\\
	0.601823638547099	0.546190796303432\\
	0.614252826045365	0.535247221794862\\
	0.629694543058503	0.519303696182809\\
	0.64364520907059	0.502330866619832\\
	0.655976139956707	0.484737079683673\\
	0.666599098749704	0.466886289884356\\
	0.675442355415143	0.449197092339295\\
	0.682487328362402	0.432016348114582\\
	0.687761130722823	0.415621187257219\\
	0.691327379177328	0.400228900958506\\
	0.693278490196869	0.386004622713047\\
	0.693729139028534	0.373067364840077\\
	0.692810616869069	0.361494930754931\\
	0.69065647442462	0.351360103642554\\
	0.687390940655143	0.342774206456733\\
	0.683158998309441	0.335789669924146\\
	0.678122118315366	0.330394278914683\\
	0.672437014729711	0.326566490591172\\
	0.66628977471344	0.324153530603084\\
	0.659838732550123	0.32303456521475\\
	0.651559286221985	0.323291880731706\\
	0.643267368351231	0.325177096084883\\
	0.635171366678374	0.328490613346634\\
	0.625986416070827	0.333986368080765\\
	0.616364570055509	0.341816264269051\\
	0.607198741443983	0.351618977055779\\
	0.600178183988437	0.361426567541041\\
	0.594845131289174	0.371567953029888\\
	0.591895395237573	0.380579090766993\\
	0.591535573081478	0.386385218945403\\
	0.592967625251817	0.388814289907905\\
	0.595141274654725	0.38810250083635\\
	0.597849043343503	0.384247117570765\\
	0.600157315704308	0.37768358048032\\
	0.601486310065694	0.369039205046158\\
	0.601444388593131	0.358839473136929\\
	0.599779737539702	0.347579402794627\\
	0.596368522036963	0.335698079108423\\
	0.591214255972861	0.323474784593147\\
	0.583630736291077	0.309812002022083\\
	0.574195659216078	0.296022743427467\\
	0.560911321660426	0.279582278080717\\
	0.543091335536665	0.260247611950474\\
	0.520096817664566	0.237683451447733\\
	0.488130334992461	0.208738949336379\\
	0.452712368696716	0.178845406782636\\
	0.412490053191097	0.147090851890615\\
	0.369789906353832	0.11558998050571\\
	0.327475226326431	0.0861525479475196\\
	0.235094193522164	0.024434835295287\\
	0.188052453500988	-0.00822310096643308\\
	0.163578722693994	-0.0265994276918839\\
};

\addplot [color=red, line width=0.8pt, forget plot]
  table[row sep=crcr]{%
-0.680061691709607	-0.342737783036121\\
-0.690511671958917	-0.304070952188808\\
-0.69984903197411	-0.266240177247466\\
-0.708088628919845	-0.229250078115713\\
-0.715248538420116	-0.193102536500252\\
-0.721349631125523	-0.157796912189121\\
-0.726415171776382	-0.123330260256119\\
-0.730470443501718	-0.0896975461331874\\
-0.733542399084977	-0.0568918559673354\\
-0.735659340069637	-0.0249046001496015\\
-0.736850623874536	0.00627429166621551\\
-0.737146398519681	0.0366561853245296\\
-0.73657736412757	0.0662535676279626\\
-0.735174560036456	0.0950798830956452\\
-0.732969176126715	0.123149382333508\\
-0.72999238681561	0.150476982117653\\
-0.726275206084142	0.177078137079073\\
-0.716742188177692	0.228164929488498\\
-0.704610682390082	0.276539979174095\\
-0.69011234534748	0.322335721557232\\
-0.673468041896967	0.365684499456778\\
-0.654886627273379	0.406716774647015\\
-0.644901450128509	0.426276953282802\\
-0.634291482361777	0.444383599025332\\
-0.622988336693353	0.460538325873349\\
-0.610997819901796	0.47452069699694\\
-0.598437448143477	0.486596617834891\\
-0.585349573206261	0.496738804509355\\
-0.571823515497743	0.505136519438068\\
-0.557943584902782	0.511997095175836\\
-0.543675695813066	0.517013142322031\\
-0.50050638117999	0.528464894280561\\
-0.28248964484573	0.582273618042502\\
0.553348547033605	0.787684529725549\\
0.561804737543913	0.787069853538882\\
0.570629246487768	0.784839423305691\\
0.579785898819793	0.781012648736895\\
0.589231770605059	0.775637813199822\\
0.598943669725949	0.768699786802269\\
0.608870636895849	0.760281734274926\\
0.618972212033903	0.75043745480337\\
0.629218945313307	0.739187168460431\\
0.650014237861946	0.712616390249056\\
0.670887030482166	0.681225128854638\\
0.6915159330466	0.645576948534071\\
0.7112610281674	0.607376372551348\\
0.73631644279203	0.554381100910941\\
0.761879408853164	0.495731493132765\\
0.774802504211566	0.464752921102672\\
0.774544816893923	0.462894141734646\\
0.774050446808442	0.461581031984421\\
0.773887960461347	0.460704348119256\\
0.773549017655055	0.459858756147703\\
0.77380598846956	0.455996813199178\\
0.799449564893518	0.266938754117071\\
0.806542493229353	0.201119256868272\\
0.812756117221851	0.123863203306483\\
0.815876046768415	0.0664346505616794\\
0.817734422147838	0.0048115602654305\\
0.817624948615251	-0.0258546885070259\\
0.816441336031286	-0.0550426541305696\\
0.814106444151693	-0.0823090059123804\\
0.810510256884775	-0.10702979237742\\
0.805749855443589	-0.129270644672484\\
0.799918252810633	-0.149097185762585\\
0.793104010565732	-0.166574393499392\\
0.78539101963761	-0.181766098761363\\
0.776858413951824	-0.194734601476481\\
0.767580592564761	-0.205540397403956\\
0.757627325646338	-0.214241994784877\\
0.747063925953563	-0.220895810605907\\
0.735951467900454	-0.225556127180979\\
0.724347045428909	-0.228275114668275\\
0.712304044523669	-0.229102859556724\\
0.699872435815146	-0.228087441937809\\
0.687099071055031	-0.225275020418626\\
0.674027991131436	-0.220709967746767\\
0.660700732876041	-0.214435022279203\\
0.647156635403386	-0.206491466994171\\
0.633433237105123	-0.19691966625402\\
0.60612034628055	-0.174991794142867\\
0.579718406655512	-0.151815394386508\\
0.542116968369474	-0.1155017144865\\
0.50727013959318	-0.0784687892570632\\
0.475372100786198	-0.0414422310857815\\
0.446570734395501	-0.00515910797204211\\
0.420878896708663	0.0300255293264801\\
0.398231386162604	0.0638842030785466\\
0.37863125226624	0.0958895982508815\\
0.361846658545676	0.126243952279653\\
0.347730476542896	0.154891310072675\\
0.336100171408089	0.181845913495897\\
0.326835379218629	0.206947115365202\\
0.319737628377054	0.230275576951578\\
0.314610433981473	0.251935684876575\\
0.311266922088324	0.27203676511103\\
0.309568230324172	0.290562384812122\\
0.309306061010971	0.307744771361057\\
0.310341596523304	0.323642745707967\\
0.313546444331863	0.342819813847796\\
0.318542220892127	0.359900947758404\\
0.325059470578551	0.375027375999519\\
0.33284975396509	0.388347564923181\\
0.341692955036857	0.39998807959558\\
0.351375997541962	0.410130714572422\\
0.361732103059045	0.418865561202438\\
0.37261795052947	0.426248436028178\\
0.386755297909437	0.433713506099745\\
0.401268174279132	0.439315542731165\\
0.41593542240251	0.443209910624339\\
0.430562231715438	0.445551460593895\\
0.444989880345297	0.446449689214934\\
0.461836808066072	0.445793616505847\\
0.477989197898967	0.443376093013682\\
0.493267786349332	0.439359148554706\\
0.507510131882958	0.433944236655833\\
0.520571831943328	0.427376007331835\\
0.532370168717041	0.419792856205233\\
0.542830265597409	0.411372414219718\\
0.55328388865508	0.400686655683056\\
0.560714843717534	0.39095552442938\\
0.567582057498692	0.379157229800088\\
0.572547510785446	0.367018248101572\\
0.575657870985186	0.354716343151516\\
0.576994227834405	0.342358881348986\\
0.576656776913722	0.330011155655041\\
0.57473839993979	0.317785656460597\\
0.570764146670949	0.303996088173723\\
0.565053386512788	0.290451738236034\\
0.557786049663446	0.277146547273076\\
0.547956237355123	0.262465619320482\\
0.535237119982785	0.246505718509201\\
0.519388629661477	0.229322748798863\\
0.500258815384077	0.210947309569707\\
0.475985341563073	0.189909160688134\\
0.446249471527274	0.166363349326249\\
0.410970193383325	0.140372241331878\\
0.359967862459181	0.10513590463419\\
0.29739163454737	0.0641410017503536\\
0.23286097716099	0.022087017835198\\
};
\addplot [color=blue, line width=0.8pt, forget plot]
  table[row sep=crcr]{%
-0.680061691709607	-0.342737783036121\\
-0.690511671959332	-0.304070952187463\\
-0.699849031974694	-0.266240177245472\\
-0.708088628920542	-0.229250078113204\\
-0.715248538420745	-0.193102536497781\\
-0.721349631126078	-0.157796912186705\\
-0.726415171776878	-0.123330260253718\\
-0.730470443502191	-0.0896975461306817\\
-0.733542399085377	-0.0568918559648829\\
-0.735659340070087	-0.02490460014679\\
-0.736850623874994	0.00627429166927285\\
-0.737146398520057	0.0366561853275368\\
-0.736577364127986	0.0662535676313466\\
-0.735174560036773	0.0950798830989462\\
-0.732969176126988	0.12314938233691\\
-0.729992386815827	0.150476982121116\\
-0.72627520608427	0.177078137082487\\
-0.716742188177655	0.228164929491824\\
-0.704610682391764	0.276539979185265\\
-0.690112345348992	0.32233572156966\\
-0.673468041898153	0.365684499470154\\
-0.654886627273868	0.406716774659789\\
-0.634564254270024	0.445559780956649\\
-0.612684084102249	0.48233654482218\\
-0.589416301830847	0.517165206309147\\
-0.564918356636694	0.550158581342644\\
-0.539335363704232	0.581423914708272\\
-0.526070864481114	0.595927313641692\\
-0.512387220605819	0.609042634911685\\
-0.498216472561121	0.62025772929128\\
-0.483601249862546	0.629556323686952\\
-0.468628350601252	0.637154913970718\\
-0.453293171295758	0.642854779697872\\
-0.437673233348126	0.646853034299636\\
-0.421794864183011	0.649126703028684\\
-0.389746750365926	0.650267027745202\\
-0.281970053404054	0.65057770880785\\
-0.0161977478230186	0.650338400751494\\
-0.00466502082098386	0.651283599922325\\
0.00633561613502531	0.653931157371862\\
0.0168490351779601	0.658238824832302\\
0.0269205361109607	0.66416385897073\\
0.0662105280579268	0.692181145112717\\
0.155906817937296	0.752251788398454\\
0.400236793708097	0.913929281485956\\
0.411108513749527	0.918979358669431\\
0.422388709213141	0.922350907055777\\
0.434046958242308	0.924033663044915\\
0.446032483530097	0.924105238700395\\
0.458318243056892	0.92256006402911\\
0.47085732203013	0.91947624419981\\
0.483623921172848	0.914858739218858\\
0.496592080192752	0.908717350769961\\
0.509732982110241	0.901076212688836\\
0.522994259838463	0.892052894304437\\
0.536345745988427	0.881689322204491\\
0.549776406854369	0.869959585343593\\
0.563267929314543	0.856862820420392\\
0.576756727058085	0.842565031865942\\
0.590269192128015	0.826913719478563\\
0.603747946269321	0.810048766967676\\
0.630491687801505	0.772824093935658\\
0.656493905399038	0.732180249321808\\
0.681474645384802	0.688689346643611\\
0.704935597671461	0.643777464242849\\
0.734448938254451	0.582926828080588\\
0.761086500843221	0.524131611913769\\
0.782054739902469	0.473699613799918\\
0.784317623777968	0.469525961747752\\
0.784699046068045	0.467358209761023\\
0.784570429073988	0.466278552235971\\
0.784130483815212	0.464828553583705\\
0.784019925051499	0.462754330392486\\
0.811089818565571	0.299614399680853\\
0.819420006949383	0.235289655498995\\
0.824628315490495	0.185929557381131\\
0.829189472106454	0.131678205366085\\
0.832784036894589	0.0727945415643708\\
0.834686251289173	0.0109236735774636\\
0.834067324217788	-0.018035749093072\\
0.83210449953273	-0.0444163299267594\\
0.828895000985056	-0.0682802341696402\\
0.824532891584209	-0.0896905792871578\\
0.819108492508338	-0.10871065045048\\
0.812707989589868	-0.125403247385546\\
0.805413214640634	-0.139830219566708\\
0.797301550418653	-0.152052099946914\\
0.788445942584601	-0.162127852504313\\
0.77891496721605	-0.17011460897208\\
0.768773014987982	-0.176067669486578\\
0.758080475387937	-0.180040395181077\\
0.746893972324788	-0.182084205518033\\
0.735266624519613	-0.182248615962962\\
0.723248322223012	-0.180581308098598\\
0.710886013429552	-0.177128225432303\\
0.698223994094597	-0.171933689158113\\
0.685304197939103	-0.165040529002339\\
0.672166482290292	-0.156490225017584\\
0.658848907083891	-0.14632305680626\\
0.63223317880003	-0.122781338527424\\
0.59427517774736	-0.085767358269938\\
0.559111420798918	-0.0481214756788183\\
0.526945116845349	-0.0105928704818709\\
0.497981558065603	0.0258878062757307\\
0.472151328685348	0.0611811054225796\\
0.449450628437847	0.0948505056662475\\
0.429806580257981	0.126591154876565\\
0.413016424499922	0.15650163420171\\
0.398939557569169	0.184439063316466\\
0.387371837912374	0.210472712380117\\
0.378115684243991	0.234668325946118\\
0.370980943309076	0.257088670283234\\
0.365785272607478	0.277794698498005\\
0.362354303430219	0.296846746013761\\
0.3605401960091	0.314229821314917\\
0.360173565568438	0.330037096672829\\
0.361074012098936	0.344433135309594\\
0.363117426269944	0.357436118066042\\
0.366135959137158	0.369233624559611\\
0.371491016887542	0.38314325974694\\
0.378112135013903	0.395184678573825\\
0.385767926474222	0.405513565179593\\
0.394255857912214	0.414260168688112\\
0.403396764236093	0.421537690764678\\
0.415491303463031	0.428753515866479\\
0.428047196769124	0.434162571152391\\
0.440831449009731	0.437963287394554\\
0.456202776807647	0.440633862011317\\
0.471351772987462	0.441450323440284\\
0.486051160527761	0.440607798199085\\
0.500093321214938	0.438361557043385\\
0.515432176329979	0.434266484712737\\
0.529483022893955	0.42882512327392\\
0.542129819668884	0.422218606421286\\
0.553280471721508	0.414644245972143\\
0.562869380214656	0.406306476072428\\
0.570856323221206	0.397436417900267\\
0.577253574517057	0.388162997520676\\
0.582100814595068	0.378572298527523\\
0.585436404914555	0.368773203277477\\
0.587322661174246	0.358831697665281\\
0.58778395607681	0.347399266192237\\
0.586552821569997	0.335939100069309\\
0.583760938209234	0.324441642231066\\
0.578896458306953	0.311492284099533\\
0.571535026580793	0.297120749923966\\
0.562325217968491	0.282748518992738\\
0.549115942870009	0.265443423057173\\
0.532464813818826	0.246629032660118\\
0.512183060719661	0.2262944424127\\
0.486490166708395	0.202956085071096\\
0.456749501979496	0.178196269835075\\
0.419267482026197	0.14924929776688\\
0.377809583448643	0.119380557811065\\
0.330902437643531	0.0873794202232855\\
0.203004968157824	0.00243014299341049\\
};

\addplot [color=black, line width=1.0pt, forget plot]
  table[row sep=crcr]{%
-0.0408598394484972	0.13466520638993\\
-0.142880355317615	0.0425116099110959\\
-0.165790007183586	0.0243363978601857\\
-0.207792446652068	-0.00619805834531528\\
-0.251309841886497	-0.0355726911620831\\
-0.379727274189688	-0.122585643275196\\
-0.421790650072576	-0.153997633105531\\
-0.454559715096812	-0.180436049038848\\
-0.483834279693127	-0.205985574362196\\
-0.509368282939933	-0.230408118394342\\
-0.531115334474125	-0.253330947980684\\
-0.549797922321808	-0.275364585524183\\
-0.565273302455181	-0.296096826582156\\
-0.577507091598239	-0.315180439378931\\
-0.586604540196093	-0.332242623443335\\
-0.59319650337228	-0.347968746384261\\
-0.597203963419761	-0.361807404254658\\
-0.598738317277314	-0.372796482947882\\
-0.598487369328993	-0.379421503086881\\
-0.601013674108382	-0.376391595304695\\
-0.603298229778452	-0.375956963721876\\
-0.605382958825144	-0.377548357263932\\
-0.606062065673226	-0.381878150143506\\
-0.604756500932058	-0.385467104119656\\
-0.600950967176329	-0.389347319813641\\
-0.593036813203388	-0.393550839721977\\
-0.581939992226893	-0.396853746022613\\
-0.57146846982738	-0.398734496993717\\
-0.549923990652537	-0.40067549765775\\
-0.520712345713422	-0.400568426157832\\
-0.48894538650299	-0.39803760523401\\
-0.465253652041426	-0.394800329267211\\
-0.426403064351405	-0.387791457136648\\
-0.391196145447763	-0.379302274861166\\
-0.359422232257052	-0.370335598480062\\
-0.319700204240819	-0.357065001344296\\
-0.288124569844594	-0.345153565668635\\
-0.247324071292396	-0.327760832059443\\
-0.214828939919369	-0.312422909768168\\
-0.176746911581933	-0.2923809485882\\
-0.14320955034361	-0.273162188194853\\
-0.108049113236798	-0.250945858736107\\
-0.0731543065995718	-0.227199976915445\\
-0.0390435046810411	-0.201710425738353\\
-0.00278326072706225	-0.172752919736408\\
0.034772879575632	-0.14018062297753\\
0.0744666055293236	-0.103946966911396\\
0.12222101782624	-0.0602793486005441\\
0.147837350724528	-0.0385532467579188\\
0.178487104403983	-0.0151200011325642\\
0.212810258784891	0.00964261761365393\\
0.275609120993059	0.0516861747166277\\
0.344679669943556	0.0978743471052701\\
0.38580933749421	0.126997562165325\\
0.426399826572588	0.157575129460162\\
0.459944653744835	0.184945781978777\\
0.487766468086906	0.209596864355854\\
0.512879924711121	0.233950222808982\\
0.534165750290406	0.256779768610305\\
0.552351719063648	0.278652838292397\\
0.567314517920431	0.299164072917108\\
0.578637708782478	0.317295628238769\\
0.587640683137094	0.334696596103494\\
0.593750547823124	0.350025366254316\\
0.597413298595257	0.363672869514385\\
0.598602969473123	0.374401703859823\\
0.598389048382896	0.37919286264729\\
0.601347697055104	0.376059034750358\\
0.604012463108397	0.376139481137542\\
0.605509629425276	0.378694277559492\\
0.605912531093133	0.382258224076919\\
0.604390648910618	0.385901578121503\\
0.601074971259114	0.389067324428908\\
0.594758932813508	0.392691231697316\\
0.582576504204521	0.396672003284468\\
0.564403335774116	0.399665331628139\\
0.543587382427411	0.40080958201004\\
0.525906967395448	0.400619845561091\\
0.493525231384274	0.398472747847316\\
0.457289818556047	0.393617338228412\\
0.431254641022287	0.388666524818272\\
0.391411175855567	0.379480134315336\\
0.356982004703935	0.369503052077386\\
0.321403811832181	0.357792774820342\\
0.282259302604925	0.342771841054673\\
0.247536568695296	0.328034523460971\\
0.210148426275164	0.310106423125453\\
0.175221557745807	0.291813367885065\\
0.140041108045668	0.271273132572285\\
0.103106199311675	0.248019373640193\\
0.0716233453926332	0.226124507137795\\
0.0342547183523489	0.198287682157441\\
0.00287794162871391	0.172714639547868\\
-0.00172868802770976	0.168826818035807\\
-0.0408598394484972	0.13466520638993\\
};
\addplot [color=black, draw=none, mark size=4.0pt, mark=asterisk, mark options={solid, black}, forget plot]
  table[row sep=crcr]{%
-0.680061691709607	-0.342737783036121\\
};

\addplot[<-, thick, color=red] coordinates
{ (-0.365,0.57) (-0.4094,0.5557)};
\addplot[->, thick, color=blue] coordinates
{ (-0.7371,0.0366) (-0.7365,0.0662)};
\addplot[<-, thick, color=blue] coordinates
{ (0.2956,0.8479) (0.0440,0.6795)};
\addplot[->, thick, color=blue] coordinates
{ (0.7840,0.4627) (0.8110,0.2996)};
\addplot[->, thick, color=green] coordinates
{ (-0.0932,0.8997) (-0.0468,0.9205)};
\addplot[->, thick, color=green] coordinates
{ (1.2187,0.5351) (1.2303,0.4959)};

\addplot[->, thick, color=black] coordinates
{ (0.4976,0.2181) (0.470,0.193)};
\addplot[->, thick, color=black] coordinates
{ (-0.3356,-0.3625) (-0.374,-0.3747)};
\addplot[->, thick, color=black] coordinates
{ (-0.1849,0.010) (-0.1492,0.037)};

\end{axis}
\end{tikzpicture}%

%% file: path_th3_eight_rev_MIQP.tex
%
%
\begin{tikzpicture}

\begin{axis}[%
width=\figurewidth,
height=\figureheight,
at={(0\figurewidth,0\figureheight)},
scale only axis,
xmin=0,
xmax=80,
xtick={0,20,40,60,80},
xlabel={$t$ [s]},
ymin=-1,
ymax=1,
xlabel style={font=\color{white!15!black},at={(axis description cs:0.5,-0.2)},anchor=north},
ylabel style={font=\color{white!15!black},at={(axis description cs:-0.24,.5)},anchor=south},
ylabel={$\tilde \theta_3$ [rad]},
axis background/.style={fill=white},
xmajorgrids,
ymajorgrids
]
\addplot [color=red, line width=0.7pt, forget plot]
  table[row sep=crcr]{%
0	-5.19380563446248e-06\\
0.340000000000003	-0.0284940976627155\\
0.780000000000001	-0.0629597477739168\\
1.42	-0.109955573775807\\
1.87	-0.14195632464785\\
1.91	-0.144875861997633\\
2.18000000000001	-0.164035889861125\\
2.22	-0.166972464536826\\
2.48	-0.185602320876953\\
2.52	-0.188572288909853\\
2.75	-0.205321671404647\\
3.04000000000001	-0.226350431350284\\
3.08	-0.229300018178122\\
3.38	-0.250523035495661\\
3.42	-0.25340319859734\\
3.72	-0.274060609157715\\
3.76000000000001	-0.276866102029757\\
3.98999999999999	-0.292277736092629\\
4.03	-0.295020395532248\\
4.25	-0.309426896302298\\
4.29000000000001	-0.31210592614174\\
4.52	-0.326743539178892\\
4.55	-0.32875343205896\\
4.78	-0.342984641184387\\
4.81	-0.344942694265512\\
5.04000000000001	-0.358748155190185\\
5.07000000000001	-0.360651440171921\\
5.3	-0.374012255003464\\
5.33	-0.375857764816615\\
5.56	-0.388753198953964\\
5.59	-0.390538598447918\\
5.82000000000001	-0.402948161721966\\
5.84999999999999	-0.404669815872012\\
6.08	-0.416572430604703\\
6.11	-0.418227566338089\\
6.33	-0.429146178889255\\
6.36	-0.430733436756256\\
6.59	-0.441578203559999\\
6.62	-0.443092234015182\\
6.84999999999999	-0.453363483189563\\
6.88	-0.454800949405438\\
7.18000000000001	-0.467331516063339\\
7.20999999999999	-0.468666798548227\\
7.5	-0.47982109113768\\
7.53	-0.481050769077484\\
8.25	-0.504645765255404\\
8.88	-0.520003956304237\\
9.47	-0.529755674399439\\
9.97	-0.53408084698134\\
10.44	-0.534825956686603\\
10.84	-0.532650540278567\\
11.24	-0.528021796246122\\
11.58	-0.52192068278562\\
11.92	-0.513940781757938\\
12.22	-0.505188447217549\\
12.49	-0.495899381740188\\
12.75	-0.485535351790077\\
13.01	-0.47387926327896\\
13.31	-0.459165318078121\\
13.63	-0.44244562127696\\
14.02	-0.42118171596303\\
14.12	-0.415561063788431\\
15.27	-0.353286875243512\\
16.21	-0.307236928502888\\
16.51	-0.29316957034824\\
16.54	-0.291633473769849\\
16.88	-0.27607011767283\\
16.91	-0.274545955145157\\
17.16	-0.263306043821672\\
17.19	-0.261785766025326\\
17.44	-0.250667553581479\\
17.47	-0.24914977532687\\
17.71	-0.238551420801983\\
17.74	-0.237034058807694\\
17.98	-0.22652068307039\\
18.01	-0.225005385179287\\
18.25	-0.214574658847496\\
18.28	-0.213059729858486\\
18.52	-0.202711114085375\\
18.55	-0.201196537236157\\
18.79	-0.190927947230776\\
18.81	-0.189814142219618\\
18.97	-0.183116600320631\\
18.99	-0.182000761750615\\
19.15	-0.175340382809026\\
19.17	-0.174222322662914\\
19.32	-0.167991247071882\\
19.34	-0.166870132633278\\
19.5	-0.160281467111687\\
19.52	-0.159158403486728\\
19.68	-0.152605551521972\\
19.7	-0.151480914699704\\
19.86	-0.144963576349312\\
19.88	-0.143837167012293\\
20.03	-0.137736800553284\\
20.05	-0.136607806594967\\
20.21	-0.130159066777452\\
20.23	-0.129028583967951\\
20.39	-0.122614430120663\\
20.41	-0.121482677184645\\
20.56	-0.115477350028769\\
20.58	-0.114343518389248\\
20.74	-0.107996904959521\\
20.76	-0.106861993379383\\
20.91	-0.100918247161559\\
20.93	-0.0997814131549291\\
21.09	-0.0935002455095884\\
21.11	-0.0923626048555093\\
21.27	-0.086114276790795\\
21.29	-0.0849759337543787\\
21.44	-0.0791216828281023\\
21.46	-0.0779817721145548\\
21.62	-0.0717957576568296\\
21.64	-0.0706553461041324\\
21.79	-0.0648579337569686\\
21.81	-0.0637162430152074\\
21.97	-0.0575905388690359\\
21.99	-0.0564486346326305\\
22.14	-0.050706066054147\\
22.16	-0.0495630895283057\\
22.32	-0.0434953325081722\\
22.34	-0.0423523758605029\\
22.49	-0.0366623598039411\\
22.51	-0.0355186053246968\\
22.67	-0.0295106603786905\\
22.69	-0.0283686789897217\\
22.84	-0.0227701608684328\\
22.86	-0.0216360714645987\\
23.02	-0.0157899944476014\\
23.04	-0.0146659877455448\\
23.2	-0.00892808226551267\\
23.22	-0.00781499566923571\\
23.37	-0.00251514815232667\\
23.39	-0.00141304579651091\\
23.55	0.00408577195138093\\
23.57	0.00517422432947967\\
23.73	0.0105426848294314\\
23.75	0.0116166865305019\\
23.91	0.016848730956383\\
23.93	0.0179075009224334\\
24.09	0.0229983975636685\\
24.11	0.0240405970476587\\
24.26	0.0287062181839985\\
24.28	0.0297327180041265\\
24.44	0.0345315642518358\\
24.46	0.035540587477044\\
24.61	0.0399295677252098\\
24.63	0.0409224259673096\\
24.79	0.0454214259623882\\
24.81	0.0463948967821892\\
24.96	0.0504958915278309\\
24.98	0.0514512794176341\\
25.13	0.0554098121173467\\
25.15	0.0563474423428261\\
25.31	0.0603804256749498\\
25.33	0.0612978320639002\\
25.48	0.0649598516944252\\
25.5	0.0658571177687577\\
25.64	0.0691737087928175\\
25.66	0.0700538694388371\\
25.81	0.0734296029872041\\
25.83	0.0742902313895542\\
25.98	0.0775184036102559\\
26	0.0783579238383254\\
26.14	0.0812680618543169\\
26.16	0.0820893546083283\\
26.3	0.0848685211050508\\
26.32	0.08567111359244\\
26.47	0.0884704232947655\\
26.49	0.0892513559244463\\
26.63	0.0917607749109095\\
26.65	0.0925225816654489\\
26.79	0.0949018577800871\\
26.81	0.095643932281348\\
26.94	0.0977690590097211\\
26.96	0.0984931383390091\\
27.1	0.100618835518631\\
27.12	0.101323112010746\\
27.26	0.103320281248372\\
27.28	0.104004139943555\\
27.41	0.10577632503788\\
27.43	0.106442051078929\\
27.56	0.108102668405465\\
27.58	0.108749732892775\\
27.72	0.110378564214926\\
27.74	0.11100540302904\\
27.87	0.112438635542091\\
27.89	0.113046762286359\\
28.02	0.11437119894552\\
28.04	0.114960898385846\\
28.17	0.116178068548336\\
28.19	0.116749099381195\\
28.31	0.117813041075948\\
28.33	0.118367451040967\\
28.46	0.119380190485401\\
28.48	0.11991605606967\\
28.61	0.120825649605095\\
28.63	0.121343170581397\\
28.75	0.122127152227407\\
28.77	0.122628094899085\\
28.9	0.123342702632812\\
28.92	0.123825478969223\\
29.04	0.124432216895173\\
29.06	0.124898652571289\\
29.18	0.125421623309791\\
29.2	0.125871827364747\\
29.33	0.126309858515427\\
29.35	0.126742259513748\\
29.47	0.127097712266774\\
29.49	0.127514192589842\\
29.61	0.127791173479679\\
29.63	0.128191861793468\\
29.75	0.128391947089199\\
29.77	0.128777030821013\\
29.89	0.128901827825075\\
29.91	0.129271473325176\\
30.02	0.129355925290525\\
30.04	0.129711785406982\\
30.16	0.129695635953851\\
30.18	0.130036424745327\\
30.3	0.129950185773566\\
30.32	0.130276084191351\\
30.44	0.130122064650905\\
30.46	0.130433406095548\\
30.57	0.130267657217303\\
30.59	0.130565989161354\\
30.71	0.130285957454319\\
30.73	0.130570173085957\\
30.84	0.130291518235595\\
30.86	0.130563205378479\\
30.98	0.130164074448459\\
31	0.130422050992919\\
31.11	0.130037170786352\\
31.13	0.130283062232564\\
31.24	0.129848810756499\\
31.26	0.130082829775475\\
31.38	0.129520499995067\\
31.4	0.129741567257881\\
31.51	0.129211123118367\\
31.53	0.129420817081282\\
31.64	0.128845841670909\\
31.66	0.129044377676152\\
31.77	0.128426533973141\\
31.79	0.128614045263006\\
31.9	0.127955175659125\\
31.92	0.128132048778312\\
32.03	0.127433407130383\\
32.05	0.127599850044575\\
32.16	0.126863149625933\\
32.18	0.127019295240203\\
32.29	0.126246257309901\\
32.31	0.126392547512083\\
32.42	0.125584605465164\\
32.44	0.125721107242455\\
32.55	0.124879359732489\\
32.57	0.125006383220338\\
32.68	0.124132370295825\\
32.7	0.124250075336079\\
32.81	0.123345738855633\\
32.83	0.123454544564154\\
32.94	0.12252097844798\\
32.96	0.122621093798571\\
33.07	0.121659917473878\\
33.09	0.121751498306267\\
33.19	0.12088441573151\\
33.21	0.12096869735359\\
33.32	0.119958110701674\\
33.34	0.120034362438005\\
33.45	0.119000123057276\\
33.47	0.119068650904453\\
33.58	0.118012199354254\\
33.6	0.11807317128104\\
33.7	0.117122706433605\\
33.72	0.117177222672709\\
33.89	0.115615725772955\\
33.92	0.115531306169359\\
34.08	0.114056002918474\\
34.11	0.113960249926222\\
34.27	0.112447804738579\\
34.3	0.11234151287033\\
34.46	0.110796034556472\\
34.49	0.110680074610784\\
34.65	0.109105617213771\\
34.68	0.108980817240649\\
34.83	0.107514852402247\\
34.86	0.107382176841384\\
35.02	0.105759757093054\\
35.05	0.105619489913963\\
35.21	0.103978514899921\\
35.24	0.103831248780466\\
35.39	0.102309972158551\\
35.42	0.102156715724163\\
35.58	0.100487042510281\\
35.61	0.10032792926431\\
35.76	0.0987838822722438\\
35.79	0.098619762555515\\
35.95	0.0969331605923713\\
35.98	0.096764556149239\\
36.14	0.0950735604570667\\
36.17	0.0949009334183302\\
36.32	0.0933423956307564\\
36.35	0.0931664365309786\\
36.51	0.0914743302471379\\
36.54	0.0912953685018039\\
36.69	0.0897388563696353\\
36.72	0.0895575076641819\\
36.87	0.0880045573322121\\
36.9	0.0878210472386343\\
37.06	0.0861430265682088\\
37.09	0.0859579250461024\\
37.24	0.0844187601907578\\
37.28	0.0841034441512249\\
37.48	0.0821319366600477\\
37.67	0.0803023934939233\\
37.71	0.0799889105817044\\
37.92	0.0779277715783593\\
37.96	0.0776164882469033\\
38.16	0.0757040037663472\\
38.2	0.0753953377557224\\
38.41	0.0733895122329216\\
38.45	0.0730844560936958\\
38.65	0.0712288265752505\\
38.69	0.0709276119239064\\
38.9	0.0689879750165403\\
38.94	0.0686913893249965\\
39.14	0.0669018630357812\\
39.18	0.0666101050772738\\
39.38	0.0648554501769638\\
39.43	0.0644595358032092\\
39.63	0.0627434551321784\\
39.68	0.0623560878719616\\
39.93	0.0603013980777831\\
39.98	0.0599245619254702\\
40.24	0.0578288723236682\\
40.29	0.0574635685326967\\
40.55	0.0554353282899456\\
40.6	0.0550818478841819\\
40.85	0.0532156234935002\\
40.9	0.0528736480171119\\
41.16	0.0509798535414632\\
41.23	0.050590240956609\\
41.53	0.0484219998720477\\
41.6	0.0480495072132499\\
41.91	0.0460069277356041\\
42.26	0.0437302504684709\\
42.33	0.043391516677957\\
42.69	0.041165099385239\\
42.77	0.0407728264859628\\
43.12	0.0387487633876589\\
43.2	0.0383801605192815\\
43.57	0.0364347172554602\\
43.98	0.0343445152234523\\
44.07	0.033959971395447\\
44.48	0.0320212060621259\\
44.57	0.0316655462515456\\
44.99	0.0298910475952567\\
45.46	0.0279783563448603\\
45.56	0.02762523032807\\
46.02	0.0259224420506996\\
46.13	0.025554874558722\\
46.64	0.0238523639741288\\
46.75	0.0235197572225303\\
47.26	0.0219772773397438\\
47.39	0.0216586049493515\\
48	0.0199719346860121\\
48.14	0.0196564240402921\\
48.8	0.018051871217466\\
49.47	0.0166224251647691\\
50.18	0.0152310889940281\\
50.35	0.0149345479156295\\
51.13	0.0136341780119977\\
51.99	0.0123106661086752\\
52.2	0.0120390121017238\\
53.17	0.0107695919351301\\
53.41	0.0104887642875582\\
54.48	0.00931309525218182\\
54.76	0.00904465758662809\\
56.04	0.0078582728130101\\
56.37	0.00759891310595151\\
57.85	0.00647007470905692\\
58.25	0.00620942651698897\\
60.03	0.00512634656811883\\
60.55	0.00485821684615928\\
62.84	0.00379213283204649\\
63.57	0.00350837984697705\\
66.49	0.00255781657675414\\
70.63	0.00163024692561464\\
75.95	0.000921859070260211\\
80	0.000605144549339798\\
};
\addplot [color=blue, line width=0.7pt, forget plot]
  table[row sep=crcr]{%
0	-5.19380563446248e-06\\
0.280000000000001	-0.0230390326984633\\
0.629999999999995	-0.0497729112613854\\
1.03	-0.077997341655859\\
1.58	-0.114063278933799\\
3.03	-0.20307572960651\\
3.23999999999999	-0.216056413311392\\
3.28	-0.218626397385236\\
3.5	-0.232247826770447\\
3.54000000000001	-0.234779760512609\\
3.82000000000001	-0.25164056694264\\
3.86	-0.254087828913086\\
4.07000000000001	-0.26630818029102\\
4.11	-0.268685886749125\\
4.32000000000001	-0.280517090186351\\
4.34999999999999	-0.282297236615108\\
4.57000000000001	-0.294251975772724\\
4.59999999999999	-0.295978440516919\\
4.81	-0.307015051471126\\
4.84	-0.308685261692233\\
5.06	-0.319784058728487\\
5.09	-0.321397241459067\\
5.3	-0.331598210061628\\
5.33	-0.333152276332868\\
5.54000000000001	-0.342928938594255\\
5.57000000000001	-0.344424383039737\\
5.78	-0.353767772530318\\
5.81	-0.355202395865859\\
6.02	-0.364102488459466\\
6.05	-0.365474144575813\\
6.2	-0.371550152012674\\
6.23	-0.37288706631935\\
6.44	-0.380984434053602\\
6.47	-0.382230261208676\\
7.12	-0.403358808318757\\
7.32000000000001	-0.408759340216434\\
7.34999999999999	-0.409655226664398\\
8.01000000000001	-0.423462856665566\\
8.61	-0.430546539949674\\
9.16	-0.432012450514165\\
9.66	-0.42885824247827\\
10.12	-0.422005644857407\\
10.5	-0.41284651412056\\
10.93	-0.399016415831127\\
11.42	-0.379589610892282\\
11.94	-0.35525799510485\\
12.68	-0.316291164884859\\
13.03	-0.297128152756798\\
13.09	-0.293904380093636\\
13.36	-0.279469358813643\\
13.42	-0.276369240534393\\
13.69	-0.262589756658599\\
13.76	-0.25912510182512\\
14.04	-0.245481052392407\\
14.11	-0.242164663361294\\
14.44	-0.226768697123291\\
14.76	-0.212339024261652\\
15.11	-0.197131501876271\\
15.48	-0.181720115756121\\
15.84	-0.167399625197746\\
15.98	-0.162029605285127\\
16.36	-0.147915676925436\\
16.76	-0.133777898942711\\
16.98	-0.126372533272018\\
17.4	-0.11287107117117\\
18.59	-0.079359896129958\\
19.07	-0.0677693352878066\\
19.56	-0.0570506737593064\\
20.09	-0.0466632698954186\\
20.66	-0.0368318292011622\\
21.27	-0.027756902574211\\
21.59	-0.0235408419672041\\
22.26	-0.0158577684212986\\
23	-0.0089959469945029\\
23.83	-0.00304152639708377\\
24.8	0.00194712557212995\\
25.96	0.00571568396262023\\
27.41	0.00799966647926453\\
29.53	0.00837996805825014\\
39.05	0.00249335622837066\\
40.16	0.00213465053670348\\
43.59	0.00140830734736141\\
48.81	0.000865006269592072\\
57.64	0.000387211358457762\\
74.28	8.76217831091708e-05\\
80	6.15479756902459e-05\\
};
\addplot [color=green, line width=0.7pt, forget plot]
  table[row sep=crcr]{%
0	-5.19380562913341e-06\\
0.24	-0.0193690167441121\\
0.45	-0.0345605261152508\\
0.65	-0.047405268715087\\
0.83	-0.0573866079936982\\
1.03	-0.0668699226748588\\
1.25	-0.075439343277738\\
1.51	-0.0838866481267919\\
1.77	-0.0905766130073822\\
2.09	-0.0971628699621139\\
2.43	-0.102400069880276\\
2.88	-0.107663192597398\\
3.46	-0.112451858453416\\
3.98	-0.114823036177876\\
4.44	-0.115113621197135\\
4.86	-0.113220627156461\\
5.35	-0.108918599550348\\
5.85	-0.102307405247094\\
6.45	-0.092304495653396\\
6.82	-0.0854093090506254\\
};
\addplot [color=black, dotted, line width=1.0pt, forget plot]
table[row sep=crcr]{%
	0	-5.19380563446248e-06\\
	0.200000000000003	-0.0173292717541926\\
	0.370000000000005	-0.0312146516742899\\
	0.609999999999999	-0.0496986322886102\\
	0.849999999999994	-0.0671252474875104\\
	0.969999999999999	-0.0756709544523915\\
	1.29000000000001	-0.0973214733456018\\
	1.43000000000001	-0.106572040128256\\
	2	-0.14253419624886\\
	3.18000000000001	-0.214723242617623\\
	3.2	-0.216088373379051\\
	3.58	-0.239619351332408\\
	3.59999999999999	-0.240981833725655\\
	3.83	-0.254845089730537\\
	3.84	-0.25560717423329\\
	4.08	-0.269632675671787\\
	4.09	-0.270381555185267\\
	4.33	-0.283962403329014\\
	4.34	-0.284696519694933\\
	4.52	-0.294529903714164\\
	4.53	-0.295252919374079\\
	4.58	-0.297820344270264\\
	4.59	-0.29854002546125\\
	4.64	-0.301082608177879\\
	4.65000000000001	-0.301798054632499\\
	4.83	-0.311192241906355\\
	4.84	-0.311896007657438\\
	4.89	-0.314336956404787\\
	4.90000000000001	-0.315037159046099\\
	4.95	-0.317451783695148\\
	4.95999999999999	-0.318147579912292\\
	5.01000000000001	-0.320538831224582\\
	5.02	-0.321230208872933\\
	5.07000000000001	-0.323595608943251\\
	5.08	-0.324283268676041\\
	5.13	-0.326623507438541\\
	5.14	-0.32730666700057\\
	5.19	-0.329622262303246\\
	5.2	-0.33030161552432\\
	5.26000000000001	-0.333045060905064\\
	5.27	-0.333720138604448\\
	5.32000000000001	-0.335978939740514\\
	5.33	-0.336649404450952\\
	5.38	-0.338882402403158\\
	5.39	-0.33954889485176\\
	5.44	-0.341755354693674\\
	5.45	-0.342417145822282\\
	5.5	-0.344598881959314\\
	5.51000000000001	-0.345255952517164\\
	5.56	-0.347410221695327\\
	5.57000000000001	-0.348063151518915\\
	5.62	-0.350191800942525\\
	5.63	-0.350839917513156\\
	5.68000000000001	-0.352942118031507\\
	5.69	-0.353585995355914\\
	5.73999999999999	-0.355661073189182\\
	5.75	-0.356300041186486\\
	5.8	-0.358349634459799\\
	5.81	-0.358983669756796\\
	5.86	-0.361005285288542\\
	5.87	-0.361634904637356\\
	5.92	-0.363632547452283\\
	5.93000000000001	-0.36425833833205\\
	5.98	-0.366232742142117\\
	5.98999999999999	-0.366854114242443\\
	6.04000000000001	-0.36880918855438\\
	6.05	-0.369429015884137\\
	6.09999999999999	-0.371376233626933\\
	6.11	-0.371994881577692\\
	6.16	-0.373933790685214\\
	6.17	-0.374551843094935\\
	6.22	-0.37648485929067\\
	6.23	-0.377102542025753\\
	6.27	-0.378645724418632\\
	6.28	-0.379262298257359\\
	6.33	-0.381181953280546\\
	6.34	-0.381794445804985\\
	6.39	-0.383680584943093\\
	6.40000000000001	-0.38428681723606\\
	6.45	-0.38613796828102\\
	6.45999999999999	-0.386737755638151\\
	6.51000000000001	-0.388553393334405\\
	6.52	-0.389146538571765\\
	6.56	-0.390573029618039\\
	6.57000000000001	-0.391159478428776\\
	6.62	-0.392909509153256\\
	6.63	-0.393488915156837\\
	6.68000000000001	-0.395201607205962\\
	6.69	-0.39577373506566\\
	6.73999999999999	-0.397448331459827\\
	6.75	-0.398012933177\\
	6.8	-0.399648615209884\\
	6.81	-0.40020542938106\\
	6.86	-0.401801321383886\\
	6.87	-0.402350073958388\\
	6.92	-0.403905243301935\\
	6.93000000000001	-0.404445647969411\\
	6.98999999999999	-0.406256963418272\\
	7	-0.406788283103168\\
	7.05	-0.408250786954881\\
	7.06	-0.408773102841735\\
	7.11	-0.410191602853729\\
	7.12	-0.410704593062931\\
	7.18000000000001	-0.412347460860602\\
	7.19	-0.412850198338774\\
	7.23999999999999	-0.414168474646701\\
	7.25	-0.414661161933296\\
	7.3	-0.415931739553088\\
	7.31	-0.416414027954715\\
	7.37	-0.417874484138153\\
	7.38	-0.418345254675032\\
	7.44	-0.419734117561475\\
	7.45	-0.420192893530853\\
	7.5	-0.421293789946475\\
	7.51000000000001	-0.421740983046604\\
	7.57000000000001	-0.422995663618607\\
	7.58	-0.42343157340369\\
	7.64	-0.42462542061682\\
	7.65000000000001	-0.425051966805157\\
	7.70999999999999	-0.426191904598284\\
	7.72	-0.426608843788316\\
	7.78	-0.427694689109018\\
	7.79000000000001	-0.428102745025043\\
	7.84	-0.428967386598501\\
	7.84999999999999	-0.4293667725413\\
	7.91	-0.430355373074505\\
	7.92	-0.430745320330246\\
	7.98	-0.431680070691868\\
	7.98999999999999	-0.432060847659471\\
	8.04000000000001	-0.432797637656435\\
	8.05	-0.433169497948313\\
	8.11	-0.434004144754539\\
	8.12	-0.434366483657158\\
	8.18000000000001	-0.4351450021218\\
	8.19	-0.435497258243998\\
	8.23999999999999	-0.43610011482204\\
	8.25	-0.436442565295593\\
	8.31	-0.437109640651258\\
	8.32000000000001	-0.437441417666861\\
	8.38	-0.438046696562765\\
	8.39	-0.438367583475056\\
	8.44	-0.438823983352506\\
	8.45	-0.439134477031232\\
	8.58	-0.440349125656965\\
	8.59	-0.440637099052481\\
	8.64	-0.440942611018968\\
	8.65000000000001	-0.441219760734668\\
	8.78	-0.442033654876909\\
	8.79000000000001	-0.442287344947758\\
	9.06	-0.44349132708858\\
	9.34	-0.44343916720436\\
	9.62	-0.442049960853936\\
	9.76000000000001	-0.440840045834435\\
	10.05	-0.43716974555575\\
	10.32	-0.432103402221699\\
	10.42	-0.429998426660504\\
	10.71	-0.422343546861654\\
	10.98	-0.413721224889869\\
	11.12	-0.408745325988235\\
	11.44	-0.396102611606267\\
	11.55	-0.39149040119824\\
	11.89	-0.375864688587683\\
	11.99	-0.371105014964584\\
	12.4	-0.350096746278609\\
	12.76	-0.330711681378915\\
	13.3	-0.30101850006092\\
	13.33	-0.299453813577273\\
	13.5	-0.290350007920168\\
	13.67	-0.281442780614654\\
	13.84	-0.272736346798965\\
	14.01	-0.264219112301618\\
	14.18	-0.25587506835312\\
	14.36	-0.247197510800248\\
	14.73	-0.229789895855376\\
	15.2	-0.208293698599832\\
	15.52	-0.194132596059035\\
	15.79	-0.182525006839029\\
	16	-0.173752108797814\\
	16.36	-0.159170265775188\\
	16.6	-0.149807286700465\\
	16.85	-0.140350554707098\\
	17.11	-0.130845042405909\\
	17.38	-0.121328601159661\\
	17.66	-0.11184638411963\\
	17.86	-0.105310334938167\\
	18.16	-0.0958799342333378\\
	18.48	-0.0863207292004233\\
	18.81	-0.0770024400131604\\
	19.14	-0.0682258826466153\\
	19.5	-0.0592461185284208\\
	19.88	-0.0504553573062196\\
	20.29	-0.041714584279049\\
	20.73	-0.0332107472295746\\
	21.21	-0.0248873527871183\\
	21.73	-0.0169906091612404\\
	22.29	-0.0096499208972034\\
	22.93	-0.00259618007405038\\
	23.67	0.00387632435030127\\
	24.53	0.00947934403409079\\
	25.14	0.0123638140009632\\
	25.18	0.0125495463115612\\
	25.76	0.0145692227307848\\
	26.28	0.0158758650740793\\
	26.79	0.0167703338154155\\
	27.55	0.0174977218662633\\
	28.37	0.0176323419423028\\
	29.32	0.0171629588932092\\
	30.57	0.01589796926541\\
	32.84	0.0128077495604373\\
	34.47	0.0105495432610923\\
	36.04	0.00862995817264789\\
	37.74	0.00691654195259162\\
	39.14	0.00580002228889498\\
	41.31	0.00446413008123159\\
	43.25	0.00362186984021662\\
	45.67	0.00284041021963333\\
	49.9	0.00192585499934239\\
	54.46	0.00126903150049884\\
	61.25	0.00067434598624061\\
	71.41	0.000263047707633746\\
	80	0.00012865400790929\\
};

\addplot [color=red, dotted, line width=0.7pt, forget plot]
table[row sep=crcr]{%
	0	-5.19380563446248e-06\\
	0.349999999999994	-0.0301340095348905\\
	0.939999999999998	-0.0785970692253528\\
	1.23999999999999	-0.10259139817083\\
	1.27	-0.105084988480286\\
	1.48	-0.121825135370614\\
	1.51000000000001	-0.124325100337231\\
	1.72	-0.141146752359674\\
	1.75	-0.143669542373502\\
	1.94	-0.159063769007517\\
	1.97	-0.161618189841221\\
	2.16	-0.177189116294215\\
	2.31	-0.189633416610533\\
	2.34	-0.192218573026992\\
	2.53	-0.207813755236387\\
	2.56	-0.210382906940126\\
	2.73999999999999	-0.225056505858589\\
	2.77	-0.227608118379123\\
	2.95999999999999	-0.24294147685832\\
	2.98999999999999	-0.245473795441811\\
	3.17	-0.259882254289849\\
	3.2	-0.262394420842909\\
	3.38	-0.276662714630262\\
	3.41	-0.279153469425978\\
	3.59	-0.293274170176048\\
	3.62	-0.295742188070292\\
	3.8	-0.30970764809436\\
	3.82000000000001	-0.311401796925338\\
	4.01000000000001	-0.325953725074726\\
	4.03	-0.327631951567881\\
	4.22	-0.342002247662478\\
	4.23999999999999	-0.343663425342797\\
	4.43000000000001	-0.357842732350633\\
	4.45	-0.359485679869081\\
	4.64	-0.37346433445552\\
	4.66	-0.3750878113227\\
	4.84999999999999	-0.388855851624953\\
	4.87	-0.390458553109696\\
	5.05	-0.40332028458387\\
	5.07000000000001	-0.404901225654157\\
	5.26000000000001	-0.41822891457052\\
	5.28	-0.419786344335691\\
	5.47	-0.432871831313079\\
	5.48999999999999	-0.434404255030543\\
	5.68000000000001	-0.44723616654511\\
	5.7	-0.448742030933545\\
	5.89	-0.461308617962231\\
	5.91	-0.462786300859008\\
	6.09999999999999	-0.475075450574721\\
	6.12	-0.476523263518004\\
	6.3	-0.487920486394444\\
	6.32000000000001	-0.489337454282264\\
	6.51000000000001	-0.501049248867702\\
	6.53	-0.502432850147272\\
	6.72	-0.513829212120569\\
	6.73999999999999	-0.51517756195372\\
	6.94	-0.526796609323\\
	6.95999999999999	-0.528106756697156\\
	7.15000000000001	-0.538814062735995\\
	7.17	-0.540084918163842\\
	7.36	-0.550434311929422\\
	7.38	-0.551663784552531\\
	7.57000000000001	-0.56164047389592\\
	7.59	-0.562826416031996\\
	7.86	-0.576375276298492\\
	7.88	-0.577498386115948\\
	8.47	-0.604586401611897\\
	8.98999999999999	-0.625149852673943\\
	9.44	-0.640213501522595\\
	9.83	-0.651051406527159\\
	10.17	-0.658708651839035\\
	10.49	-0.664375735772978\\
	10.79	-0.668262662453202\\
	11.05	-0.670433775247474\\
	11.3	-0.671472806379271\\
	11.53	-0.671489905784654\\
	11.74	-0.670698922376573\\
	11.94	-0.669200861489145\\
	12.13	-0.667120988439294\\
	12.31	-0.664534816424109\\
	12.48	-0.661528404668317\\
	12.64	-0.658244171607492\\
	12.7	-0.656824529852585\\
	12.84	-0.653218909024531\\
	13	-0.648686560247882\\
	13.04	-0.647381287436133\\
	13.18	-0.642801993284621\\
	13.22	-0.641322550015346\\
	13.36	-0.636176933954189\\
	13.39	-0.634912367233838\\
	13.53	-0.629318123835944\\
	13.56	-0.627951918281966\\
	13.71	-0.621519569290697\\
	13.74	-0.620050155024913\\
	13.89	-0.613174933345306\\
	13.92	-0.611606293252734\\
	14.07	-0.604310913333734\\
	14.09	-0.603135147862162\\
	14.25	-0.594958925168683\\
	14.35	-0.589590674904898\\
	14.37	-0.588317092720771\\
	14.53	-0.579556148482681\\
	14.55	-0.578224625695981\\
	14.71	-0.56913826065788\\
	14.82	-0.562663979791353\\
	14.84	-0.561253116727968\\
	15.01	-0.551160569521443\\
	15.03	-0.549706264914775\\
	15.21	-0.538810552992373\\
	15.23	-0.53732056345693\\
	15.4	-0.5268770284315\\
	15.42	-0.525361300818176\\
	15.61	-0.513672999098418\\
	15.63	-0.51214303420906\\
	15.81	-0.501109369733314\\
	15.83	-0.499577190424191\\
	16.02	-0.488138514525801\\
	16.04	-0.48661930043184\\
	16.23	-0.475513572864301\\
	16.24	-0.474574113170959\\
	16.44	-0.463358595673526\\
	16.45	-0.462427646030164\\
	16.66	-0.451174953617823\\
	16.67	-0.450252486221117\\
	16.89	-0.438999694219589\\
	16.9	-0.438085219304668\\
	17	-0.433333156654285\\
	17.01	-0.432421084854155\\
	17.11	-0.427770589512761\\
	17.12	-0.426860313044784\\
	17.22	-0.422305005559295\\
	17.23	-0.421395890746837\\
	17.33	-0.416929089291529\\
	17.34	-0.416020482129781\\
	17.44	-0.411635249597992\\
	17.45	-0.410726479922246\\
	17.55	-0.406415657375689\\
	17.56	-0.405506041622417\\
	17.66	-0.401262268127567\\
	17.67	-0.400351110695084\\
	17.77	-0.396166829721579\\
	17.78	-0.395253422949708\\
	17.87	-0.391532431007519\\
	17.88	-0.390615172466752\\
	17.98	-0.386523220711098\\
	17.99	-0.385602320730328\\
	18.09	-0.381546657143872\\
	18.1	-0.380621360569606\\
	18.19	-0.376995585188183\\
	18.2	-0.376064430264933\\
	18.29	-0.372455292604613\\
	18.3	-0.371517610034459\\
	18.39	-0.367918292671689\\
	18.4	-0.366973385714488\\
	18.49	-0.363376761448464\\
	18.5	-0.362423900716976\\
	18.59	-0.358822488071283\\
	18.6	-0.357860905990194\\
	18.69	-0.354246823233225\\
	18.7	-0.353275708339012\\
	18.79	-0.349640627959047\\
	18.8	-0.348659119270266\\
	18.88	-0.345403263420138\\
	18.89	-0.344410309648225\\
	18.98	-0.340711182218612\\
	18.99	-0.339706043198305\\
	19.07	-0.336377918098719\\
	19.08	-0.335359622071394\\
	19.17	-0.331560967257744\\
	19.18	-0.330528470760811\\
	19.26	-0.327094837811913\\
	19.27	-0.326047265961336\\
	19.35	-0.322559399694143\\
	19.36	-0.321498671062798\\
	19.44	-0.317976327587942\\
	19.45	-0.316904870152371\\
	19.53	-0.313368328941934\\
	19.54	-0.312287594146639\\
	19.63	-0.308296367017846\\
	19.64	-0.307207217468317\\
	19.72	-0.303651293175676\\
	19.73	-0.302553538874534\\
	19.82	-0.298547250626257\\
	19.83	-0.297441787630291\\
	19.91	-0.293878833036445\\
	19.92	-0.292765537688013\\
	20.01	-0.288759821854867\\
	20.02	-0.287639840751694\\
	20.11	-0.283642250968924\\
	20.12	-0.28251626267317\\
	20.21	-0.278533498257005\\
	20.22	-0.27740225869772\\
	20.31	-0.273441401564611\\
	20.32	-0.272305592344026\\
	20.41	-0.268372039086714\\
	20.42	-0.267232025808909\\
	20.51	-0.263326109430153\\
	20.52	-0.262181810805984\\
	20.61	-0.258308243584196\\
	20.62	-0.257160756175736\\
	20.72	-0.252899478297607\\
	20.73	-0.251749719115068\\
	20.82	-0.247946074245448\\
	20.83	-0.246791899819655\\
	20.92	-0.243014472328909\\
	20.93	-0.241855916814103\\
	21.03	-0.237688253292575\\
	21.04	-0.236526046232342\\
	21.13	-0.232799918616337\\
	21.14	-0.231633270197122\\
	21.23	-0.227931094347952\\
	21.24	-0.226759810223584\\
	21.34	-0.222670560072274\\
	21.35	-0.221495335013984\\
	21.44	-0.217834625196446\\
	21.45	-0.216654565174181\\
	21.54	-0.213012027690567\\
	21.55	-0.211826912671924\\
	21.64	-0.208199919260181\\
	21.65	-0.207009629927413\\
	21.74	-0.203396482136696\\
	21.75	-0.202201012609351\\
	21.85	-0.198201489184157\\
	21.86	-0.197001778026277\\
	21.95	-0.193415875315026\\
	21.96	-0.192211086171596\\
	22.05	-0.188638310658192\\
	22.06	-0.18742849894862\\
	22.15	-0.183868821022458\\
	22.16	-0.182654044506137\\
	22.25	-0.17910731290489\\
	22.26	-0.177887622176357\\
	22.35	-0.174353660850215\\
	22.36	-0.173129110736014\\
	22.45	-0.16960777641178\\
	22.46	-0.16837842481614\\
	22.55	-0.164869594544129\\
	22.56	-0.163635500881014\\
	22.65	-0.160139086770812\\
	22.66	-0.15890031181047\\
	22.75	-0.155416222655091\\
	22.76	-0.154172826217689\\
	22.85	-0.150700933381387\\
	22.86	-0.149452974636461\\
	22.95	-0.145993151790179\\
	22.96	-0.144740691806874\\
	23.05	-0.141292821574282\\
	23.06	-0.140035922484699\\
	23.15	-0.136599890956262\\
	23.16	-0.135338615617556\\
	23.25	-0.131914309916951\\
	23.26	-0.130648721780048\\
	23.35	-0.127236029025752\\
	23.36	-0.125966192112031\\
	23.45	-0.122564998960485\\
	23.46	-0.121290977858123\\
	23.55	-0.117901170287411\\
	23.56	-0.116623030173969\\
	23.65	-0.113244493423707\\
	23.66	-0.111962300090809\\
	23.75	-0.108594918649302\\
	23.76	-0.107308738528559\\
	23.85	-0.103952396143853\\
	23.86	-0.102662296330251\\
	23.95	-0.0993168760030727\\
	23.96	-0.0980229242711346\\
	24.05	-0.0946883083126835\\
	24.06	-0.0933905731481275\\
	24.15	-0.09006664322267\\
	24.16	-0.088765193831648\\
	24.25	-0.0854518309454591\\
	24.26	-0.0841467372688527\\
	24.35	-0.0808438217901255\\
	24.36	-0.0795351545164209\\
	24.45	-0.076242566212926\\
	24.46	-0.0749303967929791\\
	24.55	-0.071648014848364\\
	24.56	-0.0703324155007721\\
	24.64	-0.0674236021116741\\
	24.65	-0.0661037723449596\\
	24.74	-0.0628412119082213\\
	24.75	-0.0615180965776858\\
	24.84	-0.0582653927081083\\
	24.85	-0.0569390663554827\\
	24.94	-0.0536960960667727\\
	24.95	-0.0523666340310598\\
	25.04	-0.0491332737895789\\
	25.05	-0.0478007522104775\\
	25.14	-0.0445768779673443\\
	25.15	-0.0432413737938475\\
	25.24	-0.04002686108646\\
	25.25	-0.0386884520949025\\
	25.34	-0.0354831761294605\\
	25.35	-0.0341419408967312\\
	25.44	-0.030945776305515\\
	25.45	-0.0296017941394524\\
	25.54	-0.0264146146316904\\
	25.55	-0.0250679655321449\\
	25.63	-0.0222427573922772\\
	25.64	-0.0208926250984263\\
	25.73	-0.0177229681462592\\
	25.74	-0.0163703285866603\\
	25.83	-0.0132092907973345\\
	25.84	-0.011854226208655\\
	25.93	-0.00870167935579502\\
	25.94	-0.00734427272243465\\
	26.03	-0.00420008817954454\\
	26.04	-0.00284042323875155\\
	26.13	0.000295528002638434\\
	26.14	0.00165736678705741\\
	26.23	0.00478521422508038\\
	26.24	0.00614914166084191\\
	26.33	0.00926901519184753\\
	26.34	0.0106349453524075\\
	26.43	0.0137469764290756\\
	26.44	0.0151148230930289\\
	26.52	0.0178741877346056\\
	26.53	0.0192447817911727\\
	26.62	0.022341441123686\\
	26.63	0.0237137785913148\\
	26.72	0.0268029725522609\\
	26.73	0.028176964624322\\
	26.82	0.0312588246290346\\
	26.83	0.0326343818430388\\
	26.92	0.0357090399427591\\
	26.93	0.0370860721500321\\
	27.02	0.0401536607806037\\
	27.03	0.0415320772234224\\
	27.12	0.0445927296032664\\
	27.13	0.0459724388812219\\
	27.22	0.0490262886149253\\
	27.23	0.0504071987149217\\
	27.31	0.0531157119997658\\
	27.32	0.0544986672434078\\
	27.41	0.0575390954636674\\
	27.42	0.0589230680904933\\
	27.51	0.0619567284167317\\
	27.52	0.0633409164242096\\
	27.61	0.0663528493545726\\
	27.62	0.0677358750914436\\
	27.71	0.0707158370274215\\
	27.72	0.0720961616907374\\
	27.81	0.0750369778550777\\
	27.82	0.0764142508338352\\
	27.91	0.0793150828436495\\
	27.92	0.0806891847555278\\
	28.01	0.0835499813851328\\
	28.02	0.0849207830546561\\
	28.11	0.0877399421313783\\
	28.12	0.0891073932569526\\
	28.22	0.0921912242039298\\
	28.23	0.0935539532451912\\
	28.32	0.0962858878747284\\
	28.33	0.0976450744821449\\
	28.42	0.100334983297344\\
	28.43	0.101690437890326\\
	28.52	0.104337782911585\\
	28.53	0.10568936603589\\
	28.62	0.108293952693188\\
	28.63	0.109641513304709\\
	28.72	0.112203163998601\\
	28.73	0.113546530571625\\
	28.82	0.116064922049517\\
	28.83	0.117403946511601\\
	28.92	0.119879030720981\\
	28.93	0.121213618127996\\
	29.02	0.123645784034878\\
	29.03	0.124975787647159\\
	29.12	0.127364997321692\\
	29.13	0.128690257501546\\
	29.22	0.131036347784715\\
	29.23	0.132356759947811\\
	29.32	0.13466014260365\\
	29.33	0.135975536771312\\
	29.42	0.13823595723909\\
	29.43	0.139546181896634\\
	29.52	0.141763576942097\\
	29.53	0.143068520795964\\
	29.62	0.145243197775599\\
	29.63	0.14654264983632\\
	29.71	0.148438441450466\\
	29.72	0.149733758046921\\
	29.81	0.151825301566092\\
	29.82	0.15311484809888\\
	29.91	0.15516311360912\\
	29.92	0.156446714606076\\
	30	0.158230090995374\\
	30.01	0.159508926884712\\
	30.1	0.161473725158629\\
	30.11	0.162746373264227\\
	30.2	0.164667077415658\\
	30.21	0.165933310126533\\
	30.29	0.167603246778668\\
	30.3	0.168864189672917\\
	30.39	0.170693365066029\\
	30.4	0.171946721368556\\
	30.48	0.17353390561388\\
	30.49	0.174781895454586\\
	30.58	0.176523669708871\\
	30.59	0.177764594665746\\
	30.67	0.179273912512471\\
	30.68	0.180509128506714\\
	30.76	0.181980756458515\\
	30.77	0.183210107012428\\
	30.86	0.184821235015306\\
	30.87	0.186043001609022\\
	30.95	0.187435430462685\\
	30.96	0.18865100677057\\
	31.04	0.190005096943722\\
	31.05	0.191214321713929\\
	31.13	0.19252984576049\\
	31.14	0.193732562447607\\
	31.22	0.19500929541833\\
	31.23	0.196205350418168\\
	31.32	0.197595910488729\\
	31.33	0.198783557253932\\
	31.41	0.199978455827917\\
	31.42	0.201159115691567\\
	31.5	0.202314343876864\\
	31.51	0.203487820324142\\
	31.59	0.204602752819682\\
	31.6	0.205768855808586\\
	31.68	0.206839756428579\\
	31.69	0.20799787218715\\
	31.77	0.209026105158543\\
	31.78	0.2101764093307\\
	31.86	0.211161812454449\\
	31.87	0.212304166296647\\
	31.95	0.213246439710048\\
	31.96	0.214380702916415\\
	32.04	0.215279575874305\\
	32.05	0.216405609710932\\
	32.13	0.217260841405079\\
	32.14	0.218378508753617\\
	32.21	0.219090112502172\\
	32.22	0.220200939151496\\
	32.3	0.220972119273668\\
	32.31	0.22207431675352\\
	32.39	0.222801394999351\\
	32.4	0.223894751356283\\
	32.48	0.224575760008946\\
	32.49	0.225659973403467\\
	32.57	0.226295868732564\\
	32.58	0.227371021806789\\
	32.65	0.227889638843038\\
	32.66	0.228957303894546\\
	32.74	0.229506747226594\\
	32.75	0.230565116759337\\
	32.83	0.231069324070177\\
	32.84	0.232118284142317\\
	32.91	0.232521590288954\\
	32.92	0.233562727250941\\
	33	0.233980443991854\\
	33.01	0.23501195627945\\
	33.08	0.235339734395637\\
	33.09	0.236363168915631\\
	33.17	0.236693756641202\\
	33.18	0.237707230352029\\
	33.25	0.237958111623968\\
	33.26	0.238963324172701\\
	33.34	0.239206397686544\\
	33.35	0.240201548204539\\
	33.42	0.24037616422919\\
	33.43	0.241362873256961\\
	33.5	0.241501298079854\\
	33.51	0.242479480056886\\
	33.58	0.242581875824314\\
	33.59	0.243551443926549\\
	33.67	0.243625041993155\\
	33.68	0.244584158574199\\
	33.75	0.244610453416712\\
	33.76	0.245560768849145\\
	33.83	0.245550924337749\\
	33.84	0.246492367960457\\
	33.91	0.246446604961079\\
	33.92	0.247379110966804\\
	33.99	0.247297640417571\\
	34	0.24822116725403\\
	34.08	0.248085454065375\\
	34.09	0.248998203943472\\
	34.16	0.248841912943249\\
	34.17	0.249745504168388\\
	34.24	0.249553684799835\\
	34.25	0.250448067489486\\
	34.32	0.25022102324634\\
	34.33	0.251106152248454\\
	34.39	0.250883008148378\\
	34.4	0.251760569117863\\
	34.47	0.251467381053203\\
	34.48	0.252335630834182\\
	34.55	0.252007491722395\\
	34.56	0.252866309410734\\
	34.63	0.252503582592382\\
	34.64	0.253352939704087\\
	34.71	0.252955993869236\\
	34.72	0.253795866995873\\
	34.79	0.253365020257235\\
	34.8	0.25419544720836\\
	34.86	0.253799165885667\\
	34.87	0.254621796723356\\
	34.94	0.254127370378626\\
	34.95	0.254940391324709\\
	35.02	0.254412654574182\\
	35.03	0.25521605915641\\
	35.09	0.254736808121748\\
	35.1	0.255532294251609\\
	35.17	0.254942772342048\\
	35.18	0.255728715320771\\
	35.25	0.255106532291123\\
	35.26	0.25588277558262\\
	35.32	0.255323167366967\\
	35.33	0.256091340946909\\
	35.4	0.255408184830202\\
	35.41	0.256166674653016\\
	35.47	0.255555405367744\\
	35.48	0.256305887992511\\
	35.55	0.255564032207303\\
	35.56	0.256304980048725\\
	35.62	0.255644140199195\\
	35.63	0.256376959970112\\
	35.7	0.255576907902395\\
	35.71	0.256300030941432\\
	35.77	0.25558964847221\\
	35.78	0.256304726332502\\
	35.84	0.255572010101844\\
	35.85	0.256279148927192\\
	35.92	0.255396612011353\\
	35.93	0.256094048769825\\
	35.99	0.255313788742384\\
	36	0.256003160276492\\
	36.06	0.255201517494697\\
	36.07	0.255882935337993\\
	36.14	0.254920430383748\\
	36.15	0.255592169289528\\
	36.21	0.254744224916934\\
	36.22	0.25540785655916\\
	36.28	0.254538375909519\\
	36.29	0.255193968876981\\
	36.35	0.25430429412701\\
	36.36	0.254952001647965\\
	36.42	0.254041495015201\\
	36.43	0.254681111213074\\
	36.5	0.253593576422134\\
	36.51	0.254223672320123\\
	36.57	0.253270272096827\\
	36.58	0.253892393651157\\
	36.64	0.25291998981362\\
	36.65	0.253534332396228\\
	36.71	0.252542665792546\\
	36.72	0.253149018712477\\
	36.78	0.252137943349297\\
	36.79	0.252736371313148\\
	36.85	0.251707325359476\\
	36.86	0.252298031211993\\
	36.92	0.251250349989263\\
	36.93	0.251833131999362\\
	36.99	0.250766948880482\\
	37	0.251341874544849\\
	37.06	0.250258803215374\\
	37.07	0.250826082802817\\
	37.13	0.24972504876564\\
	37.14	0.25028449056154\\
	37.2	0.249165897165156\\
	37.21	0.249717569720403\\
	37.27	0.248581898539982\\
	37.28	0.249125850488923\\
	37.34	0.247974373298064\\
	37.35	0.248510862706524\\
	37.41	0.247343352871596\\
	37.42	0.247872181723082\\
	37.48	0.246688550883547\\
	37.49	0.247209768266643\\
	37.54	0.246211715515898\\
	37.55	0.246726882269257\\
	37.61	0.245514573667265\\
	37.62	0.246022190017143\\
	37.68	0.24479468747252\\
	37.69	0.245294801732172\\
	37.75	0.244053593699377\\
	37.76	0.244546477851529\\
	37.82	0.243290895684297\\
	37.83	0.243776369732487\\
	37.89	0.242506585954757\\
	37.9	0.242984697017491\\
	37.95	0.241916581444514\\
	37.96	0.242388825189622\\
	38.02	0.241093972427436\\
	38.03	0.241558945410205\\
	38.09	0.240250854213741\\
	38.1	0.24070860066638\\
	38.15	0.239609317475441\\
	38.16	0.240061301876665\\
	38.22	0.238729952716298\\
	38.23	0.239174818241153\\
	38.29	0.237831206728842\\
	38.3	0.238268993954605\\
	38.36	0.236914830438991\\
	38.37	0.237345829545319\\
	38.42	0.236207966621492\\
	38.43	0.236633151560298\\
	38.49	0.235256828720807\\
	38.5	0.235675094288055\\
	38.55	0.234520472002259\\
	38.56	0.234933215851058\\
	38.62	0.233537290577047\\
	38.63	0.233943251633605\\
	38.69	0.23253703498375\\
	38.7	0.232936245183126\\
	38.75	0.231757469577886\\
	38.76	0.232151290988256\\
	38.82	0.230727299229002\\
	38.83	0.231114519428104\\
	38.88	0.2299203146666\\
	38.89	0.230301989140699\\
	38.95	0.228861277408697\\
	38.96	0.229236674997068\\
	39.01	0.228029106365327\\
	39.02	0.22839909635475\\
	39.08	0.226941912278022\\
	39.09	0.227305495479172\\
	39.15	0.225841332215694\\
	39.16	0.226198825160864\\
	39.21	0.224972464631577\\
	39.22	0.22532470447527\\
	39.28	0.223845887509256\\
	39.29	0.224191910655847\\
	39.34	0.222954586752905\\
	39.35	0.223295649980727\\
	39.4	0.222052901075557\\
	39.41	0.22238885613632\\
	39.47	0.220891319228841\\
	39.48	0.221221235396811\\
	39.53	0.219968747952549\\
	39.54	0.220293851102525\\
	39.6	0.218785397085114\\
	39.61	0.219104656521765\\
	39.66	0.21784313350247\\
	39.67	0.218157420819779\\
	39.73	0.216638838197753\\
	39.74	0.216947602946874\\
	39.79	0.215678374673146\\
	39.8	0.215982263359635\\
	39.85	0.214709976630473\\
	39.86	0.215009259509756\\
	39.92	0.213477793264033\\
	39.93	0.213771526612206\\
	39.98	0.212491720437953\\
	39.99	0.212780727280062\\
	40.05	0.211241873310328\\
	40.06	0.211525635522975\\
	40.11	0.210240157485899\\
	40.12	0.210519364955502\\
	40.17	0.209230970507889\\
	40.18	0.209505603721112\\
	40.24	0.207956994739263\\
	40.25	0.208226580910562\\
	40.3	0.20693375180187\\
	40.31	0.207198944556524\\
	40.36	0.205903900023571\\
	40.37	0.206164672582432\\
	40.43	0.20460860699896\\
	40.44	0.204864529901656\\
	40.49	0.203566321130552\\
	40.5	0.203817927369101\\
	40.55	0.202518723952039\\
	40.56	0.202766265395155\\
	40.62	0.201205247597528\\
	40.63	0.201447945133239\\
	40.68	0.200145794194441\\
	40.69	0.200384332362432\\
	40.74	0.199081548667351\\
	40.75	0.199316184689124\\
	40.8	0.198012265206131\\
	40.81	0.198242905807916\\
	40.87	0.196677518899051\\
	40.88	0.196903510652646\\
	40.93	0.195598942039226\\
	40.94	0.195821194017128\\
	40.99	0.194516261946461\\
	41	0.194734601914718\\
	41.05	0.193429996221823\\
	41.06	0.193644657303096\\
	41.12	0.19207892948954\\
	41.13	0.192289243538738\\
	41.18	0.190984631412803\\
	41.19	0.19119119212553\\
	41.24	0.189887244116079\\
	41.25	0.190090288930293\\
	41.3	0.18878653239365\\
	41.31	0.188985984765822\\
	41.37	0.187422483418842\\
	41.38	0.187617784152948\\
	41.43	0.186316145691393\\
	41.44	0.18650809169435\\
	41.49	0.185207349310531\\
	41.5	0.185395789750103\\
	41.55	0.184096581402045\\
	41.56	0.184281729390278\\
	41.61	0.182983325094753\\
	41.62	0.183165117659485\\
	41.68	0.181608869572287\\
	41.69	0.181786810572575\\
	41.74	0.180492128169419\\
	41.75	0.180666936594335\\
	41.8	0.179373653090153\\
	41.81	0.179545267251811\\
	41.86	0.178253710130932\\
	41.87	0.17842211588173\\
	41.92	0.177132490008958\\
	41.93	0.1772979145383\\
	41.98	0.176010482109646\\
	41.99	0.176172787317412\\
	42.05	0.174631317587782\\
	42.06	0.174790269843029\\
	42.11	0.173508193371603\\
	42.12	0.173664182407592\\
	42.17	0.172384399492699\\
	42.18	0.172537421441632\\
	42.23	0.171260356592583\\
	42.24	0.171410613855727\\
	42.29	0.170136162875025\\
	42.3	0.170283540504641\\
	42.35	0.169012250403256\\
	42.36	0.169156925358365\\
	42.41	0.167888176099808\\
	42.42	0.168030105791061\\
	42.48	0.166511600876376\\
	42.49	0.166650420206821\\
	42.54	0.165388881266765\\
	42.55	0.165525148740542\\
	42.6	0.1642666412057\\
	42.61	0.16440031434972\\
	42.66	0.163145122720806\\
	42.67	0.163276197175819\\
	42.72	0.162024506888244\\
	42.74	0.16190317169729\\
	42.83	0.160035594702208\\
	42.9	0.158670430850336\\
	42.92	0.158544075381343\\
	43.03	0.156195468541171\\
	43.05	0.156066036626768\\
	43.15	0.15397548320604\\
	43.17	0.153843269026865\\
	43.27	0.151763310500201\\
	43.29	0.151628420462544\\
	43.39	0.149560341758928\\
	43.41	0.149423262804262\\
	43.51	0.147367716263972\\
	43.53	0.147228467806755\\
	43.63	0.145185718406196\\
	43.65	0.14504459453515\\
	43.75	0.143015926882185\\
	43.77	0.142873088399483\\
	43.88	0.140625794290031\\
	43.9	0.140481376324232\\
	44	0.138484119923731\\
	44.02	0.138338487988847\\
	44.12	0.136356737006849\\
	44.14	0.136210049928664\\
	44.24	0.134244977140582\\
	44.26	0.134097372998909\\
	44.36	0.132149007961232\\
	44.38	0.132000566955227\\
	44.48	0.130069707055611\\
	44.5	0.129920718892635\\
	44.6	0.128007874886961\\
	44.62	0.127858504076215\\
	44.72	0.125963663280544\\
	44.74	0.125814057729983\\
	44.84	0.123938167204173\\
	44.86	0.123788439655385\\
	44.96	0.12193146050798\\
	44.98	0.121781677529071\\
	45.08	0.119944243725485\\
	45.1	0.119794648634183\\
	45.2	0.11797712334311\\
	45.22	0.117827838966392\\
	45.32	0.116030152329557\\
	45.34	0.115881314063444\\
	45.45	0.113901945254909\\
	45.47	0.113753808203825\\
	45.57	0.111999481212933\\
	45.59	0.111851941387712\\
	45.69	0.110118615734692\\
	45.71	0.109971910656768\\
	45.81	0.108259717031459\\
	45.83	0.108113841091097\\
	45.93	0.106422761887373\\
	45.95	0.106277883579267\\
	46.05	0.104608412918438\\
	46.07	0.1044645964067\\
	46.17	0.102816547549821\\
	46.19	0.102673834275748\\
	46.29	0.101047622056669\\
	46.31	0.100906201921774\\
	46.41	0.0993018587500671\\
	46.43	0.0991617162944749\\
	46.53	0.0975791594669744\\
	46.55	0.0974404345968196\\
	46.65	0.0958800094074661\\
	46.67	0.0957427439301171\\
	46.77	0.094204206630053\\
	46.79	0.0940684257298585\\
	46.89	0.0925520636114499\\
	46.91	0.0924179148855018\\
	47.01	0.0909236693077986\\
	47.03	0.0907911267139809\\
	47.13	0.0893188540085816\\
	47.15	0.0891880316725064\\
	47.25	0.0877379504323699\\
	47.27	0.0876088702898699\\
	47.37	0.0861806881880085\\
	47.39	0.0860533605528531\\
	47.49	0.0846472563866172\\
	47.51	0.0845217923018424\\
	47.61	0.0831376313064851\\
	47.64	0.0828580123881011\\
	47.74	0.0814978960162023\\
	47.77	0.0812231558374066\\
	47.89	0.0797676520066659\\
	48.03	0.078040158654801\\
	48.16	0.0764917631262847\\
	48.19	0.0762324781167933\\
	48.34	0.0744343497802049\\
	48.37	0.0741817368271853\\
	48.52	0.0724290074385721\\
	48.55	0.072182999219649\\
	48.7	0.0704754034543953\\
	48.73	0.0702359313630723\\
	48.89	0.0684419978267954\\
	48.92	0.0682095204285389\\
	49.07	0.0665930666764893\\
	49.1	0.0663670233710008\\
	49.25	0.0647937754198296\\
	49.28	0.0645741580440671\\
	49.43	0.0630432691966121\\
	49.46	0.06282998143881\\
	49.61	0.0613409748854536\\
	49.64	0.0611339103894295\\
	49.8	0.0595718525315334\\
	49.83	0.05937133963495\\
	49.98	0.0579665145379948\\
	50.01	0.0577721455610174\\
	50.16	0.0564069851894544\\
	50.19	0.0562185990885951\\
	50.34	0.0548923481893269\\
	50.37	0.0547098761423683\\
	50.53	0.0533197350763714\\
	50.56	0.0531433761221933\\
	50.71	0.0518948741615191\\
	50.75	0.0516256413468881\\
	50.89	0.0505119188147347\\
	50.93	0.0502510114975081\\
	51.08	0.0490766011627954\\
	51.12	0.0488242740084956\\
	51.31	0.0474428331743724\\
	51.5	0.0461056750321376\\
	51.54	0.0458714866770151\\
	51.75	0.0444172830266751\\
	51.79	0.0441934849997949\\
	51.99	0.0428806707234628\\
	52.03	0.0426663755776246\\
	52.24	0.0413287658983563\\
	52.28	0.0411240476209969\\
	52.49	0.0398429023307614\\
	52.53	0.0396474183222466\\
	52.73	0.0384920463261125\\
	52.77	0.0383050481227798\\
	52.98	0.03712791610522\\
	53.02	0.0369494185883212\\
	53.22	0.0358883888187336\\
	53.27	0.0356520948102883\\
	53.47	0.0346365834257369\\
	53.52	0.0344107379964527\\
	53.75	0.033344030699638\\
	54.02	0.0320861493925975\\
	54.27	0.0309983888270438\\
	54.32	0.0308029083325039\\
	54.58	0.029717367345711\\
	54.63	0.0295325179946531\\
	54.89	0.0285038719399466\\
	54.94	0.0283289883790729\\
	55.2	0.0273537088086186\\
	55.27	0.0271665905255674\\
	55.57	0.0260609396773788\\
	55.64	0.0258866287464912\\
	55.95	0.0248726508015551\\
	56.32	0.0236661810753702\\
	56.39	0.0235150022643893\\
	56.75	0.0224302367601155\\
	56.83	0.0222530961177654\\
	57.18	0.0212783241169916\\
	57.26	0.0211140461519364\\
	57.62	0.0201699705086895\\
	57.7	0.0200177309558995\\
	58.12	0.0190028439194663\\
	58.21	0.0188320666387654\\
	58.61	0.0179505860342033\\
	58.7	0.01779240371215\\
	59.13	0.0169389735546588\\
	59.61	0.0160093553536171\\
	59.71	0.0158482821214534\\
	60.17	0.0150385774313264\\
	60.27	0.014889762501781\\
	60.79	0.0140479358732364\\
	60.9	0.0138894958461862\\
	61.42	0.0131135534944633\\
	61.55	0.0129644924846843\\
	62.17	0.0120958206170201\\
	62.31	0.0119418331392893\\
	62.98	0.0110957107781502\\
	63.13	0.0109391721392313\\
	63.85	0.0101214806008443\\
	64.01	0.00996447642272358\\
	64.79	0.00916527408992351\\
	64.96	0.00900979689080827\\
	65.79	0.00824851340509269\\
	65.99	0.00809445627614025\\
	66.97	0.00728551469563854\\
	67.2	0.00711699468187987\\
	68.28	0.00634229613123694\\
	68.55	0.00617391885056406\\
	69.85	0.00537621170987279\\
	71.47	0.00450726192129025\\
	71.87	0.00432380520392428\\
	73.71	0.00354049333245143\\
	74.31	0.00332139494793182\\
	76.65	0.00257598687281302\\
	77.57	0.00233473045246058\\
	80	0.00179965852721864\\
};

\end{axis}
\end{tikzpicture}%

%% file: path_z3_eight_rev_MIQP.tex
%
%
\begin{tikzpicture}

\begin{axis}[%
width=\figurewidth,
height=\figureheight,
at={(0\figurewidth,0\figureheight)},
scale only axis,
xmin=0,
xmax=80,
xtick={0,20,40,60,80},
xlabel={$t$ [s]},
ymin=-8,
ymax=8,
ytick={-8,-4,0,4,8},
xlabel style={font=\color{white!15!black},at={(axis description cs:0.5,-0.2)},anchor=north},
ylabel style={font=\color{white!15!black},at={(axis description cs:-0.2,.5)},anchor=south},
ylabel={$\tilde z_3$ [m]},
axis background/.style={fill=white},
xmajorgrids,
ymajorgrids
]
\addplot [color=red, line width=0.7pt, forget plot]
  table[row sep=crcr]{%
0	-3.99995877000784\\
0.379999999999995	-3.99612318303554\\
0.769999999999996	-3.98466758291386\\
1.17	-3.96522538153546\\
1.55	-3.93920775386525\\
1.92	-3.90600013599753\\
2.25	-3.86891626243816\\
2.56	-3.82678101750238\\
2.88	-3.77523957469654\\
3.22	-3.71275167440581\\
3.56	-3.6426570808667\\
3.90000000000001	-3.5650245331947\\
4.25	-3.4773373180856\\
4.59999999999999	-3.38188265501346\\
4.95999999999999	-3.27580278899002\\
5.32000000000001	-3.16195150106461\\
5.7	-3.0336723167681\\
6.09999999999999	-2.89008389019926\\
6.51000000000001	-2.73438931062938\\
6.92	-2.57076993955594\\
7.38	-2.37871312344986\\
7.89	-2.15680474035253\\
8.44	-1.90908901338869\\
9.2	-1.55745152171176\\
10.68	-0.8710179135511\\
11.22	-0.630638547310184\\
11.67	-0.438459670583256\\
12.09	-0.267665204980801\\
12.55	-0.0896371277298016\\
13.12	0.123519130941446\\
13.46	0.243297025503722\\
13.8	0.355589955941724\\
14.14	0.460146225935006\\
14.49	0.559869403157649\\
14.85	0.654453878818373\\
15.23	0.746321075992896\\
15.66	0.842329771558511\\
16.17	0.948332569327519\\
16.68	1.04669288081705\\
17.27	1.15190997285633\\
17.79	1.23722632570886\\
18.36	1.32240681940922\\
18.86	1.38966311775523\\
19.41	1.45540762061999\\
19.94	1.51058952642624\\
20.49	1.55936979168065\\
20.98	1.59559691905275\\
21.53	1.62817231811084\\
22.04	1.65087308347695\\
22.6	1.66749940642588\\
23.1	1.67505486160302\\
23.7	1.67559123645367\\
24.29	1.66769396976781\\
24.85	1.65295424419637\\
25.48	1.62831492295153\\
26.04	1.59969595653143\\
26.73	1.55620702314731\\
27.37	1.50851137767299\\
28.09	1.44741781072501\\
28.96	1.3650725964619\\
30.2	1.23655178589848\\
33.97	0.835547487268045\\
35.16	0.723021844792058\\
36.25	0.629824701425179\\
37.05	0.56781838862895\\
37.85	0.51125874637674\\
38.64	0.460631696942642\\
39.48	0.412225860181664\\
40.48	0.361365162175076\\
41.36	0.322185660537684\\
42.65	0.273028979155299\\
44.06	0.228869035351863\\
45.75	0.186544293369522\\
47.23	0.156898309125097\\
49.05	0.12768020706244\\
50.59	0.107771578848968\\
52.8	0.0849270109797402\\
57.44	0.0521400354759152\\
61.79	0.0331734200469782\\
75.36	0.00823666262405709\\
80	0.0050394351622316\\
};
\addplot [color=blue, line width=0.7pt, forget plot]
  table[row sep=crcr]{%
0	-3.99995877000784\\
0.390000000000001	-3.99610690247906\\
0.799999999999997	-3.98450771868876\\
1.23	-3.96466925766335\\
1.64	-3.93833105470382\\
2.04000000000001	-3.90496332102335\\
2.40000000000001	-3.86743736019547\\
2.73999999999999	-3.82437869218406\\
3.05	-3.77760093397167\\
3.34999999999999	-3.72443519157606\\
3.69	-3.65583129136294\\
4.04000000000001	-3.57756980862976\\
4.40000000000001	-3.48929519513544\\
4.77	-3.39065702558034\\
5.16	-3.27838823463873\\
5.55	-3.15810012417555\\
5.97	-3.02021097864851\\
6.40000000000001	-2.87048666570527\\
7.72	-2.3988893498279\\
9.39	-1.78877655983229\\
9.90000000000001	-1.61256003565431\\
11.1	-1.20738086329973\\
11.48	-1.08962311810899\\
11.86	-0.979885129612157\\
12.24	-0.878596713802082\\
12.62	-0.785935821144008\\
13	-0.701853031287726\\
13.37	-0.628022594983904\\
13.75	-0.55969718165737\\
14.21	-0.485105254525948\\
14.7	-0.413311713134988\\
15.24	-0.342488903207126\\
15.84	-0.272492259117939\\
16.47	-0.207618727646178\\
17.09	-0.151661780776493\\
17.69	-0.104391687709821\\
18.38	-0.0579748434837342\\
19.09	-0.0185015299902744\\
19.76	0.011678233764016\\
20.5	0.0377604595348231\\
21.43	0.0611916679364413\\
22.31	0.0753342379833981\\
23.32	0.0839866271536778\\
24.71	0.0863446463512219\\
25.96	0.082446368889407\\
29.36	0.0607689631492576\\
32.61	0.0409071898169913\\
35.2	0.0296328465745006\\
42.38	0.0137988178736776\\
71.39	0.00108400978110978\\
80	0.000461998043221001\\
};
\addplot [color=green, line width=0.7pt, forget plot]
  table[row sep=crcr]{%
0	-3.99995877000784\\
0.4	-3.99622336104652\\
0.98	-3.9823533852231\\
4.41	-3.88653136062994\\
5.28	-3.85677079077055\\
6.05	-3.83774587490111\\
6.82	-3.82579006250074\\
};
\addplot [color=black, dotted, line width=1.0pt, forget plot]
table[row sep=crcr]{%
	0	-3.99995877000784\\
	0.25	-3.99821449493449\\
	0.519999999999996	-3.99287448926871\\
	0.799999999999997	-3.98390730897796\\
	1.09	-3.97110617383963\\
	1.37	-3.95537789834374\\
	1.65000000000001	-3.93621787446763\\
	1.91	-3.91512998396099\\
	2.16	-3.89157702290373\\
	2.40000000000001	-3.86562307548168\\
	2.63	-3.8373422427282\\
	2.84999999999999	-3.80680893785795\\
	3.06	-3.77415075944863\\
	3.26000000000001	-3.7395586341476\\
	3.45999999999999	-3.70131612742256\\
	3.7	-3.65162737933807\\
	3.92	-3.60285253477581\\
	4.16	-3.54621885675881\\
	4.41	-3.48349363038595\\
	4.66	-3.41707328228772\\
	4.92	-3.34417752085903\\
	5.17	-3.27055223620667\\
	5.43000000000001	-3.19044658545074\\
	5.68000000000001	-3.11016681217508\\
	5.95999999999999	-3.0165863988409\\
	6.16	-2.94659834140649\\
	6.36	-2.87272437944739\\
	7.25	-2.54125213815144\\
	7.54000000000001	-2.43706835879291\\
	7.88	-2.3152542417639\\
	8.17	-2.20738768928622\\
	9.36	-1.761101629823\\
	9.72	-1.63092753346443\\
	10.09	-1.50099468185327\\
	11.19	-1.11960279103266\\
	11.44	-1.03841932603602\\
	11.7	-0.957587535141016\\
	11.94	-0.886378266493139\\
	12.19	-0.815826875415482\\
	12.45	-0.746433910934229\\
	12.68	-0.688470894643459\\
	12.91	-0.63371130966182\\
	13.14	-0.582114426265193\\
	13.38	-0.531575491268583\\
	13.65	-0.478617739468973\\
	13.91	-0.431090219679135\\
	14.19	-0.383136036635975\\
	14.5	-0.33337471557175\\
	14.86	-0.279310094235655\\
	15.19	-0.232934924632048\\
	15.56	-0.18480305824805\\
	15.89	-0.145180991655579\\
	16.26	-0.104112271904356\\
	16.63	-0.0663531430902964\\
	17.04	-0.0281037107084785\\
	17.45	0.00659332758242215\\
	17.86	0.0379195148011462\\
	18.29	0.0673001069488777\\
	18.79	0.0972173801688427\\
	19.21	0.118968765342956\\
	19.72	0.141448756227817\\
	20.21	0.15920580816109\\
	20.66	0.172453471484531\\
	21.17	0.184169697652109\\
	21.76	0.193761148952419\\
	22.24	0.198792249310102\\
	22.75	0.201751264489431\\
	23.49	0.202253802089999\\
	24.38	0.198129374589115\\
	25.1	0.191970356028122\\
	25.97	0.182166820588265\\
	29.32	0.135281482424887\\
	30.54	0.118718477747152\\
	33.15	0.0885463261373616\\
	36.4	0.0617077406898119\\
	39.39	0.0452064898734648\\
	44.22	0.0284853332483834\\
	46.46	0.0231996497529963\\
	51.14	0.0151713206196149\\
	63.76	0.00488704536266482\\
	66.73	0.00375394274368546\\
	68.31	0.0032449458755508\\
	69.73	0.00287048182468652\\
	71.38	0.00246635583471289\\
	73.55	0.00201635068872008\\
	75.87	0.00160565728316442\\
	77.61	0.00135733137305749\\
	79.28	0.00114441294101653\\
	80	0.00104402816364768\\
};

\addplot [color=red, dotted, line width=0.7pt, forget plot]
table[row sep=crcr]{%
	0	-3.99995877000784\\
	0.329999999999998	-3.99689911965417\\
	0.670000000000002	-3.98751780668286\\
	1	-3.97235044131516\\
	1.33	-3.95088884575715\\
	1.64	-3.9244821159029\\
	1.94	-3.89254191791247\\
	2.22	-3.8565508398479\\
	2.52	-3.81155895874853\\
	2.81	-3.76218095675533\\
	3.11	-3.705016063391\\
	3.42	-3.6394336558681\\
	3.72	-3.56968501634384\\
	4.02	-3.49376638062894\\
	4.33	-3.40885128602064\\
	4.64	-3.31745097005853\\
	4.95	-3.21959644292144\\
	5.26000000000001	-3.11541624517739\\
	5.57000000000001	-3.00498901541947\\
	5.89	-2.8846235690237\\
	6.22	-2.75385491252577\\
	6.55	-2.61659105162219\\
	6.89	-2.46864358455534\\
	7.25	-2.30511946706301\\
	7.62	-2.13014200015701\\
	8	-1.94371024090852\\
	8.39	-1.74595156435265\\
	8.84	-1.51073255877006\\
	9.34999999999999	-1.23672654095721\\
	10.01	-0.873996055824477\\
	11.62	0.0139891080885235\\
	12.07	0.252824812679961\\
	12.46	0.452895083759174\\
	13.77	1.10888945072361\\
	14.1	1.26410100875043\\
	14.4	1.39871922531884\\
	14.7	1.52683847019478\\
	15	1.64829068905608\\
	15.31	1.76671563163913\\
	15.59	1.86756345573869\\
	15.9	1.97248414280578\\
	16.21	2.07047853040305\\
	16.54	2.16763958209677\\
	16.92	2.27244981464361\\
	17.36	2.38725629256682\\
	17.97	2.5400602855885\\
	19.44	2.89911521350919\\
	19.77	2.97271987940758\\
	20.1	3.03995078059275\\
	20.47	3.10843313728934\\
	20.86	3.17389665753106\\
	21.29	3.24000820517006\\
	21.77	3.30737487167288\\
	22.21	3.36336053910034\\
	22.67	3.41566209022504\\
	23.09	3.45777678071642\\
	23.52	3.49520143532162\\
	23.98	3.52884205659846\\
	24.4	3.55380789760372\\
	24.84	3.57403433275218\\
	25.31	3.58906725102759\\
	25.78	3.5972622878263\\
	26.22	3.59876396855543\\
	26.72	3.59331633088431\\
	27.19	3.58131680339807\\
	27.67	3.56222100882492\\
	28.17	3.5357137526433\\
	28.65	3.50403150659356\\
	29.1	3.46869004847204\\
	29.61	3.42189651216445\\
	30.15	3.36464626803716\\
	30.64	3.30607751473774\\
	31.16	3.23720883250401\\
	31.64	3.16786782052596\\
	32.19	3.08234530069565\\
	32.78	2.98359403117608\\
	33.33	2.88522948965326\\
	34.04	2.75015475889525\\
	34.7	2.61748775744377\\
	35.66	2.41537067162562\\
	36.92	2.14158874085624\\
	38.16	1.87364195721044\\
	39.2	1.65763604343468\\
	39.79	1.54058405341453\\
	40.42	1.42102226858107\\
	40.95	1.32518244910936\\
	41.61	1.21229777776006\\
	42.13	1.12862208619313\\
	42.76	1.03355776384872\\
	43.32	0.954862591771075\\
	43.84	0.88664783505881\\
	44.47	0.810086830092303\\
	45.01	0.749610141412361\\
	45.57	0.691657129255162\\
	46.33	0.620423598824203\\
	46.96	0.567317213508346\\
	47.8	0.504241163958127\\
	48.68	0.446683532143098\\
	49.44	0.403106590006843\\
	50.48	0.351549881584361\\
	51.58	0.305591410471067\\
	52.82	0.262442672887573\\
	54.05	0.2269060604103\\
	55.48	0.192760542406049\\
	56.96	0.163711490424106\\
	58.91	0.132836451326227\\
	60.43	0.11321258684211\\
	64.1	0.0774160580043457\\
	66.89	0.0581105293928061\\
	69.17	0.0459595246834112\\
	76.18	0.0223622872346141\\
	80	0.0150343193961646\\
};
\end{axis}
\end{tikzpicture}%

%% file: Eight_suboptimalitygap.tex
%
%
\begin{tikzpicture}

\begin{axis}[%
width=\figurewidth,
height=\figureheight,
at={(0\figurewidth,0\figureheight)},
scale only axis,
unbounded coords=jump,
xmin=0,
xmax=200,
xlabel style={font=\color{white!15!black}},
xlabel={MPC iteration},
ymin=0,
ymax=0.25,
ytick={0,0.1,0.2},
yticklabels={0,10,20},
ylabel style={font=\color{white!15!black}},
ylabel={Optimality gap [\%]},
axis background/.style={fill=white},
xmajorgrids,
ymajorgrids
]
\addplot [color=red, draw=none, mark=x, mark size = 1.0pt, forget plot]
  table[row sep=crcr]{%
1	0.148720438418991\\
2	0\\
3	0\\
4	0\\
5	0\\
6	0\\
7	0\\
8	0\\
9	0\\
10	0\\
11	0.00114681712477704\\
12	0.00332354353201936\\
13	0.00293727737977179\\
14	0.00248963099173238\\
15	0.00731628902508419\\
16	0.0105994608763069\\
17	0.013681596265684\\
18	0.0174770679355731\\
19	0.021121076595449\\
20	0.0236853379852562\\
21	0.024766081190279\\
22	0.027916700508797\\
23	0.0337921390475344\\
24	0.039570517616454\\
25	0.0454224719209719\\
26	0.0509572746216804\\
27	0.0562105027298685\\
28	0.0605613649114787\\
29	0.0677809405550533\\
30	0.0747952811659047\\
31	0.0814284254217625\\
32	0.0849820623307664\\
33	0.0951890836522011\\
34	0.0858861654078282\\
35	0.0934517772532786\\
36	0.0985548910645377\\
37	0.103488981441103\\
38	0.121493734906579\\
39	0.131207432967898\\
40	0.125143356215119\\
41	0.140290211820059\\
42	0.1317753432389\\
43	0.147991950603398\\
44	0.160495984903406\\
45	0.13950361759268\\
46	0.126174606781518\\
47	0.123176004021047\\
48	0.00592963557150483\\
49	0.00231672267909744\\
50	0\\
51	0\\
52	0\\
53	0.000823149865595951\\
54	0.00556652341765584\\
55	0.00976446387664964\\
56	0.0172432379436884\\
57	0.0238271313988605\\
58	0.0320796426851473\\
59	0.0396949393901593\\
60	0.0489491016027159\\
61	0.0574562913134855\\
62	0.0674731647503108\\
63	0.0727765130451132\\
64	0\\
65	0\\
66	0\\
67	0\\
68	0\\
69	0\\
70	0\\
71	0\\
72	0\\
73	0\\
74	0\\
75	0\\
76	0\\
77	0\\
78	0\\
79	0\\
80	0\\
81	0\\
82	0\\
83	0\\
84	0\\
85	0\\
86	0\\
87	0\\
88	0\\
89	0\\
90	0\\
91	0\\
92	0\\
93	0\\
94	0\\
95	0\\
96	0\\
97	0\\
98	0\\
99	0\\
100	0\\
101	0\\
102	0\\
103	0\\
104	0\\
105	0\\
106	0\\
107	0\\
108	0\\
109	0\\
110	0\\
111	0\\
112	0\\
113	0\\
114	0\\
115	0\\
116	0\\
117	0\\
118	0\\
119	0\\
120	0\\
121	0\\
122	0\\
123	0\\
124	0\\
125	0\\
126	0\\
127	0\\
128	0\\
129	0\\
130	0\\
131	0\\
132	0\\
133	0\\
134	0\\
135	0\\
136	0\\
137	0\\
138	0\\
139	0\\
140	0\\
141	0\\
142	0\\
143	0\\
144	0\\
145	0\\
146	0\\
147	0\\
148	0\\
149	0\\
150	0\\
151	0\\
152	0\\
153	0\\
154	0\\
155	0\\
156	0\\
157	0\\
158	0\\
159	0\\
160	0\\
161	0\\
162	0\\
163	0\\
164	0\\
165	0\\
166	0\\
167	0\\
168	0\\
169	0\\
170	0\\
171	0\\
172	0\\
173	0\\
174	0\\
175	0\\
176	0\\
177	0\\
178	0\\
179	0\\
180	0\\
181	0\\
182	0\\
183	0\\
184	0\\
185	0\\
186	0\\
187	0\\
188	0\\
189	0\\
190	0\\
191	0\\
192	0\\
193	0\\
194	0\\
195	0\\
196	0\\
197	0\\
198	0\\
199	0\\
200	0\\
};
\addplot [color=red, dashed, forget plot]
  table[row sep=crcr]{%
0	0.199999999999989\\
200	0.199999999999989\\
};
\addplot [color=green, draw=none, mark=x, mark size = 1.0pt, forget plot]
  table[row sep=crcr]{%
1	0.0708285797908559\\
2	0\\
3	0\\
4	0\\
5	0\\
6	0\\
7	0\\
8	0\\
9	0\\
10	0\\
11	0\\
12	0.00455101334185315\\
13	0.006779133499208\\
14	0.00802768418938626\\
15	0.00765888726135699\\
16	0.0123536757383249\\
17	0.0129725346951375\\
18	0.0132582017043603\\
19	0.0186724134132703\\
20	0.024016256501227\\
21	0.0250473018836033\\
22	0.025877754042483\\
23	0.0384677704155649\\
24	0.0432751648237399\\
25	0.0479640464254203\\
26	0.0532818756357756\\
27	0.0578241491097629\\
28	0.0596126507801102\\
29	0.0587879362132355\\
30	0.06592511147133\\
31	0.0715604909326544\\
32	0.0613749859727193\\
33	0.0675065383410356\\
34	0.0719421627043175\\
35	0.0836352510307847\\
36	0.0884413995385955\\
37	0.101120175150697\\
38	0.031255833806938\\
39	0.0427986763218087\\
40	0.0522900215256925\\
41	0.0613327331970197\\
42	0.0698071827184208\\
43	0.0761418568608576\\
44	0.0809301666914166\\
45	0.0956320250733427\\
46	0.000139291171819877\\
47	-0\\
48	0\\
49	0.00125251422556971\\
50	0.007561172731954\\
51	0.0140277299947229\\
52	0.0199451815710745\\
53	0.0276098404990535\\
54	0.0346102319491024\\
55	0.0432100870076795\\
56	0\\
57	0.00339100513139101\\
58	0.00779794881583484\\
59	0.0138432683840506\\
60	0.0199099884815723\\
61	0.0267703600346181\\
62	0\\
63	0\\
64	0\\
65	0\\
66	0\\
67	0\\
68	0\\
69	0\\
70	0\\
71	0\\
72	0\\
73	0\\
74	0\\
75	0\\
76	0\\
77	0\\
78	0\\
79	0\\
80	0\\
81	0\\
82	0\\
83	0\\
84	0\\
85	0\\
86	0\\
87	0\\
88	0\\
89	0\\
90	0\\
91	0\\
92	0\\
93	0\\
94	0\\
95	0\\
96	0\\
97	0\\
98	0\\
99	0\\
100	0\\
101	0\\
102	0\\
103	0\\
104	0\\
105	0\\
106	0\\
107	0\\
108	0\\
109	0\\
110	0\\
111	0\\
112	0\\
113	0\\
114	0\\
115	0\\
116	0\\
117	0\\
118	0\\
119	0\\
120	0\\
121	0\\
122	0\\
123	0\\
124	0\\
125	0\\
126	0\\
127	0\\
128	0\\
129	0\\
130	0\\
131	0\\
132	0\\
133	0\\
134	0\\
135	0\\
136	0\\
137	0\\
138	0\\
139	0\\
140	0\\
141	0\\
142	0\\
143	0\\
144	0\\
145	0\\
146	0\\
147	0\\
148	0\\
149	0\\
150	0\\
151	0\\
152	0\\
153	0\\
154	0\\
155	0\\
156	0\\
157	0\\
158	0\\
159	0\\
160	0\\
161	0\\
162	0\\
163	0\\
164	0\\
165	0\\
166	0\\
167	0\\
168	0\\
169	0\\
170	0\\
171	0\\
172	0\\
173	0\\
174	0\\
175	0\\
176	0\\
177	0\\
178	0\\
179	0\\
180	0\\
181	0\\
182	0\\
183	0\\
184	0\\
185	0\\
186	0\\
187	0\\
188	0\\
189	0\\
190	0\\
191	0\\
192	0\\
193	0\\
194	0\\
195	0\\
196	0\\
197	0\\
198	0\\
199	0\\
200	0\\
};
\addplot [color=green, dashed, forget plot]
  table[row sep=crcr]{%
0	0.0999999999999943\\
200	0.0999999999999943\\
};
\addplot [color=blue, draw=none, mark=x, mark size = 1.0pt, forget plot]
  table[row sep=crcr]{%
1	0\\
2	0\\
3	0\\
4	0\\
5	0\\
6	0\\
7	0\\
8	0\\
9	0\\
10	0\\
11	0\\
12	0\\
13	0\\
14	0\\
15	0\\
16	0.00626580248916753\\
17	0.00759552304398881\\
18	0\\
19	0.00593773079631887\\
20	-0\\
21	0\\
22	0\\
23	0.00598147952334216\\
24	0.00588230639758081\\
25	0\\
26	0\\
27	0\\
28	0\\
29	0\\
30	0.00672194910836765\\
31	0.0139963435947834\\
32	0\\
33	0\\
34	0\\
35	0.00735989898228695\\
36	0.00863478891110958\\
37	0\\
38	0\\
39	0.00905691469273506\\
40	0.0159180332107383\\
41	0.00205914155841924\\
42	0.00028579281362795\\
nan	nan\\
45	0.00349259750927899\\
46	0.00023894049920159\\
47	0\\
48	0\\
49	0\\
50	0\\
51	0\\
52	0\\
53	0\\
54	0\\
55	0\\
56	0\\
57	0\\
58	0\\
59	0\\
60	0\\
61	0\\
62	0\\
63	0\\
64	0\\
65	0\\
66	0\\
67	0\\
68	0\\
69	0\\
70	0\\
71	0\\
72	0\\
73	0\\
74	0\\
75	0\\
76	0\\
77	0\\
78	0\\
79	0\\
80	0\\
81	0\\
82	0\\
83	0\\
84	0\\
85	0\\
86	0\\
87	0\\
88	0\\
89	0\\
90	0\\
91	0\\
92	0\\
93	0\\
94	0\\
95	0\\
96	0\\
97	0\\
98	0\\
99	0\\
100	0\\
101	0\\
102	0\\
103	0\\
104	0\\
105	0\\
106	0\\
107	0\\
108	0\\
109	0\\
110	0\\
111	0\\
112	0\\
113	0\\
114	0\\
115	0\\
116	0\\
117	0\\
118	0\\
119	0\\
120	0\\
121	0\\
122	0\\
123	0\\
124	0\\
125	0\\
126	0\\
127	0\\
128	0\\
129	0\\
130	0\\
131	0\\
132	0\\
133	0\\
134	0\\
135	0\\
136	0\\
137	0\\
138	0\\
139	0\\
140	0\\
141	0\\
142	0\\
143	0\\
144	0\\
145	0\\
146	0\\
147	0\\
148	0\\
149	0\\
150	0\\
151	0\\
152	0\\
153	0\\
154	0\\
155	0\\
156	0\\
157	0\\
158	0\\
159	0\\
160	0\\
161	0\\
162	0\\
163	0\\
164	0\\
165	0\\
166	0\\
167	0\\
168	0\\
169	0\\
170	0\\
171	0\\
172	0\\
173	0\\
174	0\\
175	0\\
176	0\\
177	0\\
178	0\\
179	0\\
180	0\\
181	0\\
182	0\\
183	0\\
184	0\\
185	0\\
186	0\\
187	0\\
188	0\\
189	0\\
190	0\\
191	0\\
192	0\\
193	0\\
194	0\\
195	0\\
196	0\\
197	0\\
198	0\\
199	0\\
200	0\\
};
\addplot [color=blue, dashed, forget plot]
  table[row sep=crcr]{%
0	0.0200000000000102\\
200	0.0200000000000102\\
};
\end{axis}
\end{tikzpicture}%

%% file: Eight_computation_times_new_2.tex
%
%
\begin{tikzpicture}

\begin{axis}[%
width=\figurewidth,
height=\figureheight,
at={(0\figurewidth,0\figureheight)},
scale only axis,
xmin=0,
xmax=200,
xlabel style={font=\color{white!15!black}},
xlabel={MPC iteration},
ymin=0,
ymax=0.17,
ytick={0,0.05,0.1,0.15},
yticklabels={0,50,100,150},
ylabel style={font=\color{white!15!black}},
ylabel={Time [ms]},
axis background/.style={fill=white},
xmajorgrids,
ymajorgrids,
/pgf/number format/precision=5
]
\addplot [color=blue, only marks, mark=x, mark size = 1.0pt, forget plot]
  table[row sep=crcr]{%
1	0.0923372000000029\\
2	0.0296342999999979\\
3	0.0380929000000094\\
4	0.031451699999991\\
5	0.0308957999999961\\
6	0.0310595999999919\\
7	0.0331837999999891\\
8	0.0395651999999984\\
9	0.0347912000000008\\
10	0.0335513999999932\\
11	0.0371637000000078\\
12	0.070232400000009\\
13	0.040424999999999\\
14	0.0646688999999867\\
15	0.0405202999999972\\
16	0.0579641999999865\\
17	0.0460233000000017\\
18	0.0709066999999948\\
19	0.0422821999999883\\
20	0.0768260000000112\\
21	0.0462197999999887\\
22	0.0476623000000131\\
23	0.0709483999999918\\
24	0.0471441999999911\\
25	0.0817878000000007\\
26	0.052686499999993\\
27	0.0546962999999892\\
28	0.0650914999999941\\
29	0.0882044999999891\\
30	0.0670899000000134\\
31	0.0784339999999872\\
32	0.0894815999999992\\
33	0.0613231000000098\\
34	0.054708799999986\\
35	0.0689092000000073\\
36	0.0721518000000003\\
37	0.0702963999999895\\
38	0.0622823999999866\\
39	0.0689175999999918\\
40	0.0752530999999976\\
41	0.0509351999999978\\
42	0.050460500000014\\
43	0.0609077000000013\\
44	0.0506216999999936\\
45	0.0501839000000075\\
46	0.0654327000000023\\
47	0.0587237000000016\\
48	0.0630260000000078\\
49	0.0530966000000035\\
50	0.0501099999999894\\
51	0.0443954000000133\\
52	0.059615100000002\\
53	0.0772412000000031\\
54	0.0418444999999963\\
55	0.0451408000000129\\
56	0.0582821000000138\\
57	0.0428105000000016\\
58	0.123900899999995\\
59	0.0512487000000021\\
60	0.0415474000000131\\
61	0.0403148000000044\\
62	0.038474100000002\\
63	0.0360987999999907\\
64	0.0274205000000052\\
65	0.0294336999999985\\
66	0.0305170000000032\\
67	0.0257464999999968\\
68	0.0245008999999925\\
69	0.0264117000000113\\
70	0.0355212999999992\\
71	0.0420249999999953\\
72	0.0262705000000096\\
73	0.0247629000000131\\
74	0.0268963999999983\\
75	0.0245558999999957\\
76	0.0223704999999939\\
77	0.0223182000000008\\
78	0.0255962999999895\\
79	0.0263624000000107\\
80	0.022340899999989\\
81	0.0221362999999997\\
82	0.0217105000000117\\
83	0.0237282999999877\\
84	0.0272383999999875\\
85	0.0256213000000116\\
86	0.0280300000000011\\
87	0.0300716000000136\\
88	0.0293966000000125\\
89	0.0288568999999939\\
90	0.0301561999999933\\
91	0.0269868000000031\\
92	0.0361536000000058\\
93	0.0312125000000094\\
94	0.0338954999999999\\
95	0.0291427999999883\\
96	0.0374769999999955\\
97	0.0290454000000011\\
98	0.0369812000000138\\
99	0.0401262999999972\\
100	0.0291315999999995\\
101	0.0257973999999876\\
102	0.0303471999999942\\
103	0.029538500000001\\
104	0.0277012999999897\\
105	0.0280430999999908\\
106	0.0291003000000103\\
107	0.0246061999999938\\
108	0.0250197000000014\\
109	0.0240957000000037\\
110	0.0206367000000114\\
111	0.0260431000000096\\
112	0.0219817999999918\\
113	0.0266665999999987\\
114	0.0277164000000028\\
115	0.0277086000000111\\
116	0.0245903999999939\\
117	0.0265000999999927\\
118	0.019929100000013\\
119	0.0229280000000074\\
120	0.0270413000000076\\
121	0.0294517000000099\\
122	0.026271400000013\\
123	0.0290898000000084\\
124	0.0287379000000101\\
125	0.0253666000000123\\
126	0.0295155999999963\\
127	0.0277289999999937\\
128	0.0407122999999956\\
129	0.0291397999999958\\
130	0.0308937000000071\\
131	0.0244974000000013\\
132	0.0227444999999875\\
133	0.0219611999999927\\
134	0.021897800000005\\
135	0.0230484999999874\\
136	0.0242921999999908\\
137	0.0221129999999903\\
138	0.0313816000000031\\
139	0.0276872000000026\\
140	0.0216045999999892\\
141	0.0241392000000076\\
142	0.020851399999998\\
143	0.0197023000000058\\
144	0.0284584000000052\\
145	0.0205700999999863\\
146	0.030220700000001\\
147	0.0275019000000043\\
148	0.0242904000000124\\
149	0.0288957999999866\\
150	0.0318437999999901\\
151	0.0228965000000017\\
152	0.0256914999999935\\
153	0.0231100999999967\\
154	0.0239381000000094\\
155	0.0227920000000097\\
156	0.025209799999999\\
157	0.0214584999999943\\
158	0.0304763999999977\\
159	0.0369164999999896\\
160	0.020217400000007\\
161	0.0184345000000121\\
162	0.018822699999987\\
163	0.0189147000000105\\
165	0.0183694000000116\\
166	0.0200580000000059\\
167	0.0226201000000117\\
168	0.0191418999999939\\
169	0.019929999999988\\
170	0.0178687000000082\\
171	0.0205170000000123\\
172	0.0258388000000025\\
173	0.021946899999989\\
174	0.0192781000000082\\
175	0.0209739000000013\\
176	0.0244582999999921\\
177	0.0210929000000135\\
178	0.0256053999999892\\
179	0.0295426000000134\\
180	0.0251273000000083\\
181	0.0210783999999933\\
182	0.0252562000000012\\
183	0.0299081999999942\\
184	0.0271511999999916\\
185	0.0276044999999954\\
186	0.0214178999999888\\
187	0.0198377000000107\\
188	0.0215300000000127\\
189	0.0264603999999906\\
190	0.0254957999999874\\
191	0.0246675999999866\\
192	0.0226914000000136\\
193	0.0205870000000061\\
194	0.025591200000008\\
195	0.0221267999999952\\
196	0.0291953999999919\\
197	0.029349400000001\\
198	0.0202337999999997\\
199	0.01918280000001\\
200	0.0183500999999922\\
};
\addplot [color=black, only marks, mark=x, mark size = 1.0pt, forget plot]
  table[row sep=crcr]{%
1	0.0194041000000027\\
2	0.0144498999999882\\
3	0.0159106000000122\\
4	0.0121087999999929\\
5	0.0154807000000119\\
6	0.0134357999999963\\
7	0.0152766999999869\\
8	0.017914300000001\\
9	0.0119908000000066\\
10	0.0113274999999931\\
11	0.0117205999999896\\
12	0.0124554999999873\\
13	0.0138118000000134\\
14	0.0108856999999887\\
15	0.010900399999997\\
16	0.0156267000000128\\
17	0.0136244999999917\\
18	0.0117196999999862\\
19	0.0136967999999911\\
20	0.0108968000000118\\
21	0.0113001999999938\\
22	0.0147785000000056\\
23	0.0162276000000077\\
24	0.0108898999999951\\
25	0.0171502000000032\\
26	0.0124083000000041\\
27	0.0144261000000085\\
28	0.0132571000000041\\
29	0.0111896999999885\\
30	0.0109031000000073\\
31	0.0160143999999889\\
32	0.0144947000000002\\
33	0.0116199999999935\\
34	0.0150911000000065\\
35	0.0119760999999983\\
36	0.014588300000014\\
37	0.0152347000000077\\
38	0.0156388000000049\\
39	0.0102334000000042\\
40	0.0246735000000058\\
41	0.018529199999989\\
42	0.0164383999999984\\
43	0.0195798000000025\\
44	0.0176944000000105\\
45	0.021314799999999\\
46	0.016708599999987\\
47	0.0174744999999916\\
48	0.0146503999999936\\
49	0.0152161000000035\\
50	0.0164259999999956\\
51	0.0165653999999904\\
52	0.0164498000000037\\
53	0.0113943000000063\\
54	0.0153602000000035\\
55	0.0106633999999985\\
56	0.0152266000000054\\
57	0.0155025999999907\\
58	0.0100563999999963\\
59	0.01522030000001\\
60	0.0119349\\
61	0.0160796000000119\\
62	0.0108827999999903\\
63	0.0165216999999984\\
64	0.0134338000000014\\
65	0.0164655000000096\\
66	0.0128601000000117\\
67	0.0127310000000023\\
68	0.0140615999999909\\
69	0.0191722000000141\\
70	0.0119939999999872\\
71	0.0172966000000088\\
72	0.017777499999994\\
73	0.0111269999999877\\
74	0.010822399999995\\
75	0.0108198000000073\\
76	0.0154903000000104\\
77	0.0180296000000055\\
78	0.0116256000000021\\
79	0.0114178000000038\\
80	0.0115163000000109\\
81	0.016369300000008\\
82	0.0128856999999982\\
83	0.0136151000000098\\
84	0.010917500000005\\
85	0.0122722000000124\\
86	0.0182522000000063\\
87	0.0155009000000064\\
88	0.0116471000000047\\
89	0.0118784999999946\\
90	0.0145985999999994\\
91	0.0167954000000066\\
92	0.0118452000000104\\
93	0.0169965000000047\\
94	0.0185248999999885\\
95	0.0199254000000053\\
96	0.0153936000000101\\
97	0.0117078000000106\\
99	0.0134808999999905\\
100	0.0154665000000023\\
101	0.0119281999999998\\
102	0.0180645000000084\\
103	0.0186118000000022\\
104	0.0164613999999972\\
105	0.0147024000000044\\
106	0.0138231999999903\\
107	0.016052099999996\\
108	0.0163057999999978\\
109	0.0121656999999971\\
110	0.0142012999999963\\
111	0.0142050999999981\\
112	0.0200461999999959\\
113	0.0166903999999874\\
114	0.0119545999999957\\
115	0.0165278999999998\\
116	0.0123462999999902\\
117	0.0160635000000013\\
118	0.0144620000000089\\
119	0.015038699999991\\
120	0.0192965999999899\\
121	0.0194100999999876\\
122	0.0158943000000136\\
123	0.0125601000000017\\
124	0.0177226000000132\\
125	0.0141725999999949\\
126	0.012736799999999\\
127	0.0178319999999985\\
128	0.019478399999997\\
129	0.0206196999999975\\
130	0.013280299999991\\
131	0.0194397999999865\\
132	0.0129282000000046\\
133	0.0145848999999885\\
134	0.016336200000012\\
135	0.0161546999999871\\
136	0.0163151000000141\\
137	0.0131436000000065\\
138	0.0218575000000101\\
139	0.0208394999999939\\
140	0.0234987999999987\\
141	0.019117499999993\\
142	0.0180511000000081\\
143	0.0192490000000021\\
144	0.0151564999999891\\
145	0.0188660000000027\\
146	0.0124898999999914\\
147	0.0145355999999879\\
148	0.0162671999999873\\
149	0.0151754000000039\\
150	0.0136190999999997\\
151	0.013198999999986\\
152	0.0190130999999951\\
153	0.0216743000000008\\
154	0.02071380000001\\
155	0.0135439000000019\\
156	0.0156390999999871\\
157	0.0121780999999999\\
158	0.0159127000000012\\
159	0.016282799999999\\
160	0.015668699999992\\
161	0.0184495999999967\\
162	0.0108874999999955\\
163	0.0110149999999862\\
164	0.0138412000000017\\
165	0.0165371999999877\\
166	0.0176556000000119\\
167	0.015892000000008\\
168	0.0190628000000004\\
169	0.0211496999999952\\
170	0.0139144000000044\\
171	0.0158299\\
172	0.0137593999999979\\
173	0.0156365999999935\\
174	0.0144885999999929\\
175	0.0211074999999994\\
176	0.0233299999999872\\
177	0.0221732999999915\\
178	0.0183336000000054\\
179	0.0213852999999915\\
180	0.0201246000000026\\
181	0.0235586999999953\\
182	0.0180296000000055\\
183	0.018763100000001\\
184	0.0160548000000063\\
185	0.0116074999999967\\
186	0.0179926000000137\\
187	0.0138324000000125\\
188	0.0151879999999949\\
189	0.0118629000000112\\
190	0.0177313000000083\\
191	0.0167291999999861\\
192	0.0186013000000003\\
193	0.0125419000000022\\
194	0.0115677999999946\\
195	0.0110354999999913\\
196	0.0164039000000002\\
197	0.0146900999999957\\
198	0.0112196999999981\\
199	0.0181573000000128\\
200	0.0199194999999861\\
};
\addplot [color=black, dashed, forget plot]
  table[row sep=crcr]{%
0	0.0999999999999943\\
200	0.0999999999999943\\
};
\addplot [color=green, only marks, mark=x, mark size = 1.0pt, forget plot]
  table[row sep=crcr]{%
1	0.0655415999999889\\
2	0.0350759000000096\\
3	0.0376178999999865\\
4	0.0286020000000065\\
5	0.0274001999999882\\
6	0.0274238000000082\\
7	0.0279705000000092\\
8	0.0311288999999988\\
9	0.0361078999999904\\
10	0.0277805999999998\\
11	0.0286878999999942\\
12	0.0301878999999872\\
13	0.0305822999999918\\
14	0.0399033999999858\\
15	0.0466279999999983\\
16	0.0402363000000037\\
17	0.0416588999999874\\
18	0.0420970000000125\\
19	0.044298700000013\\
20	0.0448100000000125\\
21	0.0448558999999875\\
22	0.0392433000000096\\
23	0.0516173000000038\\
24	0.0424640999999895\\
25	0.0518481999999949\\
26	0.0540511000000095\\
27	0.0464945999999884\\
28	0.0492423999999971\\
29	0.0483692999999903\\
30	0.0434906000000126\\
31	0.0541140999999925\\
32	0.100814899999989\\
33	0.0496258999999952\\
34	0.043799000000007\\
35	0.0496784000000048\\
36	0.0928613000000098\\
37	0.0748178000000053\\
38	0.0806355999999937\\
39	0.0409750999999972\\
40	0.052665200000007\\
41	0.0553524000000039\\
42	0.043215900000007\\
43	0.0398812999999905\\
44	0.0693510000000117\\
45	0.0711829999999907\\
46	0.0904223000000002\\
47	0.0289832999999931\\
48	0.0312697999999898\\
49	0.0326704999999947\\
50	0.0420815999999888\\
51	0.0368316999999934\\
52	0.0423590999999988\\
53	0.0329724000000056\\
54	0.0316775999999948\\
55	0.0330778000000009\\
56	0.0425205999999889\\
57	0.030698000000001\\
58	0.0295238999999867\\
59	0.0312629999999956\\
60	0.0361729999999909\\
61	0.0359914999999944\\
62	0.0269039000000078\\
63	0.0247377000000029\\
64	0.024354499999987\\
65	0.0248406000000045\\
66	0.0248286000000064\\
67	0.0306425999999931\\
68	0.0393486000000109\\
69	0.035589200000004\\
70	0.0273644000000104\\
71	0.0268500999999901\\
72	0.0253906000000086\\
73	0.0217379000000051\\
74	0.0238683000000037\\
75	0.0229497999999921\\
76	0.0336690000000033\\
77	0.0308982999999898\\
78	0.0233393000000035\\
79	0.0238865999999973\\
80	0.0261102999999991\\
81	0.0243070999999873\\
82	0.0228017999999963\\
83	0.0272927000000038\\
84	0.0261648000000037\\
85	0.0241417999999953\\
86	0.0232975000000124\\
87	0.0275383999999974\\
88	0.0259225000000072\\
89	0.0284002999999871\\
90	0.0364544999999907\\
91	0.0293991000000062\\
92	0.0299196999999936\\
93	0.0281491000000074\\
94	0.0286486999999909\\
95	0.0255632999999875\\
96	0.0397854999999936\\
97	0.0295897999999966\\
98	0.0292983999999876\\
99	0.0344497000000104\\
100	0.027497899999986\\
101	0.0272166999999968\\
102	0.0302660000000117\\
103	0.0258182000000033\\
104	0.0355874999999912\\
105	0.0259909000000107\\
106	0.0273740000000089\\
107	0.02719909999999\\
108	0.0252055999999925\\
109	0.0251880000000142\\
110	0.0233201000000065\\
111	0.0310891999999967\\
112	0.022894000000008\\
113	0.0264851999999962\\
114	0.0214248000000055\\
115	0.0246999999999957\\
116	0.0211014999999861\\
117	0.0205706000000134\\
118	0.0268175999999869\\
119	0.0223168999999928\\
120	0.0261173000000099\\
121	0.0242325000000108\\
122	0.0240804000000026\\
123	0.0298569999999927\\
124	0.0297510999999986\\
125	0.0257416999999975\\
126	0.0337455000000091\\
127	0.0295294999999953\\
128	0.0327372000000139\\
129	0.0247253999999941\\
130	0.027185199999991\\
131	0.0251820000000009\\
132	0.0236231000000089\\
133	0.0262347000000034\\
134	0.0227380000000039\\
135	0.0239161000000081\\
136	0.0204366000000107\\
137	0.0236178000000109\\
138	0.0249180000000138\\
139	0.0238305999999966\\
140	0.0320112999999935\\
141	0.0312078000000042\\
142	0.0245563999999945\\
143	0.0242566999999951\\
144	0.0223528999999871\\
145	0.0279932999999915\\
146	0.0260196000000121\\
147	0.0197236999999859\\
148	0.0217366999999911\\
149	0.0232660000000067\\
150	0.0225248000000136\\
151	0.0199851000000137\\
152	0.0230069000000128\\
153	0.0224707000000137\\
154	0.0206664000000103\\
155	0.0211137999999949\\
156	0.0224166000000139\\
157	0.0227843000000121\\
159	0.0196361000000138\\
160	0.021792199999993\\
161	0.0188186999999971\\
162	0.0217845999999895\\
163	0.0210429999999917\\
164	0.023753499999998\\
165	0.024347199999994\\
166	0.0204142000000047\\
167	0.0212048999999865\\
168	0.0195156999999995\\
169	0.0203721999999971\\
170	0.0193147999999894\\
171	0.0174989999999866\\
172	0.0194692000000032\\
173	0.0248438000000135\\
174	0.0204606000000069\\
175	0.0195459999999912\\
177	0.0184782999999982\\
178	0.0174088999999924\\
179	0.0171215000000018\\
180	0.0210935000000063\\
181	0.0204769000000056\\
182	0.0188730000000135\\
183	0.0195841999999971\\
184	0.0195325999999909\\
185	0.017051500000008\\
187	0.0172933999999998\\
188	0.0197420000000079\\
189	0.0208229000000131\\
190	0.0217756000000122\\
191	0.0178496999999993\\
192	0.0187186999999938\\
193	0.0187799999999925\\
194	0.0192337000000009\\
195	0.0201973000000066\\
196	0.0185403000000122\\
197	0.0177189000000055\\
198	0.0183312000000058\\
199	0.0316094000000078\\
200	0.0194903999999951\\
};
\addplot [color=red, only marks, mark=x, mark size = 1.0pt, forget plot]
  table[row sep=crcr]{%
1	0.0424649999999929\\
2	0.0277825000000007\\
3	0.0289683000000025\\
4	0.0304265000000044\\
5	0.0289793999999972\\
6	0.0294418999999948\\
7	0.0390352000000007\\
8	0.029295999999988\\
9	0.0271971999999892\\
10	0.0287375999999995\\
11	0.0296877000000109\\
12	0.0370307999999966\\
13	0.0451812000000018\\
14	0.0321375000000046\\
15	0.0305327000000091\\
16	0.0325542999999868\\
17	0.0368268999999941\\
18	0.0395753999999897\\
19	0.0367170000000101\\
20	0.031855400000012\\
21	0.0343401000000085\\
22	0.0314703999999892\\
23	0.042461099999997\\
24	0.0436645999999996\\
25	0.0320912999999905\\
26	0.0329278000000102\\
27	0.0356363999999871\\
28	0.0341460999999867\\
29	0.0349856999999929\\
30	0.0410435000000007\\
31	0.0329539999999895\\
32	0.0416476000000046\\
33	0.0407350000000122\\
34	0.0387719999999945\\
35	0.0364936999999941\\
36	0.0394737999999961\\
37	0.0412470999999925\\
38	0.0404862999999978\\
39	0.0436626000000047\\
40	0.0481163000000038\\
41	0.0508680000000084\\
42	0.0446839999999895\\
43	0.0416629999999998\\
44	0.0454780000000028\\
45	0.046467199999995\\
46	0.0617685999999935\\
47	0.038246700000002\\
48	0.0382390000000044\\
49	0.0306558000000052\\
50	0.030135400000006\\
51	0.035308299999997\\
52	0.028630800000002\\
53	0.033984299999986\\
54	0.0305381999999952\\
55	0.032334800000001\\
56	0.0368291999999997\\
57	0.0415262999999868\\
58	0.0326262000000099\\
59	0.0347562000000039\\
60	0.0324358999999959\\
61	0.0368694000000005\\
62	0.0308483999999964\\
63	0.0266201999999964\\
64	0.0222909000000016\\
65	0.0239430000000027\\
66	0.0295986000000141\\
67	0.0241865000000132\\
68	0.0274450000000002\\
69	0.0278682999999944\\
70	0.0247746000000006\\
71	0.0230335999999909\\
72	0.0237890999999877\\
73	0.0276428999999894\\
74	0.0269921000000011\\
75	0.0235954000000049\\
76	0.0233104000000139\\
77	0.0244747000000132\\
78	0.0249934999999937\\
79	0.0250091999999995\\
80	0.0229591000000084\\
81	0.0227978000000064\\
82	0.0308382000000051\\
83	0.0252241000000026\\
84	0.0228260000000091\\
85	0.0244147999999882\\
86	0.0233085000000131\\
87	0.0264582000000075\\
88	0.0257925999999884\\
89	0.0235323999999935\\
90	0.0425984999999969\\
91	0.0359058999999888\\
92	0.0402727000000027\\
93	0.0271319000000005\\
94	0.0259801999999922\\
95	0.0320356000000004\\
96	0.0278346999999997\\
97	0.0290998000000116\\
98	0.0290119999999945\\
99	0.0273836000000074\\
100	0.033771900000005\\
101	0.0318892999999889\\
102	0.02680509999999\\
103	0.0293557999999905\\
104	0.0254741999999908\\
105	0.0273809999999912\\
106	0.0285082000000045\\
107	0.026955700000002\\
108	0.0336527000000046\\
109	0.0318646000000058\\
110	0.0243720999999937\\
111	0.0243158000000108\\
112	0.0259389999999939\\
113	0.0215209000000129\\
114	0.0248047999999983\\
115	0.0212562000000105\\
116	0.0232589000000019\\
117	0.0229730000000075\\
118	0.0200673999999879\\
119	0.0236773999999969\\
120	0.0202122000000031\\
121	0.0233339000000115\\
122	0.0233519999999885\\
123	0.0248676999999873\\
124	0.0330090999999868\\
125	0.0309460999999942\\
126	0.026193400000011\\
127	0.0299797000000126\\
128	0.0280012000000056\\
129	0.0297564999999906\\
130	0.038175300000006\\
131	0.0364313999999979\\
132	0.0355982000000097\\
133	0.0366850000000056\\
134	0.0463086999999973\\
135	0.0364820000000066\\
136	0.0359400000000107\\
137	0.0303254000000095\\
138	0.0361490999999887\\
139	0.0356491999999946\\
140	0.0398755999999878\\
141	0.041027100000008\\
142	0.0391065000000026\\
143	0.0319436000000053\\
144	0.0292921999999862\\
145	0.0297581000000093\\
146	0.025014700000014\\
147	0.0318825999999888\\
148	0.0230527999999879\\
149	0.0255831999999998\\
150	0.0239108000000101\\
151	0.0235595999999987\\
152	0.0218135000000075\\
153	0.0282829999999876\\
154	0.0217997999999966\\
155	0.0255645000000015\\
156	0.0253792999999973\\
157	0.02089509999999\\
158	0.0234601999999882\\
159	0.02141850000001\\
160	0.0208261999999877\\
161	0.019644900000003\\
162	0.0208235000000059\\
163	0.0200986000000114\\
164	0.0190444000000127\\
165	0.0208083999999928\\
166	0.0191507000000115\\
167	0.0192309999999907\\
168	0.0201984999999922\\
169	0.0185491999999954\\
170	0.023902899999996\\
171	0.0220018999999922\\
172	0.0192662999999982\\
173	0.0164532000000008\\
174	0.0189592000000118\\
175	0.0180671999999902\\
176	0.019736800000004\\
177	0.0193724000000088\\
178	0.0258200000000102\\
179	0.021702299999987\\
180	0.0197288999999898\\
181	0.0179224999999974\\
182	0.0189030999999886\\
183	0.0193868000000066\\
184	0.0188685999999905\\
185	0.0193127000000004\\
186	0.0253974000000028\\
187	0.0203080999999941\\
188	0.0195851999999945\\
189	0.0179699999999912\\
190	0.0194406999999899\\
191	0.018787599999996\\
192	0.0201570000000117\\
193	0.0178452000000107\\
194	0.0237289000000089\\
195	0.0190844999999911\\
196	0.0210442000000057\\
197	0.0203926000000081\\
198	0.0194762000000139\\
199	0.021556199999992\\
200	0.0196885999999949\\
};
\end{axis}
\end{tikzpicture}%

%% file: path_xy_straight.tex
%
%
\definecolor{mycolor1}{rgb}{0.50000,0.00000,0.50000}%
\begin{tikzpicture}

\begin{axis}[%
width=\figurewidth,
height=\figureheight,
at={(0\figurewidth,0\figureheight)},
scale only axis,
xmin=50,
xmax=155,
xtick={60,80,100,120,140},
xticklabels={-80,-60,-40,-20,0},
xlabel style={font=\color{white!15!black}},
xlabel={$x$ [m]},
ymin=-16,
ymax=12,
ytick={-10,0,10},
ylabel style={font=\color{white!15!black}},
ylabel={$y$ [m]},
axis background/.style={fill=white},
xmajorgrids,
ymajorgrids
]
\addplot [color=black, line width=1.0pt, forget plot]
  table[row sep=crcr]{%
160	0\\
39.5	0\\
};
\addplot [color=blue, line width=0.7pt, forget plot]
  table[row sep=crcr]{%
139.581420898438	-0.911224365234375\\
139.565414428711	-0.900894165039063\\
138.666213989258	-0.056060791015625\\
138.638595581055	-0.0391082763671875\\
137.920669555664	0.614669799804688\\
137.886352539063	0.62689208984375\\
137.50846862793	0.949005126953125\\
137.103927612305	1.27996826171875\\
137.077606201172	1.2857666015625\\
136.676605224609	1.59269714355469\\
136.299041748047	1.88981628417969\\
136.254638671875	1.8831787109375\\
135.866119384766	2.16099548339844\\
135.842468261719	2.15742492675781\\
135.436019897461	2.43307495117188\\
135.412322998047	2.429931640625\\
135.00373840332	2.67330932617188\\
134.956268310547	2.66023254394531\\
134.563140869141	2.86894226074219\\
134.520629882813	2.84910583496094\\
134.12907409668	3.03622436523438\\
134.110656738281	3.02804565429688\\
133.732528686523	3.18656921386719\\
133.699020385742	3.16383361816406\\
133.322845458984	3.29423522949219\\
133.305557250977	3.27951049804688\\
132.899169921875	3.40174865722656\\
132.859436035156	3.36997985839844\\
132.432037353516	3.47549438476563\\
132.416366577148	3.46461486816406\\
131.985382080078	3.54380798339844\\
131.955047607422	3.51658630371094\\
131.5654296875	3.5701904296875\\
131.547271728516	3.55574035644531\\
131.166732788086	3.58708190917969\\
131.1474609375	3.57498168945313\\
130.71549987793	3.57765197753906\\
130.698257446289	3.56100463867188\\
130.256500244141	3.54212951660156\\
130.238006591797	3.52833557128906\\
129.795684814453	3.506591796875\\
129.783203125	3.48799133300781\\
129.322174072266	3.45344543457031\\
129.294616699219	3.42329406738281\\
128.822448730469	3.36984252929688\\
128.809722900391	3.353515625\\
128.311920166016	3.27009582519531\\
128.275100708008	3.23780822753906\\
127.713729858398	3.12548828125\\
127.103240966797	2.98408508300781\\
127.079818725586	2.96783447265625\\
126.435653686523	2.7933349609375\\
125.745544433594	2.58859252929688\\
125.008605957031	2.35543823242188\\
124.991607666016	2.33880615234375\\
124.229232788086	2.0872802734375\\
122.702362060547	1.58680725097656\\
122.647857666016	1.56890869140625\\
121.11376953125	1.06491088867188\\
121.053924560547	1.06741333007813\\
120.302536010742	0.848648071289063\\
119.575912475586	0.646224975585938\\
119.507659912109	0.670745849609375\\
118.785064697266	0.49114990234375\\
118.730606079102	0.507369995117188\\
118.049652099609	0.363967895507813\\
117.37109375	0.266876220703125\\
117.350830078125	0.277420043945313\\
116.766418457031	0.1956787109375\\
116.725402832031	0.215850830078125\\
116.18440246582	0.154586791992188\\
116.162002563477	0.16973876953125\\
115.618515014648	0.1417236328125\\
115.603744506836	0.146636962890625\\
115.091506958008	0.117507934570313\\
115.070358276367	0.1298828125\\
114.490158081055	0.129013061523438\\
114.473236083984	0.137832641601563\\
113.206924438477	0.165206909179688\\
113.177978515625	0.17181396484375\\
112.49870300293	0.201919555664063\\
109.337341308594	0.333175659179688\\
108.534118652344	0.36041259765625\\
107.726577758789	0.376251220703125\\
106.920684814453	0.38067626953125\\
106.892807006836	0.371536254882813\\
106.098663330078	0.356903076171875\\
104.431182861328	0.315902709960938\\
104.402816772461	0.309051513671875\\
102.753890991211	0.219390869140625\\
102.695587158203	0.214874267578125\\
99.3579864501953	0.0178680419921875\\
98.5258636474609	-0.0298919677734375\\
96.933837890625	-0.0958404541015625\\
96.9043273925781	-0.08990478515625\\
96.1209564208984	-0.116790771484375\\
94.5312805175781	-0.138336181640625\\
94.5029907226563	-0.132980346679688\\
92.2903900146484	-0.130203247070313\\
92.2387237548828	-0.124359130859375\\
89.9908599853516	-0.113784790039063\\
84.084716796875	-0.1090087890625\\
84.0599670410156	-0.11163330078125\\
81.8549652099609	-0.11407470703125\\
80.3638305664063	-0.111175537109375\\
76.482421875	-0.10845947265625\\
69.6983795166016	-0.01416015625\\
69.638916015625	-0.0056915283203125\\
68.1819458007813	0.0135955810546875\\
68.1330413818359	0.008636474609375\\
59.2561187744141	0.0607757568359375\\
59.2096557617188	0.0558013916015625\\
56.9153137207031	0.0662689208984375\\
};
\addplot [color=red, forget plot]
  table[row sep=crcr]{%
139.186660766602	1.17979431152344\\
137.479751586914	-0.5777587890625\\
145.590301513672	-8.45458984375\\
147.297210693359	-6.69705200195313\\
139.186660766602	1.17979431152344\\
};
\addplot [color=black, line width=1.0pt, forget plot]
  table[row sep=crcr]{%
139.971252441406	0.0762786865234375\\
140.557098388672	-0.492691040039063\\
};
\addplot [color=black, line width=1.0pt, forget plot]
  table[row sep=crcr]{%
138.605743408203	-1.32975769042969\\
139.191589355469	-1.89872741699219\\
};
\addplot [color=black, line width=1.0pt, forget plot]
  table[row sep=crcr]{%
140.264175415039	-0.208206176757813\\
138.898666381836	-1.61424255371094\\
};
\addplot [color=black, line width=1.0pt, forget plot]
  table[row sep=crcr]{%
139.581420898438	-0.911224365234375\\
145.320358276367	-6.48480224609375\\
};
\addplot [color=black, line width=1.0pt, forget plot]
  table[row sep=crcr]{%
145.320358276367	-6.48480224609375\\
148.007919311523	-9.26937866210938\\
};
\addplot [color=black, line width=1.0pt, forget plot]
  table[row sep=crcr]{%
145.741928100586	-5.51042175292969\\
146.309066772461	-6.09803771972656\\
};
\addplot [color=black, line width=1.0pt, forget plot]
  table[row sep=crcr]{%
144.331649780273	-6.87156677246094\\
144.898788452148	-7.45918273925781\\
};
\addplot [color=black, line width=1.0pt, forget plot]
  table[row sep=crcr]{%
146.025497436523	-5.80422973632813\\
144.615219116211	-7.16537475585938\\
};
\addplot [color=blue, forget plot]
  table[row sep=crcr]{%
148.974685668945	-8.50471496582031\\
147.233596801758	-10.2283935546875\\
151.995727539063	-15.0386047363281\\
153.73681640625	-13.3149108886719\\
148.974685668945	-8.50471496582031\\
};
\addplot [color=black, line width=1.0pt, forget plot]
  table[row sep=crcr]{%
153.098022460938	-12.6450805664063\\
153.165557861328	-13.4589538574219\\
};
\addplot [color=black, line width=1.0pt, forget plot]
  table[row sep=crcr]{%
151.705139160156	-14.0240325927734\\
151.772674560547	-14.8379058837891\\
};
\addplot [color=black, line width=1.0pt, forget plot]
  table[row sep=crcr]{%
149.590606689453	-9.47508239746094\\
150.165161132813	-10.0554504394531\\
};
\addplot [color=black, line width=1.0pt, forget plot]
  table[row sep=crcr]{%
148.197723388672	-10.8540344238281\\
148.772277832031	-11.4344024658203\\
};
\addplot [color=black, line width=1.0pt, forget plot]
  table[row sep=crcr]{%
153.131790161133	-13.0520172119141\\
151.738906860352	-14.4309692382813\\
};
\addplot [color=black, line width=1.0pt, forget plot]
  table[row sep=crcr]{%
149.877883911133	-9.7652587890625\\
148.485000610352	-11.1442108154297\\
};
\addplot [color=black, line width=1.0pt, forget plot]
  table[row sep=crcr]{%
148.007919311523	-9.26937866210938\\
152.435348510742	-13.7414855957031\\
};
\addplot [color=red, forget plot]
  table[row sep=crcr]{%
131.348373413086	4.86700439453125\\
130.942153930664	2.45091247558594\\
142.091674804688	0.576400756835938\\
142.497894287109	2.99249267578125\\
131.348373413086	4.86700439453125\\
};
\addplot [color=black, line width=1.0pt, forget plot]
  table[row sep=crcr]{%
132.620971679688	4.40460205078125\\
133.426330566406	4.26919555664063\\
};
\addplot [color=black, line width=1.0pt, forget plot]
  table[row sep=crcr]{%
132.296020507813	2.47172546386719\\
133.101379394531	2.33633422851563\\
};
\addplot [color=black, line width=1.0pt, forget plot]
  table[row sep=crcr]{%
133.023651123047	4.33689880371094\\
132.698699951172	2.40402221679688\\
};
\addplot [color=black, line width=1.0pt, forget plot]
  table[row sep=crcr]{%
132.861175537109	3.37046813964844\\
140.750457763672	2.04408264160156\\
};
\addplot [color=black, line width=1.0pt, forget plot]
  table[row sep=crcr]{%
140.750457763672	2.04408264160156\\
143.174911499023	-0.97235107421875\\
};
\addplot [color=black, line width=1.0pt, forget plot]
  table[row sep=crcr]{%
141.258499145508	2.97630310058594\\
141.770126342773	2.33975219726563\\
};
\addplot [color=black, line width=1.0pt, forget plot]
  table[row sep=crcr]{%
139.73078918457	1.7484130859375\\
140.242416381836	1.11186218261719\\
};
\addplot [color=black, line width=1.0pt, forget plot]
  table[row sep=crcr]{%
141.514312744141	2.65802001953125\\
139.986602783203	1.43013000488281\\
};
\addplot [color=blue, forget plot]
  table[row sep=crcr]{%
144.395462036133	-0.800338745117188\\
142.022354125977	-1.40928649902344\\
143.704711914063	-7.96563720703125\\
146.077819824219	-7.356689453125\\
144.395462036133	-0.800338745117188\\
};
\addplot [color=black, line width=1.0pt, forget plot]
  table[row sep=crcr]{%
145.533096313477	-6.44660949707031\\
145.843460083008	-7.20199584960938\\
};
\addplot [color=black, line width=1.0pt, forget plot]
  table[row sep=crcr]{%
143.63459777832	-6.93376159667969\\
143.944961547852	-7.68914794921875\\
};
\addplot [color=black, line width=1.0pt, forget plot]
  table[row sep=crcr]{%
144.437255859375	-1.94891357421875\\
144.640228271484	-2.73995971679688\\
};
\addplot [color=black, line width=1.0pt, forget plot]
  table[row sep=crcr]{%
142.538757324219	-2.43606567382813\\
142.741729736328	-3.22711181640625\\
};
\addplot [color=black, line width=1.0pt, forget plot]
  table[row sep=crcr]{%
145.688278198242	-6.82431030273438\\
143.789779663086	-7.31146240234375\\
};
\addplot [color=black, line width=1.0pt, forget plot]
  table[row sep=crcr]{%
144.53874206543	-2.34443664550781\\
142.640243530273	-2.83158874511719\\
};
\addplot [color=black, line width=1.0pt, forget plot]
  table[row sep=crcr]{%
143.174911499023	-0.97235107421875\\
144.739028930664	-7.06788635253906\\
};
\addplot [color=red, forget plot]
  table[row sep=crcr]{%
126.303283691406	4.03092956542969\\
126.866622924805	1.64657592773438\\
137.869689941406	4.24623107910156\\
137.306335449219	6.63059997558594\\
126.303283691406	4.03092956542969\\
};
\addplot [color=black, line width=1.0pt, forget plot]
  table[row sep=crcr]{%
127.65559387207	4.09869384765625\\
128.450378417969	4.28648376464844\\
};
\addplot [color=black, line width=1.0pt, forget plot]
  table[row sep=crcr]{%
128.106262207031	2.19120788574219\\
128.901062011719	2.37899780273438\\
};
\addplot [color=black, line width=1.0pt, forget plot]
  table[row sep=crcr]{%
128.052978515625	4.19258117675781\\
128.503662109375	2.28511047363281\\
};
\addplot [color=black, line width=1.0pt, forget plot]
  table[row sep=crcr]{%
128.2783203125	3.23884582519531\\
136.06396484375	5.07833862304688\\
};
\addplot [color=black, line width=1.0pt, forget plot]
  table[row sep=crcr]{%
136.06396484375	5.07833862304688\\
139.924697875977	4.81062316894531\\
};
\addplot [color=black, line width=1.0pt, forget plot]
  table[row sep=crcr]{%
135.724411010742	6.084228515625\\
136.539108276367	6.02774047851563\\
};
\addplot [color=black, line width=1.0pt, forget plot]
  table[row sep=crcr]{%
135.588821411133	4.12893676757813\\
136.403518676758	4.07243347167969\\
};
\addplot [color=black, line width=1.0pt, forget plot]
  table[row sep=crcr]{%
136.131759643555	6.05598449707031\\
135.996170043945	4.10067749023438\\
};
\addplot [color=blue, forget plot]
  table[row sep=crcr]{%
140.874801635742	5.59587097167969\\
139.171157836914	3.83517456054688\\
144.035522460938	-0.871612548828125\\
145.739166259766	0.889083862304688\\
140.874801635742	5.59587097167969\\
};
\addplot [color=black, line width=1.0pt, forget plot]
  table[row sep=crcr]{%
145.02978515625	1.53511047363281\\
145.227478027344	0.74273681640625\\
};
\addplot [color=black, line width=1.0pt, forget plot]
  table[row sep=crcr]{%
143.666870117188	0.126556396484375\\
143.864562988281	-0.665817260742188\\
};
\addplot [color=black, line width=1.0pt, forget plot]
  table[row sep=crcr]{%
141.511428833008	4.63896179199219\\
142.09831237793	4.07107543945313\\
};
\addplot [color=black, line width=1.0pt, forget plot]
  table[row sep=crcr]{%
140.148513793945	3.23040771484375\\
140.735397338867	2.66252136230469\\
};
\addplot [color=black, line width=1.0pt, forget plot]
  table[row sep=crcr]{%
145.128631591797	1.13893127441406\\
143.765716552734	-0.269622802734375\\
};
\addplot [color=black, line width=1.0pt, forget plot]
  table[row sep=crcr]{%
141.804870605469	4.35502624511719\\
140.441955566406	2.94645690917969\\
};
\addplot [color=black, line width=1.0pt, forget plot]
  table[row sep=crcr]{%
139.924697875977	4.81062316894531\\
144.447174072266	0.434646606445313\\
};
\addplot [color=red, forget plot]
  table[row sep=crcr]{%
119.014221191406	1.67062377929688\\
119.800857543945	-0.649658203125\\
130.508239746094	2.98043823242188\\
129.721588134766	5.30070495605469\\
119.014221191406	1.67062377929688\\
};
\addplot [color=black, line width=1.0pt, forget plot]
  table[row sep=crcr]{%
120.354042053223	1.86615753173828\\
121.127464294434	2.12837219238281\\
};
\addplot [color=black, line width=1.0pt, forget plot]
  table[row sep=crcr]{%
120.983345031738	0.0099334716796875\\
121.756767272949	0.272148132324219\\
};
\addplot [color=black, line width=1.0pt, forget plot]
  table[row sep=crcr]{%
120.740753173828	1.99726867675781\\
121.370056152344	0.141036987304688\\
};
\addplot [color=black, line width=1.0pt, forget plot]
  table[row sep=crcr]{%
121.055404663086	1.06915283203125\\
128.6318359375	3.63777160644531\\
};
\addplot [color=black, line width=1.0pt, forget plot]
  table[row sep=crcr]{%
128.6318359375	3.63777160644531\\
131.681182861328	6.02070617675781\\
};
\addplot [color=black, line width=1.0pt, forget plot]
  table[row sep=crcr]{%
127.706665039063	4.15852355957031\\
128.350143432617	4.661376953125\\
};
\addplot [color=black, line width=1.0pt, forget plot]
  table[row sep=crcr]{%
128.913528442383	2.61415100097656\\
129.557022094727	3.11701965332031\\
};
\addplot [color=black, line width=1.0pt, forget plot]
  table[row sep=crcr]{%
128.028396606445	4.40995788574219\\
129.235275268555	2.86558532714844\\
};
\addplot [color=blue, forget plot]
  table[row sep=crcr]{%
131.497940063477	7.2396240234375\\
132.128707885742	4.87222290039063\\
138.669281005859	6.61491394042969\\
138.038513183594	8.98231506347656\\
131.497940063477	7.2396240234375\\
};
\addplot [color=black, line width=1.0pt, forget plot]
  table[row sep=crcr]{%
137.113906860352	8.688232421875\\
137.905532836914	8.48753356933594\\
};
\addplot [color=black, line width=1.0pt, forget plot]
  table[row sep=crcr]{%
137.618545532227	6.7943115234375\\
138.410171508789	6.59361267089844\\
};
\addplot [color=black, line width=1.0pt, forget plot]
  table[row sep=crcr]{%
132.646072387695	7.2919921875\\
133.435195922852	7.50225830078125\\
};
\addplot [color=black, line width=1.0pt, forget plot]
  table[row sep=crcr]{%
133.15071105957	5.3980712890625\\
133.939834594727	5.60832214355469\\
};
\addplot [color=black, line width=1.0pt, forget plot]
  table[row sep=crcr]{%
137.509719848633	8.587890625\\
138.014358520508	6.69395446777344\\
};
\addplot [color=black, line width=1.0pt, forget plot]
  table[row sep=crcr]{%
133.040634155273	7.39712524414063\\
133.545272827148	5.50318908691406\\
};
\addplot [color=black, line width=1.0pt, forget plot]
  table[row sep=crcr]{%
131.681182861328	6.02070617675781\\
137.76203918457	7.64091491699219\\
};
\addplot [color=red, line width=0.7pt, forget plot]
  table[row sep=crcr]{%
100.45108795166	5.4295654296875\\
97.9331359863281	5.5452880859375\\
97.0295639038086	5.56227111816406\\
96.3128890991211	5.5570068359375\\
95.5929870605469	5.52977752685547\\
94.8654556274414	5.47464752197266\\
94.1636047363281	5.39230346679688\\
93.4794235229492	5.28035736083984\\
92.838264465332	5.14736938476563\\
92.2628021240234	5.00537109375\\
91.4987258911133	4.78146362304688\\
90.869514465332	4.56935882568359\\
90.2041931152344	4.31976318359375\\
89.5114440917969	4.03382110595703\\
88.7807998657227	3.71157836914063\\
85.1375045776367	2.07815551757813\\
84.1910552978516	1.70532989501953\\
83.5451507568359	1.47867584228516\\
82.8984222412109	1.27272033691406\\
82.2557907104492	1.09410095214844\\
81.6145172119141	0.941184997558594\\
80.9712982177734	0.809562683105469\\
80.3223648071289	0.697654724121094\\
79.6722869873047	0.606742858886719\\
79.027473449707	0.531143188476563\\
78.056640625	0.452171325683594\\
77.0922393798828	0.411285400390625\\
76.1411285400391	0.404373168945313\\
74.8957138061523	0.42022705078125\\
72.7110824584961	0.472869873046875\\
70.1604995727539	0.517555236816406\\
68.8888320922852	0.509559631347656\\
67.959846496582	0.477302551269531\\
66.7269897460938	0.400680541992188\\
65.4789962768555	0.305023193359375\\
62.0825576782227	0.0173568725585938\\
61.1525344848633	-0.0368270874023438\\
59.5917587280273	-0.08831787109375\\
57.0877380371094	-0.121665954589844\\
52.1551742553711	-0.0833587646484375\\
49.7639617919922	-0.091766357421875\\
};
\addplot [color=red, forget plot]
  table[row sep=crcr]{%
98.7652282714844	6.94143676757813\\
98.6414184570313	4.49456787109375\\
109.951950073242	3.92225646972656\\
110.075759887695	6.36912536621094\\
98.7652282714844	6.94143676757813\\
};
\addplot [color=black, line width=1.0pt, forget plot]
  table[row sep=crcr]{%
100.092796325684	6.62895202636719\\
100.908424377441	6.58767700195313\\
};
\addplot [color=black, line width=1.0pt, forget plot]
  table[row sep=crcr]{%
99.9937515258789	4.67145538330078\\
100.809379577637	4.63018035888672\\
};
\addplot [color=black, line width=1.0pt, forget plot]
  table[row sep=crcr]{%
100.500610351563	6.60831451416016\\
100.401565551758	4.65081787109375\\
};
\addplot [color=black, line width=1.0pt, forget plot]
  table[row sep=crcr]{%
100.45108795166	5.62956237792969\\
108.440864562988	5.22528839111328\\
};
\addplot [color=black, line width=1.0pt, forget plot]
  table[row sep=crcr]{%
108.440864562988	5.22528839111328\\
112.310775756836	5.19945526123047\\
};
\addplot [color=black, line width=1.0pt, forget plot]
  table[row sep=crcr]{%
108.039077758789	6.20799255371094\\
108.855728149414	6.20253753662109\\
};
\addplot [color=black, line width=1.0pt, forget plot]
  table[row sep=crcr]{%
108.026000976563	4.24803161621094\\
108.842651367188	4.24258422851563\\
};
\addplot [color=black, line width=1.0pt, forget plot]
  table[row sep=crcr]{%
108.447402954102	6.20526123046875\\
108.434326171875	4.24530792236328\\
};
\addplot [color=blue, forget plot]
  table[row sep=crcr]{%
112.434555053711	6.42583465576172\\
112.460479736328	3.97597503662109\\
119.228851318359	4.04759216308594\\
119.202926635742	6.49745941162109\\
112.434555053711	6.42583465576172\\
};
\addplot [color=black, line width=1.0pt, forget plot]
  table[row sep=crcr]{%
118.277351379395	6.50495910644531\\
118.908760070801	5.98702239990234\\
};
\addplot [color=black, line width=1.0pt, forget plot]
  table[row sep=crcr]{%
118.29808807373	4.54506683349609\\
118.929496765137	4.02713012695313\\
};
\addplot [color=black, line width=1.0pt, forget plot]
  table[row sep=crcr]{%
113.560005187988	6.19273376464844\\
114.376625061035	6.20137023925781\\
};
\addplot [color=black, line width=1.0pt, forget plot]
  table[row sep=crcr]{%
113.580741882324	4.23284149169922\\
114.397361755371	4.24148559570313\\
};
\addplot [color=black, line width=1.0pt, forget plot]
  table[row sep=crcr]{%
118.593055725098	6.24598693847656\\
118.613792419434	4.28610229492188\\
};
\addplot [color=black, line width=1.0pt, forget plot]
  table[row sep=crcr]{%
113.968315124512	6.19705200195313\\
113.989051818848	4.23715972900391\\
};
\addplot [color=black, line width=1.0pt, forget plot]
  table[row sep=crcr]{%
112.310775756836	5.19945526123047\\
118.603424072266	5.26604461669922\\
};
\addplot [color=red, forget plot]
  table[row sep=crcr]{%
88.8193130493164	5.35151672363281\\
89.6262435913086	3.03821563720703\\
100.319351196289	6.76824188232422\\
99.5124206542969	9.08153533935547\\
88.8193130493164	5.35151672363281\\
};
\addplot [color=black, line width=1.0pt, forget plot]
  table[row sep=crcr]{%
90.1668090820313	5.56208038330078\\
90.9379119873047	5.8310546875\\
};
\addplot [color=black, line width=1.0pt, forget plot]
  table[row sep=crcr]{%
90.8123626708984	3.71144104003906\\
91.5834655761719	3.98041534423828\\
};
\addplot [color=black, line width=1.0pt, forget plot]
  table[row sep=crcr]{%
90.552360534668	5.69657135009766\\
91.1979141235352	3.84593200683594\\
};
\addplot [color=black, line width=1.0pt, forget plot]
  table[row sep=crcr]{%
90.8751373291016	4.77124786376953\\
98.4287643432617	7.40614318847656\\
};
\addplot [color=black, line width=1.0pt, forget plot]
  table[row sep=crcr]{%
98.4287643432617	7.40614318847656\\
102.291809082031	7.63809204101563\\
};
\addplot [color=black, line width=1.0pt, forget plot]
  table[row sep=crcr]{%
97.9624252319336	8.35990905761719\\
98.7776260375977	8.40885162353516\\
};
\addplot [color=black, line width=1.0pt, forget plot]
  table[row sep=crcr]{%
98.0799026489258	6.40343475341797\\
98.8951034545898	6.45237731933594\\
};
\addplot [color=black, line width=1.0pt, forget plot]
  table[row sep=crcr]{%
98.3700256347656	8.38438415527344\\
98.4875030517578	6.42790222167969\\
};
\addplot [color=blue, forget plot]
  table[row sep=crcr]{%
103.056663513184	8.60469818115234\\
101.758949279785	6.52661895751953\\
107.500183105469	2.94133758544922\\
108.797897338867	5.01941680908203\\
103.056663513184	8.60469818115234\\
};
\addplot [color=black, line width=1.0pt, forget plot]
  table[row sep=crcr]{%
107.888221740723	5.45058441162109\\
108.408988952637	4.82149505615234\\
};
\addplot [color=black, line width=1.0pt, forget plot]
  table[row sep=crcr]{%
106.850044250488	3.78811645507813\\
107.370811462402	3.15902709960938\\
};
\addplot [color=black, line width=1.0pt, forget plot]
  table[row sep=crcr]{%
103.879348754883	7.80210113525391\\
104.572036743164	7.36952972412109\\
};
\addplot [color=black, line width=1.0pt, forget plot]
  table[row sep=crcr]{%
102.841171264648	6.13963317871094\\
103.53385925293	5.70706176757813\\
};
\addplot [color=black, line width=1.0pt, forget plot]
  table[row sep=crcr]{%
108.14860534668	5.13603973388672\\
107.110427856445	3.47357177734375\\
};
\addplot [color=black, line width=1.0pt, forget plot]
  table[row sep=crcr]{%
104.225692749023	7.5858154296875\\
103.187515258789	5.92334747314453\\
};
\addplot [color=black, line width=1.0pt, forget plot]
  table[row sep=crcr]{%
102.291809082031	7.63809204101563\\
107.629516601563	4.3048095703125\\
};
\addplot [color=red, forget plot]
  table[row sep=crcr]{%
78.2158432006836	1.74242401123047\\
78.6829605102539	-0.662628173828125\\
89.8002090454102	1.49662017822266\\
89.3330917358398	3.90167236328125\\
78.2158432006836	1.74242401123047\\
};
\addplot [color=black, line width=1.0pt, forget plot]
  table[row sep=crcr]{%
79.5796051025391	1.75772094726563\\
80.3812866210938	1.91342926025391\\
};
\addplot [color=black, line width=1.0pt, forget plot]
  table[row sep=crcr]{%
79.9533081054688	-0.16632080078125\\
80.7549896240234	-0.0106124877929688\\
};
\addplot [color=black, line width=1.0pt, forget plot]
  table[row sep=crcr]{%
79.9804458618164	1.83557891845703\\
80.3541488647461	-0.088470458984375\\
};
\addplot [color=black, line width=1.0pt, forget plot]
  table[row sep=crcr]{%
80.1672973632813	0.873558044433594\\
88.0205459594727	2.39884948730469\\
};
\addplot [color=black, line width=1.0pt, forget plot]
  table[row sep=crcr]{%
88.0205459594727	2.39884948730469\\
91.5146408081055	4.06264495849609\\
};
\addplot [color=black, line width=1.0pt, forget plot]
  table[row sep=crcr]{%
87.2305526733398	3.10810852050781\\
87.967887878418	3.45921325683594\\
};
\addplot [color=black, line width=1.0pt, forget plot]
  table[row sep=crcr]{%
88.0732040405273	1.33849334716797\\
88.8105392456055	1.68959045410156\\
};
\addplot [color=black, line width=1.0pt, forget plot]
  table[row sep=crcr]{%
87.5992202758789	3.28366088867188\\
88.4418716430664	1.51404571533203\\
};
\addplot [color=blue, forget plot]
  table[row sep=crcr]{%
90.8167953491211	5.07868194580078\\
92.4193649291992	3.22550964355469\\
97.5392379760742	7.65302276611328\\
95.9366683959961	9.50619506835938\\
90.8167953491211	5.07868194580078\\
};
\addplot [color=black, line width=1.0pt, forget plot]
  table[row sep=crcr]{%
95.3316345214844	8.64540100097656\\
95.9356231689453	9.19507598876953\\
};
\addplot [color=black, line width=1.0pt, forget plot]
  table[row sep=crcr]{%
96.6136932373047	7.16285705566406\\
97.2176818847656	7.71253204345703\\
};
\addplot [color=black, line width=1.0pt, forget plot]
  table[row sep=crcr]{%
91.8264236450195	5.62787628173828\\
92.4441452026367	6.16207122802734\\
};
\addplot [color=black, line width=1.0pt, forget plot]
  table[row sep=crcr]{%
93.1084823608398	4.14533233642578\\
93.726203918457	4.67952728271484\\
};
\addplot [color=black, line width=1.0pt, forget plot]
  table[row sep=crcr]{%
95.6336288452148	8.92023468017578\\
96.9156875610352	7.43769836425781\\
};
\addplot [color=black, line width=1.0pt, forget plot]
  table[row sep=crcr]{%
92.1352844238281	5.89497375488281\\
93.4173431396484	4.41242980957031\\
};
\addplot [color=black, line width=1.0pt, forget plot]
  table[row sep=crcr]{%
91.5146408081055	4.06264495849609\\
96.274658203125	8.17896270751953\\
};
\addplot [color=green, dashed, line width=0.7pt, forget plot]
  table[row sep=crcr]{%
121.478637695313	-4.05079650878906\\
121.497741699219	-4.05902862548828\\
120.236221313477	-3.51160430908203\\
119.287994384766	-3.12163543701172\\
118.306549072266	-2.74651336669922\\
117.320190429688	-2.40360260009766\\
116.547714233398	-2.15924072265625\\
115.749984741211	-1.92912292480469\\
114.920471191406	-1.71216583251953\\
114.081985473633	-1.51445007324219\\
113.247268676758	-1.33808898925781\\
112.399185180664	-1.17906951904297\\
111.537261962891	-1.03933715820313\\
110.636993408203	-0.9136962890625\\
109.427139282227	-0.776123046875\\
108.207695007324	-0.671058654785156\\
106.978370666504	-0.593055725097656\\
105.740982055664	-0.539642333984375\\
103.843322753906	-0.4873046875\\
99.9460144042969	-0.402732849121094\\
98.3853454589844	-0.336532592773438\\
96.5800552368164	-0.227249145507813\\
92.9270248413086	0.0426788330078125\\
91.1208267211914	0.173347473144531\\
89.9517669677734	0.230827331542969\\
88.4506759643555	0.275558471679688\\
86.6187973022461	0.303390502929688\\
85.0757904052734	0.302421569824219\\
82.0813598632813	0.259956359863281\\
77.9451675415039	0.187255859375\\
72.2341384887695	0.121055603027344\\
70.4274444580078	0.0992507934570313\\
65.638053894043	0.0331649780273438\\
63.2347106933594	0.00174713134765625\\
60.1975555419922	-0.0549774169921875\\
56.5218505859375	-0.140792846679688\\
49.8620223999023	-0.236923217773438\\
};
\addplot [color=red, forget plot]
  table[row sep=crcr]{%
120.371078491211	-2.22419738769531\\
119.383285522461	-4.46624755859375\\
129.747009277344	-9.03228759765625\\
130.734802246094	-6.79025268554688\\
120.371078491211	-2.22419738769531\\
};
\addplot [color=black, line width=1.0pt, forget plot]
  table[row sep=crcr]{%
121.50008392334	-2.98934936523438\\
122.247428894043	-3.31861114501953\\
};
\addplot [color=black, line width=1.0pt, forget plot]
  table[row sep=crcr]{%
120.709846496582	-4.78298187255859\\
121.457191467285	-5.11224365234375\\
};
\addplot [color=black, line width=1.0pt, forget plot]
  table[row sep=crcr]{%
121.873756408691	-3.15397644042969\\
121.083518981934	-4.94761657714844\\
};
\addplot [color=black, line width=1.0pt, forget plot]
  table[row sep=crcr]{%
121.478637695313	-4.05079650878906\\
128.799591064453	-7.27626037597656\\
};
\addplot [color=black, line width=1.0pt, forget plot]
  table[row sep=crcr]{%
128.799591064453	-7.27626037597656\\
132.300537109375	-8.92558288574219\\
};
\addplot [color=black, line width=1.0pt, forget plot]
  table[row sep=crcr]{%
128.847854614258	-6.21568298339844\\
129.586654663086	-6.56373596191406\\
};
\addplot [color=black, line width=1.0pt, forget plot]
  table[row sep=crcr]{%
128.01252746582	-7.98878479003906\\
128.751327514648	-8.33682250976563\\
};
\addplot [color=black, line width=1.0pt, forget plot]
  table[row sep=crcr]{%
129.217254638672	-6.38970947265625\\
128.381927490234	-8.16279602050781\\
};
\addplot [color=blue, forget plot]
  table[row sep=crcr]{%
132.842819213867	-7.81866455078125\\
132.015701293945	-10.1248321533203\\
138.387054443359	-12.4099273681641\\
139.214172363281	-10.1037750244141\\
132.842819213867	-7.81866455078125\\
};
\addplot [color=black, line width=1.0pt, forget plot]
  table[row sep=crcr]{%
138.333404541016	-9.78457641601563\\
138.776428222656	-10.4706420898438\\
};
\addplot [color=black, line width=1.0pt, forget plot]
  table[row sep=crcr]{%
137.671722412109	-11.6295166015625\\
138.11474609375	-12.3155670166016\\
};
\addplot [color=black, line width=1.0pt, forget plot]
  table[row sep=crcr]{%
133.817092895508	-8.42837524414063\\
134.585800170898	-8.70407104492188\\
};
\addplot [color=black, line width=1.0pt, forget plot]
  table[row sep=crcr]{%
133.155410766602	-10.2733001708984\\
133.924118041992	-10.5490112304688\\
};
\addplot [color=black, line width=1.0pt, forget plot]
  table[row sep=crcr]{%
138.554916381836	-10.1276092529297\\
137.89323425293	-11.9725341796875\\
};
\addplot [color=black, line width=1.0pt, forget plot]
  table[row sep=crcr]{%
134.201446533203	-8.56622314453125\\
133.539764404297	-10.4111480712891\\
};
\addplot [color=black, line width=1.0pt, forget plot]
  table[row sep=crcr]{%
132.300537109375	-8.92558288574219\\
138.224075317383	-11.0500793457031\\
};
\end{axis}
\end{tikzpicture}%

%% file: trajectory_beta23_straight.tex
%
%
\definecolor{mycolor1}{rgb}{0.50000,0.00000,0.50000}%
\begin{tikzpicture}

\begin{axis}[%
width=\figurewidth,
height=\figureheight,
at={(0\figurewidth,0\figureheight)},
scale only axis,
xmin=-1.2,
xmax=1.2,
xtick={-0.8,0,0.8},
xlabel style={font=\color{white!15!black},at={(axis description cs:0.5,-0.14)},anchor=north},
ylabel style={font=\color{white!15!black},at={(axis description cs:-0.2,.5)},anchor=south},
xlabel={$\beta_3$ [rad]},
ymin=-1.2,
ymax=1.2,
ytick={-0.8,0,0.8},
ylabel={$\beta_2$ [rad]},
axis background/.style={fill=white},
xmajorgrids,
ymajorgrids
]

\addplot[area legend, draw=gray, fill=gray, fill opacity=0.3, forget plot]
table[row sep=crcr] {%
x	y\\
-0.67	-0.8\\
0.87	-0.4\\
0.67	0.8\\
-0.87	0.4\\
-0.67	-0.8\\
}--cycle;
\addplot [color=blue, line width=0.7pt, forget plot]
  table[row sep=crcr]{%
-0.0550615191459656	0.0772539377212524\\
-0.0677133798599243	0.109508395195007\\
-0.0873400568962097	0.154980599880219\\
-0.104607343673706	0.192079126834869\\
-0.119558453559875	0.221760869026184\\
-0.13688987493515	0.254535913467407\\
-0.155728995800018	0.287248253822327\\
-0.168928921222687	0.308008849620819\\
-0.182198107242584	0.327187240123749\\
-0.188238024711609	0.334860801696777\\
-0.195813536643982	0.346005797386169\\
-0.215824663639069	0.370519042015076\\
-0.229985475540161	0.386908769607544\\
-0.243903636932373	0.401207625865936\\
-0.258482277393341	0.416126906871796\\
-0.273523449897766	0.430745244026184\\
-0.281228303909302	0.437602460384369\\
-0.311410367488861	0.460045874118805\\
-0.340815007686615	0.475401282310486\\
-0.369491159915924	0.481996476650238\\
-0.383066713809967	0.482492208480835\\
-0.390631318092346	0.483681261539459\\
-0.397508382797241	0.483457803726196\\
-0.405301630496979	0.484745502471924\\
-0.426492929458618	0.482856392860413\\
-0.433675408363342	0.481702864170074\\
-0.446985065937042	0.477287888526917\\
-0.453164339065552	0.474087357521057\\
-0.460002183914185	0.471516489982605\\
-0.472889482975006	0.464445769786835\\
-0.497974514961243	0.445400357246399\\
-0.503162086009979	0.439022600650787\\
-0.508879482746124	0.432917773723602\\
-0.51374626159668	0.425634682178497\\
-0.51993191242218	0.419321060180664\\
-0.524534046649933	0.411232829093933\\
-0.530507922172546	0.404449760913849\\
-0.534889280796051	0.395846843719482\\
-0.54084187746048	0.389313817024231\\
-0.545172333717346	0.381002128124237\\
-0.551214993000031	0.375247776508331\\
-0.555577754974365	0.367526650428772\\
-0.566064119338989	0.354975819587708\\
-0.572211682796478	0.350043058395386\\
-0.576587498188019	0.342711746692657\\
-0.587283372879028	0.331109583377838\\
-0.593543291091919	0.3267702460289\\
-0.597890079021454	0.319367468357086\\
-0.60407280921936	0.314590632915497\\
-0.609324336051941	0.308342814445496\\
-0.61555403470993	0.303518652915955\\
-0.625885784626007	0.290189802646637\\
-0.631932199001312	0.284605085849762\\
-0.646629393100739	0.261833190917969\\
-0.657318413257599	0.249315857887268\\
-0.661417365074158	0.238607347011566\\
-0.666857779026031	0.231689691543579\\
-0.671165704727173	0.222979187965393\\
-0.67454868555069	0.213089287281036\\
-0.679615139961243	0.205518484115601\\
-0.683716237545013	0.196719169616699\\
-0.688553512096405	0.188864171504974\\
-0.691568315029144	0.178653657436371\\
-0.700046539306641	0.16279548406601\\
-0.70393705368042	0.15219247341156\\
-0.706511735916138	0.141601741313934\\
-0.710462212562561	0.1324662566185\\
-0.71351283788681	0.122256636619568\\
-0.717165052890778	0.112724661827087\\
-0.719302952289581	0.101568818092346\\
-0.722690999507904	0.0918607711791992\\
-0.725336074829102	0.0814613699913025\\
-0.728505253791809	0.0716339945793152\\
-0.73032134771347	0.0605961680412292\\
-0.73330169916153	0.0510028600692749\\
-0.735446989536285	0.0402928590774536\\
-0.738112151622772	0.0302718877792358\\
-0.739568829536438	0.0192424058914185\\
-0.744278371334076	-0.00106751918792725\\
-0.746343314647675	-0.0116685628890991\\
-0.752495169639587	-0.0559478998184204\\
-0.754230499267578	-0.0921710729598999\\
-0.75440502166748	-0.10331666469574\\
-0.752834856510162	-0.144927084445953\\
-0.749361872673035	-0.197088062763214\\
-0.738671004772186	-0.349829256534576\\
-0.733089625835419	-0.397180199623108\\
-0.729913890361786	-0.414498209953308\\
-0.725462794303894	-0.432337999343872\\
-0.718553304672241	-0.455411911010742\\
-0.716063559055328	-0.463432192802429\\
-0.689827740192413	-0.534143626689911\\
-0.675926148891449	-0.567698776721954\\
-0.666802227497101	-0.58750194311142\\
-0.657543838024139	-0.60540497303009\\
-0.64979076385498	-0.619009017944336\\
-0.641306459903717	-0.632074594497681\\
-0.634069621562958	-0.641351640224457\\
-0.611148178577423	-0.66671895980835\\
-0.599225521087646	-0.678192377090454\\
-0.574309289455414	-0.701568067073822\\
-0.559922277927399	-0.712773263454437\\
-0.556370973587036	-0.714246094226837\\
-0.547650933265686	-0.719009399414063\\
-0.535482943058014	-0.723591446876526\\
-0.531625509262085	-0.724886059761047\\
-0.518019318580627	-0.731209218502045\\
-0.513860106468201	-0.732867062091827\\
-0.508259117603302	-0.736046850681305\\
-0.504079878330231	-0.737202227115631\\
-0.494751691818237	-0.740946173667908\\
-0.473824381828308	-0.748610019683838\\
-0.46912008523941	-0.749640166759491\\
-0.465486764907837	-0.749269843101501\\
-0.455187559127808	-0.750803768634796\\
-0.450358867645264	-0.750689327716827\\
-0.448018550872803	-0.747666656970978\\
-0.443885385990143	-0.746095061302185\\
-0.440481424331665	-0.743574857711792\\
-0.436263203620911	-0.742022097110748\\
-0.427308917045593	-0.733128905296326\\
-0.410928249359131	-0.724272668361664\\
-0.407371163368225	-0.720096409320831\\
-0.395475506782532	-0.713579118251801\\
-0.393173813819885	-0.710352599620819\\
-0.388296365737915	-0.709387302398682\\
-0.377985417842865	-0.708870232105255\\
-0.373966455459595	-0.707213401794434\\
-0.351622045040131	-0.704429149627686\\
-0.346964001655579	-0.704387068748474\\
-0.341103255748749	-0.705171406269073\\
-0.337863445281982	-0.703075528144836\\
-0.332045197486877	-0.70379650592804\\
-0.322188675403595	-0.702486753463745\\
-0.307520270347595	-0.699518859386444\\
-0.301244556903839	-0.699703216552734\\
-0.296402513980865	-0.697859704494476\\
-0.29058974981308	-0.698325276374817\\
-0.285651862621307	-0.696546375751495\\
-0.27962201833725	-0.695379137992859\\
-0.275385916233063	-0.69314432144165\\
-0.268210291862488	-0.692747592926025\\
-0.256709814071655	-0.690272986888886\\
-0.247044503688812	-0.684839129447937\\
-0.242456495761871	-0.681464970111847\\
-0.236183285713196	-0.6803337931633\\
-0.231430530548096	-0.676580846309662\\
-0.226703941822052	-0.673685312271118\\
-0.220524370670319	-0.67159104347229\\
-0.216091513633728	-0.667751431465149\\
-0.20631331205368	-0.661649465560913\\
-0.199494302272797	-0.661316812038422\\
-0.194478213787079	-0.659000158309937\\
-0.187189102172852	-0.659776329994202\\
-0.181814193725586	-0.65802127122879\\
-0.177703857421875	-0.655437588691711\\
-0.172004580497742	-0.653085470199585\\
-0.164770185947418	-0.653110504150391\\
-0.160218000411987	-0.650803208351135\\
-0.1540567278862	-0.648632168769836\\
-0.144238352775574	-0.643814265727997\\
-0.136868715286255	-0.643827795982361\\
-0.130293190479279	-0.64189201593399\\
-0.12574028968811	-0.63972544670105\\
-0.120104968547821	-0.637679100036621\\
-0.108312726020813	-0.634647607803345\\
-0.102163374423981	-0.632974147796631\\
-0.0938567519187927	-0.63437157869339\\
-0.08771812915802	-0.632646501064301\\
-0.078315794467926	-0.634101629257202\\
-0.0752522349357605	-0.628790497779846\\
-0.0656023621559143	-0.630446374416351\\
-0.0589495897293091	-0.62921929359436\\
-0.0521442294120789	-0.628596723079681\\
-0.0468756556510925	-0.625478446483612\\
-0.0389268398284912	-0.624805569648743\\
-0.0323563814163208	-0.622876107692719\\
-0.023613452911377	-0.623363733291626\\
-0.0194348692893982	-0.617536008358002\\
-0.0106091499328613	-0.617877602577209\\
-0.00476402044296265	-0.615866601467133\\
0.000522613525390625	-0.609089910984039\\
0.00825047492980957	-0.608466446399689\\
0.0161344408988953	-0.605787813663483\\
0.02195805311203	-0.603129625320435\\
0.0294866561889648	-0.599032878875732\\
0.0362712144851685	-0.596430957317352\\
0.0428162813186646	-0.593277275562286\\
0.0508379340171814	-0.593270659446716\\
0.0564551949501038	-0.587363958358765\\
0.0656419396400452	-0.587951600551605\\
0.0727830529212952	-0.586163938045502\\
0.0785746574401855	-0.583232581615448\\
0.0866708755493164	-0.580314874649048\\
0.115669012069702	-0.573058724403381\\
0.124094247817993	-0.573790907859802\\
0.131427764892578	-0.571952223777771\\
0.141313076019287	-0.573872983455658\\
0.148629546165466	-0.572055101394653\\
0.155097603797913	-0.566649794578552\\
0.161280453205109	-0.56390917301178\\
0.178066730499268	-0.560479044914246\\
0.185661494731903	-0.558424353599548\\
0.195573925971985	-0.559301674365997\\
0.202470123767853	-0.555002927780151\\
0.208243310451508	-0.550069868564606\\
0.216472327709198	-0.54626852273941\\
0.222409188747406	-0.541542410850525\\
0.230546593666077	-0.537132501602173\\
0.237373411655426	-0.531226634979248\\
0.243704497814178	-0.524963557720184\\
0.248669028282166	-0.518104493618011\\
0.258430778980255	-0.515198349952698\\
0.260695517063141	-0.502753674983978\\
0.267965614795685	-0.496478140354156\\
0.27253520488739	-0.488335907459259\\
0.281497597694397	-0.483560144901276\\
0.284490704536438	-0.470634460449219\\
0.292273104190826	-0.465714931488037\\
0.296606481075287	-0.457222759723663\\
0.302795648574829	-0.447945296764374\\
0.306236028671265	-0.43615186214447\\
0.313988149166107	-0.432112216949463\\
0.323493182659149	-0.415030539035797\\
0.329701840877533	-0.406902730464935\\
0.336468875408173	-0.402939736843109\\
0.340961515903473	-0.390480101108551\\
0.346809387207031	-0.383296251296997\\
0.354172587394714	-0.378952860832214\\
0.359916150569916	-0.371502280235291\\
0.370337665081024	-0.353648960590363\\
0.3773512840271	-0.348318457603455\\
0.382493019104004	-0.33921480178833\\
0.389492690563202	-0.333792865276337\\
0.393417239189148	-0.32336562871933\\
0.398674488067627	-0.311950623989105\\
0.411089241504669	-0.297005534172058\\
0.418789029121399	-0.295068264007568\\
0.422114610671997	-0.278292894363403\\
0.422526001930237	-0.257491946220398\\
0.427305519580841	-0.247288763523102\\
0.429563939571381	-0.231372714042664\\
0.441857218742371	-0.199748516082764\\
0.443001687526703	-0.183148801326752\\
0.446360051631927	-0.172078430652618\\
0.445310175418854	-0.151997864246368\\
0.4464071393013	-0.137146174907684\\
0.445383071899414	-0.119876325130463\\
0.446637630462646	-0.10478413105011\\
0.445816278457642	-0.0882766842842102\\
0.446892499923706	-0.0739361047744751\\
0.447460055351257	-0.0328898429870605\\
0.446924567222595	-0.00428694486618042\\
0.446218132972717	0.0103678703308105\\
0.444414377212524	0.0276671648025513\\
0.443865120410919	0.0525269508361816\\
0.44332093000412	0.0648543834686279\\
0.438569962978363	0.0977684855461121\\
0.435681343078613	0.123864710330963\\
0.429622292518616	0.156365156173706\\
0.422958672046661	0.18780529499054\\
0.411625862121582	0.235059440135956\\
0.397659420967102	0.284924924373627\\
0.388239085674286	0.314991116523743\\
0.378208875656128	0.344466388225555\\
0.367223381996155	0.374119579792023\\
0.352353513240814	0.409728944301605\\
0.340869903564453	0.43579363822937\\
0.317457020282745	0.480750322341919\\
0.311344861984253	0.491690218448639\\
0.305646002292633	0.500378012657166\\
0.2993443608284	0.511125862598419\\
0.293471276760101	0.519257307052612\\
0.281649053096771	0.537427961826324\\
0.268436968326569	0.555463910102844\\
0.262780368328094	0.562978506088257\\
0.25664746761322	0.568808138370514\\
0.250207543373108	0.577376902103424\\
0.226223886013031	0.598248541355133\\
0.214582324028015	0.604411244392395\\
0.192804038524628	0.609538614749908\\
0.181719839572906	0.608519554138184\\
0.166462540626526	0.608666241168976\\
0.15663093328476	0.605787575244904\\
0.151404201984406	0.603545308113098\\
0.146602213382721	0.603045344352722\\
0.123207747936249	0.587564051151276\\
0.119678318500519	0.583547472953796\\
0.115456402301788	0.580089330673218\\
0.103535592556	0.567815482616425\\
0.0909368991851807	0.557933986186981\\
0.0821847319602966	0.552095472812653\\
0.0592455863952637	0.539150774478912\\
0.0535385608673096	0.536912620067596\\
0.0490872859954834	0.533692002296448\\
0.0373557209968567	0.529520273208618\\
0.0325350165367126	0.527705252170563\\
0.027827262878418	0.524638473987579\\
0.022405207157135	0.521905720233917\\
0.0171611905097961	0.51832389831543\\
0.0118244886398315	0.515697956085205\\
0.00188195705413818	0.508811354637146\\
-0.0040745735168457	0.506274938583374\\
-0.00930923223495483	0.5027756690979\\
-0.0141324400901794	0.496883273124695\\
-0.0186705589294434	0.493322134017944\\
-0.0231968760490417	0.488721072673798\\
-0.0276846289634705	0.48186457157135\\
-0.032671332359314	0.477558493614197\\
-0.0414689183235168	0.466815173625946\\
-0.053435742855072	0.448371708393097\\
-0.0578374862670898	0.442348599433899\\
-0.0614123344421387	0.435574769973755\\
-0.0656155943870544	0.429324269294739\\
-0.0705159306526184	0.424580931663513\\
-0.075303316116333	0.42085063457489\\
-0.0803069472312927	0.414422750473022\\
-0.0847837924957275	0.409783780574799\\
-0.0894913077354431	0.403860688209534\\
-0.0948749780654907	0.398052990436554\\
-0.0975695848464966	0.390240907669067\\
-0.10087525844574	0.382618427276611\\
-0.112830996513367	0.362051367759705\\
-0.115980446338654	0.35492467880249\\
-0.119258880615234	0.345723271369934\\
-0.122917413711548	0.338217675685883\\
-0.128568351268768	0.32115113735199\\
-0.130513787269592	0.309725344181061\\
-0.135184764862061	0.291736781597137\\
-0.143674850463867	0.238735556602478\\
-0.147109031677246	0.157281517982483\\
-0.146630644798279	0.0693600177764893\\
-0.143949270248413	0.0243889689445496\\
-0.137313306331635	-0.0273038744926453\\
-0.128672480583191	-0.0788843035697937\\
-0.116073668003082	-0.125947833061218\\
-0.109689056873322	-0.144497036933899\\
-0.0966092944145203	-0.174093425273895\\
-0.0930946469306946	-0.180284678936005\\
-0.090232789516449	-0.186692833900452\\
-0.0781802535057068	-0.204932868480682\\
-0.074626088142395	-0.209710717201233\\
-0.072620153427124	-0.211419820785522\\
-0.0662063956260681	-0.219140768051147\\
-0.0630936622619629	-0.222255766391754\\
-0.060854434967041	-0.222978115081787\\
-0.0551717281341553	-0.227770686149597\\
-0.0437561273574829	-0.233978390693665\\
-0.0382795929908752	-0.234964489936829\\
-0.0296006202697754	-0.237184226512909\\
-0.0242535471916199	-0.236369788646698\\
-0.0156417489051819	-0.237419843673706\\
-0.00504189729690552	-0.234782695770264\\
-0.00222718715667725	-0.234597682952881\\
0.00109386444091797	-0.233389973640442\\
0.00268679857254028	-0.230804741382599\\
0.00533580780029297	-0.228502094745636\\
0.00785094499588013	-0.22819197177887\\
0.0152615904808044	-0.222720086574554\\
0.0174797177314758	-0.221092224121094\\
0.020682692527771	-0.219741821289063\\
0.0251529812812805	-0.214831054210663\\
0.027013897895813	-0.2128826379776\\
0.0320444703102112	-0.209012150764465\\
0.0363554358482361	-0.203407108783722\\
0.0393238067626953	-0.194969475269318\\
0.0415908694267273	-0.192304730415344\\
0.0429618954658508	-0.188772261142731\\
0.0450442433357239	-0.185915648937225\\
0.0465954542160034	-0.182203114032745\\
0.0485188364982605	-0.179140567779541\\
0.0514972805976868	-0.17101526260376\\
0.0521813631057739	-0.16639506816864\\
0.0537894368171692	-0.162502467632294\\
0.0543475151062012	-0.157720386981964\\
0.0557928085327148	-0.154085755348206\\
0.057618260383606	-0.14519214630127\\
0.0588939189910889	-0.142030537128448\\
0.0630176067352295	-0.126116871833801\\
0.0710303783416748	-0.0937954783439636\\
0.0739849805831909	-0.0733270049095154\\
0.0748513340950012	-0.0709710717201233\\
0.0755714178085327	-0.0672765970230103\\
0.0768711566925049	-0.055397093296051\\
0.0785080790519714	-0.0246863961219788\\
0.0779338479042053	-0.0126879811286926\\
0.0775437355041504	-0.00238192081451416\\
0.0764049887657166	0.0120610594749451\\
0.0748156309127808	0.0230692625045776\\
0.0730180740356445	0.0339486598968506\\
0.0675789713859558	0.0562463998794556\\
0.0574737191200256	0.0876300930976868\\
0.0514048933982849	0.100302278995514\\
0.0480674505233765	0.106056034564972\\
0.0355161428451538	0.121166944503784\\
0.0318346619606018	0.123702764511108\\
0.0273785591125488	0.127108752727509\\
0.0266493558883667	0.126866340637207\\
0.021942138671875	0.129179120063782\\
0.0158610939979553	0.132993042469025\\
0.0107759237289429	0.134580671787262\\
0.00906997919082642	0.135812819004059\\
0.00205975770950317	0.137312531471252\\
0.000900506973266602	0.13636726140976\\
-5.11407852172852e-05	0.134748935699463\\
-0.00147801637649536	0.134605050086975\\
-0.00505298376083374	0.131229102611542\\
-0.00704872608184814	0.128456890583038\\
-0.00955486297607422	0.126352906227112\\
-0.0121217966079712	0.124861717224121\\
-0.0156075954437256	0.120581805706024\\
-0.0174144506454468	0.120474517345428\\
-0.0197440981864929	0.117750704288483\\
-0.0213691592216492	0.117177486419678\\
-0.0213073492050171	0.113198101520538\\
-0.0221894383430481	0.111076951026917\\
-0.0267435908317566	0.0891811847686768\\
-0.0277455449104309	0.0843106508255005\\
-0.0285108685493469	0.0718218088150024\\
-0.0310059189796448	0.0586836338043213\\
-0.0328472852706909	0.0363042950630188\\
-0.0323506593704224	0.0221171975135803\\
-0.0295277833938599	0.00233054161071777\\
-0.0264608860015869	-0.0119748711585999\\
-0.0229044556617737	-0.0235023498535156\\
-0.0188992619514465	-0.0334550738334656\\
-0.0151665210723877	-0.0397601127624512\\
-0.0123856663703918	-0.044269859790802\\
-0.0103799104690552	-0.0463542342185974\\
-0.00933420658111572	-0.0468568205833435\\
-0.00763893127441406	-0.0486264824867249\\
-0.00736773014068604	-0.0478366613388062\\
-0.00642985105514526	-0.0484738945960999\\
-0.00540071725845337	-0.0491312146186829\\
-5.63859939575195e-05	-0.0523401498794556\\
0.000541269779205322	-0.0524324178695679\\
0.0021405816078186	-0.0538941621780396\\
0.00518703460693359	-0.0544916987419128\\
0.00671738386154175	-0.0541108846664429\\
0.00781089067459106	-0.0549056529998779\\
0.0114061832427979	-0.0532690286636353\\
0.0134130120277405	-0.0460516214370728\\
0.0145124793052673	-0.0450169444084167\\
0.0157235860824585	-0.039512038230896\\
0.0156576037406921	-0.0346103310585022\\
0.0160431861877441	-0.033812403678894\\
0.0166450142860413	-0.0309116840362549\\
0.0171827077865601	-0.0285425186157227\\
0.0172342658042908	-0.0247671008110046\\
0.0177299380302429	-0.0143098831176758\\
0.0184472799301147	-0.000997662544250488\\
0.0166876316070557	0.014958918094635\\
0.0119951367378235	0.0282639861106873\\
0.00907611846923828	0.0341644287109375\\
0.00556319952011108	0.0384870171546936\\
0.00289714336395264	0.0413160920143127\\
0.00214320421218872	0.0411520600318909\\
0.0022997260093689	0.0405977964401245\\
-0.000880658626556396	0.0433554649353027\\
-0.000782310962677002	0.0425795316696167\\
-0.00335586071014404	0.0433677434921265\\
-0.00435823202133179	0.0438728928565979\\
-0.00521171092987061	0.0434693694114685\\
-0.00620150566101074	0.0437347292900085\\
-0.00724297761917114	0.0414143204689026\\
-0.0077860951423645	0.0344487428665161\\
-0.00887364149093628	0.0349735617637634\\
-0.0113958120346069	0.0275351405143738\\
-0.0120913982391357	0.0272256731987\\
-0.0124817490577698	0.0246016979217529\\
-0.0127569437026978	0.0234607458114624\\
-0.0117655992507935	0.0164754390716553\\
-0.0114297866821289	0.0136699676513672\\
-0.0104579329490662	0.00952082872390747\\
-0.0104820728302002	0.00857377052307129\\
-0.0102022290229797	0.00789803266525269\\
-0.0108544230461121	0.00849276781082153\\
-0.0105845332145691	0.0078997015953064\\
-0.0106197595596313	0.00853526592254639\\
-0.010559618473053	0.00785559415817261\\
-0.0121079087257385	0.0108206272125244\\
-0.0122030377388	0.0102043151855469\\
-0.0127913355827332	0.0115228295326233\\
-0.0142612457275391	0.0133112668991089\\
-0.0149286985397339	0.0142862200737\\
-0.0150268077850342	0.0136690735816956\\
-0.0154850482940674	0.0130344033241272\\
-0.0162129402160645	0.0138693451881409\\
-0.0168038010597229	0.0132546424865723\\
-0.0167285799980164	0.0120757818222046\\
-0.0179954767227173	0.0110417604446411\\
-0.0188325643539429	0.0113244652748108\\
-0.0188717842102051	0.0105481147766113\\
-0.0194359421730042	0.0100135803222656\\
-0.0202056765556335	0.0104926824569702\\
-0.0194905996322632	0.00513380765914917\\
-0.0195968151092529	0.00269526243209839\\
-0.0195448994636536	-0.00160896778106689\\
-0.0172615051269531	-0.00991547107696533\\
-0.0170415043830872	-0.0105690360069275\\
-0.0142489075660706	-0.0158721804618835\\
-0.0136649608612061	-0.0154827833175659\\
-0.0132708549499512	-0.0143536329269409\\
-0.012062132358551	-0.0158559083938599\\
-0.0122026801109314	-0.0151822566986084\\
-0.0115958452224731	-0.0159356594085693\\
-0.0109958052635193	-0.0161673426628113\\
-0.00921857357025146	-0.0187959671020508\\
-0.00932425260543823	-0.0180933475494385\\
-0.00870734453201294	-0.0188982486724854\\
-0.00825172662734985	-0.0177714228630066\\
-0.00816881656646729	-0.0170062184333801\\
-0.00766986608505249	-0.0164629817008972\\
-0.00869965553283691	-0.0117195248603821\\
-0.00849223136901855	-0.0107996463775635\\
-0.00770270824432373	-0.0120015144348145\\
-0.00816881656646729	-0.0101170539855957\\
-0.007671058177948	-0.00889003276824951\\
-0.00758534669876099	-0.00783038139343262\\
-0.0066801905632019	-0.00935357809066772\\
-0.00722455978393555	-0.00742602348327637\\
-0.00784409046173096	-0.00554788112640381\\
-0.00797349214553833	-0.0048370361328125\\
-0.00736510753631592	-0.00579833984375\\
-0.00743067264556885	-0.00507456064224243\\
-0.00714266300201416	-0.00585830211639404\\
-0.00660926103591919	-0.00639629364013672\\
-0.006428062915802	-0.00705069303512573\\
-0.00604206323623657	-0.00747519731521606\\
-0.00594329833984375	-0.00741451978683472\\
};
\addplot [color=blue, draw=none, mark size=3.0pt, mark=asterisk, mark options={solid, blue}, forget plot]
  table[row sep=crcr]{%
-0.0550614967942238	0.0772539153695107\\
};
\addplot [color=red, line width=0.7pt, forget plot]
  table[row sep=crcr]{%
0.0438827872276306	0.0172549486160278\\
0.0165939927101135	0.101051032543182\\
0.00716161727905273	0.126268267631531\\
0.00478070974349976	0.132429838180542\\
0.00438129901885986	0.133474051952362\\
-0.0217404365539551	0.196181654930115\\
-0.0515594482421875	0.257540285587311\\
-0.0799539089202881	0.305931687355042\\
-0.0980018377304077	0.332721352577209\\
-0.122750103473663	0.364955723285675\\
-0.135898292064667	0.38080221414566\\
-0.148863434791565	0.396087467670441\\
-0.156083583831787	0.403575479984283\\
-0.189121186733246	0.432602047920227\\
-0.209821999073029	0.448207676410675\\
-0.216416597366333	0.452251195907593\\
-0.2363520860672	0.462410807609558\\
-0.241932511329651	0.463142275810242\\
-0.249022305011749	0.465927839279175\\
-0.254818081855774	0.466607689857483\\
-0.261739194393158	0.468418180942535\\
-0.267193675041199	0.467950165271759\\
-0.279845416545868	0.468983173370361\\
-0.291915237903595	0.467553555965424\\
-0.297618806362152	0.465555131435394\\
-0.308677494525909	0.459866225719452\\
-0.317880749702454	0.451260566711426\\
-0.333016335964203	0.438920199871063\\
-0.34266459941864	0.42912483215332\\
-0.346493482589722	0.421693086624146\\
-0.351454496383667	0.416345536708832\\
-0.364502489566803	0.397854864597321\\
-0.368831098079681	0.39067929983139\\
-0.372381150722504	0.382712185382843\\
-0.376532256603241	0.375236928462982\\
-0.38000363111496	0.367142796516418\\
-0.384321630001068	0.360386490821838\\
-0.387721002101898	0.352294623851776\\
-0.399511396884918	0.33092999458313\\
-0.40334278345108	0.324546635150909\\
-0.407887279987335	0.319141864776611\\
-0.41153472661972	0.312208771705627\\
-0.416026592254639	0.306648313999176\\
-0.42472916841507	0.294440448284149\\
-0.42811906337738	0.287095606327057\\
-0.432980239391327	0.281465590000153\\
-0.450168907642365	0.257122039794922\\
-0.458432614803314	0.243399024009705\\
-0.461684584617615	0.235162496566772\\
-0.465889871120453	0.227721631526947\\
-0.472415685653687	0.210580348968506\\
-0.478170692920685	0.192214727401733\\
-0.482931673526764	0.171876430511475\\
-0.484833598136902	0.160565257072449\\
-0.486125111579895	0.148750007152557\\
-0.488390922546387	0.138067364692688\\
-0.491817355155945	0.101619899272919\\
-0.493536710739136	0.0907454490661621\\
-0.498641669750214	0.0357920527458191\\
-0.499868869781494	0.0255754590034485\\
-0.501216053962708	-0.00812411308288574\\
-0.501846611499786	-0.0311551094055176\\
-0.501574873924255	-0.0552964806556702\\
-0.49956488609314	-0.0817578434944153\\
-0.497192621231079	-0.108438611030579\\
-0.490795314311981	-0.152664244174957\\
-0.482255816459656	-0.197766721248627\\
-0.475609958171844	-0.228318572044373\\
-0.467818796634674	-0.260105431079865\\
-0.459343254566193	-0.291566550731659\\
-0.450017392635345	-0.322880744934082\\
-0.436755657196045	-0.363787710666656\\
-0.42741322517395	-0.389778733253479\\
-0.408704042434692	-0.437908113002777\\
-0.394204556941986	-0.470243990421295\\
-0.379491686820984	-0.499607920646667\\
-0.369333326816559	-0.517691731452942\\
-0.349402546882629	-0.547396123409271\\
-0.333724141120911	-0.567685663700104\\
-0.328493416309357	-0.574105501174927\\
-0.31805157661438	-0.584595441818237\\
-0.30812668800354	-0.592563331127167\\
-0.30240786075592	-0.597355961799622\\
-0.298301637172699	-0.599623382091522\\
-0.292446315288544	-0.604932308197021\\
-0.275019943714142	-0.618559300899506\\
-0.239329278469086	-0.639618933200836\\
-0.21636962890625	-0.648033022880554\\
-0.205013811588287	-0.650997579097748\\
-0.199191451072693	-0.652680933475494\\
-0.175774872303009	-0.655792772769928\\
-0.16973739862442	-0.655675649642944\\
-0.164025485515594	-0.656261503696442\\
-0.146534442901611	-0.65522700548172\\
-0.141506195068359	-0.65416419506073\\
-0.135680794715881	-0.653786182403564\\
-0.124090135097504	-0.652161359786987\\
-0.118561565876007	-0.649896442890167\\
-0.101411283016205	-0.645260870456696\\
-0.0898003578186035	-0.640360176563263\\
-0.0837520956993103	-0.638803124427795\\
-0.0780026316642761	-0.635994732379913\\
-0.072634756565094	-0.632710456848145\\
-0.0667080879211426	-0.629918456077576\\
-0.060899555683136	-0.628157675266266\\
-0.0541045665740967	-0.625465631484985\\
-0.0409789085388184	-0.622503399848938\\
-0.0276447534561157	-0.619867265224457\\
-0.0137824416160584	-0.616302907466888\\
-0.00960773229598999	-0.609790742397308\\
-0.00192493200302124	-0.609605491161346\\
0.0123902559280396	-0.607495903968811\\
0.0268711447715759	-0.60304582118988\\
0.0490979552268982	-0.594985365867615\\
0.0697512030601501	-0.583909451961517\\
0.0768380761146545	-0.578841626644135\\
0.0840728878974915	-0.57468044757843\\
0.0904439091682434	-0.569140732288361\\
0.102696239948273	-0.555514276027679\\
0.119508802890778	-0.532125890254974\\
0.141766667366028	-0.497403204441071\\
0.14684009552002	-0.488373279571533\\
0.15116024017334	-0.478305339813232\\
0.156252562999725	-0.469519078731537\\
0.164228558540344	-0.449823975563049\\
0.167145967483521	-0.438825786113739\\
0.177506804466248	-0.406689345836639\\
0.184122800827026	-0.384291052818298\\
0.187644600868225	-0.374631524085999\\
0.195884823799133	-0.3564293384552\\
0.204538524150848	-0.34052973985672\\
0.211998462677002	-0.322847664356232\\
0.21680098772049	-0.315971255302429\\
0.220895707607269	-0.308304667472839\\
0.225441217422485	-0.301048517227173\\
0.228360712528229	-0.290509581565857\\
0.231664538383484	-0.281229615211487\\
0.236124515533447	-0.273779571056366\\
0.239094734191895	-0.263984501361847\\
0.246162056922913	-0.246869564056396\\
0.26236230134964	-0.198752403259277\\
0.26494574546814	-0.187591850757599\\
0.273589253425598	-0.158144533634186\\
0.275747299194336	-0.14690774679184\\
0.280496656894684	-0.12717467546463\\
0.282475352287292	-0.11683863401413\\
0.283607244491577	-0.103649497032166\\
0.285367429256439	-0.0936868190765381\\
0.287047505378723	-0.0689293146133423\\
0.289065420627594	-0.0468451976776123\\
0.290461540222168	-0.0244583487510681\\
0.292213320732117	0.029949426651001\\
0.29177987575531	0.0797027349472046\\
0.290014624595642	0.109560012817383\\
0.284176230430603	0.156437754631042\\
0.274067163467407	0.205878674983978\\
0.272119581699371	0.211174726486206\\
0.267978668212891	0.227051138877869\\
0.261752605438232	0.246772050857544\\
0.25927859544754	0.252212166786194\\
0.256765782833099	0.259453177452087\\
0.254441976547241	0.263992428779602\\
0.252480685710907	0.270014107227325\\
0.24743527173996	0.279360115528107\\
0.244888842105865	0.284416258335114\\
0.240099847316742	0.291674554347992\\
0.238022744655609	0.296009063720703\\
0.230408489704132	0.304858922958374\\
0.225119113922119	0.311834335327148\\
0.219096481800079	0.317688643932343\\
0.215526938438416	0.324034690856934\\
0.212537288665771	0.32621568441391\\
0.203498601913452	0.33585911989212\\
0.197595834732056	0.339925825595856\\
0.188874304294586	0.345736443996429\\
0.185975968837738	0.34550929069519\\
0.183489680290222	0.347129821777344\\
0.180550158023834	0.347190797328949\\
0.174564123153687	0.350351750850677\\
0.171465039253235	0.351797580718994\\
0.166246116161346	0.352518141269684\\
0.163489639759064	0.352958142757416\\
0.160964548587799	0.352743983268738\\
0.158060073852539	0.351799130439758\\
0.155359089374542	0.353329300880432\\
0.149889886379242	0.351890325546265\\
0.146766841411591	0.352017641067505\\
0.144183576107025	0.353133261203766\\
0.141293704509735	0.352161109447479\\
0.138917148113251	0.352760434150696\\
0.135972857475281	0.351689517498016\\
0.133604288101196	0.352052867412567\\
0.130908489227295	0.3517085313797\\
0.128693819046021	0.350551307201386\\
0.125535130500793	0.349876344203949\\
0.119909524917603	0.350237786769867\\
0.108659565448761	0.348316490650177\\
0.105984091758728	0.34856241941452\\
0.0996584892272949	0.347472906112671\\
0.0973824858665466	0.346582055091858\\
0.0911814570426941	0.348770081996918\\
0.0842694044113159	0.34881854057312\\
0.0812093019485474	0.348012685775757\\
0.0783490538597107	0.348442673683167\\
0.0718111991882324	0.346416890621185\\
0.0651626586914063	0.346744894981384\\
0.0618290305137634	0.345899224281311\\
0.0555377006530762	0.343192756175995\\
0.0529276132583618	0.342490375041962\\
0.0468761324882507	0.338435590267181\\
0.0439012050628662	0.336216509342194\\
0.0391405820846558	0.331159651279449\\
0.0276763439178467	0.323064148426056\\
0.0243251323699951	0.32170045375824\\
0.0217180252075195	0.318775653839111\\
0.0200269818305969	0.316009402275085\\
0.0126753449440002	0.306525230407715\\
0.0111209750175476	0.302840530872345\\
0.00986337661743164	0.29821789264679\\
0.00777751207351685	0.296199500560761\\
0.00429677963256836	0.289713799953461\\
0.000148892402648926	0.283240377902985\\
-0.00303536653518677	0.28125274181366\\
-0.00497013330459595	0.277860224246979\\
-0.00788420438766479	0.275625705718994\\
-0.0121830701828003	0.269448757171631\\
-0.0161272287368774	0.261430144309998\\
-0.018473744392395	0.257760643959045\\
-0.0197583436965942	0.253946185112\\
-0.034576416015625	0.227409482002258\\
-0.0377079844474792	0.218390703201294\\
-0.0412202477455139	0.209066212177277\\
-0.0544252991676331	0.148776590824127\\
-0.0570312142372131	0.140879034996033\\
-0.0574836134910584	0.135520696640015\\
-0.0591336488723755	0.125842750072479\\
-0.061836838722229	0.106443285942078\\
-0.0619170665740967	0.100637257099152\\
-0.0631359815597534	0.0918216109275818\\
-0.0628911256790161	0.0805474519729614\\
-0.0632287263870239	0.0714834928512573\\
-0.0627824664115906	0.0478593111038208\\
-0.0629562139511108	0.0397453308105469\\
-0.0597187280654907	0.00271546840667725\\
-0.0563098192214966	-0.0220687389373779\\
-0.0521929264068604	-0.0408896207809448\\
-0.0453609228134155	-0.061507523059845\\
-0.0337279438972473	-0.0876003503799438\\
-0.0210835337638855	-0.105912268161774\\
-0.009449303150177	-0.117344558238983\\
-0.00041353702545166	-0.12063866853714\\
0.000618040561676025	-0.120107114315033\\
0.00302660465240479	-0.120737552642822\\
0.00685995817184448	-0.120333731174469\\
0.00851505994796753	-0.120940148830414\\
0.00999873876571655	-0.120459139347076\\
0.0117058157920837	-0.121528506278992\\
0.0132834315299988	-0.120988726615906\\
0.0163581967353821	-0.12117338180542\\
0.0178297162055969	-0.120452344417572\\
0.0194555521011353	-0.12031102180481\\
0.0206473469734192	-0.118517577648163\\
0.0223628282546997	-0.118384063243866\\
0.0261257886886597	-0.1149982213974\\
0.0278939008712769	-0.114962637424469\\
0.0289706587791443	-0.113073527812958\\
0.0304802656173706	-0.112409710884094\\
0.0313948392868042	-0.110200226306915\\
0.0326723456382751	-0.109075546264648\\
0.036689281463623	-0.100546538829803\\
0.0406491160392761	-0.0917348861694336\\
0.0419303178787231	-0.0859581232070923\\
0.0511994957923889	-0.0521543622016907\\
0.0524927377700806	-0.0357394218444824\\
0.0520519614219666	-0.0279863476753235\\
0.0507994294166565	-0.00664794445037842\\
0.0467005968093872	0.0253812670707703\\
0.0408130288124084	0.0474424362182617\\
0.0345221757888794	0.0632233619689941\\
0.0290365815162659	0.0722439885139465\\
0.0257781147956848	0.0763633847236633\\
0.0199489593505859	0.0835655331611633\\
0.0156401395797729	0.0864161849021912\\
0.0125711560249329	0.0886729955673218\\
0.0118978023529053	0.0884025692939758\\
0.00997436046600342	0.0889695882797241\\
0.00945150852203369	0.0884967446327209\\
0.00571668148040771	0.0893799662590027\\
0.00445914268493652	0.0875262022018433\\
0.00349372625350952	0.0877655744552612\\
0.00316983461380005	0.086323082447052\\
0.00190085172653198	0.0868982672691345\\
0.000155806541442871	0.0858142375946045\\
-0.00132757425308228	0.0821645855903625\\
-0.00209826231002808	0.0814304351806641\\
-0.0031430721282959	0.0770972967147827\\
-0.00477522611618042	0.0726977586746216\\
-0.00495004653930664	0.0712376832962036\\
-0.00494992733001709	0.0675068497657776\\
-0.00723069906234741	0.0623626112937927\\
-0.00819540023803711	0.0607503056526184\\
-0.00974446535110474	0.0602678656578064\\
-0.0101709365844727	0.0552006959915161\\
-0.0122887492179871	0.039438784122467\\
-0.013456404209137	0.0332135558128357\\
-0.0132757425308228	0.0290068387985229\\
-0.0138134956359863	0.0230755805969238\\
-0.0129804611206055	0.0175190567970276\\
-0.0136474967002869	0.0151904821395874\\
-0.0129088759422302	0.0123206973075867\\
-0.0128718614578247	0.00978827476501465\\
-0.010972797870636	0.00216805934906006\\
-0.0102163553237915	-0.000199556350708008\\
-0.009665846824646	-0.00101238489151001\\
-0.00988560914993286	-0.000222861766815186\\
-0.00888389348983765	-0.00235670804977417\\
-0.00715458393096924	-0.00500661134719849\\
-0.00599640607833862	-0.00775504112243652\\
-0.0059126615524292	-0.00709468126296997\\
-0.00433099269866943	-0.0101618766784668\\
-0.00329047441482544	-0.0118997097015381\\
-0.00345206260681152	-0.0109632015228271\\
-0.00165677070617676	-0.0146652460098267\\
0.00334656238555908	-0.0195950865745544\\
0.00524330139160156	-0.0190824866294861\\
0.00616484880447388	-0.0192688703536987\\
0.00651454925537109	-0.0187494158744812\\
0.00661319494247437	-0.0155748724937439\\
0.00683218240737915	-0.0095982551574707\\
0.00655537843704224	-0.00571686029434204\\
0.00635850429534912	-0.00100117921829224\\
0.00657665729522705	-0.001545250415802\\
0.00657367706298828	-0.000587940216064453\\
0.00701755285263062	-0.00118529796600342\\
0.00757336616516113	-0.000580728054046631\\
0.00863975286483765	-0.0016745924949646\\
0.00834888219833374	-0.000721871852874756\\
0.00923478603363037	-0.000128686428070068\\
0.00877916812896729	0.00204223394393921\\
0.0093766450881958	0.0031164288520813\\
0.00844508409500122	0.00804108381271362\\
0.00452816486358643	0.0201197862625122\\
0.00176548957824707	0.0253611207008362\\
0.00192207098007202	0.0246936082839966\\
0.00131964683532715	0.025668203830719\\
0.00142180919647217	0.0246961712837219\\
0.00089108943939209	0.0253532528877258\\
-0.00103026628494263	0.0275965929031372\\
-0.0021016001701355	0.0286498665809631\\
-0.00230413675308228	0.027134120464325\\
-0.00295537710189819	0.0289329886436462\\
-0.00243991613388062	0.0277490019798279\\
-0.00181668996810913	0.027280867099762\\
-0.00119352340698242	0.0276530981063843\\
-0.00108122825622559	0.0277163982391357\\
};
\addplot [color=red, draw=none, mark size=3.0pt, mark=asterisk, mark options={solid, red}, forget plot]
  table[row sep=crcr]{%
0.0438827723264694	0.0172549597918987\\
};
\addplot [color=green, dashed, line width=0.7pt, forget plot]
  table[row sep=crcr]{%
-0.0252749621868134	0.0959042310714722\\
-0.0256432294845581	0.0968656241893768\\
-0.0233517289161682	0.0908872485160828\\
-0.0327918231487274	0.115267783403397\\
-0.0431816875934601	0.139828324317932\\
-0.0541383922100067	0.163332909345627\\
-0.0827577114105225	0.222252905368805\\
-0.0887050926685333	0.233150333166122\\
-0.0982786417007446	0.245844632387161\\
-0.107700139284134	0.255721151828766\\
-0.120416104793549	0.263963043689728\\
-0.124563485383987	0.26672175526619\\
-0.137028425931931	0.26905170083046\\
-0.140411347150803	0.268498182296753\\
-0.144401967525482	0.267211198806763\\
-0.14802822470665	0.266785800457001\\
-0.159061670303345	0.261287331581116\\
-0.165549069643021	0.255651235580444\\
-0.169166773557663	0.253669172525406\\
-0.172490417957306	0.250571131706238\\
-0.178332209587097	0.24208790063858\\
-0.180636435747147	0.236162096261978\\
-0.183216094970703	0.23114949464798\\
-0.189082890748978	0.22361233830452\\
-0.198192417621613	0.201144158840179\\
-0.202357113361359	0.189315319061279\\
-0.208141088485718	0.171418637037277\\
-0.215334862470627	0.147348344326019\\
-0.217126846313477	0.141405016183853\\
-0.220212697982788	0.128936290740967\\
-0.222082883119583	0.123489499092102\\
-0.223118394613266	0.117295980453491\\
-0.223837912082672	0.108701080083847\\
-0.227487504482269	0.0908943116664886\\
-0.228932529687881	0.076238214969635\\
-0.230488955974579	0.0640577375888824\\
-0.232695966959	0.0450815856456757\\
-0.236417353153229	0.011015772819519\\
-0.237411469221115	-0.0126366913318634\\
-0.23916095495224	-0.03718301653862\\
-0.240239173173904	-0.059842437505722\\
-0.240529149770737	-0.0726923048496246\\
-0.240684270858765	-0.0892626643180847\\
-0.239656567573547	-0.109096854925156\\
-0.237893313169479	-0.122491210699081\\
-0.236780017614365	-0.132896721363068\\
-0.229980438947678	-0.160500288009644\\
-0.224304497241974	-0.177370250225067\\
-0.218046933412552	-0.190053969621658\\
-0.213764846324921	-0.195905119180679\\
-0.211189836263657	-0.199841380119324\\
-0.208288669586182	-0.202750772237778\\
-0.204568266868591	-0.205664813518524\\
-0.203747779130936	-0.20601424574852\\
-0.195591509342194	-0.214320749044418\\
-0.192916750907898	-0.216092735528946\\
-0.18905770778656	-0.216615796089172\\
-0.186431169509888	-0.218980371952057\\
-0.183742940425873	-0.220634609460831\\
-0.18078139424324	-0.222364068031311\\
-0.179545998573303	-0.222429633140564\\
-0.176442533731461	-0.224093526601791\\
-0.172939360141754	-0.223701149225235\\
-0.169283390045166	-0.223678588867188\\
-0.166633397340775	-0.223694294691086\\
-0.165870696306229	-0.222895920276642\\
-0.161990940570831	-0.222635924816132\\
-0.161192357540131	-0.221736401319504\\
-0.16013777256012	-0.221418082714081\\
-0.158105552196503	-0.22292223572731\\
-0.157157629728317	-0.222283005714417\\
-0.155453085899353	-0.222608089447021\\
-0.151392966508865	-0.225333034992218\\
-0.149879157543182	-0.225728541612625\\
-0.148342937231064	-0.225511580705643\\
-0.143345713615417	-0.226698100566864\\
-0.141986578702927	-0.226154148578644\\
-0.132151931524277	-0.226907700300217\\
-0.130725294351578	-0.226704746484756\\
-0.129857391119003	-0.224800556898117\\
-0.128453969955444	-0.224039018154144\\
-0.127353221178055	-0.222785860300064\\
-0.12589231133461	-0.222698122262955\\
-0.125668078660965	-0.219554722309113\\
-0.124541133642197	-0.218812674283981\\
-0.121419727802277	-0.21860283613205\\
-0.120494574308395	-0.216825604438782\\
-0.116815537214279	-0.217881113290787\\
-0.115341454744339	-0.217566221952438\\
-0.109679490327835	-0.218930810689926\\
-0.108601301908493	-0.217472612857819\\
-0.106830537319183	-0.217551440000534\\
-0.103773772716522	-0.216406881809235\\
-0.0999604761600494	-0.217400103807449\\
-0.0982294976711273	-0.217087715864182\\
-0.0929572880268097	-0.21228700876236\\
-0.0875526964664459	-0.211958140134811\\
-0.0853914618492126	-0.212285101413727\\
-0.0775750875473022	-0.209626168012619\\
-0.0753518044948578	-0.205336004495621\\
-0.0735234022140503	-0.204689919948578\\
-0.0723465979099274	-0.203223675489426\\
-0.071138858795166	-0.200688451528549\\
-0.0696220695972443	-0.199862092733383\\
-0.0642281472682953	-0.194337427616119\\
-0.0627254843711853	-0.193367809057236\\
-0.0618900656700134	-0.190248012542725\\
-0.0591293573379517	-0.187278658151627\\
-0.0563271939754486	-0.185076236724854\\
-0.0541341006755829	-0.181769847869873\\
-0.0527581870555878	-0.180952429771423\\
-0.0488949120044708	-0.176554322242737\\
-0.0459828972816467	-0.174395352602005\\
-0.042056679725647	-0.170523405075073\\
-0.0402338206768036	-0.166050642728806\\
-0.0380259454250336	-0.162775695323944\\
-0.0372433662414551	-0.160381525754929\\
-0.0360292792320251	-0.15959358215332\\
-0.0326491296291351	-0.154980182647705\\
-0.0287406444549561	-0.152532190084457\\
-0.0278485119342804	-0.150205135345459\\
-0.0266923606395721	-0.149224579334259\\
-0.0235322117805481	-0.143560081720352\\
-0.0211219489574432	-0.141410917043686\\
-0.0203540623188019	-0.139579027891159\\
-0.0189367532730103	-0.134899646043777\\
-0.0157927274703979	-0.129214316606522\\
-0.0132290422916412	-0.123217582702637\\
-0.0121423900127411	-0.121970891952515\\
-0.011476457118988	-0.11983197927475\\
-0.0103911757469177	-0.118409335613251\\
-0.010037362575531	-0.115437507629395\\
-0.00754892826080322	-0.110761135816574\\
-0.00719565153121948	-0.107920497655869\\
-0.00486129522323608	-0.102869421243668\\
-0.00487321615219116	-0.0993093252182007\\
-0.00303980708122253	-0.0943224430084229\\
-0.00273942947387695	-0.091435968875885\\
-0.0020081102848053	-0.0902665853500366\\
0.000732064247131348	-0.0677751004695892\\
0.000860214233398438	-0.0605933666229248\\
0.00129988789558411	-0.0595673322677612\\
0.00111934542655945	-0.0359571278095245\\
0.000436633825302124	-0.0293079316616058\\
-0.00133597850799561	-0.0156858563423157\\
-0.00359970331192017	-0.00608327984809875\\
-0.00736579298973084	0.00617805123329163\\
-0.0114367306232452	0.0158825516700745\\
-0.0149502456188202	0.0219337642192841\\
-0.021232545375824	0.02994304895401\\
-0.0232946574687958	0.0316860377788544\\
-0.024800032377243	0.0322922170162201\\
-0.0268065631389618	0.0326286256313324\\
-0.027446061372757	0.0324620306491852\\
-0.0295473039150238	0.0342782735824585\\
-0.0306127071380615	0.03366419672966\\
-0.0315462052822113	0.0342504680156708\\
-0.0334846973419189	0.033557802438736\\
-0.038903146982193	0.0340533256530762\\
-0.0394730269908905	0.0334267020225525\\
-0.0422928631305695	0.0341684520244598\\
-0.0432515144348145	0.0322068333625793\\
-0.0448437035083771	0.0303352475166321\\
-0.0459261834621429	0.0271944403648376\\
-0.0466421544551849	0.0272558033466339\\
-0.052191436290741	0.0197649896144867\\
-0.0529299676418304	0.0201377272605896\\
-0.0535790622234344	0.0181850790977478\\
-0.0555709600448608	0.0146716237068176\\
-0.0567167997360229	0.0142281949520111\\
-0.0582520663738251	0.00935381650924683\\
-0.0591748058795929	0.0084560215473175\\
-0.0604304969310761	0.0064014196395874\\
-0.0608710050582886	0.00549417734146118\\
-0.0613577365875244	-0.00140511989593506\\
-0.0621914267539978	-0.00428363680839539\\
-0.0618211627006531	-0.00949719548225403\\
-0.0621364414691925	-0.0133077502250671\\
-0.0614989101886749	-0.0190859436988831\\
-0.06100994348526	-0.0240635871887207\\
-0.0577291548252106	-0.0415316522121429\\
-0.0512172877788544	-0.0590235888957977\\
-0.0490229427814484	-0.0634561479091644\\
-0.0444687902927399	-0.0687378644943237\\
-0.0434749126434326	-0.0693071186542511\\
-0.0400444865226746	-0.0734896063804626\\
-0.0391161441802979	-0.0734719336032867\\
-0.0387474298477173	-0.0730183124542236\\
-0.0369710624217987	-0.0737230479717255\\
-0.0366018414497375	-0.0732417106628418\\
-0.0350330770015717	-0.0745055675506592\\
-0.0308903157711029	-0.0760054588317871\\
-0.0305180549621582	-0.0754730403423309\\
-0.0288499891757965	-0.0767669081687927\\
-0.0276805460453033	-0.0761289298534393\\
-0.0245954394340515	-0.0758697986602783\\
-0.0240867733955383	-0.0752646327018738\\
-0.023310124874115	-0.075627475976944\\
-0.0223366022109985	-0.0745700597763062\\
-0.0192514359951019	-0.0737514197826386\\
-0.0173726081848145	-0.0733910501003265\\
-0.0165033936500549	-0.0721228122711182\\
-0.0138840079307556	-0.0715481340885162\\
-0.0108488798141479	-0.0680423080921173\\
-0.00727164745330811	-0.0648883283138275\\
-0.00621607899665833	-0.0652487277984619\\
-0.00489908456802368	-0.0627803206443787\\
-0.00373148918151855	-0.0618320405483246\\
-0.0028526782989502	-0.0562230050563812\\
-0.00200161337852478	-0.0542768239974976\\
-0.000900536775588989	-0.0530783534049988\\
0.00010189414024353	-0.04227414727211\\
0.00138711929321289	-0.034913033246994\\
0.00180384516716003	-0.0344072282314301\\
0.00244107842445374	-0.0300138294696808\\
0.00359463691711426	-0.0283816754817963\\
0.00405287742614746	-0.0237205922603607\\
0.00536927580833435	-0.0204988718032837\\
};
\addplot [color=green, draw=none, mark size=3.0pt, mark=asterisk, mark options={solid, mycolor1}, forget plot]
  table[row sep=crcr]{%
-0.0252749547362328	0.0959042310714722\\
};
\addplot [color=red, dashed, line width=0.7pt, forget plot]
  table[row sep=crcr]{%
0.0451151728630066	0.0113265514373779\\
0.0437105298042297	0.0159660577774048\\
0.043515145778656	0.0166618824005127\\
0.0454691648483276	0.0109235644340515\\
0.0417500734329224	0.0237957835197449\\
0.0300401449203491	0.0598820447921753\\
0.0170882344245911	0.0963019728660584\\
0.00304645299911499	0.132011830806732\\
-0.0113565921783447	0.165016353130341\\
-0.0240353345870972	0.192504346370697\\
-0.0442904233932495	0.233897745609283\\
-0.0683179497718811	0.279126346111298\\
-0.103476881980896	0.342864155769348\\
-0.117551684379578	0.366465628147125\\
-0.140277624130249	0.404161870479584\\
-0.164911925792694	0.444975852966309\\
-0.200263321399689	0.50238174200058\\
-0.325631976127625	0.696973204612732\\
-0.347502648830414	0.727512121200562\\
-0.358664214611053	0.742230474948883\\
-0.382029354572296	0.770465731620789\\
-0.405983150005341	0.796564877033234\\
-0.448352634906769	0.840502738952637\\
-0.52352634906769	0.90502738952637\\
};
\addplot [color=blue, dashed, line width=0.7pt, forget plot]
  table[row sep=crcr]{%
-0.0569646954536438	0.0866867899894714\\
-0.0556769967079163	0.0827354788780212\\
-0.0553459525108337	0.0816285610198975\\
-0.0567765235900879	0.0859630107879639\\
-0.0692293047904968	0.122698605060577\\
-0.0826412439346313	0.158876955509186\\
-0.098339319229126	0.197953641414642\\
-0.107480406761169	0.218829035758972\\
-0.122956156730652	0.252221763134003\\
-0.140863120555878	0.28737199306488\\
-0.167234659194946	0.336349904537201\\
-0.212378025054932	0.414183974266052\\
-0.245918035507202	0.46967750787735\\
-0.28978818655014	0.539248943328857\\
-0.347558379173279	0.628650963306427\\
-0.397712707519531	0.704970419406891\\
-0.490067720413208	0.846145510673523\\
-0.507003962993622	0.873022139072418\\
-0.547003962993622	0.933022139072418\\
};

\addplot[->, thick, color=blue] coordinates
{ (-0.4728,0.4644) (-0.4979,0.4454)};
\addplot[->, thick, color=blue] coordinates
{ (-0.7493,-0.1970) (-0.7386,-0.3498)};
\addplot[->, thick, color=blue] coordinates
{ (0.3882,0.3149) (0.3782,0.3444)};

\addplot[->, thick, color=red] coordinates
{ (-0.4861,0.1487) (-0.4883,0.1380)};
\addplot[->, thick, color=red] coordinates
{ (-0.0608,-0.6281) (-0.0541,-0.6254)};

\end{axis}
\end{tikzpicture}%

%% file: curvature_straight.tex
%
%
\definecolor{mycolor1}{rgb}{0.50000,0.00000,0.50000}%
\begin{tikzpicture}

\begin{axis}[%
width=\figurewidth,
height=\figureheight,
at={(0\figurewidth,0\figureheight)},
scale only axis,
xmin=0,
xmax=120,
xtick={0,30,60,90,120},
ylabel style={font=\color{white!15!black},at={(axis description cs:-0.2,.5)},anchor=south},
xlabel={$t$ [s]},
ymin=-0.27,
ymax=0.27,
ytick = {-0.181,0,0.181},
yticklabels={$-u_{\text{max}}$,0,$u_{\text{max}}$},
xlabel style={font=\color{white!15!black},at={(axis description cs:0.5,-0.17)},anchor=north},
ylabel={$u$ [m$^{-1}$]},
axis background/.style={fill=white},
xmajorgrids,
ymajorgrids
]
\addplot [color=blue, line width=0.7pt, forget plot]
table[row sep=crcr]{%
	0	-0.181220059660632\\
	2.09999999999999	-0.181234292889002\\
	2.40000000000001	-0.174258041457904\\
	2.69999999999999	-0.166057085728255\\
	3	-0.156712907878671\\
	3.30000000000001	-0.145269251210635\\
	3.59999999999999	-0.132780341489365\\
	3.90000000000001	-0.119014957372457\\
	4.19999999999999	-0.103193225230058\\
	4.5	-0.0863707148991466\\
	4.80000000000001	-0.0701155800699382\\
	5.09999999999999	-0.0535022130779339\\
	5.40000000000001	-0.0373002127201971\\
	5.69999999999999	-0.0213538048274984\\
	6.30000000000001	0.0144335339923316\\
	6.59999999999999	0.0333982478653923\\
	6.90000000000001	0.0510071880330543\\
	7.19999999999999	0.0689724076461857\\
	7.5	0.083810850744328\\
	7.80000000000001	0.100863393491011\\
	8.09999999999999	0.121101867932282\\
	8.40000000000001	0.142048961463189\\
	8.69999999999999	0.159270525218318\\
	9	0.176870494655134\\
	9.30000000000001	0.178210152065333\\
	9.59999999999999	0.168682309871826\\
	9.90000000000001	0.14375158008221\\
	10.2	0.134550490420366\\
	10.5	0.124512646296722\\
	10.8	0.122546887236297\\
	11.1	0.123846700577161\\
	11.4	0.129734554481104\\
	11.7	0.139582836172764\\
	12	0.145862653464349\\
	12.3	0.1512496069235\\
	12.6	0.153631739563593\\
	12.9	0.154970822612285\\
	13.2	0.150333567555151\\
	13.5	0.146825040889041\\
	13.8	0.15090471308946\\
	14.1	0.145686724564399\\
	14.4	0.142681931938824\\
	14.7	0.136894077683962\\
	15	0.146865870431697\\
	15.3	0.156332008767691\\
	15.6	0.159606250492402\\
	15.9	0.13627972235173\\
	16.2	0.115220198282429\\
	16.5	0.0811346709751319\\
	16.8	0.0585406784385611\\
	17.1	0.036772574186358\\
	17.4	0.0279327365109339\\
	17.7	0.0124265103652021\\
	18	0.00505859976871648\\
	18.3	-0.00676177746777284\\
	18.6	-0.0152126108594643\\
	18.9	-0.0261706712194609\\
	19.2	-0.0212780301542068\\
	19.5	-0.00602096798658636\\
	19.8	0.0158076326026446\\
	20.1	0.0385371608649052\\
	20.4	0.0570181848016773\\
	20.7	0.052632919756519\\
	21	0.027917295664821\\
	21.3	-0.00503670242508747\\
	21.6	-0.0335686995941273\\
	21.9	-0.0502076876220201\\
	22.2	-0.0636672738703226\\
	22.5	-0.0646471583976052\\
	22.8	-0.0673856085393538\\
	23.1	-0.0688004788346177\\
	23.4	-0.0710298737698736\\
	23.7	-0.0846739630686955\\
	24	-0.106523646803765\\
	24.6	-0.101656098669508\\
	24.9	-0.092633460007562\\
	25.2	-0.102988402397699\\
	25.5	-0.115159174899276\\
	25.8	-0.131014656411537\\
	26.1	-0.15761294908836\\
	26.4	-0.174777000005179\\
	26.7	-0.176679620712093\\
	27	-0.161853236472581\\
	27.3	-0.137217117396887\\
	27.6	-0.12894169478389\\
	27.9	-0.121313207748187\\
	28.2	-0.112893995220873\\
	28.5	-0.114170781433018\\
	28.8	-0.115754585744327\\
	29.1	-0.11893207260951\\
	29.4	-0.12608129918496\\
	29.7	-0.128470281945084\\
	30	-0.134760211870542\\
	30.3	-0.134360317236485\\
	30.6	-0.141321315340093\\
	30.9	-0.135289296316529\\
	31.2	-0.127751668931836\\
	31.5	-0.113832859444216\\
	31.8	-0.110675326004639\\
	32.1	-0.104513190729392\\
	32.4	-0.111127283739535\\
	32.7	-0.103974622286501\\
	33	-0.104733149323692\\
	33.3	-0.0999863996634076\\
	33.6	-0.0994808947463639\\
	33.9	-0.103392081011947\\
	34.2	-0.103688264777844\\
	34.5	-0.11153827469073\\
	34.8	-0.108297763662932\\
	35.1	-0.107756091644887\\
	35.4	-0.0970486448626104\\
	35.7	-0.097500230713365\\
	36	-0.0916060761453537\\
	36.3	-0.0881972546313534\\
	36.6	-0.0859183556429457\\
	36.9	-0.0957756602666962\\
	37.2	-0.103470424238651\\
	37.5	-0.112640941300839\\
	37.8	-0.12713836248102\\
	38.1	-0.131317747833549\\
	38.4	-0.134429913477362\\
	38.7	-0.131689220648781\\
	39	-0.123043269898744\\
	39.3	-0.118543238632213\\
	39.6	-0.117312343463595\\
	39.9	-0.115190129541674\\
	40.2	-0.118919843807419\\
	40.5	-0.131570575111908\\
	40.8	-0.148533561508827\\
	41.1	-0.16296294632059\\
	41.4	-0.181377977061459\\
	44.1	-0.181332841258183\\
	44.4	-0.168949306153678\\
	44.7	-0.151224676748996\\
	45	-0.131823288105238\\
	45.3	-0.11195097497415\\
	45.6	-0.0908548527969799\\
	45.9	-0.071412962688072\\
	46.2	-0.0536826307029514\\
	46.5	-0.0340267554055629\\
	46.8	-0.0158655392243077\\
	47.1	0.00314879706277793\\
	47.4	0.0226457346898883\\
	47.7	0.0436944641435559\\
	48	0.0644690610724581\\
	48.3	0.0854340682907662\\
	48.6	0.10120778964108\\
	48.9	0.115280687237401\\
	49.2	0.122530689393386\\
	49.5	0.130905410780514\\
	49.8	0.136049126355971\\
	50.1	0.133306954361615\\
	50.4	0.124573604013676\\
	50.7	0.114477770216212\\
	51	0.105896712757158\\
	51.3	0.103285030896302\\
	51.6	0.102235617208891\\
	51.9	0.104110887161397\\
	52.2	0.108876913869324\\
	52.5	0.115010935882452\\
	53.1	0.124270833331707\\
	53.4	0.124013271378345\\
	53.7	0.119731248171263\\
	54	0.114512438436805\\
	54.3	0.107986247741025\\
	54.6	0.111655338047314\\
	54.9	0.121007683866225\\
	55.5	0.133910261332375\\
	55.8	0.143426149958032\\
	56.1	0.144713796545972\\
	56.4	0.140646719992077\\
	56.7	0.134840558463594\\
	57	0.127729895329111\\
	57.3	0.114929900327581\\
	57.6	0.102838326657206\\
	57.9	0.0929261596313324\\
	58.2	0.0857131769654416\\
	58.5	0.0786778685995557\\
	58.8	0.0725715447114794\\
	59.1	0.0628513736325544\\
	59.4	0.0543391818769123\\
	59.7	0.0467581349716966\\
	60	0.0424721462195521\\
	60.3	0.0385009050605447\\
	60.6	0.0329699716056382\\
	60.9	0.0263084958120317\\
	61.2	0.0219458519483453\\
	61.5	0.013557928901065\\
	61.8	0.00354194318629197\\
	62.1	-0.00525204314516259\\
	62.4	-0.0164890344165372\\
	62.7	-0.0222315748822552\\
	63	-0.0269218012226986\\
	63.3	-0.0307044562474346\\
	63.6	-0.0349731639457502\\
	63.9	-0.0383850816605502\\
	64.2	-0.0430650788067624\\
	64.5	-0.0468496235950795\\
	64.8	-0.0484334038752081\\
	65.1	-0.0526166550246217\\
	65.4	-0.0538768157156539\\
	65.7	-0.0571307238426471\\
	66	-0.0584243646576397\\
	66.3	-0.0619574584725626\\
	66.6	-0.0634185352377301\\
	66.9	-0.0612535521917437\\
	67.2	-0.0598668865083027\\
	67.5	-0.0592990543518681\\
	67.8	-0.0553901469440916\\
	68.1	-0.0513134769746557\\
	68.4	-0.0467606988749765\\
	68.7	-0.0466711178635535\\
	69	-0.0432362938740027\\
	69.3	-0.0431092560543505\\
	69.6	-0.041139033437986\\
	69.9	-0.0396193831502103\\
	70.2	-0.0400022811789995\\
	70.8	-0.03787930053889\\
	71.4	-0.0354140031239751\\
	71.7	-0.0322339176869946\\
	72	-0.0314839495047465\\
	72.3	-0.0297623451206164\\
	72.6	-0.0272623359147133\\
	72.9	-0.0244449456177449\\
	73.2	-0.0194842153237573\\
	73.5	-0.0149686431221028\\
	73.8	-0.0111411314460668\\
	74.1	-0.00832399489291902\\
	74.4	-0.010028293546128\\
	74.7	-0.00759120923217438\\
	75.3	-0.000771152297943445\\
	75.6	0.00156885912772964\\
	76.2	0.00948164759992665\\
	76.5	0.0127523851798514\\
	76.8	0.0136624200365532\\
	77.1	0.0110718309806259\\
	77.4	0.00962776915054064\\
	77.7	0.0110963962390542\\
	78	0.0161206201451023\\
	78.6	0.0270443270769931\\
	78.9	0.0300419447813169\\
	79.5	0.0310241961154816\\
	80.1	0.0353713728044909\\
	80.4	0.0389204287321263\\
	80.7	0.0407409331935469\\
	81.3	0.0346873088382438\\
	81.9	0.0322104322831933\\
	82.2	0.0275571962217498\\
	82.5	0.0271838080362556\\
	82.8	0.0270849126771395\\
	83.1	0.0247324978389258\\
	83.4	0.0262315096056795\\
	83.7	0.0252439632121479\\
	84	0.0248294365304957\\
	84.3	0.0242527350072237\\
	84.6	0.0202689083520795\\
	84.9	0.0176132110732397\\
	85.2	0.0157598656107041\\
	85.5	0.0135070127026609\\
	85.8	0.0104208412716389\\
	86.1	0.00888807293367222\\
	86.4	0.0062276144033433\\
	86.7	0.00262431548122777\\
	87	0.00149788573713749\\
	87.6	-0.00253685512322477\\
	87.9	-0.00412590091792708\\
	88.2	-0.00468700860199078\\
	88.5	-0.00747068639043391\\
	88.8	-0.00727556333137613\\
	89.1	-0.007268296121822\\
	89.4	-0.00686801218495248\\
	89.7	-0.00626288529491603\\
	90	-0.00587184979136168\\
	90.3	-0.00829765409650918\\
	90.6	-0.00987276061442799\\
	90.9	-0.0107409279185049\\
	91.2	-0.0139709358814741\\
	91.5	-0.0166330061140343\\
	91.8	-0.0200892028789781\\
	92.1	-0.0222368101168229\\
	92.7	-0.0225145976849319\\
	93	-0.0221787969252034\\
	93.3	-0.019990737496471\\
	93.6	-0.0180648684105336\\
	93.9	-0.0185992232913748\\
	94.2	-0.016192769822851\\
	94.5	-0.0142722083887179\\
	94.8	-0.0134642160872147\\
	95.1	-0.0147180702624325\\
	95.4	-0.0127362282406978\\
	95.7	-0.0144678423188793\\
	96	-0.0120321968556141\\
	96.3	-0.00942686805240101\\
	96.6	-0.0100478762818454\\
	96.9	-0.00982651234676268\\
	97.2	-0.010669412978757\\
	97.5	-0.00617251035799882\\
	97.8	-0.0054934469194734\\
	98.1	-0.00538180536617006\\
	98.4	-0.00375281781163039\\
	98.7	-0.00428137198750278\\
	99	-0.00113536245817158\\
	99.3	-0.000280296055422014\\
	99.6	-4.18204100185449e-05\\
	99.9	0.00153254715249318\\
	100.2	0.00129920671736272\\
	100.5	0.0031325671744753\\
	100.8	0.00319157299713879\\
	101.1	0.00276742515430328\\
	101.4	0.00629756014731697\\
	101.7	0.00638259262950669\\
	102.3	0.00935633921699264\\
	102.6	0.0130263946841467\\
	103.2	0.0104693207392188\\
	103.5	0.0088723243109996\\
	103.8	0.00837654679440902\\
	104.1	0.0102136048367356\\
	104.4	0.0103948981245594\\
	104.7	0.00981179882722927\\
	105	0.0120565230944578\\
	105.3	0.0121132839025222\\
	105.9	0.0116183102753098\\
	106.2	0.00765608704995202\\
	106.5	0.00454589555397433\\
	106.8	0.00018181313811283\\
	107.1	-0.00246262662017216\\
	107.7	-0.00453235240203753\\
	108	-0.00385142805478722\\
	108.3	-0.00509116076847249\\
	108.6	-0.00606676228824199\\
	108.9	-0.00661966973470385\\
	109.2	-0.00294942301707124\\
	109.5	-0.00149431836618419\\
	109.8	-0.00255403430159618\\
	110.1	0.000817721984475384\\
	110.4	0.00186623737704394\\
	110.7	0.00245332262886677\\
	111	0.00196459478954125\\
	111.6	0.00192033000510605\\
	111.9	0.00252266078186381\\
	112.2	0.00162336950188546\\
	112.5	0.00342758839698831\\
	112.8	0.00448443871042059\\
	113.1	0.00607401258866958\\
	113.4	0.00650990236420057\\
	113.7	0.00872292869431135\\
	114	0.0059953784171114\\
	115.2	0.00476400650603637\\
	115.5	0.00509374151761222\\
	115.8	0.00204158410764421\\
	116.1	0.000141736201385356\\
	116.4	-0.0012671407249627\\
	116.7	-0.00197766738094174\\
	117	-0.000819173011052499\\
	117.3	-0.00157200536250457\\
	117.6	-0.00344126986070137\\
	117.9	-0.00370303120797644\\
	118.2	-0.00452331860694244\\
	118.5	-0.00443205450301321\\
	118.8	-0.00393903986898181\\
	119.1	-0.00302055292817727\\
	119.4	-0.00171168581852044\\
	119.7	-0.00129453882118469\\
	120	-0.00229030290833521\\
	120.3	-0.00389562500862439\\
	120.6	-0.00384115872637381\\
	120.9	-0.00464243342173631\\
	121.2	-0.00664446884835002\\
	121.5	-0.00885285773054534\\
	121.8	-0.0103173700405534\\
	122.1	-0.00936298672885982\\
	122.4	-0.00888787508080213\\
	122.7	-0.00908415940830309\\
	123	-0.00487841208970963\\
	123.3	-0.0026030030700781\\
	123.6	-0.00317054665194405\\
	123.9	-0.00242719405980552\\
	124.2	-0.00295760475890461\\
	124.5	-0.00303998758056423\\
	124.8	-0.00280249362711515\\
	125.1	-0.00298766534120887\\
	125.4	-0.0018310787758935\\
	125.7	0.000114451941897187\\
	138.3	0.000114451941897187\\
};
\addplot [color=red, line width=0.7pt, forget plot]
table[row sep=crcr]{%
	0	-0.181219799622411\\
	2.09999999999999	-0.181219799622411\\
	2.40000000000001	-0.176746777516271\\
	3	-0.160828793053128\\
	3.3	-0.160828793053128\\
	3.59999999999999	-0.134901707012887\\
	3.90000000000001	-0.121778006013926\\
	4.2	-0.106087636608606\\
	4.5	-0.0909031507766542\\
	4.8	-0.0739152281357462\\
	5.09999999999999	-0.0584724863320076\\
	5.40000000000001	-0.0422619498296939\\
	5.7	-0.0241932076787208\\
	6.3	0.0107093982798148\\
	6.59999999999999	0.0288319123620369\\
	6.90000000000001	0.0498196236282809\\
	7.2	0.0686310043239615\\
	7.5	0.0871028663421498\\
	7.8	0.103469610286439\\
	8.09999999999999	0.12528093267062\\
	8.40000000000001	0.137988947430486\\
	8.7	0.140963917683081\\
	9	0.149287358336835\\
	9.3	0.141620051844384\\
	9.59999999999999	0.139398984361094\\
	9.90000000000001	0.124976678807357\\
	10.5	0.111910625711886\\
	10.8	0.11684445761172\\
	11.1	0.119402953918438\\
	11.4	0.12732485460144\\
	11.7	0.141294052729407\\
	12	0.15858615597304\\
	12.3	0.17023699187105\\
	12.6	0.171748685277421\\
	12.9	0.157136704183372\\
	13.2	0.141209612731743\\
	13.5	0.136077713050838\\
	13.8	0.136420649392505\\
	14.1	0.146984352387818\\
	14.4	0.162710741621083\\
	14.7	0.17824587738194\\
	15.3	0.178242713345497\\
	15.6	0.172238122537891\\
	15.9	0.156257933010693\\
	16.2	0.132706507510491\\
	16.5	0.103722383082157\\
	16.8	0.0784070042335685\\
	17.1	0.0550039900531516\\
	17.4	0.0342404837858368\\
	17.7	0.0140887584728233\\
	18	-0.00521848227455735\\
	18.3	-0.0209526479984845\\
	18.6	-0.0344032745352507\\
	18.9	-0.0441292765251831\\
	19.2	-0.0561865243478366\\
	19.5	-0.0596798053826291\\
	19.8	-0.0650471649977078\\
	20.1	-0.07548866154481\\
	20.4	-0.0843359397932204\\
	20.7	-0.0902273674261096\\
	21	-0.0957150674910707\\
	21.6	-0.108752550899553\\
	21.9	-0.114612231915586\\
	22.2	-0.120280706797743\\
	22.5	-0.122725967535175\\
	22.8	-0.121965104333412\\
	23.1	-0.119484689367667\\
	23.4	-0.114216010445986\\
	23.7	-0.106411223755188\\
	24	-0.108484735927675\\
	24.3	-0.116784761066057\\
	24.6	-0.118108344928032\\
	24.9	-0.130721729893466\\
	25.2	-0.143628739310557\\
	25.5	-0.157936794359145\\
	25.8	-0.171417012006316\\
	26.1	-0.173468003052747\\
	26.4	-0.166597136099867\\
	26.7	-0.164659156100285\\
	27	-0.162933999665114\\
	27.3	-0.153264250525652\\
	27.6	-0.138676221438928\\
	27.9	-0.12325215368385\\
	28.2	-0.117199030848127\\
	28.5	-0.115425246062273\\
	28.8	-0.1179882986671\\
	29.1	-0.114448712427873\\
	29.4	-0.113778295587807\\
	29.7	-0.112665954806161\\
	30	-0.113128419925218\\
	30.3	-0.112159065545185\\
	30.6	-0.110905762017012\\
	30.9	-0.107176175577735\\
	31.2	-0.0971106788122285\\
	31.5	-0.0898431072314168\\
	32.1	-0.0742218621024335\\
	32.4	-0.0676343566336328\\
	32.7	-0.0606487611571964\\
	33	-0.0529044741161186\\
	33.3	-0.0430819806749838\\
	33.6	-0.0324284462726467\\
	33.9	-0.0240125406789957\\
	34.2	-0.0146169918933623\\
	34.5	-0.00759108242863249\\
	34.8	0.00195079847479462\\
	35.1	0.0079186625510772\\
	35.4	0.0147972115433106\\
	35.7	0.0179997131993161\\
	36	0.022749059662857\\
	36.3	0.029158316826198\\
	36.6	0.0385858037792275\\
	36.9	0.0449050809998255\\
	37.2	0.0477406961162785\\
	37.5	0.0496763003759355\\
	37.8	0.0511851933192844\\
	38.1	0.0535470497374746\\
	38.4	0.0509475916324362\\
	38.7	0.0512681895677076\\
	39	0.0505673982549553\\
	39.3	0.0483154447429541\\
	39.6	0.0458182020582285\\
	39.9	0.0472947892721294\\
	40.2	0.0462788728941774\\
	40.5	0.0516805656240535\\
	40.8	0.0516209664928482\\
	41.1	0.055710082225076\\
	41.4	0.0602043754443145\\
	41.7	0.0665486547106298\\
	42	0.0712268023587797\\
	42.3	0.0735554695185101\\
	42.6	0.073368286757713\\
	42.9	0.0761446782620823\\
	43.2	0.0793164179129349\\
	43.5	0.0774666733478711\\
	43.8	0.0707155781903737\\
	44.1	0.0644745829362705\\
	44.4	0.0707584929532459\\
	44.7	0.0703279667994678\\
	45	0.0711070083831942\\
	45.3	0.0704672121281078\\
	45.6	0.0704814617374012\\
	45.9	0.0726628240015685\\
	46.2	0.0731922852423139\\
	46.5	0.0694185189810099\\
	46.8	0.0643207757847648\\
	47.1	0.061101950824991\\
	47.4	0.0595275641449007\\
	47.7	0.0594666119967684\\
	48	0.0567465968190533\\
	48.3	0.054739529931112\\
	48.6	0.0519865198003941\\
	48.9	0.0462761680018957\\
	49.2	0.040760614348045\\
	49.5	0.0376092702419868\\
	50.1	0.0335386771514976\\
	50.4	0.0301529071665243\\
	50.7	0.0273821264138974\\
	51	0.0276257914458995\\
	51.3	0.0256459926756065\\
	51.6	0.0206589810786681\\
	51.9	0.019454505930895\\
	52.2	0.0153040537306879\\
	52.5	0.0105485347948786\\
	52.8	0.00486777298725372\\
	53.1	0.00193823231411727\\
	53.4	0.000221322520062017\\
	53.7	-0.00204196258175671\\
	54	-0.00509595343660862\\
	54.3	-0.00939043291927533\\
	54.6	-0.0120841127128273\\
	54.9	-0.0179490262274982\\
	55.2	-0.0203145541495218\\
	55.5	-0.0233714757519579\\
	55.8	-0.0219822068255553\\
	56.1	-0.0237633637462977\\
	56.4	-0.0262173183245551\\
	56.7	-0.029069762062278\\
	57	-0.0326572777700704\\
	57.3	-0.0340805650185274\\
	57.6	-0.0362449954009207\\
	57.9	-0.0373892350635145\\
	58.2	-0.0374915563587166\\
	58.5	-0.0372182114557091\\
	58.8	-0.0366248700044736\\
	59.1	-0.0369011115811304\\
	59.4	-0.0365616754684623\\
	59.7	-0.0387642911965287\\
	60	-0.0392884605471977\\
	60.3	-0.0396243586804559\\
	60.6	-0.0377345456797684\\
	60.9	-0.0344511641949055\\
	61.2	-0.0315294526103997\\
	61.5	-0.0267292651412987\\
	61.8	-0.0224850516983679\\
	62.1	-0.0196098472452491\\
	62.4	-0.0159680102740509\\
	62.7	-0.0126896230684821\\
	63	-0.00982725275582652\\
	63.6	-0.00489912857719332\\
	63.9	-0.00287276387371094\\
	64.2	-0.00287922228548609\\
	64.5	-0.00270248608472912\\
	64.8	-0.00134564630972989\\
	65.1	-0.00261505049114419\\
	65.7	-0.000133651276129854\\
	66	-0.000141248907453928\\
	66.3	0.00250605117446412\\
	66.6	0.00603812031576467\\
	66.9	0.00663788484185091\\
	67.2	0.00894395950291482\\
	67.5	0.011647854127375\\
	67.8	0.0151953587917859\\
	68.1	0.0164477122584969\\
	68.7	0.018591490893769\\
	69.3	0.0230367179861304\\
	69.6	0.0227659854910627\\
	69.9	0.0208097402800007\\
	70.2	0.0202708292537892\\
	70.5	0.0173406581576785\\
	70.8	0.013782872140979\\
	71.1	0.0142242935440891\\
	71.4	0.0151579908784214\\
	71.7	0.018303448211114\\
	72	0.0159250917603231\\
	72.3	0.0152335270812785\\
	72.6	0.0123762044835019\\
	72.9	0.012046172134589\\
	73.5	0.0123312581323773\\
	73.8	0.0108807130301756\\
	74.1	0.0109307263087004\\
	74.4	0.00987795593037788\\
	74.7	0.0078607815281373\\
	75	0.0045552849866084\\
	75.6	0.00486472426814544\\
	75.9	0.00375378409810878\\
	76.2	0.00325224390789458\\
	76.5	0.00307845719152056\\
	76.8	0.00234360619189999\\
	77.1	0.00188453008017575\\
	77.4	-0.00114249548445855\\
	77.7	-0.000190981995089601\\
	78	-0.000277182155045352\\
	78.3	0.00247127675436332\\
	78.6	0.00212689958752321\\
	78.9	0.00278486305344927\\
	79.2	0.00268093763584432\\
	79.5	0.00227426692583776\\
	79.8	0.0028164856690438\\
	80.1	0.00205636086313632\\
	80.4	0.00362452058075746\\
	80.7	0.00224104062522201\\
	81	0.00149685818891498\\
	81.3	0.002142798573189\\
	81.6	-0.000772960666665767\\
	81.9	-0.00315966505097265\\
	82.2	-0.0047777745897406\\
	82.5	-0.00729594231592046\\
	82.8	-0.0109328756565077\\
	83.4	-0.0130095186660384\\
	83.7	-0.0128917150935877\\
	84	-0.0121034182428019\\
	84.3	-0.0097068324723466\\
	84.6	-0.00830168984134616\\
	84.9	-0.00566619922840061\\
	85.2	-0.00438982642566543\\
	85.5	-0.00374752248337984\\
	86.1	-0.00220581134026077\\
	86.4	-0.00205736292599568\\
	86.7	-0.00410733843227717\\
	87	-0.0058659694026062\\
	87.3	-0.00859210799485766\\
	87.6	-0.00980792311631262\\
	87.9	-0.00998729431904621\\
	88.2	-0.0117554640070807\\
	88.5	-0.012571719991513\\
	88.8	-0.0106709632439674\\
	89.1	-0.00809897340432997\\
	89.4	-0.00710566741908281\\
	89.7	-0.00470782965412297\\
	90	-0.000988348753992341\\
	90.3	0.000261239384826695\\
	90.6	0.00220332418422231\\
	90.9	0.00207768906376771\\
	91.2	0.00330823718879003\\
	91.5	0.00226995136739561\\
	91.8	0.00610542887221754\\
	92.1	0.00739662625676374\\
	92.4	0.000114451941897187\\
	98.1	0.000114451941897187\\
};
\addplot [color=green, dashed, line width=0.7pt, forget plot]
table[row sep=crcr]{%
	0	-0.165412303349058\\
	0.300000000000011	-0.16545358322179\\
	0.900000000000006	-0.165134826279342\\
	1.19999999999999	-0.164568420724549\\
	1.80000000000001	-0.164211706168089\\
	2.09999999999999	-0.16456628186495\\
	2.40000000000001	-0.165258551508487\\
	3.30000000000001	-0.168031520085009\\
	3.59999999999999	-0.168775239195639\\
	4.5	-0.169902624119487\\
	5.09999999999999	-0.170239186391456\\
	5.69999999999999	-0.170136889734493\\
	6	-0.168890274055599\\
	6.30000000000001	-0.166696522623624\\
	6.59999999999999	-0.163556270956235\\
	6.90000000000001	-0.159041995110186\\
	7.19999999999999	-0.152463866297268\\
	7.5	-0.143766955049671\\
	7.80000000000001	-0.132426003054036\\
	8.09999999999999	-0.117965794581153\\
	8.40000000000001	-0.0983593782496257\\
	8.69999999999999	-0.0684945977454845\\
	9	-0.0316583439706051\\
	9.30000000000001	0.00155721093216243\\
	9.59999999999999	0.026463958155091\\
	9.90000000000001	0.0463764393181805\\
	10.2	0.0624954002837796\\
	10.5	0.0753266046870635\\
	10.8	0.0862721956063695\\
	11.1	0.0938827960955564\\
	11.4	0.09900717669953\\
	11.7	0.0994518267048079\\
	12	0.100501696870992\\
	12.3	0.0984366060361026\\
	12.6	0.0961823527837566\\
	12.9	0.0924410066106134\\
	13.2	0.0895920280595135\\
	13.5	0.0869683463305648\\
	13.8	0.0791984537011103\\
	14.1	0.0742394368569421\\
	14.4	0.0661502483619927\\
	14.7	0.0596203428690671\\
	15	0.0544322829259727\\
	15.3	0.0504618386820823\\
	15.6	0.043089449418801\\
	15.9	0.0362929283525091\\
	16.2	0.0334716770062755\\
	16.5	0.0298434085665065\\
	16.8	0.0279376799184661\\
	17.1	0.0243317552261146\\
	17.4	0.0217515801544153\\
	18	0.0144407090100458\\
	18.3	0.0102934420382041\\
	18.6	0.00556876221381231\\
	18.9	-0.000564836700931437\\
	19.2	-0.00713323832468404\\
	19.5	-0.0090483578086662\\
	19.8	-0.0145869335343605\\
	20.1	-0.0158099248528174\\
	20.4	-0.0194289322814996\\
	20.7	-0.0206381464732033\\
	21	-0.0214317235068222\\
	21.3	-0.0252083241179264\\
	21.6	-0.0253066353639042\\
	21.9	-0.0263979471508833\\
	22.2	-0.0271508774276299\\
	22.5	-0.0311520253465289\\
	22.8	-0.0324485750891199\\
	23.1	-0.0318729214655491\\
	23.4	-0.0316517991612386\\
	23.7	-0.0335694973483953\\
	24	-0.0347317036305981\\
	24.3	-0.0344460741143848\\
	24.6	-0.0347756281729517\\
	24.9	-0.0332234519599695\\
	25.2	-0.0305989025812892\\
	25.5	-0.0282469206002531\\
	25.8	-0.0301660704763265\\
	26.1	-0.0346443198564543\\
	26.4	-0.0367482521394606\\
	26.7	-0.0368691897795657\\
	27	-0.0394787405548698\\
	27.3	-0.0406082596688577\\
	27.6	-0.0363728297980401\\
	27.9	-0.031574199261172\\
	28.2	-0.0309740862255126\\
	28.5	-0.0340249397355024\\
	28.8	-0.0363462761960136\\
	29.1	-0.0370410960852041\\
	29.7	-0.0405694456877939\\
	30	-0.0382830182002749\\
	30.3	-0.0434234774820652\\
	30.6	-0.0433231985527698\\
	30.9	-0.0450007615323784\\
	31.2	-0.0418700143534636\\
	31.5	-0.0398255530907079\\
	31.8	-0.0397246515112499\\
	32.1	-0.0372358150438572\\
	32.4	-0.037418576306095\\
	32.7	-0.0356862375807054\\
	33	-0.0360999155887782\\
	33.3	-0.0368326014828995\\
	33.6	-0.0373375326029475\\
	33.9	-0.0350543453741352\\
	34.2	-0.0336499265133341\\
	34.5	-0.0338982274337241\\
	34.8	-0.0354806574179065\\
	35.1	-0.0337192068056709\\
	35.4	-0.0353088564761492\\
	35.7	-0.0321220317208883\\
	36	-0.0334336547273608\\
	36.3	-0.0327337383353665\\
	36.6	-0.0332960553416513\\
	36.9	-0.0334387088377923\\
	37.2	-0.0319243319102043\\
	37.5	-0.031680209818262\\
	37.8	-0.0308310268617618\\
	38.1	-0.0303557059430375\\
	38.4	-0.0284454911353009\\
	38.7	-0.028071952731608\\
	39	-0.0260879630381794\\
	39.3	-0.0245524967323547\\
	39.6	-0.0233923832370522\\
	39.9	-0.0220698201023026\\
	40.5	-0.0209929757576219\\
	40.8	-0.0198711458463379\\
	41.1	-0.0175303366898163\\
	41.4	-0.0186852297716484\\
	41.7	-0.0163758502099256\\
	42	-0.014462911006575\\
	42.3	-0.0122972643945616\\
	42.6	-0.00917213670697947\\
	42.9	-0.0085556654162815\\
	43.2	-0.00738770701707381\\
	43.5	-0.00481357900693524\\
	43.8	-0.00378813613670559\\
	44.1	-0.00154933826740944\\
	44.4	-0.000109260488898144\\
	44.7	0.00248639529058892\\
	45.3	0.00445758720442768\\
	45.6	0.00495572567811564\\
	45.9	0.00747821235094648\\
	46.8	0.00897150574382977\\
	47.1	0.00995489888026668\\
	47.4	0.0114225012289921\\
	47.7	0.0119356244468918\\
	48	0.0100909120393737\\
	48.3	0.00780562905367788\\
	48.6	0.00751518178856259\\
	48.9	0.00763975202923461\\
	49.2	0.00599006278133629\\
	49.5	0.00724203679621382\\
	49.8	0.00769707186466917\\
	50.1	0.00704739520332964\\
	50.4	0.00797303836500873\\
	50.7	0.0097445709256192\\
	51	0.00885209488913574\\
	51.6	0.0137637800445418\\
	52.2	0.0131830767178087\\
	52.5	0.0127749303457279\\
	52.8	0.0126172630138797\\
	53.1	0.0112677254885227\\
	53.4	0.0111876331118594\\
	54	0.00681051810352074\\
	54.3	0.00483688537499916\\
	54.6	0.00305335017620223\\
	54.9	0.000207308593132893\\
	55.2	-0.00125354388066512\\
	55.5	-0.000912663671385872\\
	55.8	-0.000781507768948586\\
	56.1	-0.00204849091855408\\
	56.4	-0.00508345105606622\\
	56.7	-0.0058806184841842\\
	57	-0.00808139065020441\\
	57.3	-0.00894212580021758\\
	57.6	-0.00925628897013553\\
	57.9	-0.0129710807051424\\
	58.5	-0.0117957871373164\\
	58.8	-0.0122045915659044\\
	59.1	-0.0133252419202847\\
	59.4	-0.0149739251492065\\
	59.7	-0.0153078227135097\\
	60	-0.0154012788824218\\
	60.3	-0.0160696616625557\\
	60.6	-0.0165639855525797\\
	60.9	-0.0166314799125473\\
	61.5	-0.018588144637647\\
	61.8	-0.0178348422841452\\
	62.1	-0.0191662156115626\\
	62.4	-0.0194679704394218\\
	62.7	-0.0192133198411568\\
	63	-0.0194560576628078\\
	63.3	-0.0202127597407298\\
	63.6	-0.0220968877084431\\
	63.9	-0.0222771415830039\\
	64.2	-0.0207735532683557\\
	64.5	-0.0185301563891755\\
	64.8	-0.0195880607132608\\
	65.1	-0.0161777583228968\\
	65.4	-0.0144676850587473\\
	65.7	-0.013012464460445\\
	66	-0.0123711300654179\\
	66.3	-0.0122816064254891\\
	66.6	-0.0113369290608603\\
	66.9	-0.0131247646555153\\
	67.2	-0.011716406778703\\
	67.5	-0.0126073410922061\\
	67.8	-0.012819273411651\\
	68.1	-0.012461413956629\\
	68.4	-0.0105507934909497\\
	68.7	-0.00917486038426318\\
	69	-0.00956068323330328\\
	69.6	-0.00410643916569597\\
	69.9	-0.00568485661472096\\
	70.2	-0.0063876487022867\\
	70.5	-0.006611443958775\\
	70.8	-0.00452701368413955\\
	71.1	-0.00215019315129439\\
	71.4	-0.00212642899904836\\
	71.7	-0.0023778422421401\\
	72	-0.00575003394536111\\
	72.3	-0.00666297497417645\\
	72.6	-0.00567260346099374\\
	72.9	-0.00490507992773814\\
	73.2	-0.00203617733407668\\
	73.5	-0.00075329824753112\\
	73.8	-0.00170922860894507\\
	74.1	0.00186777173925634\\
	74.7	0.000926350185892488\\
	75	0.00103164219149221\\
	75.3	0.000646584795930494\\
	75.6	-0.000482164572986221\\
	75.9	0.00136765637569169\\
	76.2	0.00107481161640521\\
	76.5	0.00216081531942791\\
	76.8	0.00202396789646286\\
	77.1	0.00313880603192729\\
	77.4	0.00160842013229967\\
	77.7	0.00112095376212551\\
	78	0.00118777062436948\\
	78.6	0.00406354024593725\\
	78.9	0.00483222095456881\\
	79.2	0.00706553973046198\\
	79.5	0.00647724053132492\\
	79.8	0.00624072407489962\\
	80.1	0.00502676786291545\\
	80.4	0.00331915657136506\\
	80.7	0.00131243915271284\\
	81.6	-0.00169949206159004\\
	81.9	-8.97672998121379e-05\\
	82.2	-0.000359742509601801\\
	82.5	-0.00100641278638136\\
	82.8	-0.00263730537906781\\
	83.1	-0.00397638409666001\\
	83.4	-0.0027080102733521\\
	83.7	-0.00333553306066392\\
	84	-0.00362792137804036\\
	84.3	-0.00446680625157114\\
	84.6	-0.00387103442886882\\
	84.9	-0.004855808393188\\
	85.2	-0.0039190309236119\\
	85.5	-0.00427351581808466\\
	85.8	-0.00397972197254148\\
	86.4	-0.00167060800859531\\
	86.7	-2.58589892609962e-07\\
	87	0.000397306150546228\\
	87.3	0.0011357147959643\\
	87.6	0.00154007812312784\\
	87.9	0.000955678851937591\\
	88.2	0.000945540019245072\\
	88.5	-1.01322176249141e-05\\
	88.8	0.00054485518154479\\
	89.1	-0.00075180801394481\\
	89.4	-0.00168483585949275\\
	89.7	-0.00512950320907635\\
	90	-0.00585123439532254\\
	90.3	-0.00486453715529933\\
	90.6	-0.00501365892475292\\
	90.9	-0.00488265948766298\\
	91.2	-0.0024744299179531\\
	91.5	0.000462941367061376\\
	91.8	0.000554070353814495\\
	92.1	0.0020895092632145\\
	92.4	0.00263752982027654\\
	92.7	0.0033609204142806\\
	93	0.00147428889971479\\
	93.3	-0.000619791402840519\\
	93.6	-0.00135100270236421\\
	93.9	-0.00308135962319511\\
	94.2	-0.00121015617773423\\
	94.5	-0.00104576403469991\\
	94.8	-0.00160905906264475\\
	95.1	-0.00146039770237394\\
	95.4	-0.00354969089954693\\
	95.7	-0.00415224750923926\\
	96	-0.00504546813550633\\
	96.3	-0.00501214470617128\\
	96.6	-0.00531566608918865\\
	96.9	-0.00659152393888007\\
	97.2	-0.00853686171504364\\
	97.5	-0.00790551732217182\\
	97.8	-0.00627882317354533\\
	98.1	-0.00697793828149429\\
	98.4	-0.0045720507345095\\
	98.7	-0.00436677573520683\\
	99	-0.00449476660799064\\
	99.3	-0.00445107364265596\\
	99.6	-0.00502491497127266\\
	99.9	-0.00533849366010486\\
	100.2	-0.00392845016628485\\
	100.5	-0.00455913393946616\\
	101.1	-0.00301142999208537\\
	101.4	-0.0012555731719317\\
	101.7	-0.0015624300902175\\
	102	-0.00330585237520609\\
	102.6	-0.000874518035914207\\
	102.9	0.000592810741864014\\
	103.2	0.00130891819267731\\
	103.5	0.00289854240239151\\
	103.8	0.00345043526397149\\
	104.1	0.00261234872559157\\
	104.4	0.00419488366429732\\
	104.7	0.00308259889422402\\
	105.3	0.00196217110502062\\
	105.6	-0.00124961330294582\\
	105.9	-0.00330181016911979\\
	106.5	-0.00190578502454741\\
	106.8	0.00136195971870734\\
	107.1	0.000272954450025509\\
	107.4	0.000959465223246525\\
	107.7	0.00288547469494915\\
	108	0.00354971181792507\\
	108.3	0.0038892331989473\\
	108.6	0.002869526703563\\
	108.9	0.00440980958219939\\
	109.2	0.00548266148902599\\
	109.5	0.0035894256571396\\
	109.8	0.00418049933105635\\
	110.1	0.000793080123798973\\
	110.4	-0.00118533100615537\\
	110.7	-0.00212083067398794\\
	111	-0.00275066368519106\\
	111.3	-0.00280581693115778\\
	111.6	-0.00243966150179631\\
	111.9	-0.00269850466105481\\
	112.2	-0.00154080384660915\\
	112.5	-0.0013158145366674\\
	112.8	-0.00129819398816267\\
	113.4	-9.58793326333307e-05\\
	113.7	0.000260824836601614\\
	114	-0.00109685130249204\\
	114.3	0.000884483165037864\\
	114.6	0.0018452435401457\\
	114.9	0.00255105089951257\\
	115.2	0.00129896057330825\\
	115.5	0.00180261314210384\\
	115.8	0.00433865616000162\\
	116.1	0.00485203689279956\\
	116.4	0.00336237721188581\\
	116.7	0.00395794264346705\\
	117	0.00430448302597597\\
	117.3	0.00354017792221839\\
	117.6	0.00439688745888134\\
	117.9	0.00549427279906922\\
	118.2	0.00505558180938692\\
	118.5	0.00567010338903629\\
	118.8	0.006494315737541\\
	119.1	0.00636935824721263\\
	119.4	0.00450406990469787\\
	120	0.000316114935685619\\
	120.3	-0.000559377936610872\\
	120.6	-0.00390125596729263\\
	120.9	-0.00445308940044242\\
	121.2	-0.00447689017204311\\
	121.5	-0.00349976199592561\\
	121.8	-0.00344034936534854\\
	122.1	-0.00382923898055765\\
	122.4	-0.00386669616099766\\
	122.7	-0.00559646085866916\\
	123	-0.0055738783719903\\
	123.3	-0.00298391770806461\\
	123.6	-0.00181672124867305\\
	123.9	-0.00142939354952887\\
	124.2	-0.00072728311050696\\
	124.5	-0.00226242751753603\\
	124.8	-0.00162624461842142\\
	125.1	-0.00124546713186646\\
	125.4	-0.00124507882367197\\
	125.7	-0.00333341837293233\\
	126	-0.003157904657769\\
	126.3	-0.00189999125095142\\
	126.6	-0.00295743034737939\\
	126.9	-0.00372505590320316\\
	127.2	-0.00371276969531209\\
	127.8	0.000347495055251557\\
	128.1	0.00155395977503758\\
	128.4	0.0018196670345958\\
	128.7	0.00109529202561021\\
	129	0.00141973879925672\\
	129.3	0.000633596796205893\\
	129.6	-0.00173217265947301\\
	129.9	-0.00201944710408952\\
	130.2	-0.000460274807380756\\
	130.5	-0.000284548934644135\\
	130.8	0.000173992096534903\\
	131.1	0.000815328377370861\\
	131.4	0.00268341941654171\\
	131.7	0.00141628909918268\\
	132	0.0011027103205663\\
	132.3	-8.9395760141997e-05\\
	132.6	-0.000440079763251333\\
	132.9	-0.00379861132134351\\
	133.2	-0.00480557048248897\\
	133.5	-0.00468591567505428\\
	133.8	-0.00496273229103394\\
	134.1	-0.00346071805981296\\
	134.4	-0.002988589598516\\
	134.7	-0.00347196902544056\\
	135	-0.00219859662007593\\
	135.3	-0.00192649621226337\\
	135.6	2.09795971954918e-05\\
	135.9	0.0027803271784137\\
	136.2	0.00298283709153679\\
	136.5	0.00464892016572094\\
	137.1	0.000425565855152854\\
	137.4	-0.000791411141108256\\
	137.7	-0.00440026674135652\\
	138	-0.00661840398444724\\
	138.3	-0.00687303095864422\\
	138.6	-0.00669599907365637\\
	138.9	-0.0075760724932934\\
	139.2	-0.00763468085617092\\
	139.5	-0.00834277165651542\\
	140.1	-0.00766577444986183\\
};
\addplot [color=blue, dashed, line width=0.7pt, forget plot]
table[row sep=crcr]{%
	0	-0.18\\
	18.3	-0.18\\
	18.6	-0.149846044608829\\
	18.9	-0.152649590223518\\
	19.2	-0.149539337209625\\
	19.5	-0.147108019642523\\
	19.8	-0.145805512231405\\
	20.1	-0.142444146373027\\
	20.4	-0.14305636743445\\
	21	-0.142938791251492\\
	21.3	-0.143712158141017\\
	21.6	-0.143161479090878\\
	21.9	-0.143791494472648\\
	22.2	-0.143789306233085\\
	23.1	-0.145971998864717\\
	23.4	-0.147136547965079\\
	23.7	-0.148707369692794\\
	24	-0.150450860859813\\
	25.2	-0.158454676597287\\
};
\addplot [color=red, dashed,line width=0.7pt, forget plot]
table[row sep=crcr]{%
	0	-0.18\\
	15.6	-0.18\\
	15.9	-0.156046753133275\\
	16.2	-0.0650033831446848\\
	16.5	0.0351387472782001\\
	16.8	0.145753167377659\\
	17.1	0.18\\
	24.9	0.18\\
};
\addplot [color=black, dotted, line width=1pt, forget plot]
table[row sep=crcr]{%
	0	0.181426947281579\\
	140	0.181426947281579\\
};
\addplot [color=black, dotted, line width=1pt, forget plot]
table[row sep=crcr]{%
	0	-0.181426947281579\\
	140	-0.181426947281579\\
};
\end{axis}
\end{tikzpicture}%

%% file: heading_error_straight.tex
%
%
\definecolor{mycolor1}{rgb}{0.50000,0.00000,0.50000}%
\begin{tikzpicture}

\begin{axis}[%
width=\figurewidth,
height=\figureheight,
at={(0\figurewidth,0\figureheight)},
scale only axis,
xmin=0,
xmax=120,
xtick={0,30,60,90,120},
xlabel={$t$ [s]},
ymin=-1,
ymax=1,
xlabel style={font=\color{white!15!black},at={(axis description cs:0.5,-0.17)},anchor=north},
ylabel style={font=\color{white!15!black},at={(axis description cs:-0.25,.5)},anchor=south},
ylabel={$\tilde \theta_3$ [rad]},
axis background/.style={fill=white},
xmajorgrids,
ymajorgrids
]
\addplot [color=blue,line width=0.7pt, forget plot]
  table[row sep=crcr]{%
0	-0.770782470703125\\
1.80000305175781	-0.76971435546875\\
2.69999694824219	-0.768051147460938\\
3.60000610351563	-0.764602661132813\\
4.19999694824219	-0.760757446289063\\
4.80000305175781	-0.755264282226563\\
5.39999389648438	-0.7484130859375\\
6	-0.7401123046875\\
6.60000610351563	-0.73004150390625\\
7.19999694824219	-0.718154907226563\\
7.80000305175781	-0.704147338867188\\
8.39999389648438	-0.688262939453125\\
9	-0.670501708984375\\
9.60000610351563	-0.651077270507813\\
10.5	-0.618820190429688\\
11.1000061035156	-0.596267700195313\\
11.6999969482422	-0.571044921875\\
12.3000030517578	-0.544052124023438\\
13.1999969482422	-0.501632690429688\\
13.8000030517578	-0.472930908203125\\
14.3999938964844	-0.4429931640625\\
15	-0.41357421875\\
15.3000030517578	-0.397979736328125\\
16.8000030517578	-0.325942993164063\\
18.3000030517578	-0.255386352539063\\
19.5	-0.19293212890625\\
20.3999938964844	-0.145217895507813\\
20.6999969482422	-0.130294799804688\\
21.6000061035156	-0.0890045166015625\\
22.8000030517578	-0.038818359375\\
24.3000030517578	0.0231170654296875\\
25.5	0.072479248046875\\
26.1000061035156	0.09527587890625\\
27.6000061035156	0.148056030273438\\
27.8999938964844	0.158294677734375\\
28.1999969482422	0.169403076171875\\
28.8000030517578	0.19024658203125\\
29.6999969482422	0.222244262695313\\
30	0.232772827148438\\
31.8000030517578	0.2879638671875\\
32.1000061035156	0.295562744140625\\
32.4000091552734	0.3043212890625\\
33	0.31988525390625\\
33.6000061035156	0.33404541015625\\
34.1999969482422	0.345138549804688\\
34.8000030517578	0.35516357421875\\
35.1000061035156	0.359039306640625\\
35.6999969482422	0.364364624023438\\
36.6000061035156	0.3682861328125\\
37.1999969482422	0.367477416992188\\
37.8000030517578	0.364334106445313\\
38.4000091552734	0.357498168945313\\
38.6999969482422	0.353546142578125\\
39.3000030517578	0.342453002929688\\
39.6000061035156	0.336532592773438\\
40.1999969482422	0.321258544921875\\
40.8000030517578	0.301910400390625\\
41.1000061035156	0.291107177734375\\
42.3000030517578	0.242904663085938\\
43.8000030517578	0.182846069335938\\
44.6999969482422	0.147598266601563\\
45.3000030517578	0.125930786132813\\
46.8000030517578	0.076751708984375\\
47.6999969482422	0.04962158203125\\
48.6000061035156	0.0240325927734375\\
49.5	0.001220703125\\
50.4000091552734	-0.0182647705078125\\
51	-0.029937744140625\\
51.6000061035156	-0.0396881103515625\\
51.9000091552734	-0.043670654296875\\
52.8000030517578	-0.05242919921875\\
53.4000091552734	-0.0561370849609375\\
54	-0.058563232421875\\
54.3000030517578	-0.05914306640625\\
54.6000061035156	-0.0584564208984375\\
55.5	-0.053436279296875\\
56.1000061035156	-0.0482635498046875\\
56.4000091552734	-0.0452880859375\\
57.9000091552734	-0.0247802734375\\
59.6999969482422	0.003570556640625\\
61.1999969482422	0.028900146484375\\
61.8000030517578	0.03759765625\\
62.4000091552734	0.0457611083984375\\
63	0.0529022216796875\\
63.6000061035156	0.058929443359375\\
64.5	0.065338134765625\\
65.1000061035156	0.0683135986328125\\
65.6999969482422	0.0701446533203125\\
66.6000061035156	0.07098388671875\\
67.5	0.0695037841796875\\
68.4000091552734	0.0666351318359375\\
69.3000030517578	0.0620880126953125\\
70.1999969482422	0.0559844970703125\\
73.1999969482422	0.0293426513671875\\
75.6000061035156	0.007904052734375\\
78	-0.007293701171875\\
79.1999969482422	-0.01080322265625\\
80.1000061035156	-0.01202392578125\\
81.3000030517578	-0.0112762451171875\\
83.1000061035156	-0.007415771484375\\
86.1000061035156	0.0025634765625\\
87.6000061035156	0.006622314453125\\
90.6000061035156	0.0130615234375\\
91.5	0.0135650634765625\\
92.6999969482422	0.012481689453125\\
93.9000091552734	0.0105743408203125\\
95.4000091552734	0.006988525390625\\
98.1000061035156	1.52587890625e-05\\
101.100006103516	-0.006317138671875\\
103.5	-0.007720947265625\\
104.699996948242	-0.00701904296875\\
105.600006103516	-0.0065460205078125\\
112.800003051758	-0.0050506591796875\\
121.800003051758	0.00408935546875\\
124.199996948242	0.0048675537109375\\
127.800003051758	0.0048675537109375\\
132	0.00469970703125\\
134.699996948242	0.0046234130859375\\
137.699996948242	0.0044708251953125\\
138.300003051758	0.0044403076171875\\
};
\addplot [color=red, line width=0.7pt, forget plot]
  table[row sep=crcr]{%
0	-0.0505599975585938\\
2.09999847412109	-0.0505599975585938\\
2.40000152587891	-0.0517730712890625\\
3.59999847412109	-0.0521774291992188\\
4.19999694824219	-0.0514450073242188\\
4.80000305175781	-0.0495376586914063\\
5.40000152587891	-0.0463180541992188\\
6	-0.0416641235351563\\
6.59999847412109	-0.035614013671875\\
7.19999694824219	-0.027984619140625\\
7.80000305175781	-0.0189132690429688\\
8.40000152587891	-0.00801849365234375\\
9	0.0039520263671875\\
9.59999847412109	0.017791748046875\\
10.1999969482422	0.033111572265625\\
10.8000030517578	0.0496292114257813\\
11.4000015258789	0.0679550170898438\\
12	0.088104248046875\\
13.5	0.141433715820313\\
15	0.197830200195313\\
15.5999984741211	0.218986511230469\\
16.2000045776367	0.238510131835938\\
16.8000030517578	0.255912780761719\\
18	0.287544250488281\\
19.5	0.3240966796875\\
20.0999984741211	0.337753295898438\\
21	0.357086181640625\\
21.5999984741211	0.368759155273438\\
22.5	0.384452819824219\\
23.7000045776367	0.400688171386719\\
24.5999984741211	0.408317565917969\\
24.9000015258789	0.409744262695313\\
25.2000045776367	0.410263061523438\\
25.5	0.409996032714844\\
26.0999984741211	0.407402038574219\\
26.4000015258789	0.405487060546875\\
27.9000015258789	0.390106201171875\\
28.5	0.380821228027344\\
29.0999984741211	0.369476318359375\\
29.7000045776367	0.355880737304688\\
30.3000030517578	0.339485168457031\\
31.2000045776367	0.312705993652344\\
31.8000030517578	0.293830871582031\\
33.9000015258789	0.225326538085938\\
34.8000030517578	0.197425842285156\\
35.7000045776367	0.1715087890625\\
36.3000030517578	0.1558837890625\\
37.5	0.125785827636719\\
38.4000015258789	0.105583190917969\\
39	0.09307861328125\\
39.3000030517578	0.08660888671875\\
40.2000045776367	0.0702362060546875\\
41.0999984741211	0.056182861328125\\
42	0.0447158813476563\\
43.2000045776367	0.03326416015625\\
43.5	0.03125\\
44.0999984741211	0.02581787109375\\
44.7000045776367	0.0232467651367188\\
45.5999984741211	0.0209808349609375\\
47.0999984741211	0.0223846435546875\\
48	0.0254669189453125\\
49.5	0.032806396484375\\
51	0.0424270629882813\\
53.7000045776367	0.0590057373046875\\
54.5999984741211	0.0629730224609375\\
56.0999984741211	0.0671539306640625\\
56.7000045776367	0.0676422119140625\\
57.9000015258789	0.0663375854492188\\
58.8000030517578	0.0634918212890625\\
60	0.0570907592773438\\
66	0.0200271606445313\\
67.8000030517578	0.01214599609375\\
69	0.00901031494140625\\
71.4000015258789	0.0059967041015625\\
74.0999984741211	0.00713348388671875\\
83.0999984741211	0.011474609375\\
87.5999984741211	0.007232666015625\\
89.6999969482422	0.00341033935546875\\
92.4000015258789	0.0006256103515625\\
95.6999969482422	-0.0001068115234375\\
98.0999984741211	3.0517578125e-05\\
};
\addplot [color=green, dashed,line width=0.7pt, forget plot]
  table[row sep=crcr]{%
0	-0.414993286132813\\
7.80000305175781	-0.414474487304688\\
8.69999694824219	-0.412704467773438\\
9.30000305175781	-0.409713745117188\\
9.89999389648438	-0.405105590820313\\
10.5	-0.398956298828125\\
11.3999938964844	-0.387252807617188\\
12.3000030517578	-0.373458862304688\\
13.1999969482422	-0.357833862304688\\
15	-0.3236083984375\\
16.1999969482422	-0.299118041992188\\
17.3999938964844	-0.2733154296875\\
19.1999969482422	-0.234222412109375\\
21	-0.198074340820313\\
23.1000061035156	-0.158981323242188\\
25.1999969482422	-0.121994018554688\\
26.3999938964844	-0.102890014648438\\
27.6000061035156	-0.0852813720703125\\
28.5	-0.07275390625\\
30	-0.0548858642578125\\
31.8000030517578	-0.036895751953125\\
33.3000030517578	-0.025604248046875\\
34.8000030517578	-0.016815185546875\\
36	-0.0117950439453125\\
37.1999969482422	-0.0091094970703125\\
38.4000091552734	-0.00787353515625\\
39.9000091552734	-0.008636474609375\\
42.3000030517578	-0.012359619140625\\
44.1000061035156	-0.014739990234375\\
45.9000091552734	-0.014892578125\\
47.6999969482422	-0.0133056640625\\
49.1999969482422	-0.010894775390625\\
51.3000030517578	-0.0047149658203125\\
54.3000030517578	0.00946044921875\\
57	0.0220184326171875\\
59.1000061035156	0.02984619140625\\
60.3000030517578	0.033294677734375\\
62.1000061035156	0.0365142822265625\\
63.6000061035156	0.0376739501953125\\
65.1000061035156	0.037628173828125\\
66.6000061035156	0.0360260009765625\\
73.1999969482422	0.0258636474609375\\
76.1999969482422	0.022308349609375\\
79.8000030517578	0.0204925537109375\\
82.5	0.0200042724609375\\
87	0.0177154541015625\\
96	0.0135040283203125\\
102.300003051758	0.0054168701171875\\
105.900009155273	0.00091552734375\\
109.5	-0.0021820068359375\\
116.400009155273	-0.0069122314453125\\
118.199996948242	-0.006805419921875\\
120.900009155273	-0.0057525634765625\\
123	-0.0061798095703125\\
130.800003051758	-0.006134033203125\\
132.899993896484	-0.0058746337890625\\
135.600006103516	-0.005859375\\
137.399993896484	-0.0056610107421875\\
139.5	-0.006866455078125\\
140.100006103516	-0.0072021484375\\
};
\addplot [color=red, dashed, line width=0.7pt, forget plot]
  table[row sep=crcr]{%
0	-0.05841064453125\\
11.3999996185303	-0.05780029296875\\
12.3000011444092	-0.056396484375\\
13.2000007629395	-0.0531654357910156\\
13.8000011444092	-0.0494766235351563\\
14.7000007629395	-0.0417232513427734\\
15.3000011444092	-0.0355358123779297\\
15.8999996185303	-0.02838134765625\\
16.5	-0.0203094482421875\\
17.1000003814697	-0.0112800598144531\\
17.7000007629395	-9.5367431640625e-05\\
18	-0.000886917114257813\\
19.7999992370605	-0.00114059448242188\\
24.8999996185303	-0.005859375\\
};
\addplot [color=blue, dashed, line width=0.7pt, forget plot]
  table[row sep=crcr]{%
0	-0.720314025878906\\
12	-0.719825744628906\\
12.6000003814697	-0.718732833862305\\
13.2000007629395	-0.716423034667969\\
13.8000011444092	-0.712709426879883\\
14.3999996185303	-0.707666397094727\\
15	-0.701549530029297\\
15.8999996185303	-0.691104888916016\\
17.3999996185303	-0.672700881958008\\
18.2999992370605	-0.672979354858398\\
21.8999996185303	-0.673591613769531\\
24	-0.674417495727539\\
};
\end{axis}
\end{tikzpicture}%

%% file: lateral_error_straight.tex
%
%
\definecolor{mycolor1}{rgb}{0.50000,0.00000,0.50000}%
\begin{tikzpicture}

\begin{axis}[%
width=\figurewidth,
height=\figureheight,
at={(0\figurewidth,0\figureheight)},
scale only axis,
xmin=0,
xmax=120,
xtick={0,30,60,90,120},
xlabel={$t$ [s]},
ymin=-8,
ymax=8,
ytick={-8,-4,0,4,8},
xlabel style={font=\color{white!15!black},at={(axis description cs:0.5,-0.17)},anchor=north},
ylabel style={font=\color{white!15!black},at={(axis description cs:-0.2,.5)},anchor=south},
ylabel={$\tilde z_3$ [m]},
axis background/.style={fill=white},
xmajorgrids,
ymajorgrids
]
\addplot [color=blue, line width=0.7pt, forget plot]
  table[row sep=crcr]{%
0	-1.40040588378906\\
0.300003051757813	-1.38803100585938\\
0.600006103515625	-1.37088012695313\\
1.19999694824219	-1.32159423828125\\
1.5	-1.28999328613281\\
1.80000305175781	-1.2529296875\\
2.10000610351563	-1.20933532714844\\
2.39999389648438	-1.15560913085938\\
2.69999694824219	-1.09561157226563\\
3	-1.02761840820313\\
3.30000305175781	-0.952178955078125\\
3.60000610351563	-0.869415283203125\\
3.89999389648438	-0.780517578125\\
4.19999694824219	-0.686386108398438\\
5.10000610351563	-0.384201049804688\\
9.30000305175781	1.03834533691406\\
10.1999969482422	1.32640075683594\\
10.8000030517578	1.506591796875\\
11.1000061035156	1.59042358398438\\
11.6999969482422	1.76362609863281\\
12.3000030517578	1.92829895019531\\
12.8999938964844	2.08395385742188\\
13.5	2.22854614257813\\
14.3999938964844	2.42448425292969\\
14.6999969482422	2.48500061035156\\
15.8999938964844	2.70121765136719\\
16.5	2.79586791992188\\
17.3999938964844	2.92169189453125\\
18.3000030517578	3.03628540039063\\
18.8999938964844	3.1075439453125\\
19.5	3.17327880859375\\
20.1000061035156	3.228759765625\\
20.6999969482422	3.2720947265625\\
21.3000030517578	3.30299377441406\\
21.8999938964844	3.32341003417969\\
22.5	3.33656311035156\\
23.1000061035156	3.34225463867188\\
24	3.33836364746094\\
24.6000061035156	3.32260131835938\\
25.1999969482422	3.29638671875\\
26.1000061035156	3.2440185546875\\
26.6999969482422	3.19944763183594\\
27.3000030517578	3.14854431152344\\
27.8999938964844	3.09007263183594\\
28.8000030517578	2.97782897949219\\
29.3999938964844	2.88629150390625\\
29.6999969482422	2.83790588378906\\
30.3000030517578	2.7322998046875\\
30.8999938964844	2.615234375\\
31.8000030517578	2.41551208496094\\
32.1000061035156	2.34809875488281\\
32.6999969482422	2.19271850585938\\
33.3000030517578	2.03096008300781\\
34.8000030517578	1.58351135253906\\
35.4000091552734	1.40072631835938\\
35.6999969482422	1.30227661132813\\
36.3000030517578	1.11531066894531\\
36.6000061035156	1.01678466796875\\
38.1000061035156	0.553558349609375\\
38.4000091552734	0.461807250976563\\
39.6000061035156	0.12457275390625\\
39.9000091552734	0.0498199462890625\\
40.5	-0.087158203125\\
40.8000030517578	-0.149124145507813\\
41.1000061035156	-0.206130981445313\\
41.4000091552734	-0.257278442382813\\
41.6999969482422	-0.303131103515625\\
42	-0.343994140625\\
42.6000061035156	-0.413742065429688\\
42.9000091552734	-0.439010620117188\\
43.1999969482422	-0.45965576171875\\
43.8000030517578	-0.4822998046875\\
44.4000091552734	-0.486541748046875\\
45	-0.473846435546875\\
45.6000061035156	-0.446685791015625\\
46.5	-0.38702392578125\\
47.1000061035156	-0.3377685546875\\
47.4000091552734	-0.311264038085938\\
48.9000091552734	-0.150497436523438\\
51.3000030517578	0.106689453125\\
51.9000091552734	0.164901733398438\\
54.6000061035156	0.397186279296875\\
55.8000030517578	0.473800659179688\\
56.4000091552734	0.505599975585938\\
57	0.527099609375\\
57.6000061035156	0.538833618164063\\
58.5	0.541397094726563\\
59.1000061035156	0.532577514648438\\
59.6999969482422	0.515151977539063\\
60.6000061035156	0.47235107421875\\
61.5	0.4144287109375\\
62.1000061035156	0.37042236328125\\
63	0.295318603515625\\
66	0.0380096435546875\\
67.1999969482422	-0.05035400390625\\
68.4000091552734	-0.12506103515625\\
69.3000030517578	-0.17236328125\\
70.1999969482422	-0.209320068359375\\
71.1000061035156	-0.231399536132813\\
72	-0.242904663085938\\
72.9000091552734	-0.244537353515625\\
73.5	-0.24066162109375\\
74.6999969482422	-0.216766357421875\\
76.5	-0.17181396484375\\
77.4000091552734	-0.136093139648438\\
78	-0.11248779296875\\
80.1000061035156	-0.048675537109375\\
81	-0.030487060546875\\
82.1999969482422	-0.016082763671875\\
83.1000061035156	-0.0137176513671875\\
84.3000030517578	-0.0234832763671875\\
85.5	-0.0437469482421875\\
89.4000091552734	-0.121902465820313\\
90.9000091552734	-0.15826416015625\\
91.8000030517578	-0.174224853515625\\
93.6000061035156	-0.195327758789063\\
94.8000030517578	-0.20166015625\\
95.4000091552734	-0.203079223632813\\
97.1999969482422	-0.201248168945313\\
98.1000061035156	-0.19158935546875\\
99.3000030517578	-0.170883178710938\\
100.5	-0.147903442382813\\
103.199996948242	-0.11676025390625\\
104.699996948242	-0.112060546875\\
105.600006103516	-0.106552124023438\\
108.900009155273	-0.07073974609375\\
111.300003051758	-0.0357666015625\\
112.199996948242	-0.02301025390625\\
115.800003051758	0.00762939453125\\
117.900009155273	0.015594482421875\\
119.100006103516	0.01385498046875\\
120.300003051758	0.01812744140625\\
122.699996948242	0.0283203125\\
127.199996948242	0.040374755859375\\
128.100006103516	0.0450286865234375\\
130.199996948242	0.0479583740234375\\
134.399993896484	0.0410308837890625\\
138.300003051758	0.0425872802734375\\
};
\addplot [color=red, line width=0.7pt, forget plot]
  table[row sep=crcr]{%
0	5.62956237792969\\
2.09999847412109	5.62956237792969\\
2.40000152587891	5.64316558837891\\
3.59999847412109	5.66373443603516\\
5.40000152587891	5.70773315429688\\
6.90000152587891	5.74009704589844\\
8.09999847412109	5.75739288330078\\
9	5.76270294189453\\
9.90000152587891	5.75814819335938\\
10.8000030517578	5.74166107177734\\
11.4000015258789	5.72258758544922\\
12	5.69626617431641\\
12.9000015258789	5.64408874511719\\
13.5	5.60108947753906\\
14.4000015258789	5.52117919921875\\
15.3000030517578	5.42667388916016\\
16.5	5.28923034667969\\
17.4000015258789	5.18016815185547\\
18	5.10050964355469\\
18.5999984741211	5.01252746582031\\
19.2000045776367	4.91542816162109\\
19.8000030517578	4.80751800537109\\
20.4000015258789	4.69075775146484\\
21	4.56399536132813\\
21.5999984741211	4.42865753173828\\
22.2000045776367	4.28468322753906\\
22.8000030517578	4.13150024414063\\
23.4000015258789	3.96855163574219\\
24.3000030517578	3.70559692382813\\
27.3000030517578	2.76931762695313\\
28.2000045776367	2.4725341796875\\
29.4000015258789	2.0877685546875\\
30	1.90532684326172\\
30.5999984741211	1.73381805419922\\
31.2000045776367	1.57279205322266\\
31.8000030517578	1.42542266845703\\
32.4000015258789	1.29410552978516\\
33	1.17668151855469\\
33.5999984741211	1.07299041748047\\
34.2000045776367	0.978469848632813\\
34.5	0.935279846191406\\
35.4000015258789	0.825965881347656\\
36.3000030517578	0.737586975097656\\
36.9000015258789	0.692848205566406\\
37.5	0.658462524414063\\
38.0999984741211	0.632484436035156\\
38.7000045776367	0.613784790039063\\
39.3000030517578	0.607383728027344\\
39.9000015258789	0.605606079101563\\
40.5	0.605987548828125\\
41.0999984741211	0.613136291503906\\
42.5999984741211	0.633827209472656\\
43.8000030517578	0.658973693847656\\
44.4000015258789	0.672866821289063\\
46.2000045776367	0.702461242675781\\
47.0999984741211	0.714195251464844\\
48	0.718254089355469\\
48.9000015258789	0.712989807128906\\
49.5	0.703628540039063\\
50.4000015258789	0.677299499511719\\
51.3000030517578	0.637924194335938\\
53.0999984741211	0.537429809570313\\
54.3000030517578	0.458297729492188\\
56.7000045776367	0.29400634765625\\
58.5	0.192291259765625\\
59.0999984741211	0.167396545410156\\
60.3000030517578	0.1297607421875\\
61.5	0.107124328613281\\
64.1999969482422	0.0782623291015625\\
65.0999984741211	0.0832138061523438\\
69	0.115325927734375\\
71.6999969482422	0.114921569824219\\
72.5999984741211	0.111434936523438\\
75.9000015258789	0.1041259765625\\
77.0999984741211	0.0932769775390625\\
78.5999984741211	0.069671630859375\\
83.0999984741211	-0.00217437744140625\\
85.1999969482422	-0.0223388671875\\
87.5999984741211	-0.0275344848632813\\
92.4000015258789	0.001739501953125\\
96	0.0068511962890625\\
98.0999984741211	0\\
};
\addplot [color=green, dashed,line width=0.7pt, forget plot]
  table[row sep=crcr]{%
0	-4.05079650878906\\
3.89999389648438	-4.05682373046875\\
5.69999694824219	-4.05879211425781\\
6.30000305175781	-4.05326843261719\\
6.89999389648438	-4.04130554199219\\
7.5	-4.01524353027344\\
7.80000305175781	-3.99525451660156\\
8.10000610351563	-3.96800231933594\\
8.39999389648438	-3.93135070800781\\
8.69999694824219	-3.88539123535156\\
9	-3.83062744140625\\
9.60000610351563	-3.699462890625\\
9.89999389648438	-3.63267517089844\\
12.8999938964844	-2.90554809570313\\
13.8000030517578	-2.69943237304688\\
14.6999969482422	-2.50534057617188\\
15.8999938964844	-2.25865173339844\\
16.8000030517578	-2.08122253417969\\
18	-1.85482788085938\\
19.1999969482422	-1.64390563964844\\
20.3999938964844	-1.45365905761719\\
21.3000030517578	-1.32418823242188\\
22.1999969482422	-1.20469665527344\\
23.1000061035156	-1.09452819824219\\
24.6000061035156	-0.933120727539063\\
25.1999969482422	-0.876663208007813\\
26.3999938964844	-0.776123046875\\
27.3000030517578	-0.713470458984375\\
27.8999938964844	-0.677276611328125\\
28.8000030517578	-0.62890625\\
29.6999969482422	-0.5889892578125\\
30.6000061035156	-0.557159423828125\\
31.8000030517578	-0.524139404296875\\
32.6999969482422	-0.504608154296875\\
34.1999969482422	-0.477859497070313\\
36.3000030517578	-0.44586181640625\\
38.1000061035156	-0.40966796875\\
39	-0.385848999023438\\
40.1999969482422	-0.34429931640625\\
41.4000091552734	-0.294036865234375\\
42.9000091552734	-0.222640991210938\\
47.4000091552734	0.03057861328125\\
48.6000061035156	0.102554321289063\\
49.8000030517578	0.164138793945313\\
51	0.212265014648438\\
51.9000091552734	0.238479614257813\\
53.1000061035156	0.266510009765625\\
54.6000061035156	0.28997802734375\\
56.1000061035156	0.303756713867188\\
57.3000030517578	0.306198120117188\\
58.5	0.29840087890625\\
61.8000030517578	0.262680053710938\\
66.9000091552734	0.194229125976563\\
68.4000091552734	0.178115844726563\\
70.8000030517578	0.14923095703125\\
72.3000030517578	0.137481689453125\\
77.1000061035156	0.105316162109375\\
82.5	0.04058837890625\\
85.5	0.0212860107421875\\
89.1000061035156	-0.0212249755859375\\
90.6000061035156	-0.0433349609375\\
96	-0.140792846679688\\
98.4000091552734	-0.173797607421875\\
101.400009155273	-0.212890625\\
109.5	-0.269180297851563\\
112.199996948242	-0.274978637695313\\
116.100006103516	-0.273635864257813\\
118.800003051758	-0.276290893554688\\
122.699996948242	-0.290145874023438\\
126.600006103516	-0.285690307617188\\
129.600006103516	-0.2947998046875\\
132	-0.298355102539063\\
134.699996948242	-0.309432983398438\\
136.800003051758	-0.316177368164063\\
139.5	-0.304122924804688\\
140.100006103516	-0.302093505859375\\
};
\addplot [color=red, dashed,line width=0.7pt, forget plot]
  table[row sep=crcr]{%
0	5.58638572692871\\
12	5.58132553100586\\
13.8000011444092	5.56874465942383\\
16.2000007629395	5.53930854797363\\
18	5.51457214355469\\
22.5	5.51769638061523\\
24.8999996185303	5.51876258850098\\
};
\addplot [color=blue, dashed,line width=0.7pt, forget plot]
  table[row sep=crcr]{%
0	-1.30915069580078\\
6.30000114440918	-1.3218936920166\\
10.8000011444092	-1.30477333068848\\
11.1000003814697	-1.2943058013916\\
11.3999996185303	-1.27466583251953\\
11.7000007629395	-1.24553680419922\\
12	-1.20686531066895\\
12.3000011444092	-1.15714645385742\\
12.6000003814697	-1.09786987304688\\
12.8999996185303	-1.02925109863281\\
14.3999996185303	-0.651948928833008\\
15	-0.517898559570313\\
15.6000003814697	-0.399887084960938\\
16.2000007629395	-0.298894882202148\\
16.7999992370605	-0.211847305297852\\
17.1000003814697	-0.173740386962891\\
17.3999996185303	-0.146949768066406\\
18.2999992370605	-0.142673492431641\\
21	-0.130451202392578\\
24	-0.124168395996094\\
};
\end{axis}
\end{tikzpicture}%

%% file: curvature_eight.tex
%
%
\definecolor{mycolor1}{rgb}{0.50000,0.00000,0.50000}%
\begin{tikzpicture}

\begin{axis}[%
width=\figurewidth,
height=\figureheight,
at={(0\figurewidth,0\figureheight)},
scale only axis,
xmin=0,
xmax=140,
xtick={0,30,60,90,120},
ylabel style={font=\color{white!15!black},at={(axis description cs:-0.2,.5)},anchor=south},
xlabel={$t$ [s]},
ymin=-0.27,
ymax=0.27,
ytick = {-0.181,0,0.181},
yticklabels={$-u_{\text{max}}$,0,$u_{\text{max}}$},
xlabel style={font=\color{white!15!black},at={(axis description cs:0.5,-0.17)},anchor=north},
ylabel={$u$ [m$^{-1}$]},
axis background/.style={fill=white},
xmajorgrids,
ymajorgrids
]
\addplot [color=blue, line width=0.7pt, forget plot]
  table[row sep=crcr]{%
0	-0.0721438883267922\\
0.599999999999994	-0.0425736446041469\\
1.80000000000001	0.0221059005528161\\
2.40000000000001	0.0564501317184636\\
3	0.088672415958456\\
3.59999999999999	0.118262033110767\\
4.19999999999999	0.145019791354372\\
4.80000000000001	0.10838255610841\\
5.40000000000001	0.103764436153483\\
6	0.106704219859097\\
6.59999999999999	0.126597541524376\\
7.19999999999999	0.13061314340078\\
7.80000000000001	0.129492113448521\\
8.40000000000001	0.128866625283933\\
9	0.12428375699821\\
9.59999999999999	0.141283430690748\\
10.2	0.178246903565679\\
10.8	0.178248100786078\\
11.4	0.151628526198266\\
12	0.134338820717232\\
12.6	0.151997298426096\\
13.2	0.107193490218094\\
13.8	0.073275158066707\\
14.4	0.0564271550746867\\
15	0.0418516081204814\\
15.6	0.0285265456878676\\
16.8	0.00592185310156879\\
17.4	-0.00156698629547236\\
18	-0.0137100147098863\\
18.6	-0.0370423040860715\\
19.2	-0.051014159501733\\
19.8	-0.056433643231685\\
20.4	-0.0595643624898798\\
21	-0.065199379157832\\
21.6	-0.0672248137576901\\
22.2	-0.069834905760132\\
22.8	-0.0782678793444518\\
23.4	-0.0801277666270437\\
24	-0.0814536598613813\\
24.6	-0.089970686181374\\
25.2	-0.0924481372299226\\
25.8	-0.0796659755207827\\
26.4	-0.0755832074591467\\
27	-0.0765566632059631\\
27.6	-0.0732893660941727\\
28.2	-0.0736388935281695\\
28.8	-0.0575577725114726\\
29.4	-0.0525835871353877\\
30	-0.0468220608573233\\
30.6	-0.0421668041073815\\
31.2	-0.0330230481773128\\
31.8	-0.0286136520067544\\
32.4	-0.023307192704408\\
33	-0.0265370343733196\\
33.6	-0.0233161752433091\\
34.2	-0.0238971897998681\\
34.8	-0.0247938697774828\\
35.4	-0.0287047331550809\\
36	-0.023001084321038\\
36.6	-0.018241969850493\\
37.2	-0.0243146951426922\\
37.8	-0.02881243947607\\
38.4	-0.0315028986177595\\
39	-0.0269219304799151\\
39.6	-0.0252188672022555\\
40.2	-0.028383717803024\\
40.8	-0.0235002653228094\\
42	-0.0189585156139742\\
42.6	-0.0244304045872639\\
43.2	-0.0227822350605607\\
43.8	-0.0256732153255541\\
44.4	-0.0185702369168723\\
45	-0.0198136398648217\\
45.6	-0.0275865236979769\\
46.2	-0.0326908877708263\\
46.8	-0.0393521059613704\\
47.4	-0.0435595771008934\\
48	-0.0511797203170659\\
49.2	-0.0613109924243247\\
49.8	-0.0647158745037757\\
50.4	-0.0723856647277614\\
51	-0.0836994769636021\\
51.6	-0.0777960284132178\\
52.2	-0.069397863728085\\
52.8	-0.0633878041809623\\
53.4	-0.0557727092860887\\
54	-0.0560553496606531\\
54.6	-0.0599610285101448\\
55.2	-0.0668037187099628\\
55.8	-0.0809799173904082\\
56.4	-0.0855419824332273\\
57	-0.0830890357624412\\
57.6	-0.083874920192244\\
58.2	-0.0887267293994967\\
58.8	-0.0857549507142039\\
59.4	-0.0855880733280969\\
60	-0.089052308522497\\
60.6	-0.0895331485226336\\
61.2	-0.0792643241778421\\
61.8	-0.0605569349921211\\
62.4	-0.0575065143899849\\
63	-0.059256482754023\\
63.6	-0.0682697385274764\\
64.2	-0.0641130640916288\\
64.8	-0.0662628150780336\\
65.4	-0.0641421449918198\\
66	-0.0570089335666069\\
66.6	-0.0572234276274344\\
67.2	-0.0554462120205415\\
67.8	-0.0553642030026822\\
68.4	-0.0601853944216373\\
69	-0.0595814643650954\\
69.6	-0.0552976844066109\\
70.2	-0.0592329407494674\\
70.8	-0.0542418411201311\\
71.4	-0.0541380387228116\\
72	-0.0617036301968028\\
72.6	-0.0648604626941562\\
73.2	-0.0626241425946716\\
73.8	-0.0631449193354854\\
74.4	-0.06498264870811\\
75	-0.0588752146441891\\
75.6	-0.0589703396461232\\
76.2	-0.0565280005849047\\
76.8	-0.0606879864783991\\
77.4	-0.0525301290574305\\
78	-0.0564140498402423\\
78.6	-0.068491538905505\\
79.2	-0.0701925059611881\\
79.8	-0.070864358955788\\
80.4	-0.0759721965617075\\
81	-0.0795391073316978\\
81.6	-0.070362559672617\\
82.2	-0.0694650579382596\\
82.8	-0.0798587637856087\\
83.4	-0.0865319002977856\\
84	-0.0883757705128687\\
84.6	-0.079428535939428\\
85.2	-0.0701949147269829\\
85.8	-0.0652540540991708\\
86.4	-0.0608888340226486\\
87	-0.054723136742922\\
87.6	-0.0620657108870546\\
88.2	-0.0676044761272294\\
88.8	-0.069375990602623\\
89.4	-0.0721768601001713\\
90	-0.08413537509756\\
90.6	-0.0789273719490211\\
91.2	-0.0697041442408306\\
91.8	-0.0670343067629631\\
92.4	-0.0651291067261752\\
93	-0.0647862782436732\\
93.6	-0.0658707012071318\\
94.2	-0.0700917740021225\\
94.8	-0.0633939800603969\\
95.4	-0.0661055465888865\\
96	-0.0699034277232329\\
96.6	-0.085774466445514\\
97.2	-0.0855500553622619\\
97.8	-0.0864675579308312\\
98.4	-0.0878644491422165\\
99	-0.0839081808998401\\
99.6	-0.0937722485908239\\
100.2	-0.0907752516348808\\
100.8	-0.0852623364573617\\
101.4	-0.0752207593090191\\
102	-0.0714132563197722\\
102.6	-0.0639053720941263\\
103.2	-0.0533868010102196\\
103.8	-0.0470795154892016\\
104.4	-0.0525299996524495\\
105	-0.0574226564272351\\
105.6	-0.0546329498893954\\
106.2	-0.0497981081775265\\
106.8	-0.0455321331195648\\
107.4	-0.0537755068814647\\
108	-0.0514611773230342\\
108.6	-0.0562339052218306\\
109.2	-0.0527141195948104\\
109.8	-0.0564382820384139\\
110.4	-0.0533181551580526\\
111	-0.0530910840159322\\
111.6	-0.048498393396045\\
112.2	-0.0485849289279372\\
112.8	-0.0604273848287846\\
113.4	-0.0619862076361528\\
114	-0.0657679607547834\\
114.6	-0.0726337571450699\\
115.2	-0.0743543523263099\\
115.8	-0.0697789662597472\\
116.4	-0.0659408738679872\\
117	-0.0680223855527231\\
117.6	-0.0624969233401202\\
118.2	-0.0551988685606943\\
118.8	-0.0483194811037038\\
119.4	-0.0462135284862484\\
120	-0.0464920763034513\\
120.6	-0.0277844040404318\\
121.2	-0.0155216130164604\\
121.8	-0.0202340557332548\\
122.4	-0.0300628026547827\\
123	-0.0307599460338395\\
123.6	-0.0247623633786418\\
124.2	-0.0161795336346131\\
124.8	-0.016905130505279\\
126	-0.0203677224763794\\
126.6	-0.0204503517131513\\
127.2	-0.0173214228227039\\
127.8	-0.00899411890736701\\
128.4	-0.00139355978629396\\
129	-0.00444635004902239\\
129.6	-0.00724210920003543\\
130.2	-0.00466193501432599\\
130.8	-0.00745717886030661\\
131.4	-0.00624877234804444\\
132	-0.0004447298001935\\
132.6	-0.00770886415918426\\
133.2	-0.0101451118796092\\
133.8	-0.013084856382477\\
134.4	-0.0183792467438479\\
135	-0.0264596126050662\\
135.6	-0.0224962503580457\\
136.2	-0.0257554471056949\\
136.8	-0.0238639628522037\\
137.4	-0.0205611304654099\\
138	-0.0189933162722866\\
138.6	-0.0191902036745546\\
139.2	-0.017380570790408\\
139.8	-0.019571860628588\\
140.4	-0.0148174912942523\\
};
\addplot [color=black, dotted, forget plot]
  table[row sep=crcr]{%
0	0.181426947281579\\
140	0.181426947281579\\
};
\addplot [color=black, dotted, forget plot]
  table[row sep=crcr]{%
0	-0.181426947281579\\
140	-0.181426947281579\\
};
\addplot [color=red, line width=0.7pt, forget plot]
  table[row sep=crcr]{%
0	0.0576794757140533\\
0.599999999999994	0.0819282565841775\\
1.19999999999999	0.118072094618043\\
1.80000000000001	0.148221643732229\\
2.40000000000001	0.148776691849434\\
3	0.0963162657689907\\
3.59999999999999	0.0841355387725287\\
4.19999999999999	0.093522517508319\\
4.80000000000001	0.0962453852301053\\
5.40000000000001	0.0870571382970127\\
6	0.0888061664761324\\
6.59999999999999	0.10520193987719\\
7.19999999999999	0.140315070361112\\
7.80000000000001	0.154584249874375\\
8.40000000000001	0.13157427219457\\
9	0.132770903691494\\
9.59999999999999	0.130085955502636\\
10.2	0.138088003356614\\
10.8	0.147096310785145\\
11.4	0.178254309333425\\
12	0.178255250037665\\
12.6	0.17105596302514\\
13.2	0.139950202047856\\
13.8	0.079303715630374\\
14.4	0.0217853668424652\\
15	    0.0305743744850986\\
15.6	-0.00734710156112328\\
16.2	-0.0217877342203792\\
16.8	-0.0306530128822089\\
17.4	-0.0163751832571393\\
18	    -0.0164287708924462\\
18.6	-0.0407050642662546\\
19.2	-0.064411795064002\\
19.8	-0.0965923757892142\\
20.4	-0.0625587323886509\\
21	    -0.0377155525038495\\
21.6	-0.0168182649457265\\
22.2	-0.0771837163119642\\
22.8	-0.080911684930074\\
23.4	-0.0511337665980159\\
24	    -0.0474495532953938\\
24.6	-0.0650055410699019\\
25.2	-0.0866585745470809\\
25.8	-0.110570870080039\\
26.4	-0.174241144487496\\
27	-0.177733929994957\\
28.2	-0.1201857899172933\\
28.8	 -0.0847826409617551\\
29.4	 -0.0751672094459457\\
30	    -0.072118872675378\\
30.6	-0.0830696654992175\\
31.2	-0.107000409997937\\
31.8	-0.17067145250914\\
32.4	-0.181317387586972\\
34.8	-0.181341288337279\\
35.4	-0.164806604965008\\
36	-0.117182773997058\\
36.6	-0.0935582033812148\\
37.2	-0.0713202342194825\\
37.8	-0.0790076966529796\\
38.4	-0.083090975684722\\
39	    -0.1066260236427161\\
39.6	-0.12165873394817595\\
40.2	-0.1515456578832811\\
41.4	-0.179085868375154\\
42	-0.181323492611057\\
42.6	-0.166601631318912\\
43.2	-0.126523838604186\\
43.8	-0.101956941533729\\
44.4	-0.121521655024424\\
45	-0.134910973213977\\
45.6	-0.157328341138395\\
46.2	-0.165154157875151\\
46.8	-0.155991580986324\\
47.4	-0.137139611581887\\
48	-0.150703765128839\\
48.6	-0.181339536452271\\
52.2	-0.181349145973655\\
52.8	-0.170884476688627\\
53.4	-0.132584282283375\\
54	-0.0959322975825216\\
54.6	-0.0842252858422512\\
55.2	-0.0665691946956031\\
55.8	-0.0616432998780567\\
56.4	-0.0464333094548408\\
57	-0.0363759134775137\\
57.6	-0.0295661189309158\\
58.2	-0.0214959746069212\\
58.8	-0.00737325174910097\\
59.4	-0.0101110659133781\\
60	-0.00176901404634577\\
60.6	0.0213037243576935\\
61.2	0.0292926647304057\\
61.8	0.0279102456632927\\
62.4	0.0285154062170818\\
63	0.0299013230353467\\
63.6	0.0353497543596291\\
64.2	0.0354393403934807\\
64.8	0.0334632496449956\\
65.4	0.0290611987377929\\
66	0.0229298735200416\\
66.6	0.0140868922391917\\
67.2	0.00448402390489377\\
67.8	0.00324189244145145\\
68.4	-0.00140311490795852\\
69	-0.00693168770942521\\
69.6	-0.00931227981490679\\
70.2	-0.0207891559822144\\
70.8	-0.0180283705551005\\
71.4	-0.0252073642541575\\
72	-0.0291355911286075\\
72.6	-0.0363891593271717\\
73.2	-0.0454224570866302\\
73.8	-0.0497785500743646\\
74.4	-0.0499642653436183\\
75	-0.0553135347083185\\
75.6	-0.0656729439164963\\
76.2	-0.079820856716168\\
76.8	-0.0824298089071078\\
77.4	-0.0840876127824686\\
78	-0.0826375003291844\\
78.6	-0.0665850028059651\\
79.2	-0.0525188817411788\\
79.8	-0.0448026160711663\\
80.4	-0.0462630115607112\\
81	-0.0650605015022165\\
81.6	-0.0751100671684242\\
82.2	-0.0805915814628975\\
82.8	-0.081172537751911\\
83.4	-0.0793500335193755\\
84	-0.0665299447864811\\
84.6	-0.0524617866087738\\
85.2	-0.0354163138776471\\
85.8	-0.0491994425456994\\
86.4	-0.0723158308203722\\
87	-0.0755587087412835\\
87.6	-0.0816590545118174\\
88.2	-0.0823960705053537\\
88.8	-0.090743788292059\\
89.4	-0.0902189983889627\\
90	-0.0869091247867004\\
90.6	-0.0818685541492243\\
91.2	-0.089728753127531\\
91.8	-0.093563534783101\\
92.4	-0.106710070595398\\
93	-0.106469263991102\\
93.6	-0.103273663393367\\
94.2	-0.0980568971864386\\
94.8	-0.0858189588139737\\
95.4	-0.0732053719854946\\
96	-0.0721770524962722\\
96.6	-0.0818276213273919\\
97.2	-0.0819014286907702\\
97.8	-0.0809388167290308\\
99	-0.058109210201053\\
99.6	-0.0489474955587355\\
100.2	-0.051122451937573\\
100.8	-0.0527718927114904\\
101.4	-0.0593544197005258\\
102	-0.0586280828315751\\
102.6	-0.0549310968112025\\
103.2	-0.0465363701685533\\
103.8	-0.040834813791804\\
104.4	-0.0269413839080528\\
105	-0.0368846854774176\\
105.6	-0.0434061601573603\\
106.2	-0.0467672648440782\\
106.8	-0.050836480397237\\
107.4	-0.0570721417778941\\
108	-0.0609374082589511\\
108.6	-0.0578176071404073\\
109.2	-0.0558319069532445\\
109.8	-0.0488549096374697\\
110.4	-0.0605878215099551\\
111	-0.0550153689869433\\
111.6	-0.0626570891724612\\
112.2	-0.0611283422870486\\
112.8	-0.0614559404603483\\
113.4	-0.0628539250477616\\
114	-0.0593246584892029\\
114.6	-0.0683634877331372\\
115.2	-0.066939652390289\\
115.8	-0.0575489019491613\\
116.4	-0.0536273495605712\\
117	-0.0454327969068231\\
117.6	-0.0469475094575671\\
118.2	-0.0499873867855456\\
118.8	-0.0476607800313502\\
119.4	-0.0403619701881723\\
120	-0.0242760871427379\\
120.6	-0.0257385287347915\\
121.2	-0.0198109906674233\\
121.8	-0.0186840769444814\\
122.4	-0.0162481580723863\\
123	-0.0113258606878048\\
123.6	-0.0134536054808052\\
124.2	-0.0165790689494827\\
124.8	-0.0237116433166307\\
125.4	-0.0260226892764308\\
126	-0.0204643138952179\\
126.6	-0.0143800743551026\\
127.2	-0.0156001265992245\\
127.8	-0.0163904822843222\\
128.4	-0.0167582498924901\\
129	-0.0189969744255905\\
129.6	-0.0171969414862474\\
130.2	-0.00948221168462737\\
130.8	2.33985275315263e-05\\
131.4	-0.000682150916674118\\
132	0.00598747693669566\\
132.6	0.00550758830348741\\
133.2	0.00426252144006867\\
133.8	-1.96169534660839e-05\\
134.4	-0.00379129788416321\\
135	-0.00396599900477668\\
135.6	-0.00827378519292665\\
136.2	-0.0110873599580543\\
136.8	-0.012357962907231\\
137.4	-0.0176359442358773\\
138	-0.01882181200898\\
138.6	-0.0160230367597194\\
139.2	-0.0178066739564144\\
139.8	-0.0245285169506531\\
};
\addplot [color=green, line width=0.7pt, forget plot]
  table[row sep=crcr]{%
0	-0.181219972981182\\
0.599999999999994	-0.181224393676331\\
1.19999999999999	-0.174125755544679\\
1.80000000000001	-0.157998797873375\\
2.40000000000001	-0.15104131224598\\
3	-0.15104131224598\\
3.59999999999999	-0.0783140311572197\\
4.19999999999999	-0.0459441733253527\\
4.80000000000001	-0.0130557692164075\\
5.40000000000001	0.0231744000743674\\
6	0.0603007223883196\\
6.59999999999999	0.0937862095312028\\
7.19999999999999	0.129656792203036\\
7.80000000000001	0.148483925375359\\
8.40000000000001	0.15686545157709\\
9	0.156167006882072\\
9.59999999999999	0.133865027236681\\
10.2	0.146048754040663\\
10.8	0.15919543757434\\
11.4	0.176415182962302\\
12	0.163826815724264\\
12.6	0.157401581728607\\
13.2	0.153861480773344\\
13.8	0.153907211225857\\
14.4	0.114170832999804\\
15	    0.0441449055587952\\
15.6	-0.00432260345391455\\
16.2	0.0127304239939576\\
16.8	-0.0286868922644317\\
17.4	-0.0646780924739403\\
18	    -0.109588672242154\\
18.6	-0.060001790452219\\
19.2	-0.0566015912780892\\
19.8	-0.0496579910121113\\
20.4	-0.0812226611288338\\
21	    -0.104136226930109\\
21.6	-0.1023814014728998\\
22.2	-0.1166223629204126\\
22.8	-0.11574638450233\\
23.4	-0.1030965195949977\\
24	    -0.0976785530127586\\
24.6	-0.0907792401793301\\
25.2	-0.112169596475119\\
25.8	-0.136437734210347\\
26.4	-0.0727115925764963\\
27	    -0.0629152082064866\\
27.6	-0.0617453265435756\\
28.2	-0.0745736016992851\\
28.8	-0.138303931130849\\
29.4	-0.181326250333854\\
30	    -0.181317404930553\\
30.6	-0.132085967843807\\
31.2	-0.113374022555654\\
31.8	-0.0838248609538698\\
32.4	-0.071969895453792\\
33	    -0.07331191846\\
33.6	-0.0982211690301904\\
34.2	-0.122030414724236\\
34.8	-0.181324810763556\\
37.2	-0.181325018893602\\
37.8	-0.171770858350328\\
38.4	-0.131293830063555\\
39	    -0.11675366640406594\\
39.6	-0.10300857446180203\\
40.2	-0.0838461395735135\\
40.8	-0.0711333080080658\\
41.4	-0.1043432702823065\\
42	    -0.128376769084183\\
42.6	-0.162053284288191\\
43.2	-0.181317439617715\\
45.6	-0.181326614563972\\
46.2	-0.160192468062974\\
46.8	-0.121519537945915\\
47.4	-0.0759638055079108\\
48	-0.0613682269407423\\
48.6	-0.0520582544309968\\
49.2	-0.0469684724915567\\
49.8	-0.0360083939597473\\
50.4	-0.03552297460586715\\
51	-0.0320489125325309\\
51.6	-0.0480466385829175\\
52.2	-0.0732211374603082\\
52.8	-0.0937341907192604\\
53.4	-0.0802647133735377\\
54	-0.0631030982593188\\
54.6	-0.075036139802279\\
55.2	-0.107249997793929\\
55.8	-0.108023439113566\\
56.4	-0.0789995021434606\\
57	-0.0752901474334351\\
57.6	-0.0942927146862189\\
58.2	-0.117210255787114\\
58.8	-0.139193223680991\\
59.4	-0.138899156305968\\
60	-0.144802397059635\\
60.6	-0.161874591870259\\
61.2	-0.181327655224806\\
68.4	-0.181363456961009\\
69	-0.17939589046199\\
69.6	-0.14434501117799\\
70.2	-0.0943822900808016\\
70.8	-0.0494115231393835\\
71.4	-0.00782261091387682\\
72	0.0279374987459562\\
72.6	0.0592156773700765\\
73.2	0.0770343094549162\\
73.8	0.0833337430831023\\
74.4	0.0809349038420351\\
75	0.0849014559954071\\
75.6	0.0918989750835806\\
76.2	0.0847651883121898\\
76.8	0.0714134309178007\\
77.4	0.0629646921777862\\
78	0.0507301368684523\\
78.6	0.0298579075432315\\
79.2	0.0195276306303356\\
79.8	0.00864872802708305\\
80.4	-0.000690920584588639\\
81	-0.00658295615428983\\
81.6	-0.0182376460129206\\
82.2	-0.0283369648366545\\
82.8	-0.0316441649221986\\
83.4	-0.0290140295942649\\
84	-0.0382103699994332\\
84.6	-0.0433096415649175\\
85.2	-0.0621954246249459\\
85.8	-0.0750547884685489\\
86.4	-0.0768912660057595\\
87	-0.0678910313744723\\
87.6	-0.0690275698646019\\
88.2	-0.0770354976076817\\
88.8	-0.0910751019304143\\
89.4	-0.102219221706264\\
90	-0.105471494949995\\
90.6	-0.115758188608197\\
91.2	-0.11208447308374\\
91.8	-0.10483657280227\\
92.4	-0.0837154128052475\\
93	-0.0712149000854936\\
93.6	-0.0547466221706259\\
94.2	-0.0439511574036544\\
94.8	-0.0257879938347969\\
95.4	-0.0172848800130225\\
96	-0.0122964203268907\\
96.6	-0.00919029078505673\\
97.2	-0.0170714186811836\\
97.8	-0.00979498490855235\\
98.4	-0.00881706184614472\\
99	-0.00822179820605129\\
99.6	-0.00603345022392432\\
100.2	-0.00981296691904276\\
100.8	-0.0147972004433257\\
101.4	-0.0142102040681493\\
102	-0.0174727336333831\\
102.6	-0.0134370009824067\\
103.2	-0.00405250985019734\\
103.8	-0.00340622619361852\\
104.4	-0.00480290989153787\\
105	-0.00364869107337995\\
106.2	-0.0077560774086578\\
107.4	-0.00839933598331299\\
108	-0.0125488171491668\\
108.6	-0.0145347265617488\\
109.2	-0.00861870529104181\\
109.8	-0.00686344190927457\\
110.4	-0.00532881197142387\\
111	-0.002412713547983\\
111.6	-0.00154555400399659\\
112.2	-0.00438712880625758\\
112.8	-0.0053753047778855\\
113.4	-0.00937133373096799\\
114	-0.0114259402375865\\
114.6	-0.0138886414613069\\
115.2	-0.0149465956467907\\
115.8	-0.0184801148826068\\
116.4	-0.0197891016427718\\
117	-0.0198109996520941\\
117.6	-0.0166886150053358\\
118.2	-0.0142559959952564\\
119.4	-0.0197355731582434\\
120	-0.0187755100633069\\
120.6	-0.0110620482448383\\
121.2	-0.00294111174841305\\
121.8	-0.000942418650851096\\
122.4	0.00526985519763912\\
123.6	0.0115647561246135\\
124.2	0.0193175586341852\\
124.8	0.0239224446750086\\
125.4	0.0267858120097912\\
126	0.0322627872598389\\
126.6	0.0345130386430696\\
127.2	0.0371611861443455\\
127.8	0.0363122555326925\\
128.4	0.0276162213156113\\
129	0.0265501428107768\\
129.6	0.0319443530005969\\
130.2	0.0288456522490605\\
130.8	0.0280036836338127\\
131.4	0.0216860746720897\\
132	0.021228057201256\\
132.6	0.0347799394052402\\
133.2	0.0380267097198157\\
133.8	0.0426136682046376\\
134.4	0.0391358836499478\\
135	0.0368976949174566\\
135.6	0.034026024545625\\
136.2	0.0396507231957912\\
136.8	0.0405848066353371\\
137.4	0.0454333015163968\\
138	0.0496164958930194\\
138.6	0.056131072370647\\
139.8	0.0540389709223916\\
140.4	0.0499626023149631\\
};
\addplot [color=blue, dashed, line width=0.7pt, forget plot]
  table[row sep=crcr]{%
0	-0.000213546949648702\\
1.8	-0.000213546949648702\\
2.4	-0.180\\
10.8	-0.180\\
11.4	-0.0821897045731923\\
12	0.113788242112573\\
12.6	0.0895069953791818\\
13.2	0.0882497907448254\\
13.8	0.0874947528700147\\
14.4	0.0856907018130961\\
15	0.0827841043272883\\
15.6	0.0787089115528126\\
16.8	0.06902175160468\\
17.4	0.0650434612833131\\
};
\addplot [color=green, dashed, line width=0.7pt, forget plot]
  table[row sep=crcr]{%
0	-0.180\\
11.4	-0.180\\
};
\addplot [color=red, dashed, line width=0.7pt, forget plot]
  table[row sep=crcr]{%
0	-0.180\\
10.8	-0.180\\
};
\end{axis}
\end{tikzpicture}%

%% file: heading_error_eight.tex
%
%
\definecolor{mycolor1}{rgb}{0.50000,0.00000,0.50000}%
\begin{tikzpicture}

\begin{axis}[%
width=\figurewidth,
height=\figureheight,
at={(0\figurewidth,0\figureheight)},
scale only axis,
xmin=0,
xmax=140,
xtick={0,30,60,90,120},
xlabel={$t$ [s]},
ymin=-1,
ymax=1,
xlabel style={font=\color{white!15!black},at={(axis description cs:0.5,-0.17)},anchor=north},
ylabel style={font=\color{white!15!black},at={(axis description cs:-0.25,.5)},anchor=south},
ylabel={$\tilde \theta_3$ [rad]},
axis background/.style={fill=white},
xmajorgrids,
ymajorgrids
]
\addplot [color=blue,line width=0.7pt, forget plot]
  table[row sep=crcr]{%
0	0.0325880733132351\\
0.900000000000006	0.0174706701934326\\
2.40000000000001	-0.00253824397921676\\
3.30000000000001	-0.0125048856437218\\
3.59999999999999	-0.0155232341215026\\
4.80000000000001	-0.0238205652870249\\
5.40000000000001	-0.0263932954706263\\
6	-0.0285701320692908\\
6.90000000000001	-0.0282534764567401\\
8.09999999999999	-0.0233891678228986\\
8.40000000000001	-0.021549606435002\\
9.30000000000001	-0.0129340986907494\\
9.90000000000001	-0.00402910143137092\\
10.2	0.0023538705706585\\
10.5	0.00738365694880372\\
10.8	0.0139539588987816\\
11.4	0.0248087172210205\\
11.7	0.0323330995440472\\
12.3	0.0434165161848057\\
12.6	0.0507517617940891\\
13.2	0.0593384367227543\\
13.5	0.0653494310379017\\
13.8	0.06997048228979\\
14.1	0.0763076853752125\\
14.7	0.0853979077935207\\
15	0.091766723692416\\
15.3	0.096119092106818\\
15.6	0.102605120539664\\
16.2	0.11120978116989\\
16.5	0.117411540150641\\
16.8	0.121715397834777\\
17.1	0.127966658473014\\
17.7	0.135898784995078\\
18	0.141460211277007\\
18.3	0.144721270799636\\
18.6	0.150024683475493\\
19.2	0.155096934437751\\
19.5	0.15960116803646\\
19.8	0.161830649971961\\
20.1	0.16623790681362\\
20.4	0.168182255029677\\
20.7	0.172469781041144\\
21.3	0.175650880932807\\
21.6	0.179332957863807\\
21.9	0.180341811180114\\
22.2	0.183561728596686\\
22.5	0.184197620749472\\
22.8	0.186909453272818\\
23.4	0.186823458671569\\
23.7	0.188678236007689\\
24	0.187739865183829\\
24.3	0.18898229420185\\
24.6	0.187449545860289\\
24.9	0.187963322997092\\
25.2	0.185958595275878\\
25.5	0.186287821531295\\
25.8	0.183939070701598\\
26.1	0.183327288627623\\
26.4	0.180353419184684\\
26.7	0.179167852997779\\
27	0.175582156777381\\
27.3	0.17452645421028\\
27.6	0.170494229793547\\
27.9	0.168590993881224\\
28.2	0.16399732351303\\
28.5	0.162064806818961\\
28.8	0.157748923897742\\
29.1	0.155154855251311\\
29.4	0.150515483021735\\
29.7	0.148056761622428\\
30	0.143477798700332\\
30.3	0.140911222696303\\
30.6	0.136318148374556\\
30.9	0.133918763399123\\
31.2	0.129317642450332\\
31.5	0.126757876276969\\
31.8	0.122397781610488\\
32.1	0.119793341755866\\
32.4	0.115443871617316\\
33	0.110302628278731\\
33.3	0.105962590575217\\
33.6	0.103540660142897\\
33.9	0.0991027516126621\\
34.2	0.0966999691724766\\
34.5	0.0926102933287609\\
34.8	0.0901687678694714\\
35.1	0.0860725504159916\\
35.4	0.0836180460453022\\
35.7	0.0798269224166859\\
36	0.0772766259312618\\
36.3	0.0731764969229687\\
36.6	0.070572146475314\\
36.9	0.0667862084507931\\
37.5	0.061932676732539\\
37.8	0.0580019888281811\\
38.1	0.0559995782375324\\
38.4	0.0524771597981442\\
38.7	0.0506839987635601\\
39	0.0472731080651272\\
39.3	0.0457450684905041\\
39.6	0.0425947469472874\\
39.9	0.0412797865271557\\
40.2	0.0384191063046444\\
40.5	0.0374796998500813\\
40.8	0.0351025474071491\\
41.1	0.0344753053784359\\
41.4	0.0323333230614651\\
41.7	0.0322244253754604\\
42	0.030633212327956\\
42.3	0.0306767237186421\\
42.6	0.0293503303825844\\
42.9	0.0300479319691647\\
43.2	0.0289784421026695\\
43.5	0.0299583871662605\\
44.1	0.0292930205166329\\
44.4	0.0309133504331101\\
44.7	0.0307587359845627\\
45	0.0327306327223766\\
45.3	0.0326235008239735\\
45.6	0.0346762996911991\\
45.9	0.0345558089017857\\
46.2	0.0368181926012028\\
46.5	0.0370527741312969\\
46.8	0.0392084059119213\\
47.1	0.0397263032197941\\
47.4	0.0422598358988751\\
48	0.0430637088417996\\
48.3	0.0455658522248257\\
48.6	0.0457584273815144\\
48.9	0.0480385732650745\\
49.5	0.0485434842109669\\
49.8	0.050421671271323\\
50.1	0.0503038626909245\\
50.4	0.0520026919245709\\
51	0.0514142823219288\\
51.3	0.0530229446291912\\
51.6	0.0526864391565312\\
51.9	0.0544650119543064\\
52.2	0.0542595621943462\\
52.5	0.0558458504080761\\
53.1	0.0556543853878964\\
53.4	0.0575990837812412\\
53.7	0.0576424238085735\\
54	0.0595698964595783\\
54.3	0.0594400700926769\\
54.6	0.0610394266247738\\
55.2	0.060122878551482\\
55.5	0.0613533121347416\\
55.8	0.0604866904020298\\
56.1	0.0614393514394749\\
56.4	0.0605304849147785\\
56.7	0.061285966336726\\
57	0.0602223289012898\\
57.3	0.0608423885703075\\
57.9	0.0583386582136143\\
58.2	0.0586253342032421\\
58.5	0.057277404963969\\
58.8	0.0574571725726116\\
59.1	0.0560336422920216\\
59.4	0.0560940665006626\\
59.7	0.0544650045037258\\
60	0.054225647151469\\
60.3	0.0523772847652424\\
60.6	0.0520920914411533\\
60.9	0.0502570730447758\\
61.2	0.0502188590168942\\
61.5	0.0486678120493877\\
61.8	0.0487032321095455\\
62.1	0.0472672817111004\\
62.4	0.0472854536771763\\
62.7	0.045978361070155\\
63	0.0459156495332707\\
63.3	0.0445270699262608\\
63.6	0.0444261518120754\\
64.2	0.0417560127377499\\
64.5	0.041859911084174\\
64.8	0.0404111158847797\\
65.1	0.0403892409801472\\
65.4	0.0390223649144161\\
65.7	0.0391754966974247\\
66	0.038001627922057\\
66.3	0.0382108774781216\\
66.6	0.0371512633562077\\
66.9	0.0374749389290798\\
67.2	0.0363969740271557\\
67.5	0.0364759948849667\\
67.8	0.0355321404337872\\
68.1	0.0356808614730824\\
68.4	0.0347457093000401\\
69.3	0.0345759627223003\\
69.9	0.0327912878990162\\
70.2	0.0330063116550434\\
70.5	0.0323006671667088\\
70.8	0.032711291015147\\
71.1	0.0322315183281887\\
71.4	0.0327019330859173\\
71.7	0.0320738677680481\\
72	0.0325599473714817\\
72.3	0.0319866847991932\\
72.6	0.0324619796872128\\
72.9	0.0318991963565338\\
73.2	0.032408935278653\\
73.8	0.0312993240356434\\
74.1	0.0317864727973927\\
74.4	0.0309869398176659\\
74.7	0.0312864494323719\\
75	0.0308064271509636\\
75.9	0.0314507794380177\\
76.8	0.0308463697135437\\
77.7	0.0318267506360996\\
82.8	0.0268710193037975\\
84.6	0.0230103370547283\\
86.1	0.0221491616964329\\
87	0.0217077222466457\\
88.8	0.0209881451725948\\
89.4	0.0201506961882103\\
90.9	0.0177420836687077\\
96.3	0.0185970936715592\\
97.8	0.016165496110915\\
98.7	0.0143011932075012\\
99.6	0.0110329289734352\\
101.7	0.00226002305745965\\
103.8	-0.0035390562564146\\
105	-0.00585055790841693\\
106.2	-0.00621556185186023\\
108	-0.00483517415821666\\
109.2	-0.00406278923153991\\
110.1	-0.00307020984590167\\
113.4	0.00159703314304238\\
114.9	0.000752347782253082\\
117.3	-0.00537532448768729\\
118.5	-0.00914188914001102\\
119.4	-0.0117387592047464\\
120.3	-0.0142002335935842\\
121.8	-0.0142603821307432\\
123.6	-0.0142607062309992\\
124.8	-0.0129768799245369\\
126.6	-0.0104852455854427\\
127.5	-0.00843531891703719\\
129.9	-0.000748290419579689\\
131.4	0.00656767681240922\\
132.6	0.0113148961961258\\
134.4	0.0160386276245106\\
137.4	0.0187098172307003\\
139.5	0.015747824013232\\
140.1	0.0143999804556358\\
};
\addplot [color=red,line width=0.7pt, forget plot]
  table[row sep=crcr]{%
0	-0.46318107843399\\
1	-0.4702688827754\\
3	-0.502834614410119\\
4	-0.517106194780098\\
5	-0.527499647728888\\
6	-0.533626255439145\\
7	-0.534303927812886\\
8	-0.528854107674846\\
9	-0.518147651149064\\
10	-0.500803231563083\\
11	-0.47879697848407\\
12	-0.452815718847688\\
15	-0.367095761530408\\
17	-0.308261312948531\\
18	-0.280002991457081\\
19	-0.253717254618465\\
23	-0.128015394400933\\
25	-0.0558612816818709\\
26	-0.0215374829470818\\
27	0.00857189936348846\\
28	0.0426364249451865\\
29	0.0803770430712234\\
30	0.120007271238535\\
32	0.185972103950604\\
33	0.216192649847756\\
34	0.242636239996898\\
35	0.271234534759799\\
36	0.307909449105239\\
37	0.347061011935978\\
38	0.389258295365522\\
39	0.426601213960765\\
40	0.466465547389021\\
41	0.494210809831003\\
43	0.537323203544076\\
44	0.555424817486454\\
45	0.569035993042689\\
46	0.578511077830058\\
47	0.585420860900825\\
48	0.587398526690436\\
49	0.582480739265407\\
50	0.572729340298338\\
51	0.556882364471704\\
52	0.535334197627975\\
53	0.509332388183594\\
54	0.479199985135864\\
55	0.444935624371226\\
56	0.407370622233714\\
58	0.328623349539726\\
59	0.290624225685463\\
60	0.254411334820901\\
61	0.221772228614498\\
62	0.193223197160478\\
63	0.1687387677039\\
64	0.149448658534652\\
65	0.135457399522721\\
66	0.124690987966488\\
67	0.116863716156246\\
68	0.110863324314181\\
70	0.104150374988336\\
72	0.100945879557827\\
76	0.0962103859732224\\
78	0.0897093153104436\\
79	0.0865378620394495\\
80	0.0846626880115195\\
81	0.0840369254084408\\
83	0.0796697086120162\\
84	0.0771684089111773\\
85	0.0758309541803897\\
86	0.0766663661591735\\
87	0.0755184689234341\\
88	0.0735745796133926\\
92	0.0611011391758325\\
93	0.0558232983573816\\
95	0.0419142529424619\\
97	0.0267444423864163\\
99	0.0106742824794992\\
100	0.00392067815067776\\
101	-0.00154408194200073\\
102	-0.00604812043724223\\
103	-0.0097424122923826\\
104	-0.0122452823678429\\
105	-0.0130701085419389\\
106	-0.0126448757250728\\
112	-0.00377827026858313\\
113	-0.00320432336903309\\
115	-0.00452421438376405\\
119	-0.00787280669982238\\
121	-0.00971644365012025\\
123	-0.00951018537503501\\
126	-0.00557522758876416\\
131	0.00328685823518526\\
133	0.00987118583995539\\
137	0.0239434174926316\\
139	0.0275679671329669\\
140	0.0275480191177166\\
};
\addplot [color=green,line width=0.7pt, forget plot]
  table[row sep=crcr]{%
0	-0.806800484657288\\
1	-0.806216216791114\\
3	-0.818428089723483\\
4	-0.828339254660449\\
5	-0.843513161395208\\
6	-0.850698454356603\\
7	-0.852293082359097\\
8	-0.850086165590568\\
9	-0.841314121391946\\
10	-0.826535338991647\\
11	-0.809506020210108\\
12	-0.785853933230158\\
14	-0.729225281931349\\
15	-0.696591908407356\\
16	-0.665608940726344\\
17	-0.638344556954394\\
18	-0.604218789504642\\
19	-0.572040176041781\\
20	-0.542769437771057\\
21	-0.506619729583264\\
22	-0.473292723510639\\
23	-0.444958909743775\\
24	-0.408514281426818\\
25	-0.374003588201845\\
26	-0.344974559037411\\
27	-0.307782370525189\\
28	-0.271523766609135\\
29	-0.240860632507861\\
32	-0.138826937437301\\
33	-0.102930943674153\\
34	-0.0650139953858684\\
36	-0.0078667546733584\\
37	0.0211661555151466\\
38	0.0467602226030408\\
39	0.0754638566048129\\
40	0.109917530953282\\
41	0.142440404199078\\
42	0.169834562967878\\
43	0.200366832645415\\
44	0.224530880260829\\
45	0.245401593595403\\
46	0.27235967444858\\
48	0.324987256312426\\
49	0.352952434851318\\
50	0.386649608554563\\
51	0.419202403318963\\
52	0.447239907394135\\
53	0.470362086117206\\
54	0.496629493128609\\
55	0.51943287114338\\
57	0.553790061594214\\
58	0.570069006759894\\
59	0.582370290602228\\
60	0.58991416209426\\
61	0.592394809127342\\
62	0.588930033968438\\
63	0.580147224520573\\
64	0.566755274734646\\
65	0.547736630385884\\
66	0.521419408911385\\
67	0.488678367362212\\
68	0.450993377176445\\
69	0.408339732610074\\
71	0.316844298730331\\
72	0.273132490003661\\
73	0.234819396881079\\
74	0.203313394900619\\
75	0.179808013751511\\
76	0.16319468672512\\
77	0.151747936303337\\
78	0.144693588046778\\
79	0.140886472652426\\
80	0.139559636203074\\
82	0.140497915180333\\
85	0.143778256640474\\
86	0.143009083781408\\
87	0.141398057140407\\
88	0.138790674750879\\
89	0.133984175392754\\
90	0.126312727511021\\
91	0.11554991984579\\
92	0.102364554325447\\
94	0.074624028622992\\
95	0.0623520446148973\\
96	0.0521656329376867\\
97	0.0434884084325802\\
98	0.0362639073395883\\
99	0.0307299704467425\\
100	0.0268910079870466\\
101	0.0240192254432543\\
103	0.0209274583168053\\
104	0.0206094411349227\\
106	0.0222577941553652\\
111	0.0279929549583358\\
113	0.0311824474896127\\
115	0.0332721759566539\\
117	0.0336442977711329\\
118	0.032913323685392\\
119	0.0311336543945515\\
121	0.0254187512378508\\
124	0.0162647043337358\\
126	0.0114522609189578\\
129	0.00717916276187225\\
131	0.00408303082610928\\
133	-0.000997755866990246\\
135	-0.00579888552530861\\
141	-0.02432317520109\\
};
\addplot [color=blue, dashed,line width=0.7pt, forget plot]
  table[row sep=crcr]{%
0	0.00774223636835814\\
6	0.00673526494391652\\
7	0.00349650888002628\\
8	-0.00298318201083347\\
10	-0.0216141064443427\\
11	-0.028867816855513\\
12	-0.0333344151868857\\
14	-0.0373414654065414\\
17	-0.0403738110133531\\
};
\addplot [color=green, dashed,line width=0.7pt, forget plot]
  table[row sep=crcr]{%
0	-0.814496695995331\\
1	-0.811864267767033\\
2	-0.812529460605324\\
3	-0.81721170148764\\
5	-0.839048698817351\\
6	-0.846761293287491\\
7	-0.851628328028088\\
8	-0.848816330951545\\
9	-0.847166982137054\\
11	-0.847113353044742\\
};
\addplot [color=red, dashed,line width=0.7pt, forget plot]
  table[row sep=crcr]{%
0	-0.461065459251405\\
1	-0.467431692162418\\
2	-0.483179826358457\\
4	-0.52335439405911\\
5	-0.538043557920746\\
6	-0.549027086481104\\
7	-0.55324840994475\\
9	-0.554943627983716\\
10	-0.556199549622589\\
};
\end{axis}
\end{tikzpicture}%

%% file: lateral_error_eight.tex
%
%
\definecolor{mycolor1}{rgb}{0.50000,0.00000,0.50000}%
\begin{tikzpicture}

\begin{axis}[%
width=\figurewidth,
height=\figureheight,
at={(0\figurewidth,0\figureheight)},
scale only axis,
xmin=0,
xmax=140,
xtick={0,30,60,90,120},
xlabel={$t$ [s]},
ymin=-10,
ymax=10,
xlabel style={font=\color{white!15!black},at={(axis description cs:0.5,-0.14)},anchor=north},
ylabel style={font=\color{white!15!black},at={(axis description cs:-0.23,.5)},anchor=south},
ylabel={$\tilde z_3$ [m]},
axis background/.style={fill=white},
xmajorgrids,
ymajorgrids
]
\addplot [color=blue, line width=0.7pt, forget plot]
  table[row sep=crcr]{%
0	2.95311791065689\\
1.19999999999999	2.92410987283782\\
2.69999999999999	2.90262995344449\\
4.5	2.89779692848947\\
8.09999999999999	2.90584034809066\\
9.59999999999999	2.89200288306523\\
10.8	2.85632999410961\\
11.7	2.81903873297094\\
12.9	2.75944147499987\\
14.1	2.68617416537177\\
15	2.61710696202255\\
15.9	2.53583952575212\\
17.4	2.37365383848163\\
19.8	2.07419165873046\\
21	1.90638243936689\\
26.1	1.15090041005939\\
27.6	0.94411175975111\\
28.5	0.830083248913553\\
30	0.66403522320249\\
31.5	0.5202115807692\\
32.4	0.442251896978547\\
34.5	0.291839126092697\\
35.7	0.232461489484621\\
36.9	0.202562620191543\\
39.9	0.20325064524215\\
43.2	0.242745764664591\\
51.3	0.294067165954402\\
53.4	0.283392514653485\\
56.7	0.229995947346254\\
59.4	0.175410875570265\\
63.3	0.140537044320041\\
64.8	0.11795086372797\\
68.1	0.0782770153333558\\
71.1	0.0484960603013462\\
72.6	0.0354230295123159\\
74.1	0.0154706801089901\\
77.1	-0.0423873267198758\\
78.3	-0.0674787971566957\\
79.5	-0.0872781553266577\\
82.2	-0.121859467202171\\
93.6	-0.177927508968821\\
98.7	-0.244536470920423\\
100.2	-0.249026783081263\\
103.2	-0.250395752897873\\
108	-0.278721847222783\\
109.2	-0.274729728804203\\
110.1	-0.27093156460225\\
113.4	-0.277028752858143\\
114.9	-0.25792863323494\\
116.4	-0.218436756112226\\
117.9	-0.171017265075449\\
119.1	-0.118135532462048\\
121.5	0.0265593032909521\\
123.9	0.160592866289619\\
124.8	0.205420116332789\\
126	0.261552497650513\\
127.2	0.301617754159537\\
129	0.335415813708693\\
129.9	0.337498323524926\\
131.4	0.31955936580772\\
138.6	0.166575041072804\\
140.1	0.126746626708467\\
};
\addplot [color=red, line width=0.7pt, forget plot]
  table[row sep=crcr]{%
0	3.43948693687247\\
0.599999999999994	3.58256711149724\\
1.19999999999999	3.7387969112695\\
2.09999999999999	3.98629760582742\\
3.30000000000001	4.33154113394474\\
3.59999999999999	4.41221146380917\\
5.09999999999999	4.85634356306264\\
5.69999999999999	5.03160752564159\\
7.19999999999999	5.48691221600109\\
8.09999999999999	5.74373423850258\\
8.69999999999999	5.90418258491383\\
10.2	6.27725034656328\\
10.8	6.42199455858082\\
11.7	6.6094674705908\\
12.9	6.84310443432145\\
13.5	6.94898517766597\\
15.9	7.32957597463829\\
17.4	7.53366255910004\\
18.9	7.7003136520282\\
19.8	7.78353677541858\\
21.3	7.90054172927543\\
22.8	7.97157700553652\\
24.6	8.01659773798289\\
25.2	8.0212464835729\\
25.8	8.01960863346679\\
26.4	8.0171929400685\\
27.9	7.9821803264123\\
28.5	7.95661010544237\\
29.1	7.92416485608925\\
29.7	7.88639554230096\\
30.3	7.84927115347014\\
32.1	7.7023945336704\\
32.4	7.66989610625137\\
33	7.61693386238883\\
33.6	7.54712715004175\\
34.5	7.42992371794304\\
39.3	6.66916985880562\\
39.9	6.57708203675872\\
41.4	6.31799869880956\\
42	6.19219320883363\\
43.2	5.91569423020823\\
43.5	5.83378419732443\\
44.4	5.6093897227538\\
46.5	5.04396364597122\\
46.8	4.95329167127048\\
50.4	3.75510369547689\\
50.7	3.65196818153174\\
51	3.55924720594965\\
51.6	3.35316618031169\\
52.5	3.04723301120566\\
54.3	2.47558442140482\\
54.9	2.3058517968025\\
55.5	2.1488626269182\\
56.1	2.00991484386881\\
56.7	1.88236904562174\\
57.3	1.76424906327838\\
57.9	1.66054143680751\\
58.8	1.51776121062676\\
59.1	1.4715892738559\\
60.9	1.24985735158677\\
61.5	1.19122004254632\\
63	1.05857487117649\\
64.5	0.94413325472766\\
66	0.842244890231342\\
69.3	0.652065988964893\\
70.5	0.589216801888227\\
72	0.503953800680762\\
73.8	0.412304303341472\\
74.7	0.364945546925014\\
77.1	0.27385919605419\\
81.3	0.119312629450633\\
84.3	0.0182676353674651\\
84.9	-0.0148780773639032\\
85.8	-0.0643560182656415\\
86.7	-0.106420063358826\\
88.5	-0.19314717774057\\
89.7	-0.24459716249973\\
90.9	-0.302457620442851\\
92.1	-0.346834147866161\\
96.9	-0.480611529731704\\
98.4	-0.495588790373347\\
99.6	-0.498899913786858\\
101.1	-0.489780086012445\\
102.6	-0.47572905551371\\
103.8	-0.46551647111886\\
105.9	-0.461100329716658\\
109.8	-0.442102138493709\\
111.9	-0.405468208848305\\
112.5	-0.388254186091501\\
117.3	-0.212917875065727\\
118.2	-0.163337738116212\\
120	-0.0461957727160609\\
120.9	0.0167330135227814\\
122.1	0.0890520463992175\\
124.2	0.199322559485381\\
126	0.28322216475658\\
127.5	0.357429042870166\\
129	0.43522504561426\\
131.1	0.508178352286279\\
132.3	0.535153774168066\\
134.7	0.557983660462128\\
137.7	0.550328974281655\\
139.2	0.547814608934118\\
140.1	0.561773781914411\\
};
\addplot [color=green, line width=0.7pt, forget plot]
  table[row sep=crcr]{%
0	1.21713251938155\\
0.300000000000011	1.23374495193835\\
0.599999999999994	1.25960541021237\\
0.900000000000006	1.29523285053668\\
1.19999999999999	1.33868343791886\\
1.5	1.39054175299469\\
1.80000000000001	1.45818634649351\\
2.09999999999999	1.49726506011322\\
3.30000000000001	1.49726506011322\\
3.59999999999999	1.99380292924104\\
6.30000000000001	3.02557674199937\\
7.19999999999999	3.37942250098186\\
9.30000000000001	4.157531701316\\
9.59999999999999	4.27625139277873\\
10.2	4.48495742061633\\
11.7	4.99547543727104\\
13.5	5.55274456902924\\
14.4	5.80482422205202\\
15	5.96184402671361\\
15.6	6.11353551907177\\
18.3	6.8086379023388\\
19.5	7.10333225533392\\
21	7.42472878201085\\
21.6	7.53788041652311\\
22.2	7.65554681151886\\
24	7.97048196708252\\
24.9	8.11228093984977\\
25.5	8.20187175618346\\
26.1	8.28186891258079\\
26.7	8.36640613050832\\
27.9	8.51425288289261\\
28.5	8.57790860082548\\
29.1	8.63429921484786\\
30	8.71775937341985\\
30.6	8.76520039632231\\
31.2	8.81348900130649\\
32.7	8.90084944267076\\
34.5	8.95253619772151\\
35.1	8.96299172587777\\
35.7	8.97399426268512\\
37.5	8.97739579883122\\
37.8	8.97262105935866\\
38.4	8.97678216918968\\
39	8.96493612559439\\
41.4	8.89444200624848\\
42.6	8.82267729379197\\
43.5	8.75342051041267\\
44.7	8.66181282357408\\
45.3	8.60656328538437\\
45.9	8.54919619353771\\
46.8	8.45507492769929\\
47.4	8.38374856907424\\
47.7	8.34037899419121\\
48.3	8.2649614951421\\
49.5	8.07945269975204\\
50.7	7.8786822761775\\
51	7.82682879255148\\
51.9	7.65190924628209\\
52.5	7.53678092095214\\
53.4	7.34533788869214\\
55.5	6.85951176197838\\
55.8	6.77901296742107\\
56.1	6.71171326183665\\
57.6	6.32669662562296\\
59.4	5.80888547576234\\
60.9	5.33746285626168\\
61.2	5.25352397677779\\
62.1	4.98150259741936\\
63.3	4.62229023461845\\
65.1	4.04661443476596\\
65.7	3.85582271860517\\
66	3.76421502022939\\
66.3	3.65982459496115\\
67.5	3.28608994585244\\
68.4	3.03318567017138\\
69	2.88222398116744\\
69.3	2.80670417125739\\
69.9	2.68004550664335\\
70.5	2.57149805025418\\
71.7	2.40295017944933\\
72.6	2.30452824774468\\
75.9	1.96342483728591\\
76.8	1.86160325199219\\
78.6	1.64903571666267\\
80.4	1.44631709652566\\
82.2	1.2476218510615\\
83.1	1.13834315967259\\
83.7	1.06286789222239\\
84.6	0.951695789105088\\
87.3	0.639351737971623\\
88.2	0.549202791389121\\
89.4	0.447941403569416\\
90.3	0.392951614012873\\
91.2	0.356509666510306\\
93.6	0.320177149267806\\
95.4	0.331109915916187\\
96.9	0.363402572863237\\
98.1	0.388837430486717\\
102	0.46774650784829\\
103.5	0.496414646932038\\
104.7	0.507224798464136\\
108	0.526291545365893\\
108.9	0.533078368305155\\
110.4	0.522484941228129\\
111.6	0.50571534486474\\
113.4	0.486758529880291\\
120	0.428016380106044\\
122.1	0.432078598649781\\
125.4	0.452983743600697\\
127.5	0.436433087333739\\
131.1	0.402711407295953\\
132.3	0.392939472752204\\
133.5	0.380009795209702\\
136.5	0.358956681234531\\
138.3	0.360985172529752\\
139.2	0.365364639205353\\
140.1	0.36959821577102\\
};
\addplot [color=blue, dashed, line width=0.7pt, forget plot]
  table[row sep=crcr]{%
0	3.1297425829964\\
7.8	3.13689267112844\\
9.3	3.15167904501424\\
11.1	3.18966766384028\\
12.3	3.20747281252773\\
16.5	3.23686559249344\\
17.4	3.24115297455769\\
};
\addplot [color=green, dashed, line width=0.7pt, forget plot]
  table[row sep=crcr]{%
0	1.21012879626413\\
1.2	1.21653825236139\\
1.5	1.23661952552109\\
1.8	1.26737069969766\\
2.1	1.30951313493285\\
2.4	1.36333021264165\\
2.7	1.43106716678163\\
3	1.50796774096931\\
3.6	1.68997109470038\\
3.9	1.79035471214742\\
4.5	2.02592311672043\\
5.1	2.2553299004086\\
5.7	2.46971991304617\\
6.3	2.66121144423252\\
6.6	2.74504627504315\\
7.2	2.89755689758036\\
7.8	3.02472026944679\\
8.1	3.07928101460648\\
8.4	3.11012485850182\\
9.3	3.12967378420568\\
10.5	3.15460623084098\\
11.4	3.16691580856026\\
};
\addplot [color=red, dashed, line width=0.7pt, forget plot]
  table[row sep=crcr]{%
0	3.4219914187788\\
0.6	3.49562559900748\\
1.2	3.58229815136923\\
1.8	3.68286164090129\\
2.1	3.73767849019179\\
3.3	3.98257083415782\\
3.9	4.10043485224033\\
4.8	4.25760213102724\\
5.4	4.34172923381355\\
6	4.40514490832422\\
6.6	4.45021623954125\\
6.9	4.46309007005276\\
10.2	4.4949197600386\\
10.8	4.49815187050953\\
};
\end{axis}
\end{tikzpicture}%